\documentclass[english]{smfart}
\usepackage{smfthm,smfenum,bull}
\usepackage{amsfonts,amssymb,amscd,amsmath,latexsym}
\usepackage{enumerate}
\usepackage[english]{babel}
\usepackage[dvips]{graphicx}

\theoremstyle{plain}

\def\surl#1_#2{\mathrel{\mathop{\kern 0pt #1}\limits_{#2}}}

\newcounter{definition}
\author{Franck LESIEUR}
\address{LMNO\\
         Universit\'e de Caen\\
         BP 5186\\
         F-14032 Caen Cedex}
\email{franck.lesieur@math.unicaen.fr}
\urladdr{http://www.math.unicaen.fr/~lesieur}
\title[MEASURED QUANTUM GROUPOIDS]{Measured quantum groupoids}
\alttitle{Groupo\"{\i}des quantiques mesur\'es}

\begin{document}
\frontmatter

\begin{abstract}
In this article, we give a definition for measured quantum
groupoids. We want to get objects with duality extending both
quantum groups and groupoids. We base ourselves on J. Kustermans and
S. Vaes' works about locally compact quantum groups that we
generalize thanks to formalism introduced by M. Enock and J.M.
Vallin in the case of inclusion of von Neumann algebras. From a
structure of Hopf-bimodule with left and right invariant
operator-valued weights, we define a fundamental
pseudo-multiplicative unitary. To get a satisfying duality in the
general case, we assume the existence of an antipode given by its
polar decomposition. This theory is illustrated with many examples
among others inclusion of von Neumann algebras (M. Enock) and a sub
family of measured quantum groupoids with easier axiomatic.
\end{abstract}

\begin{altabstract}
Dans cet article, on d\'efinit une notion de groupo\"{\i}des
quantiques mesur\'es. On cherche \`a obtenir des objets munis d'une
dualit\'e qui \'etend celle des groupo\"{\i}des et des groupes
quantiques. On s'appuie sur les travaux de J. Kustermans et S. Vaes
concernant les groupes quantiques localement compacts qu'on
g\'en\'eralise gr\^{a}ce au formalisme introduit par M. Enock et
J.M. Vallin \`a propos des inclusions d'alg\`ebres de von Neumann.
\`A partir d'un bimodule de Hopf muni de poids op\'eratoriels
invariants \`a gauche et \`a droite, on d\'efinit un unitaire
pseudo-multiplicatif fondamental. Pour obtenir une dualit\'e
satisfaisante dans le cas g\'en\'eral, on suppose l'existence d'une
antipode d\'efinie par sa d\'ecomposition polaire. Cette th\'eorie
est illustr\'ee dans une derni\`ere partie par de nombreux exemples
notamment les inclusions d'alg\`ebres de von Neumann (M. Enock) et
une sous famille de groupo\"{\i}des quantiques mesur\'es \`a
l'axiomatique plus simple.
\end{altabstract}

 \altkeywords{Groupo\"{\i}des
quantiques ; antipode ; unitaire pseudo-multiplicatif}
\thanks{The author is mostly indebted to Michel Enock, Stefaan Vaes,
Leonid Va\u{\i}nerman and Jean-Michel Vallin for many fruitful
conversations.}

\maketitle \mainmatter

\tableofcontents
\section{Introduction}

\subsection{Historic}
Theory of quantum groups has lot of developments in operator
algebras setting. Many contributions are given by \cite{KaV},
\cite{W6}, \cite{ES}, \cite{MN}, \cite{BS}, \cite{W1}, \cite{W},
\cite{VD2}, \cite{KV1}. In particular, J. Kustermans and S. Vaes'
work is crucial: in \cite{KV1}, they propose a simple definition for
locally compact quantum groups which gathers all known examples
(locally compact groups, quantum compacts groupe \cite{W1}, quantum
group $ax+b$ \cite{W3}, \cite{W4}, Woronowicz' algebra \cite{MN}...)
and they find a general framework for duality of theses objects. The
very few number of axioms gives the theory a high manageability
which is proved with recent developments in many directions (actions
of locally compact quantum groups \cite{Vae2}, induced
co-representations \cite{Ku2}, cocycle bi-crossed products
\cite{VaV}). They complete their work with a theory of locally
compact quantum groups in the von Neumann setting \cite{KV2}.\\

In geometry, groups are rather defined by their actions. Groupoids
category contains groups, group actions and equivalence relation. It
is used by G.W Mackey and P. Hahn (\cite{Ma}, \cite{Ha1} and
\cite{Ha2}), in a measured version, to link theory of groups and
ergodic theory. Locally compact groupoids and the operator theory
point of view are introduced and studied by J. Renault in \cite{R1}
and \cite{R2}. It covers many interesting examples in differential
geometry \cite{C2} e.g holonomy groupoid of a foliation.\\

In \cite{V1}, J.M Vallin introduces the notion of Hopf bimodule from
which he is able to prove a duality for groupoids. Then, a natural
question is to construct a category, containing quantum groups and
groupoids, with a duality theory.\\

In the quantum group case, duality is essentially based on a
multiplicative unitary \cite{BS}. To generalize the notion up to the
groupoid case, J.M Vallin introduces pseudo-multiplicative
unitaries. In \cite{V2}, he exhibits such an object coming from Hopf
bimodule structures for groupoids. Technically speaking,
Connes-Sauvageot's theory of relative tensor products is intensively
used.\\

In the case of depth $2$ inclusions of von Neumann algebras, M.
Enock and J.M Vallin, and then, M. Enock underline two "quantum
groupoids" in duality. They also use Hopf bimodules and
pseudo-multiplicative unitaries. At this stage, a non trivial
modular theory on the basis (the equivalent for units of a groupoid)
is revealed to be necessary and a simple generalization of axioms
quantum groups is not sufficient to construct quantum groupoids
category: we have to add an axiom on the basis \cite{E1} i.e we use
a special weight to do the construction. The results are improved in
\cite{E3}.\\

In \cite{E2}, M. Enock studies in detail pseudo-multiplicative
unitaries and introduces an analogous notion of S. Baaj and G.
Skandalis' regularity. In quantum groups, the fundamental
multiplicative unitary is weakly regular and manageable in the sense
of Woronowicz. Such properties have to be satisfied in quantum
groupoids. Moreover, M. Enock defines and studies compact (resp.
discrete) quantum groupoids which have to enter into the general
theory.\\

Lot of works have been led about quantum groupoids but essentially
in finite dimension. We have to quote weak Hopf
$\text{C}^*$-algebras introduced by G. B\"{o}hm, F. Nill and K.
Szlach\'{a}nyi \cite{BNSz}, \cite{BSz}, and then studied by F. Nill
and L. Va\u{\i}nerman \cite{Ni}, \cite{N}, \cite{NV1}, \cite{NV2}.
J.M Vallin develops a quantum groupoids theory in finite dimension
thanks to multiplicative partial isometries \cite{V3}, \cite{V4}. He
proves that his theory coincide exactly with weak Hopf
$\text{C}^*$-algebras.

\subsection{Aims and Methods}
In this article, we propose a definition for measured quantum
groupoids in any dimensions. "Measured" means we are in the von
Neumann setting and we assume existence of the analogous of a
measure. We use a similar approach as J. Kustermans and S. Vaes'
theory with the formalism of Hopf bimodules and
pseudo-multiplicative unitaries. The notion has to recover all known examples and shall extend
their duality if already exist.\\

In our setting, we assume the existence of a scaling group and a
coinvolution so that we are much more closer to \cite{MNW}. Then, we
are able to construct a dual structure for theses
objects and we prove a duality theorem. We also get uniqueness of the equivalent of Haar measure.\\

We want to give many examples. First of all, we present a family of
measured quantum groupoids of a particular interest: the axiomatic
of them is easier than the general measured quantum groupoids and
very similar to J. Kustermans and S. Vaes axiomatic of locally
compact quantum groups because we can construct the antipode.
However, this new category is not self dual but we can characterize
their dual objects. Then we are interested in depth 2 inclusions of
von Neumann algebras of Enock's type which are included in our
theory and for which we can compute the dual structure. In a
forthcoming article, we will study an example of the type $G=G_1G_2$
where $G_1$ and $G_2$ are
two groupoids such that $G_1\cap G_2=G^{(0)}$.\\

We are inspired by technics develop by J. Kustermans and S. Vaes
about locally compact quantum groups in the von Neumann setting
\cite{KV1}, by M. Enock \cite{E3} as far as the density theorem is
concerned which is a key tool for duality and by author's thesis
\cite{Les}.

\subsection{Study plan}
After brief recalls about tools and technical points, we define
objects we will use. We start by associating a fundamental
pseudo-multiplicative unitary to every Hopf bimodule with invariant
operator-valued weights. In fact, we shall define several isometries
depending on which operator-valued weight we use. Each of them are
useful, especially as far as the proof of unitarity of the
fundamental isometry is concerned. This point can be also noticed in
the crucial paper of S. Baaj and G. Skandalis \cite{BS} where they
need a notion of irreducible unitary that means there exist another
unitary. Also, in \cite{KV1}, they need to introduce several
unitaries. The fundamental unitary gathers all informations on the
structure so that we can re-construct von Neumann algebra and
co-product.\\

In the first part, we give axioms of measured quantum groupoids. In
this setting, we construct a modulus, which corresponds to modulus
of groupoids and a scaling operator which corresponds to scaling
factor in locally compact quantum groups. They come from
Radon-Nikodym's cocycle of right invariant operator-valued weight
with respect to left invariant one thanks to proposition 5.2 of
\cite{Vae}. Then, we prove uniqueness of
invariant operator-valued weight up to an element of basis center.\\

Also, we prove a "manageability" property of the fundamental
pseudo-multiplicative unitary. A density result concerning bounded
elements can be handled. These are sections \ref{dtheo} and
\ref{mtheo}. They give interesting results on the structure and a
necessary preparation step for duality.\\

Then, we can proceed to the construction of the dual structure and
get a duality theorem. \\

The second part is devoted to examples. We have a lot of examples
for locally compact quantum groups thanks to Woronowicz \cite{W2},
\cite{W3}, \cite{W4}, \cite{W5} and the cocycle bi-crossed product
due to S. Vaes and L. Va\u{\i}nerman \cite{VaV}. Theory of measured
quantum groupoids has also a lot of examples.\\

First, we lay stress on Hopf bimodule with invariant operator-valued
weights which are "adapted" in a certain sense. This hypothesis
corresponds to the choice of a special weight on the basis to do the
constructions (like in the groupoid case with a quasi-invariant
measure on $G^{\{0\}}$). For them, we are able to construct the
antipode $S$, the polar decomposition of which is given by a
co-involution $R$ and a one-parameter group of automorphisms called
scaling group $\tau$. In particular, we show that $S,R$ and $\tau$
are independent of operator-valued weights.\\

Then we explain how these so-called adapted measured quantum
groupoids fit into our measured quantum groupoids. In this setting,
we develop informations about modulus, scaling operator and
uniqueness. We also characterize them and their dual in the general
theory. Groupoids, weak Hopf $\text{C}^*$-algebras, quantum groups,
quantum groupoids of compact
(resp. discrete) type\ldots are of this type.\\

Depth $2$ inclusions of von Neumann algebras also enter into our
general setting (but not in adapted measured quantum groupoids
unless the
basis is semi-finite) and we compute their dual.\\

Finally, we state stability of the category by direct sum (which
reflects the stability of groupoids under disjoint unions), finite
tensor product and direct integrals. Then, we are able to construct
new examples: in particular we can exhibit quantum
groupoids with non scalar scaling operator.\\

\section{Recalls}
\subsection{Weights and operator-valued weights \cite{St}, \cite{T}}

Let $N$ be a von Neumann and $\psi$ a normal, semi-finite faithful
(n.s.f) weight on $N$; we denote by ${\mathcal N}_{\psi}$,
${\mathcal M}_{\psi}$, $H_{\psi}$, $\pi_{\psi}$, $\Lambda_{\psi}$,
$J_{\psi}$, $\Delta_{\psi}\ldots$ canonical objects of Tomita's
theory with respect to (w.r.t) $\psi$.

\begin{defi}
Let denote by ${\mathcal T}_{\psi}$ {\bf Tomita's algebra} w.r.t
$\psi$ defined by:
$$\{ x \in {\mathcal N}_{\psi } \cap {\mathcal N}_{\psi }^*|\ x \text{
analytic w.r.t }\sigma^{\psi}\text{ such that } \sigma_z^{\psi}(x)
\in {\mathcal N}_{\psi } \cap {\mathcal N}_{\psi }^* \text{ for all
} z \in {\mathbb C} \}$$
\end{defi}

By (\cite{St}, 2.12), we have the following approximating result:

\begin{lemm}\label{analyse}
For all $x\in {\mathcal N}_{\psi}$, there exists a sequence
$(x_n)_{n\in \mathbb{N}}$ of ${\mathcal T}_{\psi}$ such that:

\begin{center}
\begin{minipage}{13cm}
\begin{enumerate}[i)]
\item $||x_n||\leq ||x||$ for all $n\in\mathbb{N}$;
\item $(x_n)_{n\in\mathbb{N}}$ converges to $x$ in the strong topology;
\item $(\Lambda_{\psi}(x_n))_{n\in\mathbb{N}}$ converges to
$\Lambda_{\psi}(x)$ in the norm topology of $H_{\psi}$.
\end{enumerate}
\end{minipage}
\end{center}
Moreover, if $x\in {\mathcal N}_{\psi}\cap {\mathcal N}_{\psi}^*$,
then we have:

\begin{center}
\begin{minipage}{13cm}
\begin{enumerate}[iv)]
\item $(x_n)_{n\in\mathbb{N}}$ converges to $x$ in the *-strong topology;
\item $(\Lambda_{\psi}(x_n^*))_{n\in\mathbb{N}}$ converges to
$\Lambda_{\psi}(x^*)$ in the norm topology of $H_{\psi}$.
\end{enumerate}
\end{minipage}
\end{center}
\end{lemm}

Let $N\subset M$ be an inclusion of von Neumann algebras and $T$ a
normal, semi-finite, faithful (n.s.f) operator-valued weight from
$M$ to $N$. We put:
$${\mathcal N}_T=\left\lbrace x\in M\,
/\, T(x^*x)\in N^+ \right\rbrace  \text{ and } {\mathcal
M}_T={\mathcal N}_T^*{\mathcal N}_T$$ We can define a n.s.f weight
$\psi\circ T$ on $M$ in a natural way.  Let us recall theorem 10.6
of \cite{EN}:
\begin{prop}\label{preintro}
Let $N\subset M$ be an inclusion of von Neumann algebras and $T$ be
a normal, semi-finite, faithful (n.s.f) operator-valued weight from
$M$ to $N$ and $\psi$ a n.s.f weight on $N$. Then we have:

\begin{enumerate}[i)]
\item for all
$x\in {\mathcal N}_T$ and $a\in {\mathcal N}_{\psi}$, $xa$ belongs
to ${\mathcal N}_T\cap {\mathcal N}_{\psi\circ T}$, there exists
$\Lambda_T(x)\in\text{Hom}_{N^o}(H_{\psi},H_{\psi\circ T})$ such
that: $$\Lambda_T(x)\Lambda_{\psi}(a)=\Lambda_{\psi\circ T}(xa)$$
and $\Lambda_T$ is a morphism of $M-N$-bimodules from ${\mathcal
N}_T$ to $\text{Hom}_{N^o}(H_{\psi},H_{\psi\circ T})$;
\item ${\mathcal N}_T\cap {\mathcal N}_{\psi\circ T}$ is a weakly
dense ideal of $M$ and $\Lambda_{\psi\circ T}({\mathcal N}_T\cap
{\mathcal N}_{\psi\circ T})$ is dense in $H_{\psi\circ T}$,
$\Lambda_{\psi\circ T}({\mathcal N}_T\cap {\mathcal N}_{\psi\circ
T}\cap {\mathcal N}^*_T\cap {\mathcal N}^*_{\psi\circ T})$ is a core
for $\Delta_{\psi\circ T}^{1/2}$ and $\Lambda_T({\mathcal N}_T)$ is
dense in $\text{Hom}_{N^o}(H_{\psi},H_{\psi\circ T})$ for the
s-topology defined by (\cite{BDH}, 1.3);
\item for all $x\in {\mathcal N}_T$ and
$z\in {\mathcal N}_T\cap {\mathcal N}_{\psi\circ T}$, $T(x^*z)$
belongs to ${\mathcal N}_{\psi}$ and:
$$\Lambda_T(x)^*\Lambda_{\psi\circ T}(z)=\Lambda_{\psi}(T(x^*z))$$
\item for all $x,y\in {\mathcal
N}_T$: $$\Lambda_T(y)^*\Lambda_T(x)=\pi_{\psi}(T(x^*y)) \text{ and }
||\Lambda_T(x)||=||T(x^*x)||^{1/2}$$ and $\Lambda_T$ is injective.
\end{enumerate}
\end{prop}

Let us also recall lemma 10.12 of \cite{EN}:

\begin{prop}
Let $N\subseteq M$ be an inclusion of von Neumann algebras, $T$ a
n.s.f operator-valued weight from $M$ to $N$, $\psi$ a n.s.f weight
on $N$ and $x\in {\mathcal M}_T\cap {\mathcal M}_{\psi\circ T}$. If
we put:
$$x_n=\sqrt\frac{n}{\pi}\int^{+\infty}_{-\infty}\!
e^{-nt^2}\sigma_t^{\psi\circ T}(x)\ dt$$ then $x_n$ belongs to
${\mathcal M}_T\cap {\mathcal M}_{\psi\circ T}$ and is analytic
w.r.t $\psi\circ T$. The sequence converges to $x$ and is bounded by
$||x||$. Moreover, $(\Lambda_{\psi\circ T}(x_n))_{n\in\mathbb{N}}$
converges to $\Lambda_{\psi\circ T}(x)$ and  $\sigma_z^{\psi\circ
T}(x_n)\in {\mathcal M}_T\cap {\mathcal M}_{\psi\circ T}$ for all
$z\in\mathbb{C}$.
\end{prop}

\begin{defi}
The set of $x\in {\mathcal N}_{\Phi } \cap {\mathcal N}_{\Phi }^*
\cap {\mathcal N}_T \cap {\mathcal N}_T^*$, analytic w.r.t
$\sigma^{\Phi}$ such that $\sigma_z^{\Phi}(x)\in {\mathcal N}_{\Phi
} \cap {\mathcal N}_{\Phi }^* \cap {\mathcal N}_T \cap {\mathcal
N}_T^*$ for all $z\in {\mathbb C}$ is denoted by ${\mathcal
T}_{\Phi}$ and is called {\bf Tomita's algebra} w.r.t $\psi\circ
T=\Phi$ and $T$.
\end{defi}
Lemma \ref{analyse} is still satisfied with Tomita's algebra w.r.t
$\Phi$ and $T$.

\subsection{Spatial theory \cite{C1},
\cite{S2}, \cite{T}}\label{intre}

Let $\alpha$ be a normal, non-degenerated representation of $N$ on a
Hilbert space $H$. So, $H$ becomes a left $N$-module and we write $\
_{\alpha}H$.

\begin{defi}\cite{C1}
An element $\xi$ of $\ _{\alpha}H$ is said to be bounded w.r.t
$\psi$ if there exists $C\in\mathbb{R}^+$ such that, for all $y\in
{\mathcal N}_{\psi}$, we have $||\alpha(y)\xi||\leq C||
\Lambda_{\psi}(y) ||$. The set of {\bf bounded elements} w.r.t
$\psi$ is denoted by $D(_{\alpha}H,\psi)$.
\end{defi}

By \cite{C1} (lemma 2), $D(_{\alpha}H,\psi)$ is dense in $H$ and
$\alpha(N)'$-stable. An element $\xi$ of $D(_{\alpha}H,\psi)$ gives
rise to a bounded operator $R^{\alpha,\psi}(\xi)$ of
$Hom_N(H_{\psi},H)$ such that, for all $y\in {\mathcal N}_{\psi}$:
$$R^{\alpha,\psi}(\xi)\Lambda_{\psi}(y)=\alpha(y)\xi$$
For all $\xi,\eta \in D(_{\alpha}H,\psi)$, we put:
$$\theta^{\alpha,\psi}(\xi,\eta)=R^{\alpha,\psi}(\xi)R^{\alpha,\psi}(\eta)^*\text{ and }
<\xi,\eta>_{\alpha,\psi}=R^{\alpha,\psi}(\eta)^*R^{\alpha,\psi}(\xi)^*$$
By \cite{C1} (lemma 2), the linear span of
$\theta^{\alpha,\psi}(\xi,\eta)$ is a weakly dense ideal of
$\alpha(N)'$. $<\xi,\eta>_{\alpha,\psi}$ belongs to
$\pi_{\psi}(N)'=J_{\psi}\pi_{\psi}(N)J_{\psi}$ which is identified
with the opposite von Neumann algebra $N^o$. The linear span of
$<\xi,\eta>_{\alpha,\psi}$ is weakly dense in $N^o$.

\noindent By \cite{C1} (proposition 3), there exists a net
$(\eta_i)_{i\in I}$ of $D(_{\alpha}H,\psi)$ such that:
$$\sum_{i\in I}\theta^{\alpha,\psi}(\eta_i,\eta_i)=1$$ Such a net is called a
$(N,\psi)$-\textbf{basis} of $\ _{\alpha}H$. By \cite{EN}
(proposition 2.2), we can choose $\eta_i$ such that
$R^{\alpha,\psi}(\eta_i)$ is a partial isometry with two-by-two
orthogonal final supports and such that
$<\eta_i,\eta_j>_{\alpha,\psi}=0$ unless $i=j$. In the following, we
assume these properties hold for all $(N,\psi)$-basis of $\
_{\alpha}H$.

Now, let $\beta$ be a normal, non-degenerated anti-representation
from $N$ on $H$. So $H$ becomes a right $N$-module and we write
$H_{\beta}$. But $\beta$ is also a representation of $N^o$. If
$\psi^o$ is the n.s.f weight on $N^o$ coming from $\psi$ then
${\mathcal N}_{\psi^o}={\mathcal N}_{\psi}^*$ and we identify
$H_{\psi^o}$ with $H_{\psi}$ thanks to:
$$(\Lambda_{\psi^o}(x^*)\mapsto J_{\psi}\Lambda_{\psi}(x))$$

\begin{defi}\cite{C1}
An element $\xi$ of $H_{\beta}$ is said to be bounded w.r.t $\psi^o$
if there exists $C\in\mathbb{R}^+$ such that, for all $y\in
{\mathcal N}_{\psi}$, we have $||\beta(y^*)\xi||\leq C||
\Lambda_{\psi}(y) ||$. The set of {\bf bounded elements} w.r.t
$\psi^o$ is denoted by $D(H_{\beta},\psi^o)$.
\end{defi}

$D(_{\alpha}H,\psi)$ is dense in $H$ and $\beta(N)'$-stable. An
element $\xi$ of $D(H_{\beta},\psi^o)$ gives rise to a bounded
operator $R^{\beta,\psi^o}(\xi)$ of $Hom_{N^o}(H_{\psi},H)$ such
that, for all $y\in {\mathcal N}_{\psi}$:
$$R^{\beta,\psi^o}(\xi)\Lambda_{\psi}(y)=\beta(y^*)\xi$$
For all $\xi,\eta \in D(H_{\beta},\psi^o)$, we put:
$$\theta^{\beta,\psi^o}(\xi,\eta)=R^{\beta,\psi^o}(\xi)R^{\beta,\psi^o}(\eta)^*\text{ and }
<\xi,\eta>_{\beta,\psi^o}=R^{\beta,\psi^o}(\eta)^*R^{\beta,\psi^o}(\xi)^*$$
The linear span of $\theta^{\beta,\psi^o}(\xi,\eta)$ is a weakly
dense ideal of $\beta(N)'$. $<\xi,\eta>_{\beta,\psi^o}$ belongs to
$\pi_{\psi}(N)$ which is identified with $N$. The linear span of
$<\xi,\eta>_{\beta,\psi^o}$ is weakly dense in $N$. In fact, we know
that $<\xi,\eta>_{\beta,\psi^o}\in {\mathcal M}_{\psi}$ by \cite{C1}
(lemma 4) and by \cite{S2} (lemma 1.5), we have
$$\Lambda_{\psi}(<\xi,\eta>_{\beta,\psi^o})=R^{\beta,\psi^o}(\eta)^*\xi$$
A net $(\xi_i)_{i\in I}$ of $\psi^o$-bounded elements of is said to
be a $(N^o,\psi^o)$-basis of $H_{\beta}$ if:
$$\sum_{i\in I}\theta^{\beta,\psi^o}(\xi_i,\xi_i)=1$$
and if $\xi_i$ such that $R^{\beta,\psi^o}(\xi_i)$ is a partial
isometry with two-by-two orthogonal final supports and such that
$<\xi_i,\xi_j>_{\alpha,\psi}=0$ unless $i=j$. Therefore, we have:
$$R^{\beta,\psi^o}(\xi_i)=\theta^{\beta,\psi^o}(\xi_i,\xi_i)R^{\beta,\psi^o}(\xi_i)=
R^{\beta,\psi^o}(\xi_i)<\xi_i,\xi_i>_{\beta,\psi^o}$$ and, for all
$\xi\in D(H_{\beta},\psi^o)$:
$$\xi=\sum_{i\in I}
R^{\beta,\psi^o}(\xi_i)\Lambda_{\psi}(<\xi,\xi_i>_{\beta,\psi^o})$$

\begin{prop}(\cite{E2}, proposition 2.10)
Let $N\subseteq M$ be an inclusion of von Neumann algebras and $T$
be a n.s.f operator-valued weight from $M$ to $N$. There exists a
net $(e_i)_{i\in I}$ of ${\mathcal N}_{T}\cap {\mathcal N}_{T}^*\cap
{\mathcal N}_{\psi\circ T}\cap {\mathcal N}_{\psi\circ T}^*$ such
that $\Lambda_{T}(e_i)$ is a partial isometry, $T(e_j^*e_i)=0$
unless $i=j$ and with orthogonal final supports of sum $1$.
Moreover, we have $e_i=e_iT(e_i^*e_i)$ for all $i\in I$, and, for
all $x\in {\mathcal N}_{T}$:
$$\Lambda_{T}(x)=\sum_{i\in I}\Lambda_{T}(e_i)T(e_i^*x)\quad\text{ and }\quad
x=\sum_{i\in I} e_iT(e_i^*x)$$ in the weak topology. Such a net is
called a basis for $(T,\psi^o)$. Finally, the net
$(\Lambda_{\psi\circ T}(e_i))_{i\in I}$ is a $(N^o,\psi^o)$-basis of
$(H_{\psi\circ T})_s$ where $s$ is the anti-representation which
sends $y\in N$ to $J_{\psi\circ T}y^*J_{\psi\circ T}$.\label{partis}
\end{prop}

\subsection{Relative tensor product \cite{C1}, \cite{S2}, \cite{T}}
Let $H$ and $K$ be Hilbert space. Let $\alpha$ (resp. $\beta$) be a
normal and non-degenerated (resp. anti-) representation of $N$ on
$K$ (resp. $H$). Let $\psi$ be a n.s.f weight on $N$. Following
\cite{S2}, we put on $D(H_{\beta},\psi^o)\odot K$ a scalar product
defined by:
$$(\xi_1 \odot \eta_1|\xi_2 \odot \eta_2)=(\alpha(<\xi_1,\xi_2>_{\beta,\psi^o})\eta_1|\eta_2)$$
for all $\xi_1,\xi_2\in D(H_{\beta},\psi^o)$ and $\eta_1,\eta_2\in
K$. We have identified $\pi_{\psi}(N)$ with $N$.

\begin{defi}
The completion of $D(H_{\beta},\psi^o)\odot K$ is called {\bf
relative tensor product} and is denoted by $H\surl{\ _{\beta}
\otimes_{\alpha}}_{\ \psi}K$.
\end{defi}

The image of $\xi\odot\eta$ in $H\surl{\ _{\beta}
\otimes_{\alpha}}_{\ \psi} K$ is denoted by $\xi\surl{\ _{\beta}
\otimes_{\alpha}}_{\ \psi}\eta$. One should bear in mind that, if we
start from another n.s.f weight $\psi'$ on $N$, we get another
Hilbert space which is canonically isomorphic to $H\surl{\ _{\beta}
\otimes_{\alpha}}_{\ \psi}K$ by (\cite{S2}, proposition 2.6).
However this isomorphism does not send $\xi\surl{\ _{\beta}
\otimes_{\alpha}}_{\ \psi}\eta$ on $\xi\surl{\ _{\beta}
\otimes_{\alpha}}_{\ \ \psi '}\eta$.

By \cite{S2} (definition 2.1), relative tensor product can be
defined from the scalar product:
$$(\xi_1 \odot \eta_1| \xi_2 \odot \eta_2)=(\beta(<\eta_1,\eta_2>_{\alpha,\psi})\xi_1|\xi_2)$$
for all $\xi_1,\xi_2 \in H$ and $\eta_1,\eta_2\in
D(_{\alpha}K,\psi)$ that's why we can define a one-to-one flip from
$H\surl{\ _{\beta} \otimes_{\alpha}}_{\ \psi}K$ onto $K \surl{\
_{\alpha} \otimes_{\beta}}_{\ \ \psi^o}H$ such that:
$$\sigma_{\psi}(\xi\surl{\ _{\beta} \otimes_{\alpha}}_{\ \psi}\eta)=\eta \surl{\
_{\alpha} \otimes_{\beta}}_{\ \ \psi^o} \xi$$ for all $\xi\in
D(H_{\beta},\psi)$ (resp. $\xi\in H$) and $\eta\in K$ (resp.
$\eta\in D(_{\alpha}K,\psi)$). The flip gives rise at the operator
level to $\varsigma_{\psi}$ from $\mathcal{L}(H \surl{\ _{\beta}
\otimes_{\alpha}}_{\ \psi} K)$ onto $\mathcal{L}(K \surl{\ _{\alpha}
\otimes_{\beta}}_{\ \ \psi^o} H)$ such that:
$$\varsigma_{\psi}(X)=\sigma_{\psi}X\sigma_{\psi}^*$$
Canonical isomorphisms of change of weights send $\varsigma_{\psi}$
on $\varsigma_{\psi '}$ so that we write $\varsigma_N$ without any
reference to the weight on $N$.

For all $\xi\in D(H_{\beta},\psi^o)$ and $\eta\in
D(_{\alpha}K,\psi)$, we define bounded operators:

$$
\begin{aligned}
\lambda^{\beta,\alpha}_{\xi}: K &\rightarrow H\surl{\ _{\beta}
\otimes_{\alpha}}_{\ \psi} K& \text{ and }\quad
\rho^{\beta,\alpha}_{\eta}: H &\rightarrow H\surl{\ _{\beta}
\otimes_{\alpha}}_{\ \psi} K \\
\eta &\mapsto \xi\surl{\ _{\beta} \otimes_{\alpha}}_{\ \psi}\eta&
\xi &\mapsto \xi\surl{\ _{\beta} \otimes_{\alpha}}_{\ \psi}\eta
\end{aligned}$$
Then, we have:
$$(\lambda^{\beta,\alpha}_{\xi})^*\lambda^{\beta,\alpha}_{\xi}=
\alpha(<\xi,\xi>_{\beta,\psi^o})\text{ and }
(\rho^{\beta,\alpha}_{\eta})^*\rho^{\beta,\alpha}_{\eta}=
\beta(<\eta,\eta>_{\alpha,\psi})$$

By \cite{S2} (remark 2.2), we know that $D(_{\alpha}K,\psi)$ is
$\alpha(\sigma_{-i/2}^{\psi}({\mathcal
D}(\sigma_{-i/2}^{\psi})))$-stable and for all $\xi\in H$, $\eta\in
D(_{\alpha}K,\psi)$ and $y\in {\mathcal D}(\sigma_{-i/2}^{\psi})$,
we have:
$$\beta(y)\xi\surl{\ _{\beta} \otimes_{\alpha}}_{\ \psi}\eta=
\xi\surl{\ _{\beta} \otimes_{\alpha}}_{\
\psi}\alpha(\sigma_{-i/2}^{\psi}(y))\eta$$

\begin{lemm}\label{aidintre}
If $\xi'\surl{\ _{\beta} \otimes_{\alpha}}_{\ \psi}\eta=0$ for all
$\xi'\in D(H_{\beta},\psi^o)$ then $\eta=0$.
\end{lemm}

\begin{proof}
For all $\xi,\xi'\in D(H_{\beta},\psi^o)$, we have:
$$\alpha(<\xi',\xi>_{\beta,\psi^o})\eta=(\lambda_{\xi}^{\beta,\alpha})^*
\lambda_{\xi'}^{\beta,\alpha}\eta=(\lambda_{\xi}^{\beta,\alpha})^*(\xi'
\surl{\ _{\beta} \otimes_{\alpha}}_{\ \psi} \eta)=0$$ Since the
linear span of $<\xi',\xi>_{\beta,\psi^o}$ is dense in $N$, we get
$\eta=0$.
\end{proof}

\begin{prop}\label{preli}
Assume $H\neq\{0\}$. Let $K'$ be a closed subspace of $K$ such that
$\alpha(N)K'\subseteq K'$. Then:
$$H\surl{\ _{\beta}
\otimes_{\alpha}}_{\ \psi} K=H\surl{\ _{\beta} \otimes_{\alpha}}_{\
\psi} K'\quad\Rightarrow\quad K=K'$$
\end{prop}

\begin{proof}
Let $\eta \in K'^{\bot}$. For all $\xi,\xi'\in D(H_{\beta},\psi^o)$
and $k \in K'$, we have:
$$(\xi \surl{\
_{\beta} \otimes_{\alpha}}_{\ \psi} k|\xi' \surl{\ _{\beta}
\otimes_{\alpha}}_{\ \psi}
\eta)=(\alpha(<\xi,\xi'>_{\beta,\nu^o})k|\eta)=0$$ Therefore, for
all $\xi'\in D(H_{\beta},\psi^o)$, we have:
$$\xi'\surl{\ _{\beta}
\otimes_{\alpha}}_{\ \psi} \eta \in (H\surl{\ _{\beta}
\otimes_{\alpha}}_{\ \psi} K')^{\bot}=(H\surl{\ _{\beta}
\otimes_{\alpha}}_{\ \psi} K)^{\bot}=\{0\}$$ By the previous lemma,
we get $\eta=0$ and $K=K'$.
\end{proof}

Let $H'$, $K'$, $\alpha'$ and $\beta'$ like $H$, $K$, $\alpha$ and
$\beta$. Let $A\in \mathcal{L}(H,H')$ and $B\in \mathcal{L}(K,K')$
such that: $$\forall n\in N,\quad A\beta(n)=\beta'(n)A\quad\text{
and}\quad B\alpha(n)=\alpha'(n)B$$ Then we can define an operator
$A\surl{\ _{\beta}\otimes_{\alpha}}_{\ \psi}B\in\mathcal{L}(H\surl{\
_{\beta}\otimes_{\alpha}}_{\ \psi}K,H'\surl{\
_{\beta}\otimes_{\alpha}}_{\ \psi}K')$ which naturally acts on
elementary tensor products. In particular, if
$x\in\beta(N)'\cap\mathcal{L}(H)$ and $y\in\alpha(N)'\cap
\mathcal{L}(K)$, we get an operator $x\surl{\
_{\beta}\otimes_{\alpha}}_{\ \psi}y$ on $H\surl{\
_{\beta}\otimes_{\alpha}}_{\ \psi} K$. Canonical isomorphism of
change of weights sends $x\surl{\ _{\beta}\otimes_{\alpha}}_{\ \psi}
y$ on $x\surl{\ _{\beta}\otimes_{\alpha}}_{\ \ \psi '}y$ so that we
write $x\surl{\ _{\beta}\otimes_{\alpha}}_{\ N} y$ without any
reference to the weight.

Let $P$ be a von Neumann algebra and $\epsilon$ a normal and
non-degenerated anti-representation of $P$ on $K$ such that
$\epsilon(P)'\subseteq \alpha(N)$. $K$ is equipped with a
$N-P$-bimodule structure denoted by $\ _{\alpha}K_{\epsilon}$. For
all $y\in P$, $1_H \surl{\ _{\beta}\otimes_{\alpha}}_{\ \psi}
\epsilon(y)$ is an operator on $H\surl{\
_{\beta}\otimes_{\alpha}}_{\ \psi}K$ so that we define a
representation of $P$ on $H\surl{\ _{\beta}\otimes_{\alpha}}_{\
\psi} K$ still denoted by $\epsilon$. If $H$ is a $Q-N$-bimodule,
then $H\surl{\ _{\beta}\otimes_{\alpha}}_{\ \psi} K$ becomes a
$Q-P$-bimodule (Connes' fusion of bimodules). If $\nu$ is a n.s.f
weight on $P$ and $\ _{\zeta}L$ a left $P$-module. It is possible to
define two Hilbert spaces $(H\surl{\ _{\beta}\otimes_{\alpha}}_{\
\psi} K) \surl{\ _{\epsilon}\otimes_{\zeta}}_{\ \nu} L$ and
$H\surl{\ _{\beta}\otimes_{\alpha}}_{\ \psi} (K \surl{\
_{\epsilon}\otimes_{\zeta}}_{\ \nu} L)$. These two
$\beta(N)'-\zeta(P)'^o$-bimodules are isomorphic. (The proof of
\cite{V1}, lemme 2.1.3, in the case of commutative $N=P$ is still
valid). We speak about associativity of relative tensor product and
we write $H\surl{\ _{\beta}\otimes_{\alpha}}_{\ \psi} K \surl{\
_{\epsilon}\otimes_{\zeta}}_{\ \nu} L$ without parenthesis.

We identify $H_{\psi} \surl{\ _{\beta} \otimes_{\alpha}}_{\ \psi} K$
and $K$ as left $N$-modules by $\Lambda_{\psi}(y) \surl{\ _{\beta}
\otimes_{\alpha}}_{\ \psi}\eta\mapsto\alpha(y)\eta$ for all $y\in
{\mathcal N}_{\psi}$. By \cite{EN}, 3.10, we have:
$$\lambda^{\beta,\alpha}_{\xi}=R^{\beta,\psi^o}(\xi) \surl{\
_{\beta} \otimes_{\alpha}}_{\ \psi} 1_K$$

We recall proposition 2.3 of \cite{E2}:

\begin{prop}\label{decomposition}
Let $(\xi_i)_{i\in I}$ be a $(N^o,\psi^o)$-basis of $H_{\beta}$.
Then:
\begin{enumerate}[i)]
\item for all $\xi\in D(H_{\beta},\psi^o)$ and $\eta\in K$, we have:
$$\xi\surl{\
_{\beta}\otimes_{\alpha}}_{\ \psi} \eta=\sum_{i\in I}\xi_i\surl{\
_{\beta}\otimes_{\alpha}}_{\
\psi}\alpha(<\xi,\xi_i>_{\beta,\psi^o})\eta$$
\item we have the following decomposition:
$$ H\surl{\
_{\beta}\otimes_{\alpha}}_{\ \psi}K=\bigoplus_{i\in I}(\xi_i \surl{\
_{\beta}\otimes_{\alpha}}_{\ \psi}
\alpha(<\xi_i,\xi_i>_{\beta,\psi^o})K)$$
\end{enumerate}
\end{prop}

We here add a proposition we will use several times.

\begin{prop}\label{tenscomp}
Let $\gamma$ a *-automorphism from $N$ such that
$\psi\circ\nu=\psi$. Then :
$$H\surl{\
_{\beta\circ\gamma}\otimes_{\alpha\circ\gamma}}_{\ \psi}K=H\surl{\
_{\beta}\otimes_{\alpha}}_{\ \psi}K
$$
\end{prop}

\begin{proof}
Because of invariance of $\psi$ with respect to $\gamma$, we have a
unitary $I$ from $H_{\psi}$ such that
$I\Lambda_{\psi}(y)=\Lambda_{\psi}(\gamma(y))$ for all
$y\in\mathcal{N}_{\psi}$. Moreover $IJ_{\psi}=J_{\psi}I$ and
$I^*nI=\gamma^{-1}(n)$ for all $n\in N$. For all $\xi \in
D(H_{\beta},\psi^o)$ and $y\in {\mathcal N}_{\psi}$, we compute:
$$
\begin{aligned}
\beta\circ\gamma(y^*))\xi &
=\beta(\gamma(y)^*)\xi=R^{\beta,\psi^o}(\xi)J_{\psi}\Lambda_{\psi}(\gamma(y))\\
&=R^{\beta,\psi^o}(\xi)J_{\psi}I\Lambda_{\psi}(y)=R^{\beta,\nu^o}(\xi)IJ_{\psi}\Lambda_{\psi}(y)
\end{aligned}$$
that's why we get:
$$D(H_{\beta\circ\gamma},\psi^o)=D(H_{\beta},\psi^o) \text{ and }
\forall \xi \in D(H_{\beta},\psi^o),\
R^{\beta\circ\gamma,\psi^o}(\xi)=R^{\beta,\psi^o}(\xi)I$$ To
conclude, we show that scalar products on $D(H_{\beta,\psi^o}) \odot
K$ used to define $H \surl{\ _{\beta} \otimes_{\alpha}}_{\ \psi}K$
and $H\surl{\ _{\beta\circ\gamma} \otimes_{\alpha\circ\gamma}}_{\
\psi} K$ are equal. If $\xi,\xi' \in D(H_{\beta},\nu^o)$ and
$\eta,\eta' \in K$, we have:
$$
\begin{aligned}
(\xi \surl{\ _{\beta\circ\gamma} \otimes_{\alpha\circ\gamma}}_{\
\psi} \eta|\xi' \surl{\ _{\beta\circ\gamma}
\otimes_{\alpha\circ\gamma}}_{\ \psi} \eta')
&=(\alpha(\gamma(<\xi,\xi'>_{\beta\circ\gamma,\psi^o}))\eta|\eta')\\
&=(\alpha(\gamma(I^*<\xi,\xi'>_{\beta,\psi^o}I))\eta|\eta')\\
&=(\alpha(<\xi,\xi'>_{\beta,\psi^o})\eta|\eta')=(\xi\surl{\ _{\beta}
\otimes_{\alpha}}_{\ \psi}\xi'|\eta \surl{\ _{\beta}
\otimes_{\alpha}}_{\ \psi} \eta')
\end{aligned}$$
\end{proof}

To end the paragraph, we detail finite dimension case. We assume
that $N$, $H$ and $K$ are of finite dimensions. $H\surl{\
_{\beta}\otimes_{\alpha}}_{\ \psi} K$ can be identified with a
subspace of $H \otimes K$. We denote by $\text{Tr}$ the normalized
canonical trace on $K$ and $\tau=\text{Tr}\circ\alpha$. There exist
a projection $e_{\beta,\alpha}\in\beta(N)\otimes\alpha(N)$ and
$n_o\in Z(N)^+$ such that
$(id\otimes\text{Tr})(e_{\beta,\alpha})=\beta(n_o)$. Let $d$ be the
Radon-Nikodym derivative of $\psi$ w.r.t $\tau$. By \cite{EV}, 2.4,
and proposition 2.7 of \cite{S2}, for all $\xi,\eta\in H$:
$$I_{\beta,\alpha}^{\psi}:\xi\surl{\ _{\beta}\otimes_{\alpha}}_{\ \psi}\eta
\mapsto\xi\surl{\ _{\beta}\otimes_{\alpha}}_{\
\tau}\alpha(d)^{1/2}\eta \mapsto
e_{\beta,\alpha}(\beta(n_o)^{-1/2}\xi\otimes\alpha(d)^{1/2}\eta)$$
defines an isometric isomorphism of
$\beta(N)'-\alpha(N)'^o$-bimodules from $H\surl{\
_{\beta}\otimes_{\alpha}}_{\ \psi}K$ onto a subspace of $H\otimes
K$, the final support of which is $e_{\beta,\alpha}$.\label{iden}

\subsection{Fiber product \cite{V1}, \cite{EV}}
We use previous notations. Let $M_1$ (resp. $M_2$) be a von Neumann
algebra on $H$ (resp. $K$) such that $\beta(N)\subseteq M_1$ (resp.
$\alpha(N)\subseteq M_2$). We denote by $M_1'\surl{\
_{\beta}\otimes_{\alpha}}_{\ N} M_2'$ the von Neumann algebra
generated by $x\surl{\ _{\beta}\otimes_{\alpha}}_{\ N}y$ with $x\in
M_1'$ and $y\in M_2'$.

\begin{defi}
The commutant of $M_1'\surl{\ _{\beta}\otimes_{\alpha}}_{\ N} M_2'$
in $\mathcal{L}(H\surl{\ _{\beta}\otimes_{\alpha}}_{\ \psi}K)$ is
denoted by $M_1\surl{\ _{\beta}\star_{\alpha}}_{\ N} M_2$ and is
called \textbf{fiber product}.
\end{defi}

If $P_1$ and $ P_2$ are von Neumann algebras like $M_1$ and $M_2$,
we have:
$$
\begin{aligned}
i) &\quad (M_1\surl{\ _{\beta}\star_{\alpha}}_{\ N} M_2)\cap
(P_1\surl{\ _{\beta}\star_{\alpha}}_{\ N} P_2)=(M_1\cap
P_1)\surl{\ _{\beta}\star_{\alpha}}_{\ N} (M_2\cap P_2)\\
ii)&\quad \varsigma_N(M_1\surl{\ _{\beta}\star_{\alpha}}_{\ N}
M_2)=M_2\surl{\ _{\alpha}\star_{\beta}}_{\ \ N^o} M_1 \\
iii) &\quad (M_1\cap\beta(N)')\surl{\ _{\beta}\otimes_{\alpha}}_{\
N} (M_2\cap\alpha(N)')\subseteq M_1\surl{\
_{\beta}\star_{\alpha}}_{\ N} M_2 \\
iv) &\quad M_1\surl{\ _{\beta}\star_{\alpha}}_{\ N}
\alpha(N)=(M_1\cap\beta(N)')\surl{\ _{\beta}\otimes_{\alpha}}_{\ N}
1
\end{aligned}$$
More generally, if $\beta$ (resp. $\alpha$) is a normal,
non-degenerated *-anti-homomorphism (resp. homomorphism) from $N$ to
a von Neumann algebra $M_1$ (resp. $M_2$), it is possible to define
a von Neumann algebra $M_1\surl{\ _{\beta}\star_{\alpha}}_{\ N} M_2$
without any reference to a specific Hilbert space. If $P_1$, $P_2$,
$\alpha'$ and $\beta'$ are like $M_1$, $M_2$, $\alpha$ and $\beta$
and if $\Phi$ (resp. $\Psi$) is a normal *-homomorphism from $M_1$
(resp. $M_2$) to $P_1$ (resp. $P_2$) such that
$\Phi\circ\beta=\beta'$ (resp. $\Psi\circ\alpha=\alpha'$), then we
define a normal *-homomorphism by \cite{S1}, 1.2.4: $$\Phi\surl{\
_{\beta}\star_{\alpha}}_{\ N}\Psi: M_1\surl{\
_{\beta}\star_{\alpha}}_{\ N}M_2\rightarrow P_1\surl{\
_{\beta'}\star_{\alpha'}}_{\ N}P_2$$

Assume $\ _{\alpha}K_{\epsilon}$ is a $N-P^o$-bimodule and $\
_{\zeta}L$ a left $P$-module. If $\alpha(N)\subseteq M_2$,
$\epsilon(P)\subseteq M_2$ and if $\zeta(P)\subseteq M_3$ where
$M_3$ is a von Neumann algebra on $L$, then we can construct
$M_1\surl{\ _{\beta}\star_{\alpha}}_{\ N}(M_2\surl{\
_{\epsilon}\star_{\zeta}}_{\ N}M_3)$ and $(M_1\surl{\
_{\beta}\star_{\alpha}}_{\ N}M_2)\surl{\
_{\epsilon}\star_{\zeta}}_{\ N}M_3$. Associativity of relative
tensor product induces an isomorphism between these fiber products
and we write $M_1\surl{\ _{\beta}\star_{\alpha}}_{\ N}M_2\surl{\
_{\epsilon}\star_{\zeta}}_{\ N}M_3$ without parenthesis.

Finally, if $M_1$ and $M_2$ are of finite dimensions, then we have:
$$M_1'\surl{\ _{\beta}\otimes_{\alpha}}_{\
N}M_2'=(I_{\beta,\alpha}^{\psi})^*(M_1'\otimes
M_2')I_{\beta,\alpha}^{\psi}\text{ and } M_1\surl{\
_{\beta}\star_{\alpha}}_{\
N}M_2=(I_{\beta,\alpha}^{\psi})^*(M_1\otimes
M_2)I_{\beta,\alpha}^{\psi}$$ Therefore the fiber product can be
identified with a reduction of $M_1\otimes M_2$ by
$e_{\beta,\alpha}$ by \cite{EV}, 2.4.

\subsection{Slice map \cite{E1}}
\subsubsection{For normal forms}
Let $A\in M_1\surl{\ _{\beta}\star_{\alpha}}_{\ N}M_2$ and
$\xi_1,\xi_2\in D(H_{\beta},\psi^o)$. We define an element of $M_2$
by:
$$(\omega_{\xi_1,\xi_2}\surl{\ _{\beta}\star_{\alpha}}_{\ \psi}
id)(A)=(\lambda^{\beta,\alpha}_{\xi_2})^*A\lambda^{\beta,\alpha}_{\xi_1}$$
so that we have $((\omega_{\xi_1,\xi_2}\surl{\
_{\beta}\star_{\alpha}}_{\ \psi}
id)(A)\eta_1|\eta_2)=(A(\xi_1\surl{\ _{\beta}\otimes_{\alpha}}_{\
\psi} \eta_1)|\xi_2\surl{\ _{\beta}\otimes_{\alpha}}_{\
\psi}\eta_2)$ for all $\eta_1,\eta_2\in K$. Also, we define an
operator of $M_1$ by:
$$(id\surl{\ _{\beta}\star_{\alpha}}_{\
\psi}\omega_{\eta_1,\eta_2})(A)=
(\rho^{\beta,\alpha}_{\eta_2})^*A\rho^{\beta,\alpha}_{\eta_1}$$ for
all $\eta_1,\eta_2\in D(_{\alpha}K,\psi)$. We have a Fubini's
formula:
$$\omega_{\eta_1,\eta_2}((\omega_{\xi_1,\xi_2}\surl{\ _{\beta}\star_{\alpha}}_{\ \psi}
id)(A))=\omega_{\xi_1,\xi_2}((id\surl{\ _{\beta}\star_{\alpha}}_{\
\psi}\omega_{\eta_1,\eta_2})(A))$$ for all $\xi_1,\xi_2\in
D(H_{\beta},\psi^o)$ and $\eta_1,\eta_2\in D(_{\alpha}K,\psi)$.

Equivalently, by (\cite{E1}, proposition 3.3), for all $\omega_1\in
M_{1*}^+$ and $k_1\in\mathbb{R}^+$ such that $\omega_1\circ\beta\leq
k_1\psi$ and for all $\omega_2\in M_{2*}^+$ and $k_2\in\mathbb{R}^+$
such $\omega_2\circ\alpha\leq k_2\psi$, we have:
$$\omega_2((\omega_1\surl{\ _{\beta}\star_{\alpha}}_{\ \psi}
id)(A))=\omega_1((id\surl{\ _{\beta}\star_{\alpha}}_{\
\psi}\omega_2)(A))$$

\subsubsection{For conditional expectations}
If $P_2$ is a von Neumann algebra such that $\alpha(N)\subseteq P_2
\subseteq M_2$ and if $E$ is a normal, faithful conditional
expectation from $M_2$ onto $P_2$, we can define a normal, faithful
conditional expectation $(id\surl{\ _{\beta}\star_{\alpha}}_{\ N}E)$
from $M_1\surl{\ _{\beta}\star_{\alpha}}_{\ N} M_2$ onto $M_1\surl{\
_{\beta}\star_{\alpha}}_{\ N} P_2$ such that:
$$(\omega\surl{\ _{\beta}\star_{\alpha}}_{\
\psi}id) (id\surl{\ _{\beta}\star_{\alpha}}_{\
N}E)(A)=E((\omega\surl{\ _{\beta}\star_{\alpha}}_{\ \psi}id)(A))$$
for all $A\in M_1\surl{\ _{\beta}\star_{\alpha}}_{\ N} M_2$,
$\omega\in M_{1*}^+$ and $k_1\in\mathbb{R}^+$ such that
$\omega\circ\beta\leq k_1\psi$.

\subsubsection{For weights}
If $\phi_1$ is n.s.f weight on $M_1$ and if $A$ is a positive
element of $M_1\surl{\ _{\beta}\star_{\alpha}}_{\ N} M_2$, we can
define an element of the extended positive part of $M_2$, denoted by
$(\phi_1\surl{\ _{\beta}\star_{\alpha}}_{\ \psi}id)(A)$, such that,
for all $\eta\in D(_{\alpha}L^2(M_2),\psi)$, we have:
$$||((\phi_1\surl{\ _{\beta}\star_{\alpha}}_{\ \psi}
id)(A))^{1/2}\eta||^2 =\phi_1((id\surl{\ _{\beta}\star_{\alpha}}_{\
\psi}\omega_{\eta})(A))$$ Moreover, if $\phi_2$ is a n.s.f weight on
$M_2$, we have:
$$\phi_2((\phi_1\surl{\ _{\beta}\star_{\alpha}}_{\ \psi} id)(A))=
\phi_1((id\surl{\ _{\beta}\star_{\alpha}}_{\ \psi}\phi_2)(A))$$

Let $(\omega_i)_{i\in I}$ be an increasing net of normal forms such
that $\phi_1=\sup_{i\in I}\omega_i$. Then we have $(\phi_1\surl{\
_{\beta}\star_{\alpha}}_{\ \psi} id)(A)=\sup_i(\omega_i\surl{\
_{\beta}\star_{\alpha}}_{\ \psi} id)(A)$.

\subsubsection{For operator-valued weights}
Let $P_1$ be a von Neumann algebra such that $\beta(N)\subseteq
P_1\subseteq M_1$ and $\Phi_i$ $(i=1,2)$ be operator-valued n.s.f
weights from $M_i$ to $P_i$. By \cite{E1}, for all positive operator
$A\in M_1\surl{\ _{\beta}\star_{\alpha}}_{\ N} M_2$, there exists an
element $(\Phi_1\surl{\ _{\beta}\star_{\alpha}}_{\ N} id)(A)$
belonging to $P_1\surl{\ _{\beta}\star_{\alpha}}_{\ N} M_2$ such
that, for all $\xi\in L^2(P_1)$ and $\eta \in D(_{\alpha}K,\psi)$,
we have:
$$||((\Phi_1\surl{\
_{\beta}\star_{\alpha}}_{\ N} id)(A))^{1/2}(\xi\surl{\
_{\beta}\otimes_{\alpha}}_{\ \psi}\eta)||^2=||[\Phi_1((id\surl{\
_{\beta}\star_{\alpha}}_{\
\psi}\omega_{\eta,\eta})(A))]^{1/2}\xi||^2$$ If $\phi_1$ is a n.s.f
weight on $P_1$, we have:
$$(\phi_1\circ\Phi_1\surl{\
_{\beta}\star_{\alpha}}_{\ N} id)(A)=(\phi_1\surl{\
_{\beta}\star_{\alpha}}_{\ \psi} id)(\Phi_1\surl{\
_{\beta}\star_{\alpha}}_{\ N} id)(A)$$ Also, we define an element
$(id\surl{\ _{\beta}\star_{\alpha}}_{\ N} \Phi_2)(A)$ of the
extended positive part of $M_1\surl{\ _{\beta}\star_{\alpha}}_{\ N}
P_2$ and we have:
$$(id\surl{\
_{\beta}\star_{\alpha}}_{\ N} \Phi_2)((\Phi_1\surl{\
_{\beta}\star_{\alpha}}_{\ N} id)(A))=(\Phi_1\surl{\
_{\beta}\star_{\alpha}}_{\ N} id)((id\surl{\
_{\beta}\star_{\alpha}}_{\ N} \Phi_2)(A))$$

\begin{rema}
We have seen that we can identify $M_1\surl{\
_{\beta}\star_{\alpha}}_{\ N} \alpha(N)$ with $M_1\cap\beta(N)'$.
Then, it is easy to check that the slice map $id\surl{\
_{\beta}\star_{\alpha}}_{\ \psi}\psi\circ\alpha^{-1}$ (if $\alpha$
is injective) is just the injection of $M_1\surl{\
_{\beta}\star_{\alpha}}_{\ N} \alpha(N)$ into $M_1$. Also we see on
that example that, if $\phi_1$ is a n.s.f weight on $M_1$, then
$\phi_1\surl{\ _{\beta}\star_{\alpha}}_{\ N}id$ (which is equal to
$\phi_{1|M_1\cap\beta(N)'}$) needs not to be semi-finite.
\end{rema}

\section{Fundamental pseudo-multiplicative unitary} In this section,
we construct a fundamental pseudo-multiplicative unitary from a Hopf
bimodule with a left invariant operator-valued weight and a right
invariant operator-valued weight. Let $N$ and $M$ be von Neumann
algebras, $\alpha$ (resp. $\beta$) be a faithful, non-degenerate,
normal (resp. anti-) representation from $N$ to $M$. We suppose that
$\alpha(N) \subseteq \beta(N)'$.

\subsection{Definitions}

\begin{defi}
A quintuplet $(N,M,\alpha,\beta,\Gamma)$ is said to be a {\bf Hopf
bimodule} of basis $N$ if $\Gamma$ is a normal *-homomorphism from
$M$ into $M \surl{\nonumber_{\beta}\star_{\alpha}}_{N} M$ such that,
for all $n,m\in N$, we have:

\begin{center}
\begin{minipage}{10cm}
\begin{enumerate}[i)]
\item $\Gamma(\alpha(n)\beta(m))=\alpha(n)\surl{\nonumber_{\beta} \otimes_{\alpha}}_{\ N}
  \beta(m)$
\item $\Gamma$ is co-associative: $(\Gamma \surl{\nonumber _{\beta}
\star_{\alpha}}_{\ N} id)\circ\Gamma=(id \surl{\ _{\beta}
\star_{\alpha}}_{\ N}\Gamma)\circ\Gamma$
\end{enumerate}
\end{minipage}
\end{center}

\end{defi}

One should notice that property i) is necessary in order to write
down the formula given in ii).
$(N^o,M,\beta,\alpha,\varsigma_N\circ\Gamma)$ is a Hopf bimodule
called opposite Hopf bimodule. If $N$ is commutative, $\alpha=\beta$
and $\Gamma=\varsigma_N\circ\Gamma$, then
$(N,M,\alpha,\alpha,\Gamma)$ is equal to its opposite: we shall
speak about a symmetric Hopf bimodule.

\begin{defi} Let
$(N,M,\alpha,\beta,\Gamma)$ be a Hopf bimodule. A normal,
semi-finite, faithful operator-valued weight from $M$ to $\alpha(N)$
is said to be {\bf left invariant} if:
$$(id\surl{\ _{\beta}\star_{\alpha}}_{\ N} T_L)\Gamma(x)=T_L(x)\surl{\ _{\beta}
\otimes_{\alpha}}_{\ N} 1\quad\quad\text{for all }x\in
\mathcal{M}_{T_L}^+$$ In the same way, a normal, semi-finite,
faithful operator-valued weight from $M$ to $\beta(N)$ is said to be
{\bf right invariant} if:
$$(T_R \surl{\
_{\beta}
  \star_{\alpha}}_{\ N} id)\Gamma(x)=1 \surl{\ _{\beta}
  \otimes_{\alpha}}_{\ N}T_R(x)\quad\quad\text{for all }x\in \mathcal{M}_{T_R}^+$$
\end{defi}

We give several examples in the last section. In this section,
$(N,M,\alpha,\beta,\Gamma)$ is a Hopf bimodule with a left
operator-valued weight $T_L$ and a right operator-valued weight
$T_R$.

\begin{defi}
A *-anti-automorphism $R$ of $M$ is said to be a {\bf co-involution}
if $R\circ\alpha=\beta$, $R^2=id$ and $\varsigma_{N^o}\circ(R\surl{\
_{\beta} \star_{\alpha}}_{\ N} R)\circ\Gamma=\Gamma \circ R$.
\end{defi}

\begin{rema}
With the previous notations, let us notice that $R\circ T_L\circ R$
is a right invariant operator-valued weight from $M$ to $\beta(N)$.
Also, let us say that $R$ is an anti-isomorphism of Hopf bimodule
from the bimodule and its symmetric.
\end{rema}

Let $\mu$ be a normal, semi-finite, faithful weight of $N$. We put:
$$\Phi =\mu\circ\alpha^{-1}\circ T_L \text{ and } \Psi
=\mu\circ\beta^{-1}\circ T_R$$ so that, for all $x\in M^+$, we have:
$$(id \surl{\ _{\beta}
  \star_{\alpha}}_{\ \mu} \Phi)\Gamma(x)=T_L(x) \text{ and }
(\Psi \surl{\ _{\beta} \star_{\alpha}}_{\ \mu} id)\Gamma(x)=T_R(x)$$

If $H$ denote a Hilbert space on which $M$ acts, then $N$ acts on
$H$, also, by way of $\alpha$ and $\beta$. We shall denote again
$\alpha$ (resp. $\beta$) for (resp. anti-) the representation of $N$
on $H$.

\subsection{Construction of the fundamental isometry}

\begin{defi}
Let define $\hat{\beta}$ and $\hat{\alpha}$ by:

$$
\begin{aligned}
\hat{\beta}: N &\rightarrow \mathcal{L}(H_{\Phi}) &\quad\text{ and }\quad\quad\hat{\alpha}: N &\rightarrow \mathcal{L}(H_{\Psi}) \\
x &\mapsto J_{\Phi}\alpha(x^*)J_{\Phi}& x &\mapsto J_{\Psi}\beta
(x^*)J_{\Psi}
\end{aligned}$$

Then $\hat{\beta}$ (resp. $\hat{\alpha}$) is a normal,
non-degenerate and faithful anti-representation (resp.
representation) from $N$ to $\mathcal{L}(H_{\Phi})$ (resp.
$\mathcal{L}(H_{\Psi})$). \label{defutil}
\end{defi}

\begin{prop} \label{prem}
We have $\Lambda_{\Phi}({\mathcal N}_{T_L} \cap {\mathcal N}_{\Phi})
\subseteq D((H_{\Phi})_{\hat{\beta}},\mu^o)$ and for all $a\in
{\mathcal N}_{T_L} \cap {\mathcal N}_{\Phi}$, we have:
$$R^{\hat{\beta},\mu^o}(\Lambda_{\Phi}(a))=\Lambda_{T_L}(a)$$ Also,
we have $\Lambda_{\Psi}({\mathcal N}_{T_R} \cap {\mathcal N}_{\Psi})
\subseteq D(_{\hat{\alpha}}(H_{\Psi}),\mu)$ and for all $b\in
{\mathcal N}_{T_R} \cap {\mathcal N}_{\Psi}$, then:
$$R^{\hat{\alpha},\mu}(\Lambda_{\Psi}(b))=\Lambda_{T_R}(b)$$
\end{prop}

\begin{rema}
We identify $H_{\mu}$ with $H_{\mu\circ\alpha^{-1}}$ and $H_{\mu}$
with $H_{\mu\circ\beta^{-1}}$.
\end{rema}

\begin{proof}
Let $y\in {\mathcal N}_{\mu}$ analytic w.r.t $\mu$. We have:
$$
\begin{aligned}
\hat{\beta}(y^*)\Lambda_{\Phi}(a)&=\Lambda_{\Phi}(a\sigma_{-i/2}^{\Phi}(\alpha(y^*)))
=\Lambda_{\Phi}(a\sigma_{-i/2}^{\mu\circ\alpha^{-1}}(\alpha(y^*)))\\
&=\Lambda_{\Phi}(a\alpha(\sigma_{-i/2}^{\mu}(y^*)))
=\Lambda_{T_L}(a)\Lambda_{\mu}(\sigma_{-i/2}^{\mu}(y^*))
=\Lambda_{T_L}(a)J_{\mu}\Lambda_{\mu}(y)
\end{aligned}$$
Thanks to lemma \ref{analyse}, we get
$\hat{\beta}(y^*)\Lambda_{\Phi}(a)
=\Lambda_{T_L}(a)J_{\mu}\Lambda_{\mu}(y)$, for all $y\in {\mathcal
N}_{\mu}$, which gives the first part of the proposition. The end of
the proof is very similar.
\end{proof}

\begin{prop}\label{evi}
We have $J_{\Phi}D((H_{\Phi})_{\hat{\beta}},
\mu^o)=D(_{\alpha}(H_{\Phi}),\mu)$ and for all $\eta\in
D((H_{\Phi})_{\hat{\beta}}, \mu^o)$, we have:
$$R^{\alpha,\mu}(J_{\Phi}\eta)=J_{\Phi}R^{\hat{\beta},\mu^o}(\eta)J_{\mu}$$
Also, we have
$J_{\Psi}D(_{\hat{\alpha}}(H_{\Psi}),\mu)=D((H_{\Phi})_{\beta},\mu^o)$
and for all $\xi\in D((H_{\Phi})_{\beta},\mu^o)$, we have:
$$R^{\beta,\mu^o}(J_{\Psi}\xi)=J_{\Psi}R^{\hat{\alpha},\mu}(\xi)J_{\mu}$$
\end{prop}

\begin{proof}
Straightforward.
\end{proof}

\begin{coro}
We have $\Lambda_{\Phi}({\mathcal T}_{\Phi,T_L}) \subseteq
D((H_{\Phi})_{\hat{\beta}}, \mu^o) \cap D(_{\alpha}(H_{\Phi}),\mu)$
and $\Lambda_{\Psi}({\mathcal T}_{\Psi,T_R}) \subseteq
D(_{\hat{\alpha}}(H_{\Psi}),\mu) \cap D((H_{\Psi})_{\beta},\mu^o)$.
\end{coro}

\begin{proof}
This is a corollary of the two previous propositions.
\end{proof}

\begin{rema}
The invariance of operator-valued weights does not play a part in
the previous propositions.
\end{rema}

\begin{prop}\label{semi}
We have $(\omega_{v,\xi}\surl{\ _{\beta} \star_{\alpha}}_{\ \mu}
id)(\Gamma(a)) \in {\mathcal N}_{T_L}\cap {\mathcal N}_{\Phi}$ for
all elements $a\in {\mathcal N}_{T_L}\cap {\mathcal N}_{\Phi}$ and
$v,\xi\in D(H_{\beta},\mu^o)$.
\end{prop}

\begin{proof} By definition of the slice maps, we have:
$$\begin{aligned}
(\omega_{v,\xi}\surl{\ _{\beta} \star_{\alpha}}_{\ \mu}
  id)(\Gamma(a))^*(\omega_{v,\xi}\surl{\ _{\beta}
  \star_{\alpha}}_{\ \mu} id)(\Gamma(a))
&=(\lambda^{\beta,
  \alpha}_v)^*\Gamma(a^*)\lambda^{\beta,
  \alpha}_{\xi}(\lambda^{\beta,
  \alpha}_{\xi})^*\Gamma(a)\lambda^{\beta, \alpha}_v\\
&\leq \|\lambda^{\beta, \alpha}_{\xi}\|^2(\omega_{v,v}\surl{\
_{\beta}
  \star_{\alpha}}_{\ \mu} id)(\Gamma(a^*a))\\
&\leq \|R^{\beta, \mu^o}(\xi)\|^2(\omega_{v,v}\surl{\ _{\beta}
  \star_{\alpha}}_{\ \mu} id)(\Gamma(a^*a))
  \end{aligned}$$
Then, on one hand, we get, thanks to left invariance of $T_L$:
$$\begin{aligned}
&\ \quad {T_L}((\omega_{v,\xi}\surl{\ _{\beta} \star_{\alpha}}_{\
\mu}
  id)(\Gamma(a))^*(\omega_{v,\xi}\surl{\ _{\beta}
  \star_{\alpha}}_{\ \mu} id)(\Gamma(a)))\\
&\leq \|R^{\beta, \mu^o}(\xi)\|^2{T_L}((\omega_{v,v}\surl{\
_{\beta}\star_{\alpha}}_{\ \mu} id)(\Gamma(a^*a)))\\
&=\|R^{\beta, \mu^o}(\xi)\|^2 (\omega_{v,v}\surl{\ _{\beta}
  \star_{\alpha}}_{\ \mu} id)(id \surl{\ _{\beta}
  \star_{\alpha}}_{\ \mu} {T_L})(\Gamma(a^*a))\\
&\leq \|R^{\beta,\mu^o}(\xi)\|^2(\lambda^{\beta,\alpha}_{v})^*
({T_L}(a^*a)\surl{\ _{\beta}\otimes_{\alpha}}_{\ \mu} 1)
\lambda^{\beta,\alpha}_{v}\\
&\leq\|R^{\beta,\mu^o}(\xi)\|^2||{T_L}(a^*a)||
\|\alpha(<v,v>_{\beta,\mu^o})\|1\\
&\leq \|R^{\beta,\mu^o}(\xi)\|^2\|{T_L}(a^*a)\|\|R^{\beta,
\mu^o}(v)\|^21 \end{aligned}$$ So, we get that
$(\omega_{v,\xi}\surl{\ _{\beta}
  \star_{\alpha}}_{\ \mu} id)(\Gamma(a)) \in {\mathcal N}_{T_L}$.
On the other hand, thanks to left invariance of $T_L$, we know that:
$$\Phi(((\omega_{v,\xi}\surl{\ _{\beta}
  \star_{\alpha}}_{\ \mu} id)(\Gamma(a)))^*(\omega_{v,\xi}\surl{\
  _{\beta} \star_{\alpha}}_{\ \mu} id)(\Gamma(a)))$$ is less or
equal to:
$$
\begin{aligned}
&\ \quad\|R^{\beta,\mu^o}(\xi)\|^2 \Phi((\omega_{v,v}\surl{\
_{\beta}\star_{\alpha}}_{\ \mu} id)(\Gamma(a^*a)))\\
&=\|R^{\beta, \mu^o}(\xi)\|^2\omega_{v,v}((id \surl{\ _{\beta}
  \star_{\alpha}}_{\ \mu}\Phi)(\Gamma(a^*a)))\\
&=\|R^{\beta,
\mu^o}(\xi)\|^2({T_L}(a^*a)v|v)\leq\|R^{\beta,\mu^o}(\xi)\|^2\|{T_L}(a^*a)\|\|v\|^2
< + \infty
\end{aligned}$$  So, we get that
$(\omega_{v,\xi}\surl{\ _{\beta} \star_{\alpha}}_{\ \mu}
id)(\Gamma(a)) \in {\mathcal N}_{\Phi}$.
\end{proof}

\begin{prop}\label{rapide}
For all $v,w\in H$ and $a,b\in {\mathcal N}_{\Phi}\cap {\mathcal
N}_{T_L}$, we have:
$$(v\surl{\ _{\alpha}
\otimes_{\hat{\beta}}}_{\ \ \mu^o}\Lambda_{\Phi}(a)|w\surl{\
_{\alpha} \otimes_{\hat{\beta}}}_{\ \
\mu^o}\Lambda_{\Phi}(b))=(T_L(b^*a)v|w)$$

For all $v,w\in H$ and $c,d\in {\mathcal N}_{\Psi}\cap {\mathcal
N}_{T_R}$, we have:
$$(\Lambda_{\Psi}(c)\surl{\ _{\hat{\alpha}} \otimes_{\beta}}_{\ \
\mu^o}v|\Lambda_{\Psi}(d)\surl{\ _{\hat{\alpha}} \otimes_{\beta}}_{\
\ \mu^o}w)=(T_R(d^*c)v|w)$$
\end{prop}

\begin{proof}
Using \ref{prem} and \ref{preintro}, we get that:
$$\begin{aligned}
(v \surl{\ _{\alpha} \otimes_{\hat{\beta}}}_{\ \ \mu^o}
\Lambda_{\Phi}(a)|w\surl{\ _{\alpha} \otimes_{\hat{\beta}}}_{\ \
\mu^o}\Lambda_{\Phi}(b))
&=(\alpha(<\Lambda_{\Phi}(a),\Lambda_{\Phi}(b)>_{\hat{\beta},\mu^o})v|w)\\
&=(\alpha(\Lambda_{T_L}(b)^*\Lambda_{T_L}(a))v|w)\\
&=(\alpha(\pi_{\mu}(\alpha^{-1}({T_L}(b^*a))))v|w)
\end{aligned}$$
which gives the result after the identification of $\pi_{\mu}(N)$
with $N$. The second point is very similar.
\end{proof}

\begin{lemm}\label{inde}
Let $a\in {\mathcal N}_{\Phi}\cap {\mathcal N}_{T_L}$ and $v\in
D(H_{\beta},\mu^o)$. The following sum: $$\sum_{i\in I} \xi_{i}
\surl{\ _{\beta} \otimes_{\alpha}}_{\ \mu}
  \Lambda_{\Phi} ((\omega_{v,\xi_i} \surl{\
      _{\beta} \star_{\alpha}}_{\ \mu} id)(\Gamma(a)))$$ converges
in $H \surl{\ _{\beta} \otimes_{\alpha}}_{\ \mu} H_{\Phi}$ for all
$(N^o,\mu^o)$-basis $(\xi_i)_{i\in I}$ of $H_{\beta}$ and it does
not depend on the $(N^o,\mu^o)$-basis of $H_{\beta}$.
\end{lemm}

\begin{proof}
By \ref{semi}, we have $(\omega_{v,\xi_i} \surl{\ _{\beta}
\star_{\alpha}}_{\ \mu} id)(\Gamma(a))\in {\mathcal N}_{\Phi}\cap
{\mathcal N}_{T_L}$ for all $i\in I$, and the vectors $\xi_{i}
\surl{\ _{\beta} \otimes_{\alpha}}_{\ \mu} \Lambda_{\Phi}
((\omega_{v,\xi_i} \surl{\ _{\beta} \star_{\alpha}}_{\ \mu}
id)(\Gamma(a)))$ are two-by-two orthogonal. Normality and left
invariance of $\Phi$ imply:
$$
\begin{aligned}
&\quad\sum_{i
    \in I} ||\xi_{i} \surl{\ _{\beta} \otimes_{\alpha}}_{\ \mu}
  \Lambda_{\Phi} ((\omega_{v,\xi_i} \surl{\
      _{\beta} \star_{\alpha}}_{\ \mu} id)(\Gamma(a)))||^2\\
&=\sum_{i\in I} (\alpha(<\xi_i,\xi_i>_{\beta,\mu^o}) \Lambda_{\Phi}
((\omega_{v,\xi_i} \surl{\
      _{\beta} \star_{\alpha}}_{\ \mu} id)(\Gamma(a)))|
\Lambda_{\Phi} ((\omega_{v,\xi_i} \surl{\
      _{\beta} \star_{\alpha}}_{\ \mu} id)(\Gamma(a))))\\
&=\Phi((\lambda_v^{\beta,\alpha})^*\Gamma(a^*)[\sum_{i\in I}
\lambda_{\xi_i}^{\beta,\alpha}(\lambda_{\xi_i}^{\beta,\alpha})^*
\lambda_{\xi_i}^{\beta,\alpha}(\lambda_{\xi_i}^{\beta,\alpha})^*]
\Gamma(a)\lambda_v^{\beta,\alpha})\\
&=\Phi((\omega_{v,v}\surl{\ _{\beta} \star_{\alpha}}_{\ \mu}
id)(\Gamma(a^*a)))=((id \surl{\ _{\beta} \star_{\alpha}}_{\ \mu}
\Phi)(\Gamma(a^*a))v|v)=(T_L(a^*a)v|v)<\infty
\end{aligned}$$ We deduce that the sum $\sum_{i\in I} \xi_{i} \surl{\ _{\beta}
\otimes_{\alpha}}_{\ \mu}\Lambda_{\Phi} ((\omega_{v,\xi_i} \surl{\
_{\beta} \star_{\alpha}}_{\ \mu} id)(\Gamma(a)))$ converges in $H
\surl{\ _{\beta} \otimes_{\alpha}}_{\ \mu} H_{\Phi}$. To prove that
the sum does not depend on the $(N^o, \mu^{o})$-basis, we compute
for all $b \in {\mathcal N}_{T_L} \cap {\mathcal N}_{\Phi}$ and $w
\in D(H_{\beta},\mu^o)$:
$$
\begin{aligned}
&\quad(\sum_{i\in I} \xi_{i} \surl{\ _{\beta} \otimes_{\alpha}}_{\
\mu} \Lambda_{\Phi} ((\omega_{v,\xi_i} \surl{\ _{\beta}
\star_{\alpha}}_{\ \mu} id)(\Gamma(a)))) | w
\surl{\ _{\beta} \otimes_{\alpha}}_{\ \mu} \Lambda_{\Phi}(b))\\
&=\sum_{i \in I}
(\alpha(<\xi_i,w>_{\beta,\mu^o})\Lambda_{\Phi}((\omega_{v,\xi_i}
\surl{\ _{\beta}\star_{\alpha}}_{\ \mu}id)
(\Gamma(a)))|\Lambda_{\Phi}(b))\\
&=\sum_{i\in I}
\Phi(b^*\alpha(<\xi_i,w>_{\beta,\mu^o})(\omega_{v,\xi_i} \surl{\
_{\beta} \star_{\alpha}}_{\ \mu} id)(\Gamma(a)))\\
&=\Phi(b^*\lambda_w^{\beta,\alpha}[\sum_{i \in I}
\lambda^{\beta,\alpha}_{\xi_i}(\lambda^{\beta,\alpha}_{\xi_i})^*]
\Gamma(a)\lambda_v^{\beta,\alpha})=\Phi(b^*(\omega_{v,w}\surl{\
_{\beta} \star_{\alpha}}_{\ \mu} id)(\Gamma(a))).
\end{aligned}$$ As $D(H_{\beta}, \mu^{o})
\odot \Lambda_{\Phi}({\mathcal N}_{T_L} \cap {\mathcal N}_{\Phi})$
is dense in $H \surl{\ _{\beta} \otimes_{\alpha}}_{\ \mu} H_{\Phi}$
and the last expression is independent of the $(N^o,
\mu^{o})$-basis, we can conclude.
\end{proof}

\begin{theo}\label{isom}
Let $H$ be a Hilbert space on which $M$ acts. There exists a unique
isometry $U_H$, called (left) \textbf{fundamental isometry}, from $
H \surl{\ _{\alpha} \otimes_{\hat{\beta}}}_{\ \ \mu^o} H_{\Phi}$ to
$H \surl{\ _{\beta} \otimes_{\alpha}}_{\ \mu} H_{\Phi}$ such that,
for all $(N^o, \mu^{o})$-basis $(\xi_i)_{i \in I}$ of $H_{\beta}$,
$a \in {\mathcal N}_{T_L} \cap {\mathcal N}_{\Phi}$ and $v \in
D(H_{\beta}, \mu^{o})$:
$$U_H(v \surl{\ _{\alpha} \otimes_{\hat{\beta}}}_{\ \ \mu^o}
\Lambda_{\Phi}(a))=\sum_{i\in I} \xi_{i} \surl{\ _{\beta}
\otimes_{\alpha}}_{\ \mu}\Lambda_{\Phi} ((\omega_{v,\xi_i} \surl{\
_{\beta} \star_{\alpha}}_{\ \mu} id)(\Gamma(a))))$$
\end{theo}

\begin{proof}
By \ref{inde}, we can define the following application:
$$\begin{aligned}
\tilde{U}:\  D(H_{\beta}, \mu^{o}) \times \Lambda_{\Phi}({\mathcal
N}_{T} \cap {\mathcal N}_{\Phi}) &\rightarrow H \surl{\ _{\beta}
\otimes_{\alpha}}_{\ \mu}
H_{\Phi} \\
(v,\Lambda_{\Phi}(a)) & \mapsto \sum_{i \in I} \xi_{i} \surl{\
_{\beta} \otimes_{\alpha}}_{\ \mu} \Lambda_{\Phi} ((\omega_{v,\xi_i}
\surl{\ _{\beta} \star_{\alpha}}_{\ \mu} id)(\Gamma(a))))
\end{aligned}$$
Let $b \in {\mathcal N}_{T_L} \cap {\mathcal N}_{\Phi}$ and $w \in
D(H_{\beta},\mu^o)$. Then, by normality and left invariance of
$\Phi$, we have:
$$
\begin{aligned}
&\ \quad(\tilde{U}(v,\Lambda_{\Phi}(a))|\tilde{U}(w,\Lambda_{\Phi}(b)))\\
&=\sum_{i,j\in I}
(\alpha(<\xi_i,\xi_j>_{\beta,\mu^o})\Lambda_{\Phi}((\omega_{v,\xi_i}
\surl{\ _{\beta} \star_{\alpha}}_{\ \mu} id)(\Gamma(a)))
|\Lambda_{\Phi}((\omega_{w,\xi_i} \surl{\ _{\beta}
\star_{\alpha}}_{\ \mu} id)(\Gamma(b))))\\
&=\sum_{i\in I} (\Lambda_{\Phi}(\alpha(<\xi_i,\xi_i>_{\beta,\mu^o})
(\omega_{v,\xi_i} \surl{\ _{\beta} \star_{\alpha}}_{\ \mu}
id)(\Gamma(a))) | \Lambda_{\Phi}((\omega_{w,\xi_i} \surl{\
_{\beta} \star_{\alpha}}_{\ \mu} id)(\Gamma(b))))\\
&=\sum_{i \in I}\Phi((\lambda_w^{\beta,\alpha})^*\Gamma(b^*)
\lambda_{\xi_i}^{\beta,\alpha}\alpha(<\xi_i,\xi_i>_{\beta,\mu^o})
(\lambda_{\xi_i}^{\beta,\alpha})^*\Gamma(a)\lambda_v^{\beta,\alpha})\\
&=\Phi ((\lambda_w^{\beta,\alpha})^*\Gamma(b^*)[\sum_{i \in I}
\lambda_{\xi_i}^{\beta,\alpha}(\lambda_{\xi_i}^{\beta,\alpha})^*
\lambda_{\xi_i}^{\beta,\alpha}(\lambda_{\xi_i}^{\beta,\alpha})^*]
\Gamma(a)\lambda_v^{\beta,\alpha}) \end{aligned}$$ Then, properties
of $(N^o, \mu^{o})$-basis $(\xi_i)_{i \in I}$ of $H_{\beta}$ imply
that:
$$\begin{aligned}
\Phi((\omega_{v,w}\surl{\ _{\beta}\star_{\alpha}}_{\ \mu}
id)(\Gamma(b^*a)))&=\omega_{v,w}((id \surl{\ _{\beta}
\star_{\alpha}}_{\ \mu}\Phi)(\Gamma(b^*a)))\\
&=\omega_{v,w}(T_L(b^*a))=(T_L(b^*a)v|w)
\end{aligned}$$
By \ref{rapide}, we get:
$$(\tilde{U}((v,\Lambda_{\Phi}(a))|\tilde{U}((w,\Lambda_{\Phi}(b))))
=(v\surl{\ _{\alpha}\otimes_{\hat{\beta}}}_{\ \ \mu^o}
\Lambda_{\Phi}(a)|w\surl{\ _{\alpha} \otimes_{\hat{\beta}}}_{\ \
\mu^o} \Lambda_{\Phi}(b))$$ so that, from $\tilde{U}$, we can easily
define a suitable application $U_H$ which is independent of the
$(N^o,\mu^{o})$-basis by \ref{inde}.
\end{proof}

One can define a right version of $U_H$ from the right invariant
weight:

\begin{theo}
Let $H$ be a Hilbert space on which $M$ acts. There exists a unique
isometry $U'_H$, called \textbf{right fundamental isometry}, from
$H_{\Psi} \surl{\ _{\hat{\alpha}} \otimes_{\beta}}_{\ \ \mu^o} H$ to
$H_{\Psi} \surl{\ _{\beta} \otimes_{\alpha}}_{\ \mu} H$ such that,
for all $(N,\mu)$-basis $(\eta_i)_{i \in I}$ of $\ _{\alpha}H$, $a
\in {\mathcal N}_{T_R} \cap {\mathcal N}_{\Psi}$ and $v \in
D(_{\alpha}H, \mu)$:
$$U'_H(\Lambda_{\Psi}(a) \surl{\ _{\hat{\alpha}} \otimes_{\beta}}_{\ \ \mu^o} v
)=\sum_{i\in I} \Lambda_{\Psi} ((id\surl{\
_{\beta}\star_{\alpha}}_{\ \mu} \omega_{v,\eta_i})(\Gamma(a)))
\surl{\ _{\beta} \otimes_{\alpha}}_{\ \mu} \eta_{i}$$
\end{theo}

\subsection{Fundamental isometry and co-product}
In this paragraph, we establish several links between fundamental
isometry and co-product. In fact, many of the following relations
are more or less equivalent to definition of fundamental unitary
and, depending of the situation, we will give priority to one or the
other relations in our demonstrations.

\begin{prop}\label{raccourci}
We have $(1\surl{\ _{\beta} \otimes_{\alpha}}_{\
N}J_{\Phi}eJ_{\Phi})U_H
\rho^{\alpha,\hat{\beta}}_{\Lambda_{\Phi}(x)}=\Gamma(x)
\rho^{\beta,\alpha}_{J_{\Phi}\Lambda_{\Phi}(e)}$ for all $e,x\in
{\mathcal N}_{\Phi}\cap {\mathcal N}_{T_L}$ and
$(J_{\Psi}fJ_{\Psi}\surl{\ _{\beta} \otimes_{\alpha}}_{\ N}1)U'_H
\lambda^{\hat{\alpha},\beta}_{\Lambda_{\Psi}(y)}=\Gamma(y)
\lambda^{\beta,\alpha}_{J_{\Psi}\Lambda_{\Psi}(f)}$\! for all
$f,y\in {\mathcal N}_{\Psi}\cap {\mathcal N}_{T_R}$.
\end{prop}

\begin{proof}
Let $v\in D(H_{\beta},\mu^o)$ and $(\xi_i)_{i \in I}$ a
$(N^o,\mu^o)$-basis of $H_{\beta}$. We have:

$$
\begin{aligned}
&\ \quad (1\surl{\ _{\beta} \otimes_{\alpha}}_{\
N}J_{\Phi}eJ_{\Phi})U_H(v\surl{\ _{\alpha}
\otimes_{\hat{\beta}}}_{\ \ \mu^o} \Lambda_{\Phi}(x))\\
&=\sum_{i\in I}\xi_i\surl{\ _{\beta} \otimes_{\alpha}}_{\ \mu}
J_{\Phi}eJ_{\Phi}\Lambda_{\Phi}((\omega_{v,\xi_i}\surl{\ _{\beta}
\star_{\alpha}}_{\ \mu} id) (\Gamma(x)))\\
&=\sum_{i\in I}\xi_i\surl{\ _{\beta} \otimes_{\alpha}}_{\ \mu}
(\omega_{v,\xi_i}\surl{\ _{\beta} \star_{\alpha}}_{\ \mu} id)
(\Gamma(x))J_{\Phi}\Lambda_{\Phi}(e)=\Gamma(x)(v \surl{\ _{\beta}
\otimes_{\alpha}}_{\ \mu} J_{\Phi}\Lambda_{\Phi}(e))
\end{aligned}$$

By \ref{prem} and \ref{evi}, we have $\Lambda_{\Phi}(x)\in
D((H_{\Phi})_{\hat{\beta}},\mu^o)$ and $J_{\Phi}\Lambda_{\Phi}(e)\in
D( _{\alpha}(H_{\Phi}),\mu)$ so that each term of the previous
equality is continuous in $v$. Density of $D(H_{\beta},\mu^o)$ in
$H$ finishes the proof. The last part is very similar.
\end{proof}

\begin{prop}\label{rap}
For all $v,w\in D(H_{\beta},\mu^o)$ and $a\in {\mathcal
N}_{\Phi}\cap {\mathcal N}_{T_L}$, we have:
$$(\lambda_w^{\beta,\alpha})^*U_H(v\surl{\
_{\alpha}\otimes_{\hat{\beta}}}_{\ \ \mu^o}\Lambda_{\Phi}(a))=
\Lambda_{\Phi}((\omega_{v,w}\surl{\ _{\beta}\star_{\alpha}}_{\ \mu}
id)(\Gamma(a)))$$ Also, for all $v',w'\in D(_{\alpha}H,\mu)$ and
$b\in {\mathcal N}_{\Psi}\cap {\mathcal N}_{T_R}$, we have:
$$(\rho_{w'}^{\beta,\alpha})^*U'_H(\Lambda_{\Psi}(b)\surl{\
_{\alpha}\otimes_{\hat{\beta}}}_{\ \ \mu^o}v')=
\Lambda_{\Psi}((id\surl{\ _{\beta}\star_{\alpha}}_{\ \mu}
\omega_{v',w'})(\Gamma(b)))$$
\end{prop}

\begin{proof}
Let $e\in {\mathcal N}_{\Phi}\cap {\mathcal N}_{T_L}$. By
\ref{raccourci}, we can compute:

$$
\begin{aligned}
J_{\Phi}eJ_{\Phi}(\lambda_w^{\beta,\alpha})^*U_H(v\surl{\
_{\alpha}\otimes_{\hat{\beta}}}_{\ \ \mu^o}\Lambda_{\Phi}(a)) &=
(\lambda_w^{\beta,\alpha})^*(1\surl{\ _{\beta}\otimes_{\alpha}}_{\
N}J_{\Phi}eJ_{\Phi})U_H\rho_{\Lambda_{\Phi}(a)}^{\alpha,\hat{\beta}}v\\
&=(\lambda_w^{\beta,\alpha})^*\Gamma(a)\rho^{\beta,\alpha}_{J_{\Phi}\Lambda_{\Phi}(e)}v\\
&=(\omega_{v,w}\surl{\ _{\beta}\star_{\alpha}}_{\ \mu}
id)(\Gamma(a))J_{\Phi}\Lambda_{\Phi}(e)\\
&=J_{\Phi}eJ_{\Phi}\Lambda_{\Phi}((\omega_{v,w}\surl{\
_{\beta}\star_{\alpha}}_{\ \mu} id)(\Gamma(a)))
\end{aligned}$$
Density of ${\mathcal N}_{\Phi}\cap {\mathcal N}_{T_L}$ in $N$
finishes the proof. The second part is very similar.
\end{proof}

\begin{coro}\label{lienGV}
For all $a \in {\mathcal N}_{T_L}\cap {\mathcal N}_{\Phi}$, $v\in D(
_{\alpha}H,\mu) \cap D(H_{\beta},\mu^o)$ and $w\in
D(H_{\beta},\mu^o)$, we have: $$(\omega_{v,w} *
id)(U_H)\Lambda_{\Phi}(a)= \Lambda_{\Phi}((\omega_{v,w} \surl{\
_{\beta} \star_{\alpha}}_{\ \mu} id)(\Gamma(a)))$$ where we denote
by $(\omega_{v,w}*id)(U_H)$ the operator
$(\lambda_w^{\beta,\alpha})^*U_H\lambda_v^{\alpha,\hat{\beta}}$ of
$\mathcal{L}(H_{\Phi})$.
\end{coro}

\begin{proof}
Straightforward.
\end{proof}

\begin{coro}\label{corres}
For all $e,x\in {\mathcal N}_{\Phi}\cap {\mathcal N}_{T_L}$ and
$\eta\in D(_{\alpha}H_{\Phi},\mu^o)$, we have: $$(id\surl{\ _{\beta}
\star_{\alpha}}_{\ \mu} \omega_{J_{\Phi}\Lambda_{\Phi}(e),\eta}
)(\Gamma(x))=(id*
\omega_{\Lambda_{\Phi}(x),J_{\Phi}e^*J_{\Phi}\eta})(U_H)$$ Also, for
all $f,y\in {\mathcal N}_{\Psi}\cap {\mathcal N}_{T_R}$ and $\xi\in
D((H_{\Psi})_{\beta},\mu^o)$, we have:
$$(\omega_{J_{\Psi}\Lambda_{\Psi}(f),\xi} \surl{\ _{\beta}
\star_{\alpha}}_{\
\mu}id)(\Gamma(y))=(\omega_{\Lambda_{\Psi}(y),J_{\Psi}f^*J_{\Psi}\xi}*
id)(U'_H)$$
\end{coro}

\begin{proof}
Straightforward by \ref{raccourci}.
\end{proof}

\begin{coro}\label{switch}
For all $a,b \in {\mathcal N}_{\Psi}\cap {\mathcal N}_{T_R} \cap
{\mathcal N}_{\Psi}^*\cap {\mathcal N}_{T_R}^*$, we have:
$$(\omega_{\Lambda_{\Psi}(a),J_{\Psi}\Lambda_{\Psi}(b)}*id)(U'_H)^*
=(\omega_{\Lambda_{\Psi}(a^*),J_{\Psi}\Lambda_{\Psi}(b^*)}*id)(U'_H)$$
\end{coro}

\begin{proof}
By \ref{corres}, we have for all $e\in {\mathcal N}_{\Psi}\cap
{\mathcal N}_{T_R}$:
$$
\begin{aligned}
(\omega_{\Lambda_{\Psi}(a),J_{\Psi}\Lambda_{\Psi}(e^*b)}*
id)(U'_H)^* &=
(\omega_{J_{\Psi}\Lambda_{\Psi}(e),J_{\Psi}\Lambda_{\Psi}(b)}
\surl{\ _{\beta} \star_{\alpha}}_{\ \mu}
id)(\Gamma(a))^* \\
&= (\omega_{J_{\Psi}\Lambda_{\Psi}(b),J_{\Psi}\Lambda_{\Psi}(e)}
\surl{\ _{\beta} \star_{\alpha}}_{\ \mu}
id)(\Gamma(a^*)) \\
&= (\omega_{\Lambda_{\Psi}(a^*),J_{\Psi}\Lambda_{\Psi}(b^*e)}*
id)(U'_H).
\end{aligned}$$

Let $(u_k)_{k\in K}$ be a family in ${\mathcal N}_{\Psi} \cap
{\mathcal N}_{\Psi}^*$ such that $u_k\rightarrow 1$ in the *-strong
topology. We denote:
$$e_k=\frac{1}{\sqrt{\pi}} \int\! e^{-t^2}\sigma_{t}^{\Psi}(u_k)\ dt$$
For all $k\in K$, $e_k$ and $\sigma_{-i/2}^{\Psi}(e^*_k)$ are
bounded and belong to ${\mathcal N}_{\Psi}$ and converge to $1$ in
the *-strong topology so that $J_{\Psi}\Lambda_{\Psi}(b^*e_k)
=\sigma_{-i/2}^{\Psi}(e^*_k)J_{\Psi}\Lambda_{\Psi}(b^*)$ converge to
$J_{\Psi}\Lambda_{\Psi}(b^*)$ in norm of $H_{\Psi}$. Let $\xi,\eta
\in D( _{\alpha}H,\mu)$ and we compute:
$$
\begin{aligned}
((\omega_{\Lambda_{\Psi}(a),J_{\Psi}\Lambda_{\Psi}(b)}*
id)(U'_H)^*\xi|\eta)&=(J_{\Psi}\Lambda_{\Psi}(b) \surl{\ _{\beta}
\otimes_{\alpha}}_{\ \mu}\xi|U'_H(\Lambda_{\Psi}(a)\surl{\
_{\hat{\alpha}} \otimes_{\beta}}_{\ \ \mu^o} \eta))\\
&=\lim_{k\in K}(J_{\Psi}\Lambda_{\Psi}(e_k^*b) \surl{\ _{\beta}
\otimes_{\alpha}}_{\ \mu}\xi|U'_H(\Lambda_{\Psi}(a)\surl{\
_{\hat{\alpha}} \otimes_{\beta}}_{\ \ \mu^o} \eta))\\
&=\lim_{k\in
K}((\omega_{\Lambda_{\Psi}(a),J_{\Psi}\Lambda_{\Psi}(e_k^*b)}*
id)(U'_H)^*\xi|\eta)\end{aligned}$$ By the previous computation,
this last expression is equal to:
$$
\begin{aligned}
&\ \quad\lim_{k\in
K}((\omega_{\Lambda_{\Psi}(a^*),J_{\Psi}\Lambda_{\Psi}(b^*e_k)}*
id)(U'_H)\xi|\eta)\\
&=\lim_{k\in K}(U'_H(\Lambda_{\Psi}(a)\surl{\ _{\hat{\alpha}}
\otimes_{\beta}}_{\ \ \mu^o} \xi)|J_{\Psi}\Lambda_{\Psi}(b^*e_k)
\surl{\ _{\beta} \otimes_{\alpha}}_{\ \mu}\eta)\\
&=(U'_H(\Lambda_{\Psi}(a^*)\surl{\ _{\hat{\alpha}}
\otimes_{\beta}}_{\ \ \mu^o} \xi)|J_{\Psi}\Lambda_{\Psi}(b^*)
\surl{\ _{\beta}\otimes_{\alpha}}_{\
\mu}\eta)=((\omega_{\Lambda_{\Psi}(a^*),J_{\Psi}\Lambda_{\Psi}(b^*)}*
id)(U'_H)\xi|\eta)\\
\end{aligned}$$
By density of $D( _{\alpha}H,\mu)$ in $H$, the result holds.
\end{proof}

\subsection{Commutation relations}\label{rcom}
In this section, we verify commutation relations which are necessary
for $U_H$ to be a pseudo-multiplicative unitary and we establish a
link between $U_H$ and $\Gamma$. We also have similar formulas for
$U'_H$.

\begin{lemm}\label{base}
Let $\xi\in D(H_{\beta},\mu^o)$ and $\eta\in D( _{\alpha}H,\mu)$.

\begin{center}
\begin{minipage}{11cm}
\begin{enumerate}[i)]
\item For all $a \in \alpha(N)'$, we have
$\lambda_\xi^{\beta,\alpha}\circ a=(1 \surl{\ _{\beta}
\otimes_{\alpha}}_{\ N} a)\lambda_{\xi}^{\beta,\alpha}$.
\item For all $b \in \beta(N)'$, we have
$\lambda_{b\xi}^{\beta,\alpha}=(b \surl{\ _{\beta}
  \otimes_{\alpha}}_{\ N} 1)\lambda_{\xi}^{\beta,\alpha}$.
\item For all $x \in {\mathcal D}(\sigma_{-i/2}^{\mu})$, we have
$\lambda_{\beta(x)\xi}^{\beta,\alpha}=\lambda_{\xi}^{\beta,\alpha}\circ
\alpha(\sigma_{-i/2}^{\mu}(x))$.
\item For all $x \in {\mathcal D}(\sigma_{i/2}^{\mu})$, we have
$\rho_{\alpha(x)\eta}^{\beta,\alpha}=\rho_{\eta}^{\beta,\alpha}\circ
\beta(\sigma_{i/2}^{\mu}(x))$.
\end{enumerate}
\end{minipage}
\end{center}
\end{lemm}

\begin{proof}
Straightforward.
\end{proof}

We recall that $\alpha(N)$ and $\beta(N)$ commute with
$\hat{\beta}(N)'$.

\begin{prop}\label{comm}
For all $n\in N$, we have:

\begin{center}
\begin{minipage}{8cm}
\begin{enumerate}[i)]
\item $U_H(1\surl{\ _{\alpha}\otimes_{\hat{\beta}}}_{\ \ N^o}
  \alpha(n))=(\alpha(n)\surl{\ _{\beta}\otimes_{\alpha}}_{\ N}
  1)U_H$;
\item $U_H(1\surl{\ _{\alpha}\otimes_{\hat{\beta}}}_{\ \ N^o}
  \beta(n))=(1\surl{\ _{\beta}\otimes_{\alpha}}_{\ N}
  \beta(n))U_H$;
\item $U_H(\beta(n)\surl{\ _{\alpha}\otimes_{\hat{\beta}}}_{\ \ N^o}
  1)=(1\surl{\ _{\beta}\otimes_{\alpha}}_{\ N}
  \hat{\beta}(n))U_H$.
\end{enumerate}
\end{minipage}
\end{center}

\end{prop}

\begin{proof}
By \ref{raccourci}, we can compute for all $n\in N $ and $e,x \in
{\mathcal N}_{T_L} \cap {\mathcal N}_{\Phi}$:
$$\begin{aligned}
(\alpha(n) \surl{\ _{\beta} \otimes_{\alpha}}_{\
N}J_{\Phi}eJ_{\Phi})U_H\rho^{\alpha,\hat{\beta}}_{\Lambda_{\Phi}(x)}
&=(\alpha(n) \surl{\ _{\beta} \otimes_{\alpha}}_{\
N}1)\Gamma(x)\rho_{J_{\Phi}\Lambda_{\Phi}(e)}^{\beta,\alpha}\\
&=\Gamma(\alpha(n)x)\rho_{J_{\Phi}\Lambda_{\Phi}(e)}^{\beta,\alpha}\\
&=(1\surl{\ _{\beta} \otimes_{\alpha}}_{\
N}J_{\Phi}eJ_{\Phi})U_H\rho^{\alpha,\hat{\beta}}_{\Lambda_{\Phi}(\alpha(n)x)}\\
&=(1\surl{\ _{\beta} \otimes_{\alpha}}_{\ N}J_{\Phi}eJ_{\Phi})U_H(1
\surl{\ _{\alpha} \otimes_{\hat{\beta}}}_{\ \ N^o}
\alpha(n))\rho^{\alpha,\hat{\beta}}_{\Lambda_{\Phi}(x)}
\end{aligned}$$
Usual arguments of density imply the first equality. The second one
can be proved in a very similar way. By \ref{raccourci} and
\ref{base}, we can compute for all $n\in {\mathcal T}_{\mu}$ and
$e,x \in {\mathcal N}_{T_L} \cap {\mathcal N}_{\Phi}$:
$$
\begin{aligned}
(1\surl{\ _{\beta} \otimes_{\alpha}}_{\
N}J_{\Phi}eJ_{\Phi}\hat{\beta}(n))U_H\rho^{\alpha,\hat{\beta}}_{\Lambda_{\Phi}(x)}
&=\Gamma(x)\rho_{J_{\Phi}\Lambda_{\Phi}(e\alpha(n^*))}^{\beta,\alpha}\\
&=\Gamma(x)\rho_{\alpha(\sigma_{-i/2}^{\mu}(n))J_{\Phi}\Lambda_{\Phi}(e)}^{\beta,\alpha}\\
&=\Gamma(x)\rho_{J_{\Phi}\Lambda_{\Phi}(e)}^{\beta,\alpha}\beta(n)\\
&=(1\surl{\ _{\beta} \otimes_{\alpha}}_{\
N}J_{\Phi}eJ_{\Phi})U_H\rho^{\alpha,\hat{\beta}}_{\Lambda_{\Phi}(x)}\beta(n)\\
&=(1\surl{\ _{\beta} \otimes_{\alpha}}_{\
N}J_{\Phi}eJ_{\Phi})U_H(\beta(n)\surl{\ _{\alpha}
\otimes_{\hat{\beta}}}_{\ \ N^o}
1)\rho^{\alpha,\hat{\beta}}_{\Lambda_{\Phi}(x)}
\end{aligned}$$
Density of ${\mathcal T}_{\mu}$ in $N$ and normality of $\beta$ and
$\hat{\beta}$ finish the proof.
\end{proof}

\begin{prop}\label{appartenance}
For all $x\in M'\cap\mathcal{L}(H)$, we have: $$U_H(x \surl{\
_{\alpha} \otimes_{\hat{\beta}}}_{\ \ N^o} 1)= (x \surl{\ _{\beta}
\otimes_{\alpha}}_{\ N} 1)U_H$$
\end{prop}

\begin{proof}
For all $e,y \in {\mathcal N}_{T_L} \cap {\mathcal N}_{\Phi}$ and
$x\in M'\cap\mathcal{L}(H)\subseteq
\alpha(N)'\cap\beta(N)'\cap\mathcal{L}(H)$, we have by
\ref{raccourci}:
$$\begin{aligned}
(x\surl{\ _{\beta} \otimes_{\alpha}}_{\
N}J_{\Phi}eJ_{\Phi})U_H\rho^{\alpha,\hat{\beta}}_{\Lambda_{\Phi}(y)}
&=(x\surl{\ _{\beta} \otimes_{\alpha}}_{\
N}1)\Gamma(y)\rho_{J_{\Phi}\Lambda_{\Phi}(e)}^{\beta,\alpha}\\
&=\Gamma(y)\rho_{J_{\Phi}\Lambda_{\Phi}(e)}^{\beta,\alpha}x\\
&=(1\surl{\ _{\beta} \otimes_{\alpha}}_{\
N}J_{\Phi}eJ_{\Phi})U_H\rho^{\alpha,\hat{\beta}}_{\Lambda_{\Phi}(y)}x\\
&=(1\surl{\ _{\beta} \otimes_{\alpha}}_{\ N}J_{\Phi}eJ_{\Phi})U_H(x
\surl{\ _{\alpha} \otimes_{\hat{\beta}}}_{\ \ N^o}
1)\rho^{\alpha,\hat{\beta}}_{\Lambda_{\Phi}(y)}
\end{aligned}$$ Usual arguments of density imply the result.
\end{proof}

\begin{coro}
For all $n \in N$, we have:

\begin{center}
\begin{minipage}{8cm}
\begin{enumerate}[i)]
\item $U_{H_{\Phi}}(\hat{\beta}(n) \surl{\ _{\alpha}
\otimes_{\hat{\beta}}}_{\ \ N^o} 1)=(\hat{\beta}(n) \surl{\ _{\beta}
\otimes_{\alpha}}_{\ N} 1)U_{H_{\Phi}}$
\item $U_{H_{\Psi}}(\hat{\alpha}(n) \surl{\ _{\alpha}
\otimes_{\hat{\beta}}}_{\ \ N^o} 1)=(\hat{\alpha}(n) \surl{\
_{\beta} \otimes_{\alpha}}_{\ N} 1)U_{H_{\Psi}}$
\end{enumerate}
\end{minipage}
\end{center}

\end{coro}

\begin{prop}\label{impl}
We have $\Gamma(m)U_H=U_H(1 \surl{\ _{\alpha}
\otimes_{\hat{\beta}}}_{\ \ N^o} m)$ for all $m \in M$.
\end{prop}

\begin{proof}
By \ref{raccourci}, we can compute for all $e,x\in {\mathcal
N}_{T_L} \cap {\mathcal N}_{\Phi}$:
$$\begin{aligned}
(1\surl{\ _{\beta} \otimes_{\alpha}}_{\
N}J_{\Phi}eJ_{\Phi})\Gamma(m)U_H\rho^{\alpha,\hat{\beta}}_{\Lambda_{\Phi}(x)}
&=\Gamma(m)(1\surl{\ _{\beta} \otimes_{\alpha}}_{\
N}J_{\Phi}eJ_{\Phi})U_H\rho^{\alpha,\hat{\beta}}_{\Lambda_{\Phi}(x)}\\
&=\Gamma(mx)\rho_{J_{\Phi}\Lambda_{\Phi}(e)}^{\beta,\alpha}\\
&=(1\surl{\ _{\beta} \otimes_{\alpha}}_{\
N}J_{\Phi}eJ_{\Phi})U_H\rho^{\alpha,\hat{\beta}}_{\Lambda_{\Phi}(mx)}\\
&=(1\surl{\ _{\beta} \otimes_{\alpha}}_{\
N}J_{\Phi}eJ_{\Phi})U_H(1\surl{\ _{\alpha} \otimes_{\hat{\beta}}}_{\
\ N^o} m)\rho^{\alpha,\hat{\beta}}_{\Lambda_{\Phi}(x)}
\end{aligned}$$ Usual arguments of density imply the result.
\end{proof}

\subsection{Unitarity of the fundamental isometry}
This is a key part of the theory and certainly one of the most
difficult. To prove unitary of $U_H$ (resp. $U'_H$), we establish a
reciprocity law where both left and right operator-valued weights
are at stake.

\subsubsection{First technical result}

We establish results needed for \ref{reciprocite}. In the following
proposition, we compute some functions $\theta$ defined in section
\ref{intre}.

\begin{prop}\label{inov}
We have for all $c\in {\mathcal N}_{\Psi}\cap {\mathcal N}_{T_R}$,
$m\in ({\mathcal N}_{\Psi}\cap {\mathcal N}_{T_R})^*$ and $v\in
D(H_{\beta},\mu^o)$:
$$\theta^{\beta,\mu^o}(v,J_{\Psi}\Lambda_{\Psi}(c))m=(\lambda_{\Lambda_{\Psi}
(m^*)}^{\hat{\alpha},\beta})^*\rho_v^{\hat{\alpha},\beta}J_{\Psi}c^*J_{\Psi}$$
\end{prop}

\begin{proof}
Let $x\in {\mathcal N}_{\Psi}\cap {\mathcal N}_{T_R}$. On one hand,
we get by \ref{prem} and \ref{evi}:
$$
\begin{aligned}
\theta^{\beta,\mu^o}(v,J_{\Psi}\Lambda_{\Psi}(c))m\Lambda_{\Psi}(x)
&=
R^{\beta,\mu^o}(v)R^{\beta,\mu^o}(J_{\Psi}\Lambda_{\Psi}(c))^*\Lambda_{\Psi}(mx)\\
&=
R^{\beta,\mu^o}(v)J_{\mu}\Lambda_{T_R}(c)^*J_{\Psi}\Lambda_{\Psi}(mx).
\end{aligned}$$

On the other hand, if $c\in {\mathcal T}_{\Psi,T_R}$, then we have
by \ref{rapide}:
$$
\begin{aligned}
(\lambda_{\Lambda_{\Psi}(m^*)}^{\hat{\alpha},\beta})^*\rho_v^{\hat{\alpha},\beta}
J_{\Psi}c^*J_{\Psi}\Lambda_{\Psi}(x)
&=(\lambda_{\Lambda_{\Psi}(m^*)}^{\hat{\alpha},\beta})^*
(J_{\Psi}c^*J_{\Psi}\Lambda_{\Psi}(x) \surl{\
_{\hat{\alpha}}\otimes_{\beta}}_{\ \mu^o} v) \\
&=T_R(mx\sigma_{-i/2}^{\Psi}(c))v\\
&=R^{\beta,\mu^o}(v)J_{\mu}
\Lambda_{\mu}(\beta^{-1}(T_R(\sigma_{i/2}^{\Psi}(c^*)x^*m^*)))\\
&=R^{\beta,\mu^o}(v)J_{\mu}\Lambda_{T_R}(c)^*J_{\Psi}\Lambda_{\Psi}(mx)
\end{aligned}$$
We obtain:
$$(\lambda_{\Lambda_{\Psi}(m^*)}^{\hat{\alpha},\beta})^*\rho_v^{\hat{\alpha},\beta}
J_{\Psi}c^*J_{\Psi}\Lambda_{\Psi}(x)=R^{\beta,\mu^o}(v)J_{\mu}\Lambda_{T_R}(c)^*J_{\Psi}\Lambda_{\Psi}(mx)$$
for all $c\in {\mathcal N}_{\Psi}\cap {\mathcal N}_{T_R}$ by
normality which finishes the proof.
\end{proof}

\begin{coro}\label{prepa}
Let $a\in ({\mathcal N}_{\Psi}\cap {\mathcal N}_{T_R})^* ({\mathcal
N}_{\Phi}\cap {\mathcal N}_{T_L})$. If $c\in {\mathcal N}_{\Psi}\cap
{\mathcal N}_{T_R}$, $e\in {\mathcal N}_{\Phi}\cap {\mathcal
N}_{T_L}$ and $\xi\in H_{\Psi},\eta \in D(_{\alpha}
(H_{\Phi}),\mu)$, $u\in H$, $v\in D(H_{\beta},\mu^o)$, then we have:
$$(v\surl{\ _{\beta}\otimes_{\alpha}}_{\ \mu}
(\lambda^{\beta,\alpha}_{J_{\Psi}\Lambda_{\Psi}(c)})^*
U_{H_{\Psi}}(\xi\surl{\ _{\alpha}\otimes_{\hat{\beta}}}_{\ \
\mu^o}\Lambda_{\Phi}(a))|u\surl{\ _{\beta}\otimes_{\alpha}}_{\ \mu}
J_{\Phi}e^*J_{\Phi}\eta)$$
$$=(J_{\Psi}c^*J_{\Psi}\xi\surl{\ _{\hat{\alpha}}\otimes_{\beta}}_{\ \ \mu^o} v|
\Lambda_{\Psi}((id\surl{\ _{\beta}\star_{\alpha}}_{\ \mu}
\omega_{\eta, J_{\Phi}\Lambda_{\Phi}(e)})(\Gamma(a^*)))\surl{\
_{\hat{\alpha}}\otimes_{\beta}}_{\ \ \mu^o}u)$$
\end{coro}

\begin{proof}
By \ref{raccourci} and \ref{inov}, we can compute:
$$\begin{aligned}
&\ \quad(v\surl{\ _{\beta}
  \otimes_{\alpha}}_{\ \mu}
  (\lambda^{\beta,\alpha}_{J_{\Psi}\Lambda_{\Psi}(c)})^*
  U_{H_{\Psi}}(\xi\surl{\ _{\alpha}\otimes_{\hat{\beta}}}_{\ \ \mu^o}\Lambda_{\Phi}(a))
  |u\surl{\ _{\beta}
  \otimes_{\alpha}}_{\ \mu} J_{\Phi}e^*J_{\Phi}\eta)\\
&=((\rho^{\beta,\alpha}_{\eta})^*
\lambda_v^{\beta,\alpha}(\lambda^{\beta,\alpha}_{J_{\Psi}\Lambda_{\Psi}(c)})^*
(1\surl{\ _{\beta}
  \otimes_{\alpha}}_{\ N}J_{\Phi}e^*J_{\Phi})
  U_{H_{\Psi}}(\xi\surl{\ _{\alpha}\otimes_{\hat{\beta}}}_{\ \
  \mu^o}\Lambda_{\Phi}(a))|u)\\
&=((\rho^{\beta,\alpha}_{\eta})^*
\lambda_v^{\beta,\alpha}(\lambda^{\beta,\alpha}_{J_{\Psi}\Lambda_{\Psi}(c)})^*
\Gamma(a)\rho^{\beta,\alpha}_{J_{\Phi}\Lambda_{\Phi}(e)}\xi|u)\\
&=\theta^{\beta,\mu^o}(v,J_{\Psi}\Lambda_{\Psi}(c))(\rho^{\beta,\alpha}_{\eta})^*
\Gamma(a)\rho^{\beta,\alpha}_{J_{\Phi}\Lambda_{\Phi}(e)}\xi|u)\\
&=((\lambda^{\hat{\alpha},\beta}_{\Lambda_{\Psi}((id\surl{\ _{\beta}
\star_{\alpha}}_{\ \mu} \omega_{\eta,
J_{\Phi}\Lambda_{\Phi}(e))})(\Gamma(a^*))})^*\rho^{\hat{\alpha},\beta}_v
J_{\Psi}c^*J_{\Psi}\xi|u)\\
&=(J_{\Psi}c^*J_{\Psi}\xi\surl{\ _{\hat{\alpha}}
  \otimes_{\beta}}_{\ \ \mu^o} v|\Lambda_{\Psi}((id\surl{\
_{\beta} \star_{\alpha}}_{\ \mu} \omega_{\eta,
J_{\Phi}\Lambda_{\Phi}(e))})(\Gamma(a^*)))\surl{\ _{\hat{\alpha}}
  \otimes_{\beta}}_{\ \ \mu^o} u)
\end{aligned}$$

\end{proof}

\subsubsection{Second technical result}

In this section, results only depend on \ref{raccourci} and
co-product relation but not on the previous technical result. Let
${\mathcal H}$ be an other Hilbert space on which $M$ acts.

\begin{lemm}\label{simple}
Let $a,e\in {\mathcal N}_{\Phi}\cap {\mathcal N}_{T_L}$, $\xi\in
D({\mathcal H}_{\beta},\mu^o)$, $\eta\in D(_{\alpha}H,\mu)$, and
$\zeta\in {\mathcal H}$. We have:
$$(1\surl{\ _{\beta} \otimes_{\alpha}}_{\ N}J_{\Phi}eJ_{\Phi})
U_H(\eta\surl{\ _{\alpha} \otimes_{\hat{\beta}}}_{\ \ \mu^o}
[(\lambda_{\xi}^ {\beta,\alpha})^*U_{\mathcal H}(\zeta\surl{\
_{\alpha} \otimes_{\hat{\beta}}}_{\ \ \mu^o}\Lambda_{\Phi} (a))])$$
$$=(\lambda_{\xi}^ {\beta,\alpha}\surl{\ _{\beta} \otimes_{\alpha}}_{\
N}1)^*(id \surl{\ _{\beta}
  \star_{\alpha}}_{\ N} \Gamma)(\Gamma(a))(\zeta
  \surl{\ _{\beta} \otimes_{\alpha}}_{\ \mu}
\eta \surl{\ _{\beta} \otimes_{\alpha}}_{\ \mu}
J_{\Phi}\Lambda_{\Phi}(e))$$
\end{lemm}

\begin{proof}
First let assume $\zeta\in D({\mathcal H}_{\beta},\mu^o)$. By
\ref{rap} and \ref{raccourci}, we can compute:
$$
\begin{aligned} &\quad(1\surl{\ _{\beta} \otimes_{\alpha}}_{\
N}J_{\Phi}eJ_{\Phi}) U_H(\eta\surl{\ _{\alpha}
\otimes_{\hat{\beta}}}_{\ \ \mu^o} [(\lambda_{\xi}^
{\beta,\alpha})^*U_{\mathcal H}(\zeta\surl{\ _{\alpha}
\otimes_{\hat{\beta}}}_{\ \ \mu^o}\Lambda_{\Phi} (a))])\\
&=(1\surl{\ _{\beta} \otimes_{\alpha}}_{\ N}J_{\Phi}eJ_{\Phi})
  U_H(\eta\surl{\ _{\alpha} \otimes_{\hat{\beta}}}_{\ \ \mu^o} \Lambda_{\Phi}
((\omega_{\zeta,\xi}\surl{\ _{\beta}
  \star_{\alpha}}_{\ \mu} id)(\Gamma(a)))\\
&=\Gamma((\omega_{\zeta,\xi}\surl{\ _{\beta}
  \star_{\alpha}}_{\ \mu} id)(\Gamma(a)))
  (\eta \surl{\ _{\beta} \otimes_{\alpha}}_{\ \mu}
  J_{\Phi}\Lambda_{\Phi}(e))\\
&=(\lambda_{\xi}^ {\beta,\alpha}\surl{\ _{\beta}
\otimes_{\alpha}}_{\ N}1)^*(id \surl{\ _{\beta}
  \star_{\alpha}}_{\ N} \Gamma)(\Gamma(a))(\zeta
  \surl{\ _{\beta} \otimes_{\alpha}}_{\ \mu}
\eta \surl{\ _{\beta} \otimes_{\alpha}}_{\ \mu}
J_{\Phi}\Lambda_{\Phi}(e))
\end{aligned}$$
So, we get the result for all $\zeta\in D({\mathcal
H}_{\beta},\mu^o)$. The first term of the equality is continuous in
$\zeta$ because $\eta\in D(_{\alpha}H,\mu)$ and
$\Lambda_{\Phi}(a)\in D((H_{\Phi})_{\hat{\beta}},\mu^o)$. Also,
since $\eta\in D(_{\alpha}H,\mu)$ and $\Lambda_{\Phi}(a)\in
D((H_{\Phi})_{\hat{\beta}},\mu^o)$, the last term of the equality is
continuous in $\zeta$. Density of $D({\mathcal H}_{\beta},\mu^o)$ in
${\mathcal H}$ finishes the proof.
\end{proof}

\begin{lemm}\label{simple2}
The sum $\sum_{i\in I} \eta_i\surl{\ _{\alpha}
\otimes_{\hat{\beta}}}_{\ \ \mu^o}
[(\lambda_{\xi}^{\beta,\alpha})^*U_{\mathcal
H}((\rho_{\eta_i}^{\beta,\alpha})^*\Xi\surl{\ _{\alpha}
\otimes_{\hat{\beta}}}_{\ \ \mu^o} \Lambda_{\Phi}(a))]$ converges
for all $\xi\in D({\mathcal H}_{\beta},\mu^o)$, $\Xi\in {\mathcal
H}\surl{\ _{\beta}\otimes_{\alpha}}_{\ \mu} H$, $a\in {\mathcal
N}_{\Phi}\cap {\mathcal N}_{T_L}$ and $(N,\mu)$-basis
$(\eta_i)_{i\in I}$ of $\ _{\alpha}H$.
\end{lemm}

\begin{proof}
First, observe that $\eta_i\surl{\
_{\alpha}\otimes_{\hat{\beta}}}_{\ \ \mu^o}
[(\lambda_{\xi}^{\beta,\alpha})^*U_{\mathcal
H}((\rho_{\eta_i}^{\beta,\alpha})^*\Xi\surl{\ _{\hat{\alpha}}
\otimes_{\beta}}_{\ \ \mu^o} \Lambda_{\Phi}(a))]$ are orthogonal. To
compute, we put: $\Omega_i=\rho_{\eta_i}^{\beta,\alpha})^*\Xi\surl{\
_{\alpha}
  \otimes_{\hat{\beta}}}_{\ \ \mu^o} \Lambda_{\Phi}(a)$. By \ref{base} and \ref{comm}, we have:
$$
\begin{aligned}
&\ \quad ||\eta_i\surl{\ _{\alpha}
  \otimes_{\hat{\beta}}}_{\ \ \mu^o} [(\lambda_{\xi}^{\beta,\alpha})^*
  U_{\mathcal H}(\Omega_i)]||^2\\
&=(\hat{\beta}(<\eta_i,\eta_i>_{\alpha,\mu})
(\lambda_{\xi}^{\beta,\alpha})^*
  U_{\mathcal H}(\Omega_i)|(\lambda_{\xi}^{\beta,\alpha})^*
  U_{\mathcal H}(\Omega_i))\\
&=((\lambda_{\xi}^{\beta,\alpha})^*(1\surl{\ _{\beta}
  \otimes_{\alpha}}_{\ \mu}\hat{\beta}(<\eta_i,\eta_i>_{\alpha,\mu}))
  U_{\mathcal H}(\Omega_i)|(\lambda_{\xi}^{\beta,\alpha})^*
  U_{\mathcal H}(\Omega_i))\\
&=((\lambda_{\xi}^{\beta,\alpha})^*
  U_{\mathcal H}(\beta(<\eta_i,\eta_i>_{\alpha,\mu})
  (\Omega_i)|(\lambda_{\xi}^{\beta,\alpha})^*U_{\mathcal H}(\Omega_i))\\
&=(\lambda_{\xi}^{\beta,\alpha}(\lambda_{\xi}^{\beta,\alpha})^*
  U_{\mathcal H}(\Omega_i)|U_{\mathcal H}(\Omega_i))
\end{aligned}$$
By \ref{rapide}, it follows that we have, for all $i\in I$:
$$
\begin{aligned}
&\ \quad ||\eta_i\surl{\ _{\alpha}
  \otimes_{\hat{\beta}}}_{\ \ \mu^o} [(\lambda_{\xi}^{\beta,\alpha})^*
  U_{\mathcal H}((\rho_{\eta_i}^{\beta,\alpha})^*\Xi\surl{\ _{\alpha}
  \otimes_{\hat{\beta}}}_{\ \ \mu^o} \Lambda_{\Phi}(a))]||^2\\
&\leq ||R^{\beta,\alpha}(\xi)||^2
  ((\rho_{\eta_i}^{\beta,\alpha})^*\Xi\surl{\ _{\alpha}
  \otimes_{\hat{\beta}}}_{\ \ \mu^o} \Lambda_{\Phi}(a)|
  (\rho_{\eta_i}^{\beta,\alpha})^*\Xi\surl{\ _{\alpha}
  \otimes_{\hat{\beta}}}_{\ \ \mu^o} \Lambda_{\Phi}(a))\\
&\leq ||R^{\beta,\alpha}(\xi)||^2(T_L(a^*a)
  (\rho_{\eta_i}^{\beta,\alpha})^*\Xi
  |(\rho_{\eta_i}^{\beta,\alpha})^*\Xi)\\
&\leq ||R^{\beta,\alpha}(\xi)||^2||T(a^*a)||(
  (\rho_{\eta_i}^{\beta,\alpha})^*\Xi
  |(\rho_{\eta_i}^{\beta,\alpha})^*\Xi)
\end{aligned}$$
So, we can sum over $i\in I$ to get that:
$$\sum_{i\in I}||\eta_i\surl{\ _{\alpha}\otimes_{\hat{\beta}}}_{\ \ \mu^o}
[(\lambda_{\xi}^{\beta,\alpha})^*U_{\mathcal
H}((\rho_{\eta_i}^{\beta,\alpha})^*\Xi\surl{\ _{\alpha}
\otimes_{\hat{\beta}}}_{\ \ \mu^o} \Lambda_{\Phi}(a))]||^2$$ is less
or equal to
$||R^{\beta,\alpha}(\xi)||^2||T(a^*a)||||\Xi||^2<\infty$. That's why
the sum converges.
\end{proof}

\begin{prop}\label{techinter}
Let $a,e\in {\mathcal N}_{\Phi}\cap {\mathcal N}_{T_L}$, $\Xi\in
{\mathcal H}\surl{\ _{\beta}\otimes_{\alpha}}_{\ \mu} H$, $\xi\in
D({\mathcal H}_{\beta},\mu^o)$, $\eta\in D(_{\alpha}(H_{\Phi}),\mu)$
and $(\eta_i)_{i\in I}$ a $(N,\mu)$-basis of $\ _{\alpha}H$. We
have:
$$(\rho_{J_{\Phi}eJ_{\Phi}\eta}^{\beta,\alpha})^*U_H(\sum_{i\in I} \eta_i\surl{\ _{\alpha}
  \otimes_{\hat{\beta}}}_{\ \ \mu^o} [(\lambda_{\xi}^{\beta,\alpha})^*
  U_{\mathcal H}((\rho_{\eta_i}^{\beta,\alpha})^*\Xi\surl{\ _{\alpha}
  \otimes_{\hat{\beta}}}_{\ \ \mu^o} \Lambda_{\Phi}(a))])$$
$$=(\lambda_{\xi}^{\beta,\alpha})^*\Gamma((id\surl{\ _{\beta}
  \star_{\alpha}}_{\ \mu}\omega_{J_{\Phi}
  \Lambda_{\Phi}(e),\eta})(\Gamma(a)))\Xi$$
\end{prop}

\begin{proof}
The existence of the first term comes from the previous lemma. By
\ref{simple} and the co-product relation, we can compute:
$$
\begin{aligned}
&\quad\sum_{i\in I}(\rho_{\eta}^{\beta,\alpha})^*(1\surl{\ _{\beta}
  \otimes_{\alpha}}_{\ N}
J_{\Phi}eJ_{\Phi})U_H( \eta_i\surl{\ _{\alpha}
  \otimes_{\hat{\beta}}}_{\ \ \mu^o} [(\lambda_{\xi}^{\beta,\alpha})^*
  U_{\mathcal H}((\rho_{\eta_i}^{\beta,\alpha})^*\Xi\surl{\ _{\alpha}
  \otimes_{\hat{\beta}}}_{\ \ \mu^o} \Lambda_{\Phi}(a))])\\
&=\sum_{i\in I}(\rho_{\eta}^{\beta,\alpha})^*
  (\lambda_{\xi}^ {\beta,\alpha}\surl{\ _{\beta} \otimes_{\alpha}}_{\
  N}1)^*(id \surl{\ _{\beta}
  \star_{\alpha}}_{\ N} \Gamma)(\Gamma(a))((\rho_{\eta_i}^{\beta,\alpha})^*
  \Xi \surl{\ _{\beta} \otimes_{\alpha}}_{\ \mu}
  \eta_i \surl{\ _{\beta} \otimes_{\alpha}}_{\ \mu}
  J_{\Phi}\Lambda_{\Phi}(e))\\
&=(\rho_{\eta}^{\beta,\alpha})^*
  (\lambda_{\xi}^ {\beta,\alpha}\surl{\ _{\beta} \otimes_{\alpha}}_{\
  N}1)^*(\Gamma\surl{\ _{\beta}
  \star_{\alpha}}_{\ N} id)(\Gamma(a))([\sum_{i\in I}
  \rho_{\eta_i}^{\beta,\alpha}
  (\rho_{\eta_i}^{\beta,\alpha})^*]\Xi\surl{\ _{\beta}
  \otimes_{\alpha}}_{\ \mu} J_{\Phi}\Lambda_{\Phi}(e)))\\
&=(\rho_{\eta}^{\beta,\alpha})^*
  (\lambda_{\xi}^ {\beta,\alpha}\surl{\ _{\beta} \otimes_{\alpha}}_{\
  N}1)^*(\Gamma\surl{\ _{\beta}
  \star_{\alpha}}_{\ N} id)(\Gamma(a))(\Xi\surl{\ _{\beta}
  \otimes_{\alpha}}_{\ \mu} J_{\Phi}\Lambda_{\Phi}(e)))\\
&=(\lambda_{\xi}^ {\beta,\alpha})^*
  (1\surl{\ _{\beta} \otimes_{\alpha}}_{\
  N}\rho_{\eta}^{\beta,\alpha})^*(\Gamma\surl{\ _{\beta}
  \star_{\alpha}}_{\ N} id)(\Gamma(a))(\Xi\surl{\ _{\beta}
  \otimes_{\alpha}}_{\ \mu} J_{\Phi}\Lambda_{\Phi}(e)))\\
&=(\lambda_{\xi}^{\beta,\alpha})^*\Gamma((id\surl{\ _{\beta}
  \star_{\alpha}}_{\ \mu}
  \omega_{J_{\Phi}\Lambda_{\Phi}(e),
  \eta})(\Gamma(a)))\Xi
\end{aligned}$$

\end{proof}

With results of the two last sections in hand, we can prove now a
reciprocity law where ${\mathcal H}$ will be equal to $H_{\Psi}$.

\subsubsection{Reciprocity law}\label{reciprocite}

For all monotone increasing net $(e_k)_{k\in K}$ in ${\mathcal
N}_{\Psi}\cap {\mathcal N}_{T_R}$ of limit equal to $1$, the
following  $(\omega_{J_{\Psi}\Lambda_{\Psi}(e_k)})_{k\in K}$ is
monotone increasing and converges to $\Psi$. So, for all $x\in
{\mathcal N}_{\Psi}\cap {\mathcal N}_{T_R}$,
$(\omega_{J_{\Psi}\Lambda_{\Psi}(e_k)}\surl{\ _{\beta}
\star_{\alpha}}_{\ \mu} id)(\Gamma(x))$ converges to $(\Psi\surl{\
_{\beta}\star_{\alpha}}_{\ \mu} id)(\Gamma(x))$ in the weak
topology. We denote $\zeta_k=J_{\Psi}\Lambda_{\Psi}(e^*_ke_k)\in
D((H_{\Psi})_{\beta},\mu^o)$ for all $k\in K$.

\begin{prop}
For all $a\in ({\mathcal N}_{\Psi}\cap {\mathcal N}_{T_R})^*
({\mathcal N}_{\Phi}\cap {\mathcal N}_{T_L}))$, $e\in {\mathcal
N}_{\Phi}\cap {\mathcal N}_{T_L},b\in {\mathcal N}_{\Psi}\cap
{\mathcal N}_{T_R}, c\in {\mathcal T}_{\Psi,T_R}$, $v \in
D(H_{\beta},\mu^o), \eta \in D(_{\alpha} (H_{\Phi}),\mu)$ and
$(N,\mu)$-basis of $\ _{\alpha}H$, $(\eta_i)_{i\in I}$ , we have
that the image of:
$$\sum_{i\in I} \eta_i\!\!\surl{\ _{\alpha}\otimes_{\hat{\beta}}}_{\ \ \mu^o}
[(\lambda_{\zeta_k}^{\beta,\alpha})^*U_{H_{\Psi}}([(\rho_{\eta_i}^{\beta,\alpha})^*
U'_H(J_{\Psi}c^*J_{\Psi}\Lambda_{\Psi}(b)\!\!\surl{\
_{\hat{\alpha}}\otimes_{\beta}}_{\ \ \mu^o}v)]\!\!\surl{\
_{\alpha}\otimes_{\hat{\beta}}}_{\ \ \mu^o} \Lambda_{\Phi}(a))]$$ by
$(\rho^{\beta,\alpha}_{J_{\Phi}e^*J_{\Phi}\eta})^*U_H$ converges, in
the weak topology, to:
$$(\rho^{\beta,\alpha}_{J_{\Phi}e^*J_{\Phi}\eta})^*
(v\surl{\ _{\beta}\otimes_{\alpha}}_{\
\mu}(\lambda^{\beta,\alpha}_{J_{\Psi}\Lambda_{\Psi}(c)})^*
U_{H_{\Psi}}(\Lambda_{\Psi}(b)\surl{\
_{\alpha}\otimes_{\hat{\beta}}}_{\ \ \mu^o}\Lambda_{\Phi}(a)))$$
\end{prop}

\begin{proof}
Let $u\in H$. We compute the value of the scalar product of:
$$U_H(\sum_{i\in I} \eta_i\surl{\ _{\alpha}
  \otimes_{\hat{\beta}}}_{\ \ \mu^o} [(\lambda_{\zeta_k}^{\beta,\alpha})^*
  U_{H_{\Psi}}([(\rho_{\eta_i}^{\beta,\alpha})^*U'_H(\Lambda_{\Psi}(bc)
  \surl{\ _{\hat{\alpha}}
  \otimes_{\beta}}_{\ \ \mu^o}v)]\surl{\ _{\alpha}
  \otimes_{\hat{\beta}}}_{\ \ \mu^o} \Lambda_{\Phi}(a))])$$ by $u\surl{\ _{\beta}
  \otimes_{\alpha}}_{\ \mu} J_{\Phi}e^*J_{\Phi}\eta$. By \ref{techinter},
we get that it is equal to:
$$(\Gamma((id\surl{\ _{\beta}
  \star_{\alpha}}_{\ \mu}
  \omega_{J_{\Phi}\Lambda_{\Phi}(e),
  \eta})(\Gamma(a)))U'_H(\Lambda_{\Psi}(bc)
  \surl{\ _{\hat{\alpha}}
  \otimes_{\beta}}_{\ \ \mu^o}v)|\zeta_k\surl{\ _{\beta}
  \otimes_{\alpha}}_{\ \mu} u)$$
By the right version of \ref{impl}, this is equal to:
$$(U'_H(\Lambda_{\Psi}((id\surl{\ _{\beta}
  \star_{\alpha}}_{\ \mu}
  \omega_{J_{\Phi}\Lambda_{\Phi}(e),
  \eta})(\Gamma(a))bc)
  \surl{\ _{\hat{\alpha}}
  \otimes_{\beta}}_{\ \ \mu^o}v)|\zeta_k\surl{\ _{\beta}
  \otimes_{\alpha}}_{\ \mu} u)$$
By \ref{raccourci}, we obtain:
$$((\omega_{J_{\Psi}\Lambda_{\Psi}(e_k)}\surl{\ _{\beta}
  \star_{\alpha}}_{\ \mu}id)(\Gamma((id\surl{\ _{\beta}
  \star_{\alpha}}_{\ \mu}
  \omega_{J_{\Phi}\Lambda_{\Phi}(e),
  \eta})(\Gamma(a))bc))v|u)$$ which converges to:
$$((\Psi\surl{\ _{\beta}
  \star_{\alpha}}_{\ \mu} id)(\Gamma((id\surl{\ _{\beta}
  \star_{\alpha}}_{\ \mu}
  \omega_{J_{\Phi}\Lambda_{\Phi}(e),
  \eta})(\Gamma(a)))bc)v|u)$$
Now, by right invariance of $T_R$, \ref{rapide} and \ref{prepa}, we
can compute this last expression:
$$
\begin{aligned}
&\ \quad ((\Psi\surl{\ _{\beta}
  \star_{\alpha}}_{\ \mu} id)(\Gamma((id\surl{\ _{\beta}
  \star_{\alpha}}_{\ \mu}
  \omega_{J_{\Phi}\Lambda_{\Phi}(e),
  \eta})(\Gamma(a)))bc)v|u)\\
&=(T_R((id\surl{\ _{\beta} \star_{\alpha}}_{\ \mu}
\omega_{J_{\Phi}\Lambda_{\Phi}(e),\eta})(\Gamma(a))bc)v|u)\\
&=(\Lambda_{\Psi}(bc) \surl{\ _{\hat{\alpha}}
  \otimes_{\beta}}_{\ \ \mu^o} v|\Lambda_{\Psi}((id\surl{\
_{\beta} \star_{\alpha}}_{\ \mu} \omega_{\eta,
J_{\Phi}\Lambda_{\Phi}(e)})(\Gamma(a^*)))) \surl{\
_{\hat{\alpha}}\otimes_{\beta}}_{\ \ \mu^o} u)\\
&=(v\surl{\ _{\beta}
  \otimes_{\alpha}}_{\ \mu}
  (\lambda^{\beta,\alpha}_{\Lambda_{\Psi}(\sigma_{-i}^{\Psi}(c^*))})^*
  U_{H_{\Psi}}(\Lambda_{\Psi}(b)\surl{\ _{\alpha}\otimes_{\hat{\beta}}}_{\ \ \mu^o}\Lambda_{\Phi}(a))
  |u\surl{\ _{\beta}
  \otimes_{\alpha}}_{\ \mu} J_{\Phi}e^*J_{\Phi}\eta)
\end{aligned}$$
which finishes the proof.
\end{proof}

Let $(\eta_i)_{i\in I}$be a $(N,\mu)$-basis of $\ _{\alpha}H$. For
all finite subset $J$ of $I$, we denote by $P_J$ the projection
$\sum_{i\in J}\theta^{\alpha,\mu}(\eta_i,\eta_i)\in \alpha(N)'$ so
that:
$$\sum_{i\in J}\rho_{\eta_i}^{\beta,\alpha}(\rho_{\eta_i}^{\beta,\alpha})^*=
1\surl{\ _{\beta} \otimes_{\alpha}}_{\ N} P_J$$ For all $e\in
{\mathcal N}_{\Phi}\cap {\mathcal N}_{T_L}$, we also denote by
$P_J^e$:
$$1\surl{\ _{\beta} \otimes_{\alpha}}_{\ N}
J_{\Phi}e^*J_{\Phi}P_JJ_{\Phi}eJ_{\Phi}=\sum_{i\in
J}\rho_{J_{\Phi}e^*J_{\Phi}\eta_i}^{\beta,\alpha}
(\rho_{J_{\Phi}e^*J_{\Phi}\eta_i}^{\beta,\alpha})^*$$

\begin{coro}\label{mef}
For all $a\in ({\mathcal N}_{\Psi}\cap {\mathcal N}_{T_R})^*
({\mathcal N}_{\Phi}\cap {\mathcal N}_{T_L})$, $b\in {\mathcal
N}_{\Psi}\cap {\mathcal N}_{T_R}$, and $c\in {\mathcal
T}_{\Psi,T_R}$, $v\in D(H_{\beta},\mu^o)$, $e\in {\mathcal
N}_{\Phi}\cap {\mathcal N}_{T_L}$ and $J$ finite subset of $I$, we
have:
$$P_J^eU_H(\sum_{i\in I} \eta_i\!\!\surl{\ _{\alpha}
  \otimes_{\hat{\beta}}}_{\ \ \mu^o}\! [(\lambda_{\zeta_k}^{\beta,\alpha})^*
  U_{H_{\Psi}}([(\rho_{\eta_i}^{\beta,\alpha})^*U'_H(J_{\Psi}c^*J_{\Psi}\Lambda_{\Psi}(b)\surl{\ _{\hat{\alpha}}
  \otimes_{\beta}}_{\ \ \mu^o}v)]\!\!\surl{\ _{\alpha}
  \otimes_{\hat{\beta}}}_{\ \ \mu^o}\! \Lambda_{\Phi}(a))])$$
converges, in the weak topology, to:
$$P_J^e(v\surl{\ _{\beta}\otimes_{\alpha}}_{\ \mu}
(\lambda^{\beta,\alpha}_{J_{\Psi}\Lambda_{\Psi}(c)})^*
U_{H_{\Psi}}(\Lambda_{\Psi}(b)\surl{\
_{\alpha}\otimes_{\hat{\beta}}}_{\ \ \mu^o}\Lambda_{\Phi}(a)))$$
\end{coro}

\begin{proof}
We apply to the reciprocity law
$\rho_{J_{\Phi}e^*J_{\Phi}\eta}^{\beta,\alpha}$ which is a
continuous linear operator of $H$ in $H\surl{\ _{\beta}
\otimes_{\alpha}}_{\ \mu}H_{\Phi}$, and also a continuous linear
operator of $H$ with weak topology in $H\surl{\ _{\beta}
\otimes_{\alpha}}_{\ \mu}H_{\Phi}$ with weak topology. Then, we take
finite sums.
\end{proof}

Until the end of the section, we denote by ${\mathcal H}_{\Phi}$ the
closed linear span in $H_{\Phi}$ of
$(\lambda_w^{\beta,\alpha})^*U_{H_{\Psi}}(v\surl{\ _{\alpha}
\otimes_{\hat{\beta}}}_{\ \ \mu^o} \Lambda_{\Phi}(a))$ where $v\in
H_{\Psi}$, $w\in J_{\Psi}\Lambda_{\Psi}({\mathcal N}_{\Psi}\cap
{\mathcal N}_{T_R})$, and $a\in ({\mathcal N}_{\Psi}\cap {\mathcal
N}_{T_R})^*{\mathcal N}_{\Phi}\cap {\mathcal N}_{T_L}$. By the third
relation of lemma \ref{base} (resp. proposition \ref{comm}),
$\alpha$ (resp. $\hat{\beta}$) is a non-degenerated (resp. anti-)
representation of $N$ on ${\mathcal H}_{\Phi}$.

\begin{lemm}
Let $a\in ({\mathcal N}_{\Psi}\cap {\mathcal N}_{T_R})^*({\mathcal
N}_{\Phi}\cap {\mathcal N}_T)$, $b\in {\mathcal N}_{\Psi}\cap
{\mathcal N}_{T_R}$, $c\in {\mathcal T}_{\Psi,T_R}$, $v\in
D(H_{\beta},\mu^o)$ and $(\eta_i)_{i\in I}$ a $(N,\mu)$-basis of $\
_{\alpha}H$. We put, for all $k\in K$:
$$\Xi_k=(\sum_{i\in I} \eta_i\surl{\ _{\alpha}
\otimes_{\hat{\beta}}}_{\ \ \mu^o}
[(\lambda_{\zeta_k}^{\beta,\alpha})^*
U_{H_{\Psi}}([(\rho_{\eta_i}^{\beta,\alpha})^*U'_H(J_{\Psi}c^*J_{\Psi}
\Lambda_{\Psi}(b)\surl{\ _{\hat{\alpha}}\otimes_{\beta}}_{\ \
\mu^o}v)]\surl{\ _{\alpha}\otimes_{\hat{\beta}}}_{\ \ \mu^o}
\Lambda_{\Phi}(a))]$$ Then the net $(\Xi_k)_{k\in K}$ is bounded.
\end{lemm}

\begin{proof}
Let $\Xi=v\surl{\ _{\beta}\otimes_{\alpha}}_{\ \mu}
(\lambda^{\beta,\alpha}_{J_{\Psi}\Lambda_{\Psi}(c)})^*
U_{H_{\Psi}}(\Lambda_{\Psi}(b)\surl{\
_{\alpha}\otimes_{\hat{\beta}}}_{\ \ \mu^o}\Lambda_{\Phi}(a))$. By
the previous corollary, we know that $P_J^eU_H\Xi_k$ weakly
converges to $P_J^e\Xi$, so that:
$$\lim_{J,||e||\leq 1}\lim_k P_J^eU_H\Xi_k =\Xi$$

Consequently, there exists $C\in\mathbb{R}^+$ such that:
$$\sup_{J,||e||\leq 1}\sup_k ||P_J^eU_H\Xi_k||\leq C$$
and, the interversion of the supremum gives:
$$C\geq\sup_k \sup_{J,||e||\leq 1}||P_J^eU_H\Xi_k||=\sup_k||U_H\Xi_k||=\sup_k||\Xi_k||$$
\end{proof}

\begin{coro}\label{precau}
For all $a\in ({\mathcal N}_{\Psi}\cap {\mathcal
N}_{T_R})^*({\mathcal N}_{\Phi}\cap {\mathcal N}_T)$, $b\in
{\mathcal N}_{\Psi}\cap {\mathcal N}_{T_R}$, $c\in {\mathcal
T}_{\Psi,T_R}$, $v\in D(H_{\beta},\mu^o)$ and $(\eta_i)_{i\in I}$ a
$(N,\mu)$-basis of $\ _{\alpha}H$, we put: $$\Xi_k=(\sum_{i\in I}
\eta_i\surl{\ _{\alpha} \otimes_{\hat{\beta}}}_{\ \ \mu^o}
[(\lambda_{\zeta_k}^{\beta,\alpha})^*U_{H_{\Psi}}([(\rho_{\eta_i}^{\beta,\alpha})^*
U'_H(J_{\Psi}c^*J_{\Psi}\Lambda_{\Psi}(b)\surl{\ _{\hat{\alpha}}
\otimes_{\beta}}_{\ \ \mu^o}v)]\surl{\ _{\alpha}
\otimes_{\hat{\beta}}}_{\ \ \mu^o} \Lambda_{\Phi}(a))]$$ for all
$k\in K$, and:
$$\Xi=v\surl{\ _{\beta}\otimes_{\alpha}}_{\ \mu}
(\lambda^{\beta,\alpha}_{J_{\Psi}\Lambda_{\Psi}(c)})^*
U_{H_{\Psi}}(\Lambda_{\Psi}(b)\surl{\
_{\alpha}\otimes_{\hat{\beta}}}_{\ \ \mu^o}\Lambda_{\Phi}(a))$$ Then
$U_H\Xi_k$ converges to $\Xi$ in the weak topology.
\end{coro}

\begin{proof}
Let $\Theta\in H\surl{\ _{\beta} \otimes_{\alpha}}_{\ \mu}H_{\Phi}$
and $\epsilon>0$. Then, there exists $e\in {\mathcal N}_{\Phi}\cap
{\mathcal N}_{T_L}$ of norm less than equal to $1$ and a finite
subset $J$ of $I$ such that $||(1-P_J^e)\Theta||\leq\epsilon$. By
\ref{mef}, there also exists $k_0$ such that
$|(P_J^eU_H\Xi_k-P_J^e\Xi|\Theta)|\leq\epsilon$ for all $k\geq k_0$.
Then, we get:
$$
\begin{aligned}
&\ \quad |(U_H\Xi_k-\Xi|\Theta)|\\
&\leq|(U_H\Xi_k-P_J^eU_H\Xi_k|\Theta)|+|(P_J^eU_H\Xi_k-P_J^e\Xi|\Theta)|+|(P_J^e\Xi-\Xi|\Theta)|\\
&\leq|(U_H\Xi_k|(1-P_J^e)\Theta)|+\epsilon+|(\Xi|(1-P_J^e)\Theta)|\\
&\leq|(U_H\Xi_k|(1-P_J^e)\Theta)|+\epsilon+|(\Xi|(1-P_J^e)\Theta)|
\leq (sup_{k\in k}||\Xi_k||+||\Xi||+1)\epsilon
\end{aligned}$$

\end{proof}

\begin{coro}
We have the following inclusion:
$$H\surl{\ _{\beta}\otimes_{\alpha}}_{\ \mu} {\mathcal H}_{\Phi} \subseteq
U_H(H\surl{\ _{\alpha}\otimes_{\hat{\beta}}}_{\ \ \mu^o} {\mathcal
H}_{\Phi})$$
\end{coro}

\begin{proof}
By the previous corollary, we know that $\Xi$ belongs to the weak
closure of $U_H(H\surl{\ _{\alpha}\otimes_{\hat{\beta}}}_{\ \ \mu^o}
{\mathcal H}_{\Phi})$ which is also the norm closure. Now, $U_H$ is
an isometry, that's why $U_H(H\surl{\
_{\alpha}\otimes_{\hat{\beta}}}_{\ \ \mu^o} {\mathcal H}_{\Phi})$ is
equal to $U_H(H\surl{\ _{\alpha}\otimes_{\hat{\beta}}}_{\ \ \mu^o}
{\mathcal H}_{\Phi})$.
\end{proof}

\begin{theo}
$U_H: H \surl{\ _{\alpha} \otimes_{\hat{\beta}}}_{\ \ \mu^o}
H_{\Phi} \rightarrow H\surl{\ _{\beta} \otimes_{\alpha}}_{\ \mu}
H_{\Phi}$ is a unitary.
\end{theo}

\begin{proof}
By the previous corollary, we have:
\begin{equation}\label{inc}
H\surl{\ _{\beta}\otimes_{\alpha}}_{\ \mu} {\mathcal H}_{\Phi}
\subseteq U_H(H\surl{\ _{\alpha}\otimes_{\hat{\beta}}}_{\ \
\mu^o}{\mathcal H}_{\Phi}) \subseteq U_H(H \surl{\
_{\alpha}\otimes_{\hat{\beta}}}_{\ \ \mu^o} H_{\Phi}) \subseteq H
\surl{\ _{\beta} \otimes_{\alpha}}_{\ \mu} H_{\Phi}.
\end{equation}

Also, using a $(N^o,\mu^o)$-basis, we have, for all $v\in H_{\Psi}$
and $a\in {\mathcal N}_{T_L}\cap {\mathcal N}_{\Phi}$:
$$U_{H_{\Psi}}(v\surl{\ _{\alpha} \otimes_{\hat{\beta}}}_{\ \
\mu^o}\Lambda_{\phi}(a))=\sum_i\xi_i\surl{\ _{\beta}
\otimes_{\alpha}}_{\
\mu}(\lambda_{\xi_i}^{\beta,\alpha})^*U_{H_{\Psi}}(v\surl{\
_{\alpha} \otimes_{\hat{\beta}}}_{\ \ \mu^o}\Lambda_{\phi}(a))$$ so
that $U_{H_{\Psi}}(H_{\Psi}\surl{\ _{\alpha}
\otimes_{\hat{\beta}}}_{\ \ \mu^o}H_{\Phi}) \subseteq
H_{\Psi}\surl{\ _{\beta} \otimes_{\alpha}}_{\ \mu} {\mathcal
H}_{\Phi}$. The reverse inclusion is the relation \eqref{inc}
applied to $H_{\Psi}$. Consequently, we get:
$$U_{H_{\Psi}}(H_{\Psi}\surl{\ _{\alpha}
\otimes_{\hat{\beta}}}_{\ \
\mu^o}H_{\Phi})=U_{H_{\Psi}}(H_{\Psi}\surl{\ _{\alpha}
\otimes_{\hat{\beta}}}_{\ \ \mu^o}{\mathcal H}_{\Phi})$$ Since
$U_{H_{\Psi}}$ is an isometry, $H_{\Psi}\surl{\ _{\alpha}
\otimes_{\hat{\beta}}}_{\ \ \mu^o}H_{\Phi}=H_{\Psi}\surl{\ _{\alpha}
\otimes_{\hat{\beta}}}_{\ \ \mu^o}{\mathcal H}_{\Phi}$ and, so
${\mathcal H}_{\Phi}=H_{\Phi}$. Finally, by inclusion \eqref{inc},
we obtain $U_H(H \surl{\ _{\hat{\alpha}} \otimes_{\beta}}_{\ \
\mu^o} H_{\Phi})=H \surl{\ _{\beta} \otimes_{\alpha}}_{\ \mu}
H_{\Phi}$.
\end{proof}

\begin{defi}
Fundamental isometry $U_H$ is now called (left) \textbf{fundamental
unitary}. Right version $U'_H$ is called \textbf{right fundamental
unitary}.
\end{defi}


\begin{coro}\label{dense} If $[F]$ denote the linear span of a
subset $F$ of a vector space $E$, we have:

$$
\begin{aligned}
H_{\Phi}&=[\Lambda_{\Phi}((\omega_{v,w}\surl{\ _{\beta}
  \otimes_{\alpha}}_{\ \mu}id)(\Gamma(a)))|v,w\in
D(H_{\beta},\mu^o), a\in {\mathcal N}_{\Phi}\cap {\mathcal N}_{T_L}]\\
&=[(\lambda_w^{\beta,\alpha})^*U_H(v\surl{\ _{\alpha}
\otimes_{\hat{\beta}}}_{\ \ \mu^o} \Lambda_{\Phi}(a))|v\in H,w\in
D(H_{\beta},\mu^o), a\in {\mathcal N}_{\Phi}\cap {\mathcal N}_{T_L}]\\
&=[(\omega_{v,w}*id)(U_H)\xi|v\in D(_{\alpha}H,\mu),w\in
D(H_{\beta},\mu^o),\xi\in H_{\Phi}]
\end{aligned}$$
\end{coro}

\begin{proof}
The second equality comes from \ref{lienGV}. The last one is clear.
It's sufficient to prove that the last subspace is equal to
$H_{\Phi}$. Let $\eta\in H_{\Phi}$ in the orthogonal of:
$$[(\omega_{v,w}*id)(U_H)\xi|v\in D(_{\alpha}H,\mu),w\in
D(H_{\beta},\mu^o),\xi\in H_{\Phi}]$$ Then, for all $v\in
D(_{\alpha}H,\mu),w\in D(H_{\beta},\mu^o)$ and $\xi\in H_{\Phi}$, we
have:
$$(U_H(v\surl{\ _{\alpha}
\otimes_{\hat{\beta}}}_{\ \ \mu^o}\xi)|w\surl{\ _{\beta}
\otimes_{\alpha}}_{\ \mu}\eta)=((\omega_{v,w}*id)(U_H)\xi|\eta)=0$$
Since $U_H$ is a unitary, $w\surl{\ _{\beta} \otimes_{\alpha}}_{\
\mu}\eta=0$ for all $w\in D(H_{\beta},\mu^o)$ from which we easily
deduce that $\eta=0$ (by \ref{aidintre} for example).
\end{proof}

\begin{coro}\label{implementation}
We have $\Gamma(m)=U_H(1 \surl{\ _{\alpha} \otimes_{\hat{\beta}}}_{\
\ N^o} m)U^*_H$ for all $m\in M$.
\end{coro}

\begin{proof}
Straightforward thanks to unitary of $U_H$ and \ref{impl}.
\end{proof}

\subsection{Pseudo-multiplicativity}
Let put $W=U^*_{H_{\Phi}}$. We have already proved commutation
relations of section \ref{rcom} and, now the aim is to prove that
$W$ is a pseudo-multiplicative unitary in the sense of M. Enock and
J.M Vallin (\cite{EV}, definition 5.6):

\begin{defi}
We call {\bf pseudo-multiplicative unitary} over $N$ w.r.t
$\alpha,\hat{\beta},\beta$ each unitary $V$ from $H \surl{\ _{\beta}
\otimes_{\alpha}}_{\ \mu} H$ onto $H \surl{\ _{\alpha}
\otimes_{\hat{\beta}}}_{\ \ \mu^o} H$ which satisfies the following
commutation relations, for all $n,m\in N$:
$$(\beta(n)\surl{\ _{\alpha} \otimes_{\hat{\beta}}}_{\ \ N^o}
\alpha(m))V=V(\alpha(m) \surl{\ _{\beta} \otimes_{\alpha}}_{\
N}\hat{\beta}(n))$$ and $$(\hat{\beta}(n)\surl{\ _{\alpha}
\otimes_{\hat{\beta}}}_{\ \ N^o}\beta(m))V= V(\hat{\beta}(n)\surl{\
_{\beta} \otimes_{\alpha}}_{\ N} \beta(m))$$ and the formula:
$$(V \surl{\ _{\alpha}
\otimes_{\hat{\beta}}}_{\ \ N^o} 1)(\sigma_{\mu^o}\surl{\
_{\alpha}\otimes_{\hat{\beta}}}_{\ \ N^o} 1)(1\surl{\ _{\alpha}
\otimes_{\hat{\beta}}}_{\ \ N^o} V)\sigma_{2\mu}(1\surl{\ _{\beta}
\otimes_{\alpha}}_{\ N} \sigma_{\mu^o})(1\surl{\ _{\beta}
\otimes_{\alpha}}_{\ N} V)=$$
$$(1\surl{\ _{\alpha} \otimes_{\hat{\beta}}}_{\ \ N^o} V)
(V\surl{\ _{\beta} \otimes_{\alpha}}_{\ N}1)$$ where the first
$\sigma_{\mu^o}$ is the flip from $H \surl{\ _{\alpha}
\otimes_{\hat{\beta}}}_{\ \ \mu^o} H$ onto $H \surl{\ _{\hat{\beta}}
\otimes_{\alpha}}_{\ \mu} H$, the second is the flip from $H \surl{\
_{\alpha} \otimes_{\beta}}_{\ \ \mu^o} H$ onto $H \surl{\ _{\beta}
\otimes_{\alpha}}_{\ \mu} H$ and $\sigma_{2\mu}$ is the flip from $H
\surl{\ _{\beta} \otimes_{\alpha}}_{\ \mu} H \surl{\ _{\hat{\beta}}
\otimes_{\alpha}}_{\ \mu} H$ onto $H \surl{\ _{\alpha}
\otimes_{\hat{\beta}}}_{\ \ \mu^o} (H \surl{\ _{\beta}
\otimes_{\alpha}}_{\ \mu} H)$. This last flip turns around the
second tensor product. Moreover, parenthesis underline the fact that
the representation acts on the furthest leg.
\end{defi}

We recall, following (\cite{E2}, 3.5), if we use an other n.s.f
weight for the construction of relative tensor product, then
canonical isomorphisms of bimodules change the pseudo-multiplicative
unitary into another pseudo-multiplicative unitary. The pentagonal
relation is essentially the expression of the co-product relation.
So, we compute $(id\surl{\ _{\beta}\star_{\alpha}}_{\
N}\Gamma)\circ\Gamma $ and $(\Gamma\surl{\
_{\beta}\star_{\alpha}}_{\ N}id)\circ\Gamma $ in terms of $U_H$ with
the following propositions \ref{igamma} and \ref{gammai}. Until the
end of the section, ${\mathcal H}$ is an other Hilbert space on
which $M$ acts.

\begin{lemm}
We have, for all $\xi_1\in D(_{\alpha}{\mathcal H},\mu)$ and
$\xi_2'\in D(H_{\beta},\mu^o)$:
$$\lambda_{\xi_1}^{\alpha,\hat{\beta}}(\lambda_{\xi_2'}^{\beta,\alpha})^*
=(\lambda_{\xi_2'}^{\beta,\alpha})^*\sigma_{2\mu^o}(1\surl{\
_{\alpha}\otimes_{\hat{\beta}}}_{\ \
N^o}\sigma_{\mu})\lambda_{\xi_1}^{\alpha,\hat{\beta}}$$ and:
$$U_{\mathcal
H}\lambda_{\xi_1}^{\alpha,\hat{\beta}}(\lambda_{\xi_2'}^{\beta,\alpha})^*U_H
=(\lambda_{\xi_2'}^{\beta,\alpha})^*(1\surl{\ _{\beta}
\otimes_{\alpha}}_{\ \ N}U_{\mathcal H})\sigma_{2\mu^o}(1\surl{\
_{\alpha}\otimes_{\hat{\beta}}}_{\ \ N^o}\sigma_{\mu})(1\surl{\
_{\alpha}\otimes_{\beta}}_{\ \
N^o}U_H)\lambda_{\xi_1}^{\alpha,\beta}$$
\end{lemm}

\begin{proof}
The first equality is easy to verify and the second one comes from
the first one.
\end{proof}

\begin{prop} The two following equations hold:
\begin{enumerate}[i)]\label{calcule}
\item for all $\xi_1\in D(
_{\alpha}{\mathcal H},\mu),\xi'_1\in D( _{\alpha}H,\mu),\xi_2\in
D({\mathcal H}_{\beta},\mu^o),\xi'_2\in D(H_{\beta},\mu^o)$ and
$\eta_1,\eta_2\in H_{\Phi}$, the scalar product of:
$$(1\surl{\ _{\beta}\otimes_{\alpha}}_{\ N}U_{\mathcal H})\sigma_{2\mu^o}(1\surl{\
_{\alpha}\otimes_{\hat{\beta}}}_{\ \ N^o}\sigma_{\mu})(1\surl{\
_{\alpha}\otimes_{\beta}}_{\ \ N^o}U_H)(\sigma_{\mu}\surl{\
_{\alpha}\otimes_{\hat{\beta}}}_{\ \ N^o}1)([\xi'_1\surl{\
_{\beta}\otimes_{\alpha}}_{\ \ \mu}\xi_1]\surl{\
_{\alpha}\otimes_{\hat{\beta}}}_{\ \ \mu^o} \eta_1)$$ by
$\xi'_2\surl{\ _{\beta}\otimes_{\alpha}}_{\ \mu}\xi_2\surl{\
_{\beta}\otimes_{\alpha}}_{\ \mu}\eta_2$ is equal to
$((\omega_{\xi_1,\xi_2}*id)(U_{\mathcal
H})(\omega_{\xi'_1,\xi'_2}*id)(U_H)\eta_1|\eta_2)$.
\item for all $a\in {\mathcal N}_{\Phi}\cap {\mathcal N}_{T_L}$,
$\xi_1\in {\mathcal H}$ and $\xi'_1,\xi'_2\in D(H_{\beta},\mu^o)$,
the value of:
$$(\lambda_{\xi_2'}^{\beta,\alpha})^*(1\surl{\ _{\beta}
\otimes_{\alpha}}_{\ \ N}U_{\mathcal H})\sigma_{2\mu^o}(1\surl{\
_{\alpha} \otimes_{\hat{\beta}}}_{\ \ N^o}\sigma_{\mu})(1\surl{\
_{\alpha}\otimes_{\beta}}_{\ \ N^o}U_H)(\sigma_{\mu}\surl{\
_{\alpha}\otimes_{\hat{\beta}}}_{\ \ N^o}1)$$ on $[\xi'_1\surl{\
_{\beta}\otimes_{\alpha}}_{\ \ \mu}\xi_1]\surl{\ _{\alpha}
\otimes_{\hat{\beta}}}_{\ \ \mu^o}\Lambda_{\Phi}(a)$ is equal to:
$$U_{\mathcal H}(\xi_1\surl{\ _{\alpha}\otimes_{\hat{\beta}}}_{\ \
N^o}\Lambda_{\Phi}((\omega_{\xi'_1,\xi'_2}\surl{\ _{\beta}
\star_{\alpha}}_{\ \mu}id)(\Gamma(a))))$$
\end{enumerate}
\end{prop}

\begin{proof}
By the previous lemma, we can compute the scalar product of i) in
the following way:
$$
\begin{aligned}
&\ \quad((\lambda_{\xi_2'}^{\beta,\alpha})^*(1\surl{\ _{\beta}
  \otimes_{\alpha}}_{\ \ N}U_{\mathcal H})\sigma_{2\mu^o}(1\surl{\ _{\alpha}
  \otimes_{\hat{\beta}}}_{\ \
  N^o}\sigma_{\mu})(1\surl{\ _{\alpha}
  \otimes_{\beta}}_{\ \ N^o}U_H)\lambda_{\xi_1}^{\alpha,\beta}(\xi'_1\surl{\ _{\alpha}
  \otimes_{\hat{\beta}}}_{\ \ \mu^o}\eta_1)|\xi_2\surl{\ _{\beta}
  \otimes_{\alpha}}_{\ \ \mu}\eta_2)\\
&=(U_{\mathcal
H}\lambda_{\xi_1}^{\alpha,\hat{\beta}}(\lambda_{\xi_2'}^{\beta,\alpha})^*U_H
(\xi'_1\surl{\ _{\alpha}
  \otimes_{\hat{\beta}}}_{\ \ \mu^o}\eta_1)|\xi_2\surl{\ _{\beta}
  \otimes_{\alpha}}_{\ \ \mu}\eta_2)\\
&=((\lambda_{\xi_2}^{\beta,\alpha})^*U_{\mathcal H}(\xi_1\surl{\
_{\alpha}
  \otimes_{\hat{\beta}}}_{\ \
  \mu^o}(\omega_{\xi'_1,\xi'_2}*id)(U_H)\eta_1 |\eta_2)\\
&=((\omega_{\xi_1,\xi_2}*id)(U_{\mathcal
H})(\omega_{\xi'_1,\xi'_2}*id)(U_H)\eta_1|\eta_2)
\end{aligned}$$

Also, the second assertion comes from the previous lemma and
\ref{rap}. Let's first assume that $\xi_1\in D(_{\alpha}{\mathcal
H},\mu)$. Then, we compute the vector in demand:
$$
\begin{aligned}
&\ \quad(\lambda_{\xi_2'}^{\beta,\alpha})^*(1\surl{\ _{\beta}
\otimes_{\alpha}}_{\ \ N}U_{\mathcal H})\sigma_{2\mu^o}(1\surl{\
_{\alpha}\otimes_{\hat{\beta}}}_{\ \ N^o}\sigma_{\mu})(1\surl{\
_{\alpha}\otimes_{\beta}}_{\ \
N^o}U_H)\lambda_{\xi_1}^{\alpha,\beta}(\xi'_1\surl{\ _{\alpha}
\otimes_{\hat{\beta}}}_{\ \ \mu^o}\Lambda_{\Phi}(a))\\
&=U_{\mathcal
H}\lambda_{\xi_1}^{\alpha,\hat{\beta}}(\lambda_{\xi_2'}^{\beta,\alpha})^*U_H
(\xi'_1\surl{\ _{\alpha}\otimes_{\hat{\beta}}}_{\ \
\mu^o}\Lambda_{\Phi}(a))\\
&=U_{\mathcal H}(\xi_1\surl{\ _{\alpha}\otimes_{\hat{\beta}}}_{\ \
N^o}\Lambda_{\Phi}((\omega_{\xi'_1,\xi'_2}\surl{\ _{\beta}
\star_{\alpha}}_{\ \mu}id)(\Gamma(a))))
\end{aligned}$$

So, we obtain the expected equality for all $\xi_1\in
D(_{\alpha}{\mathcal H},\mu)$. Since the two expressions are
continuous in $\xi_1$, density of $D(_{\alpha}{\mathcal H},\mu)$ in
${\mathcal H}$ implies that the equality is still true for all
$\xi_1\in {\mathcal H}$.
\end{proof}

\begin{prop}\label{igamma}
For all $a,b\in {\mathcal N}_{\Phi}\cap {\mathcal N}_{T_L}$, we
have:
$$(id\surl{\
_{\beta}\star_{\alpha}}_{\
\mu}\Gamma)(\Gamma(a))\rho_{J_{\Phi}\Lambda_{\Phi}(b)}^{\beta,\alpha}$$
$$=(1\surl{\ _{\beta}\otimes_{\alpha}}_{\ \ N}(1\surl{\
_{\beta}\otimes_{\alpha}}_{\ N}J_{\Phi}bJ_{\Phi})U_{\mathcal
H})\sigma_{2\mu^o}(1\!\!\surl{\ _{\alpha}\otimes_{\hat{\beta}}}_{\ \
N^o}\sigma_{\mu})(1\surl{\ _{\alpha}\otimes_{\beta}}_{\ \
N^o}U_H)(\sigma_{\mu}\surl{\ _{\alpha}\otimes_{\hat{\beta}}}_{\ \
N^o}1)\rho_{\Lambda_{\Phi}(a)}^{\alpha,\hat{\beta}}$$
\end{prop}

\begin{proof}
Let $\xi_1\in {\mathcal H}$ and $\xi'_1,\xi'_2\in
D(H_{\beta},\mu^o)$. We compose the second term of the equality on
the left by $(\lambda_{\xi_2'}^{\beta,\alpha})^*$ and we get:
$$(1\!\surl{\ _{\beta}\otimes_{\alpha}}_{\ N}\! J_{\Phi}bJ_{\Phi})
(\lambda_{\xi_2'}^{\beta,\alpha})^* (1\!\surl{\ _{\beta}
\otimes_{\alpha}}_{\ \ N}U_{\mathcal H})\sigma_{2\mu^o}(1\!\surl{\
_{\alpha} \otimes_{\hat{\beta}}}_{\ \ N^o}\sigma_{\mu})(1\surl{\
_{\alpha}\otimes_{\beta}}_{\ \ N^o}U_H)(\sigma_{\mu}\surl{\
_{\alpha}\otimes_{\hat{\beta}}}_{\ \ N^o}\!
1)\rho_{\Lambda_{\Phi}(a)}^{\alpha,\hat{\beta}}$$ which we evaluate
on $\xi'_1\!\surl{\ _{\beta}\otimes_{\alpha}}_{\ \ \mu}\xi_1$, to
get, by the previous proposition and \ref{raccourci}:
$$
\begin{aligned}
&\ \quad(1\surl{\ _{\beta}\otimes_{\alpha}}_{\
N}J_{\Phi}bJ_{\Phi})U_{\mathcal H}(\xi_1\surl{\
_{\alpha}\otimes_{\hat{\beta}}}_{\ \
N^o}\Lambda_{\Phi}((\omega_{\xi'_1,\xi'_2}\surl{\ _{\beta}
\star_{\alpha}}_{\ \mu}id)(\Gamma(a))))\\
&=\Gamma((\omega_{\xi'_1,\xi'_2}\surl{\ _{\beta}\star_{\alpha}}_{\ \
\mu}id)(\Gamma(a)))\rho_{J_{\Phi}\Lambda_{\Phi}(b)}^{\beta,\alpha}\xi_1\\
&=(\lambda_{\xi_2'}^{\beta,\alpha})^*(id\surl{\
_{\beta}\star_{\alpha}}_{\
\mu}\Gamma)(\Gamma(a))\rho_{J_{\Phi}\Lambda_{\Phi}(b)}^{\beta,\alpha}(\xi'_1\surl{\
_{\beta}\otimes_{\alpha}}_{\ \ \mu}\xi_1)
\end{aligned}$$
So, the proposition holds.
\end{proof}

\begin{lemm}
For all $X\in M\surl{\ _{\beta}\star_{\alpha}}_{\ N}M\subset
(1\surl{\ _{\beta}\otimes_{\alpha}}_{\ N}\hat{\beta}(N))'$, we have:
$$(\Gamma\surl{\ _{\beta}\star_{\alpha}}_{\
N}id)(X)=(U_H\surl{\ _{\beta}\otimes_{\alpha}}_{\ \ N}1)(1\surl{\
_{\alpha}\otimes_{\hat{\beta}}}_{\ \ N^o}X)(U_H^*\surl{\
_{\beta}\otimes_{\alpha}}_{\ \ N}1)$$
\end{lemm}

\begin{proof}
By \ref{implementation}, $\Gamma$ is implemented by $U_H$ so that we
easily deduce the lemma.
\end{proof}

\begin{prop}\label{gammai}
For all $a,b\in {\mathcal N}_{\Phi}\cap {\mathcal N}_{T_L}$, we
have:
$$(\Gamma\surl{\ _{\beta}\star_{\alpha}}_{\ N}id)(\Gamma(a))
\rho_{J_{\Phi}\Lambda_{\Phi}(b)}^{\beta,\alpha}$$
$$=(1\surl{\
_{\beta}\otimes_{\alpha}}_{\ N} 1\surl{\
_{\beta}\otimes_{\alpha}}_{\ N}J_{\Phi}bJ_{\Phi})(U_H\surl{\
_{\beta}\otimes_{\alpha}}_{\ N}1)(1\surl{\
_{\alpha}\otimes_{\hat{\beta}}}_{\ \ N^o}W^*)(U_H^*\surl{\
_{\alpha}\otimes_{\hat{\beta}}}_{\ \
N^o}1)\rho_{\Lambda_{\Phi}(a)}^{\alpha,\hat{\beta}}$$
\end{prop}

\begin{proof}
By the previous lemma and \ref{raccourci}, we can compute:
$$
\begin{aligned}
&\ \quad(1\surl{\ _{\beta}\otimes_{\alpha}}_{\ N} 1\surl{\
_{\beta}\otimes_{\alpha}}_{\ N}J_{\Phi}bJ_{\Phi})(U_H\surl{\
_{\beta}\otimes_{\alpha}}_{\ N}1)(1\surl{\
_{\alpha}\otimes_{\hat{\beta}}}_{\ \ N^o}W^*)(U_H^*\surl{\
_{\alpha}\otimes_{\hat{\beta}}}_{\ \
N^o}1)\rho_{\Lambda_{\Phi}(a)}^{\alpha,\hat{\beta}}\\
&=(U_H\surl{\ _{\beta}\otimes_{\alpha}}_{\ N}1)(1\surl{\
_{\alpha}\otimes_{\hat{\beta}}}_{\ \ N^o} 1\surl{\
_{\beta}\otimes_{\alpha}}_{\ N}J_{\Phi}bJ_{\Phi})(1\surl{\
_{\alpha}\otimes_{\hat{\beta}}}_{\ \
N^o}W^*)\rho_{\Lambda_{\Phi}(a)}^{\alpha,\hat{\beta}}U_H^*\\
&=(U_H\surl{\ _{\beta}\otimes_{\alpha}}_{\ N}1)(1\surl{\
_{\alpha}\otimes_{\hat{\beta}}}_{\ \ N^o}(1\surl{\
_{\beta}\otimes_{\alpha}}_{\
N}J_{\Phi}bJ_{\Phi})W^*\rho_{\Lambda_{\Phi}(a)}^{\alpha,\hat{\beta}})U_H^*\\
&=(U_H\surl{\ _{\beta}\otimes_{\alpha}}_{\ N}1)(1\surl{\
_{\alpha}\otimes_{\hat{\beta}}}_{\ \
N^o}\Gamma(a)\rho_{J_{\Phi}\Lambda_{\Phi}(b)}^{\beta,\alpha})U_H^*\\
&=(U_H\surl{\ _{\beta}\otimes_{\alpha}}_{\ N}1)(1\surl{\
_{\alpha}\otimes_{\hat{\beta}}}_{\ \ N^o}\Gamma(a))(U_H^*\surl{\
_{\alpha}\otimes_{\hat{\beta}}}_{\ \
N^o}1)\rho_{J_{\Phi}\Lambda_{\Phi}(b)}^{\beta,\alpha}=(\Gamma\surl{\
_{\beta}\star_{\alpha}}_{\ N}id)(\Gamma(a))
\rho_{J_{\Phi}\Lambda_{\Phi}(b)}^{\beta,\alpha}
\end{aligned}$$
\end{proof}

\begin{coro}The following relation is satisfied:
$$(U_H^*\!\!\surl{\ _{\alpha}
\otimes_{\hat{\beta}}}_{\ \ N^o} 1)(\sigma_{\mu^o}\!\! \surl{\
_{\alpha}\otimes_{\hat{\beta}}}_{\ \ N^o} 1)(1\!\! \surl{\
_{\alpha}\otimes_{\hat{\beta}}}_{\ \ N^o}U_H^*)\sigma_{2\mu}(1\!\!
\surl{\ _{\beta} \otimes_{\alpha}}_{\ N} \sigma_{\mu^o})(1\!\!
\surl{\ _{\beta} \otimes_{\alpha}}_{\ N}W)$$
$$=(1\!\!\surl{\ _{\alpha}\otimes_{\hat{\beta}}}_{\ \ N^o}W)(U_H^*\!\!\surl{\ _{\beta}
\otimes_{\alpha}}_{\ N} 1)$$
\end{coro}

\begin{proof}
We put together \ref{igamma} (with ${\mathcal H}=H_{\Phi}$) and
\ref{gammai} thanks to the co-product relation. We get:
$$(1 \surl{\ _{\beta} \otimes_{\alpha}}_{\ N} W^*)\sigma_{2\mu^o}(1
\surl{\ _{\alpha} \otimes_{\hat{\beta}}}_{\ \ N^o} \sigma_{\mu})(1
\surl{\ _{\alpha}
  \otimes_{\beta}}_{\ \ N^o}U_H)$$
$$=(U_H\surl{\ _{\beta}
\otimes_{\alpha}}_{\ N} 1)(1\surl{\ _{\alpha}
\otimes_{\hat{\beta}}}_{\ \ N^o} W^*)(U_H^*\surl{\
  _{\alpha} \otimes_{\hat{\beta}}}_{\ \ N^o} 1)(\sigma_{\mu^o}
  \surl{\ _{\alpha}
  \otimes_{\hat{\beta}}}_{\ \ N^o} 1)$$
Take adjoint and we are.
\end{proof}

\begin{theo}
$W$ is a pseudo-multiplicative unitary over $N$ w.r.t
$\alpha,\hat{\beta},\beta $.
\end{theo}

\begin{proof}
$W$ is a unitary from $H_{\Phi}\surl{\ _{\beta} \otimes_{\alpha}}_{\
\mu}H_{\Phi}$ onto $H_{\Phi}\surl{\ _{\alpha}
\otimes_{\hat{\beta}}}_{\ \ \mu^o}H_{\Phi}$ which satisfies the four
required commutation relations. The previous corollary, with
$H=H_{\Phi}$, finishes the proof.
\end{proof}

Similar results hold for the right version:
\begin{theo}
If $W'=U'_{H_{\Psi}}$, then the following relation makes sense and
holds:
$$(W'\surl{\ _{\beta}
\otimes_{\alpha}}_{\ N} 1)(\sigma_{\mu}\surl{\ _{\beta}
\otimes_{\alpha}}_{\ N} 1)(1\surl{\ _{\beta} \otimes_{\alpha}}_{\
N}U'_H)\sigma_{2\mu^o}(1\surl{\ _{\hat{\alpha}} \otimes_{\beta}}_{\
\ N^o}\sigma_{\mu})(1\surl{\ _{\hat{\alpha}} \otimes_{\beta}}_{\ \
N^o}U'_H)$$
$$=(1\surl{\ _{\beta} \otimes_{\alpha}}_{\ N}U'_H)(W'\surl{\
  _{\hat{\alpha}} \otimes_{\beta}}_{\ \ N^o} 1)$$
If $H=H_{\Psi}$, then $W'$ is a pseudo-multiplicative unitary over
$N^o$ w.r.t $\beta,\alpha, \hat{\alpha}$.
\end{theo}

\begin{proof}
For example, it is sufficient to apply the previous results with the
opposite Hopf bimodule.
\end{proof}

\subsection{Right leg
of the fundamental unitary}

In the von Neumann setting of the theory of locally compact quantum
groups, it is well-known (see \cite{KV2}) that we can recover $M$
from the right leg of the fundamental unitary. In this paragraph, we
prove the first result in that direction in our setting.

\begin{defi}
We call $A(U'_H)$ (resp. ${\mathcal A}(U'_H)$) the weak closure in
${\mathcal L}(H)$ of the vector space (resp. von Neumann algebra)
generated by $(\omega_{v,w}* id)(U'_H)$ with $v \in
D(_{\hat{\alpha}}(H_{\Psi}),\mu)$ and $w \in
D((H_{\Psi})_{\beta},\mu^o)$.
\end{defi}

\begin{prop}\label{situa}
$A(U'_H)$ is a non-degenerate involutive algebra i.e
$A(U'_H)={\mathcal A}(U'_H)$ such that:
$$\alpha(N)\cup\beta(N)\subseteq A(U'_H)={\mathcal A}(U'_H)\subseteq
M\subseteq\hat{\alpha}(N)'$$ Moreover, we have:
$$x\in {\mathcal A}(U'_H)'\cap\mathcal{L}(H)\Longleftrightarrow
U'_H(1 \surl{\ _{\hat{\alpha}} \otimes_{\beta}}_{\ \ N^o}
x)=(1\surl{\ _{\beta} \otimes_{\alpha}}_{\ N} x)U'_H$$
\end{prop}

In fact, we will see later that $A(U'_H)={\mathcal A}(U'_H)=M$.

\begin{proof}
The second and third points are obtained in \cite{EV} (theorem 6.1).
As far as the first point is concerned, it comes from \cite{E2}
(proposition 3.6) and \ref{switch} which proves that $A(U'_H)$ is
involutive.
\end{proof}

To summarize the results of this section, we state the following
theorem:

\begin{theo}\label{transition}
Let $(N,M,\alpha,\beta,\Gamma)$ be a Hopf bimodule, $T_L$ (resp.
$T_R$) be a left (resp. right) invariant n.s.f operator-valued
weight. Then, for all n.s.f weight $\mu$ on $N$, if
$\Phi=\mu\circ\alpha^{-1}\circ T_L$, then the application:
$$v \surl{\ _{\alpha} \otimes_{\hat{\beta}}}_{\ \ \mu^o}
\Lambda_{\Phi}(a)\mapsto\sum_{i
    \in I} \xi_{i} \surl{\ _{\beta} \otimes_{\alpha}}_{\ \mu}
  \Lambda_{\Phi} ((\omega_{v,\xi_i} \surl{\
      _{\beta} \star_{\alpha}}_{\ \mu} id)(\Gamma(a)))$$
for all $v\in D((H_{\Phi})_{\beta},\mu^o)$, $a\in {\mathcal
N}_{T_L}\cap {\mathcal N}_{\Phi}$, $(N^o,\mu^o)$-basis
$(\xi_i)_{i\in I}$ of $(H_{\Phi})_{\beta}$ and where
$\hat{\beta}(n)=J_{\Phi}\alpha(n^*)J_{\Phi}$, extends to a unitary
$W$, the adjoint of which $W^*$ is a pseudo-multiplicative unitary
over $N$ w.r.t $\alpha,\hat{\beta},\beta$ from $H_{\Phi}\surl{\
_{\alpha} \otimes_{\hat{\beta}}}_{\ \ \mu^o}H_{\Phi}$ onto
$H_{\Phi}\surl{\ _{\beta} \otimes_{\alpha}}_{\ \mu}H_{\Phi}$.
Moreover, for all $m\in M$, we have: $$\Gamma(m)=W^*(1\surl{\
_{\alpha} \otimes_{\hat{\beta}}}_{\ \ N^o}m)W$$

Also, we have similar results from $T_R$.
\end{theo}

We also add a key relation between $\Gamma$ and the fundamental
unitary proved in corollary \ref{corres}:

\begin{theo}\label{319}
For all $e,x\in {\mathcal N}_{\Phi}\cap {\mathcal N}_{T_L}$ and
$\eta\in D(_{\alpha}H_{\Phi},\mu^o)$, we have: $$(id\surl{\ _{\beta}
\star_{\alpha}}_{\ \mu} \omega_{J_{\Phi}\Lambda_{\Phi}(e),\eta}
)(\Gamma(x))=(id*
\omega_{\Lambda_{\Phi}(x),J_{\Phi}e^*J_{\Phi}\eta})(U_H)$$ Also, for
all $f,y\in {\mathcal N}_{\Psi}\cap {\mathcal N}_{T_R}$ and $\xi\in
D((H_{\Psi})_{\beta},\mu^o)$, we have:
$$(\omega_{J_{\Psi}\Lambda_{\Psi}(f),\xi} \surl{\ _{\beta}
\star_{\alpha}}_{\
\mu}id)(\Gamma(y))=(\omega_{\Lambda_{\Psi}(y),J_{\Psi}f^*J_{\Psi}\xi}*
id)(U'_H)$$
\end{theo}

\part{Measured quantum groupoids}

In this part, we propose a definition for measured quantum groupoids
from which we can develop a full theory that is we construct all
expected natural objects, then we perform a dual structure within
the category and we also get a duality theorem which extends duality
for locally compact quantum groups. Two main ideas are used in this
theory. First of all, we use axioms of Masuda-Nakagami-Woronowicz's
type: we assume the existence of the antipode defined by its polar
decomposition. On the other hand, we introduce a rather weak
condition on the modular group of the invariant operator-valued
weight. From this, we can proceed and we get all known examples as
we will see in the second part.

\section{Definition}
In the following, $(N,M,\alpha,\beta,\Gamma)$ denotes a
Hopf-bimodule. Like in the quantum group case (for example
\cite{KV1} or \cite{MNW}), we assume that there exist a normal
semi-finite and faithful (nsf) left invariant operator-valued weight
$T_L$. We also assume that we have an antipode. Precisely, like in
\cite{MNW}, we require the existence of a co-involution $R$ of $M$
and a scaling operator $\tau$ (deformation operator) which will lead
to polar decomposition of the antipode. Axioms we choose for them
are well known properties at the quantum groups level. They are
quite symmetric, easy to express and adapted to our developments.
They give a link between $R$, $\tau$ and the co-product $\Gamma$.
They stand for strong invariance and relative invariance of the
weight in \cite{MNW}. Finally, we add a modular condition on the
basis coming from inclusions of von Neumann algebras. The idea is
that we have to choose a weight on the basis $N$ to proceed
constructions. That is also the case for usual groupoids (see
\cite{R1}, \cite{V1} and also section \ref{groupoids}).

\begin{defi}
We call $(N,M,\alpha,\beta,\Gamma,T_L,R,\tau,\nu)$ a
\textbf{measured quantum groupoid} if $(N,M,\alpha,\beta,\Gamma)$ is
a Hopf-bimodule equipped with a nsf left invariant operator-valued
weight $T_L$ from $M$ to $\alpha(N)$, a co-involution $R$ of $M$, a
one-parameter group of automorphisms $\tau$ of $M$ and a nsf weight
$\nu$ on $N$ such that, for all $t\in\mathbb R$ and
$a,b\in\mathcal{N}_{T_L}\cap\mathcal{N}_{\Phi}$ :
$$
\begin{aligned}
R((id\surl{\ _{\beta} \star_{\alpha}}_{\
\nu}\omega_{J_{\Phi}\Lambda_{\Phi}(a)})\Gamma(b^*b)) &=(id\surl{\
_{\beta} \star_{\alpha}}_{\
\nu}\omega_{J_{\Phi}\Lambda_{\Phi}(b)})\Gamma(a^*a)\\
\text{ and } \tau_t((id\surl{\ _{\beta} \star_{\alpha}}_{\
\nu}\omega_{J_{\Phi}\Lambda_{\Phi}(a)})\Gamma(b^*b)) &=(id\surl{\
_{\beta} \star_{\alpha}}_{\
\nu}\omega_{J_{\Phi}\Lambda_{\Phi}(\sigma_t^{\Phi}(a))})\Gamma(\sigma_t^{\Phi}(b^*b))
\end{aligned}$$ where $\Phi=\nu\circ\alpha^{-1}\circ T_L$ and such that:
$$\nu\circ\gamma_t=\nu$$ where $\gamma$ is the unique one-parameter group of automorphisms
$\gamma$ of $N$ satisfying for all $n\in N,t\in \mathbb R$:
$$\sigma_t^{T_L}(\beta(n))=\beta(\gamma_t(n))$$
\end{defi}

We recall that the Hopf-bimodule does also admit a nsf right
invariant operator-valued weight $T_R=R\circ T_L\circ R$. The rest
of the section is devoted to develop several points of the
definition and clarify from where $\gamma$ comes from. Thanks to
relation concerning $\tau$, we easily get that:
$$\tau_t\circ\beta=\beta\circ\sigma_t^{\nu}\quad\text{ and }\quad
(\tau_t \surl{\ _{\beta} \star_{\alpha}}_{\ N} \sigma^{\Phi}_t)\circ
\Gamma=\Gamma\circ\sigma^{\Phi}_t$$ for all $n\in N$ and
$t\in\mathbb R$ (For the first one, make $b$ goes to $1$). The first
equality give the behavior $\tau$ should have on the basis. In fact,
it is necessary, if we want to give a meaning to $\tau_t \surl{\
_{\beta} \star_{\alpha}}_{\ N} \sigma^{\Phi}_t$. The last relation
is usual in the theory of locally compact quantum groups. Then, we
can explain how to recover $M$ from $\Gamma$:

\begin{theo}\label{densevn1}
If $<F>^{-\textsc{w}}$ is the weakly closed linear span of $F$ in
$M$, then we have:
$$\begin{aligned}
M&=<(\omega\surl{\ _{\beta} \star_{\alpha}}_{\ \nu} id)(\Gamma(m))\
| \ m\in M,\omega\in M^+_*, k\in {\mathbb
R}^+\text{ s.t } \omega\circ\beta \leq k\nu>^{-\textsc{w}}\\
&=<(id\surl{\ _{\beta}\star_{\alpha}}_{\ \nu} \omega)(\Gamma(m))\ |
\ m\in M, \omega\in M^+_*, k\in {\mathbb
R}^+ \ \text{s.t } \omega\circ\alpha \leq k\nu>^{-\textsc{w}}\\
\end{aligned}$$
\end{theo}

\begin{proof}
Let call $M_R$ the first subspace of $M$ and $M_L$ the second one.
Since $\tau_t(\beta(n))=\beta(\sigma_t^{\nu}(n))$ for all
$t\in\mathbb R$, we have:
$$M_R=<(\omega\circ\tau_t\surl{\ _{\beta} \star_{\alpha}}_{\ \nu}id)(\Gamma(m))
\ | \ m\in M,\omega \in (M_R)^+_* , k\in {\mathbb R}^+ \ \text{s.t }
\omega\circ\beta \leq k\nu>^{-\textsc{w}}$$ Moreover we have
$\sigma_t^{\Phi}((\omega\surl{\ _{\beta} \star_{\alpha}}_{\ \nu}id
)\Gamma(m))=(\omega\circ\tau_t\surl{\ _{\beta} \star_{\alpha}}_{\
\nu}id )\Gamma(\sigma_t^{\Phi}(m))$ so that
$\sigma_t^{\Phi}(M_R)=M_R$ for all $t\in\mathbb R$. On the other
hand, by proposition \ref{semi}, restriction of $\Phi$ to $M_R$ is
semi-finite. By Takesaki's theorem (\cite{St}, theorem 10.1), there
exists a unique normal and faithful conditional expectation $E$ from
$M$ to $M_R$ such that $\Phi(m)=\Phi(E(m))$ for all $m\in M^+$.
Moreover, if $P$ is the orthogonal projection on the closure of
$\Lambda_{\Phi}({\mathcal N}_{\Phi} \cap M_R)$ then $E(m)P=PmP$.

So the range of $P$ contains $\Lambda_{\Phi}((\omega \surl{\
_{\beta} \star_{\alpha}}_{\ \nu}id)\Gamma(x))$ for all $\omega$ and
$x\in {\mathcal N}_{\Phi}$. By proposition \ref{dense} implies that
$P=1$ so that $E$ is the identity and $M=M_R$. Now, it is clear that
$R(M_R)=M_L$ thanks to co-involution property what completes the
proof.
\end{proof}

The theorem enables us to understand that formulas satisfied by $R$
and $\tau$ in the definition are sufficient to determine them. For
example, we can be ensured of the commutation between $R$ and $\tau$
which can be tested on elements of the form $(id\surl{\ _{\beta}
\star_{\alpha}}_{\
\nu}\omega_{J_{\Phi}\Lambda_{\Phi}(a)})\Gamma(b^*b)$. Also, if we
put $\Psi=\nu\circ\beta^{-1}\circ T_R=\Phi\circ R$, we get, for all
$t\in\mathbb R$:
$$\sigma^{\Psi}_t=R\circ\sigma^{\Phi}_{-t}\circ R\quad\text{ and }\quad\tau_t\circ\alpha=\alpha\circ\sigma_t^{\nu}\quad\text{ and }\quad
(\sigma^{\Psi}_t\surl{\ _{\beta} \star_{\alpha}}_{\ N} \tau_{-t}
)\circ \Gamma=\Gamma\circ\sigma^{\Psi}_t$$

Then, we can precise the behavior of $\tau$ with respect to the
Hopf-bimodule structure:

\begin{prop}\label{easy}
We have $\Gamma \circ \tau_t=(\tau_t \surl{\ _{\beta}
\star_{\alpha}}_{\ N}\tau_t)\circ\Gamma$ for all $t\in \mathbb R$.
\end{prop}

\begin{proof}
Because of the behavior of $\tau$ on the basis, it is possible to
define a normal *-automorphism $\tau_t\surl{\ _{\beta}
\star_{\alpha}}_{\ N}\tau_t$ of $M \surl{\ _{\beta}
\star_{\alpha}}_{\ N}M$ which naturally acts for all $t\in\mathbb
R$. By co-product relation, we have for all $t\in\mathbb R$:
$$\begin{aligned}
(id \surl{\ _{\beta} \star_{\alpha}}_{\ \nu} \Gamma)(\sigma_t^{\Psi}
\surl{\ _{\beta} \star_{\alpha}}_{\ \nu} \tau_{-t})\circ\Gamma &=(id
\surl{\ _{\beta} \star_{\alpha}}_{\
\nu} \Gamma)\Gamma \circ \sigma_t^{\Psi}\\
&=(\Gamma \surl{\ _{\beta} \star_{\alpha}}_{\ \nu} id)\Gamma \circ
\sigma_t^{\Psi}=(\Gamma \circ\sigma_t^{\Psi}\surl{\ _{\beta}
\star_{\alpha}}_{\
\nu} \tau_{-t})\Gamma\\
&=(\sigma_t^{\Psi}\surl{\ _{\beta} \star_{\alpha}}_{\ \nu}
\tau_{-t}\surl{\ _{\beta} \star_{\alpha}}_{\ \nu}
\tau_{-t})(\Gamma\surl{\ _{\beta} \star_{\alpha}}_{\ \nu}
id)\Gamma\\
&=(\sigma_t^{\Psi}\surl{\ _{\beta}\star_{\alpha}}_{\ \nu}
[(\tau_{-t}\surl{\ _{\beta}\star_{\alpha}}_{\ \nu}
\tau_{-t})\circ\Gamma])\circ\Gamma
\end{aligned}$$
Consequently, for all $m\in M$, $\omega\in M_*^+$,
$k\in\mathbb{R}^+$ such that $\omega\circ\beta\leq k\nu$, we have:

$$
\begin{aligned}
\Gamma\circ\tau_{-t}\circ((\omega\circ\sigma_t^{\Psi})\surl{\
_{\beta} \star_{\alpha}}_{\ \nu} id)\Gamma&=(\omega\surl{\ _{\beta}
\star_{\alpha}}_{\ \nu} id\surl{\ _{\beta} \star_{\alpha}}_{\ \nu}
id)(\sigma_t^{\Psi}\!\!\surl{\ _{\beta}
\star_{\alpha}}_{\ \nu}\!(\Gamma\circ\tau_{-t}))\circ\Gamma \\
&=(\omega\surl{\ _{\beta}\star_{\alpha}}_{\ \nu} id\surl{\ _{\beta}
\star_{\alpha}}_{\ \nu} id)(\sigma_t^{\Psi}\!\!\surl{\ _{\beta}
\star_{\alpha}}_{\ \nu}\![(\tau_{-t}\surl{\
_{\beta}\star_{\alpha}}_{\ \nu}\tau_{-t})\circ\Gamma])\\
&=[(\tau_{-t}\surl{\ _{\beta}\star_{\alpha}}_{\ \nu}
\tau_{-t})\circ\Gamma]\circ((\omega\circ\sigma_t^{\Psi})\surl{\
_{\beta} \star_{\alpha}}_{\ \nu} id)\Gamma
\end{aligned}$$

The theorem \ref{densevn1} allows us to conclude.
\end{proof}

Then, we get a nice and useful characterization of elements of the
basis thanks to $\Gamma$:

\begin{prop}\label{clef3}
For all $x\in M\cap\alpha(N)'$, we have $\Gamma(x)=1 \!\surl{\
_{\beta} \otimes_{\alpha}}_{\ N} x \Leftrightarrow x\in\beta(N)$.
Also we have, for all $x\in M\cap\beta(N)'$,  $\Gamma(x)=x\! \surl{\
_{\beta} \otimes_{\alpha}}_{\ N} 1 \Leftrightarrow x\in\alpha(N)$.
\end{prop}

\begin{proof}
Let $x\in M\cap\alpha(N)'$ such that $\Gamma(x)=1\surl{\ _{\beta}
\otimes_{\alpha}}_{\ N} x$. For all $n \in \mathbb N$, we define in
the strong topology:
$$x_n=\frac{n}{\sqrt{\pi}}\int\!
exp(-n^2t^2)\sigma_t^{\Psi}(x)\ dt \quad \text{analytic w.r.t }
\sigma^{\Psi},$$ and:
$$y_n=\frac{n}{\sqrt{\pi}}\int\!
exp(-n^2t^2)\tau_{-t}(x)\ dt \quad \text{belongs to } \alpha(N)'.$$
Then we have $\Gamma(x_n)=1\!\surl{\ _{\beta} \otimes_{\alpha}}_{\
N}y_n$. If $d \in ({\mathcal M}_{\Psi} \cap {\mathcal M}_{T_R})^+$,
then, for all $n\in\mathbb N$, we have $dx_n \in {\mathcal M}_{\Psi}
\cap {\mathcal M}_{T_R}$. Let $\omega \in M^+_*$ and $k\in {\mathbb
R}^+$ such that $\omega\circ\alpha\leq k\nu$. By right invariance,
we get:

$$
\begin{aligned}
\omega\circ T_R(dx_n)&=\omega((\Psi \surl{\ _{\beta}
\star_{\alpha}}_{\ \nu} id)(\Gamma(dx_n)))\\
&=\Psi((id \surl{\ _{\beta} \star_{\alpha}}_{\ \nu}
\omega)(\Gamma(dx_n)))=\Psi((id \surl{\ _{\beta}
\star_{\alpha}}_{\ \nu} (y_n\omega))(\Gamma(d)))\\
&=\omega((\Psi \surl{\ _{\beta} \star_{\alpha}}_{\ \nu}
id)(\Gamma(d))y_n)=\omega(T_R(d)y_n)
\end{aligned}$$

Take the limit over $n\in\mathbb{N}$ to obtain $T_R(dx)=T_R(d)x$ for
all $d\in {\mathcal M}_{\Psi} \cap {\mathcal M}_{T_R}$ and, by
semi-finiteness of $T_R$, we conclude that $x$ belongs to
$\beta(N)$. Reverse inclusion comes from axioms. If we apply this
result to the opposite Hopf-bimodule, then we get the second point.
\end{proof}

Finally, we are able to explain existence and uniqueness of $\gamma$
for the definition:

\begin{prop}
There exists a unique one-parameter group of automorphisms $\gamma$
of $N$ such that:
$$\sigma_t^{T_L}(\beta(n))=\beta(\gamma_t(n))$$ for all $n\in N$
and $t\in \mathbb R$.
\end{prop}

\begin{proof}
For all $n\in N$ and $t\in \mathbb R$, we have
$\sigma_t^{\Phi}(\beta(n))$ belongs to $M\cap\alpha(N)'$. Then, we
can compute:
$$
\begin{aligned}
\Gamma\circ\sigma_t^{\Phi}(\beta(n))&=(\tau_t \surl{\ _{\beta}
\star_{\alpha}}_{\ N}
\sigma^{\Phi}_t)\circ\Gamma(\beta(n))\\
&=(\tau_t \surl{\ _{\beta} \star_{\alpha}}_{\ N}
\sigma^{\Phi}_t)(1\surl{\ _{\beta} \otimes_{\alpha}}_{\
N}\beta(n))=1\surl{\ _{\beta} \otimes_{\alpha}}_{\
N}\sigma_t^{\Phi}(\beta(n))
\end{aligned}$$ By the previous
proposition, we deduce that $\sigma_t^{\Phi}(\beta(n))$ belongs to
$\beta(N)$ i.e there exists a unique element $\gamma_t(n)$ in $N$
such that $\sigma_t^{\Phi}(\beta(n))=\beta(\gamma_t(n))$. The rest
of the proof is straightforward.
\end{proof}

In our definition, we ask $\gamma$ to leave invariant $\nu$. Just
before investigating the structure of these objects, we re-formulate
at the Hilbert level relations for $R$ and $\tau$ with $U_H$ (or
$W$) coming from theorem \ref{transition}. Depending on the
situation, we will use one or the other expression.

\begin{prop}\label{PUJ}
Let $I$ be a unitary anti-linear operator which implements $R$ that
is $R(m)=Im^*I$ for all $m\in M$ and $P$ be a strictly positive
operator which implements $\tau$ that is $\tau_t(m)=P^{-it}mP^{it}$
for all $m\in M$ and $t\in\mathbb{R}$. For all $t\in\mathbb R$ and
$v,w\in D(_{\alpha}H_{\Phi},\nu)$, we have :
$$\begin{aligned} R((id*\omega_{J_{\Phi}v,w})(U_H))
&=(id*\omega_{J_{\Phi}w,v})(U_H)\\
\tau_t((id*\omega_{J_{\Phi}v,w})(U_H))
&=(id*\omega_{\Delta_{\Phi}^{-it}J_{\Phi}v,\Delta_{\Phi}^{-it}w})(U_H)
\end{aligned}$$
$$(I\surl{\ _{\alpha}\otimes_{\hat{\beta}}}_{\
\nu^o}J_{\Phi})U_H^*=U_H(I\surl{\ _{\beta}\otimes_{\alpha}}_{\
\nu}J_{\Phi})\text{ and }(P^{it}\surl{\ _{\beta}\otimes_{\alpha}}_{\
\nu}\Delta^{it}_{\Phi})U_H=U_H(P^{it}\surl{\
_{\alpha}\otimes_{\hat{\beta}}}_{\ \nu^o}\Delta^{it}_{\Phi})$$
$$\varsigma_{N^o}\circ(R\surl{\ _{\beta} \star_{\alpha}}_{\ N}
R)\circ\Gamma=\Gamma \circ R\text{ and }(\tau_t \surl{\ _{\beta}
\star_{\alpha}}_{\ N} \sigma^{\Phi}_t)\circ
\Gamma=\Gamma\circ\sigma^{\Phi}_t$$
\end{prop}

\begin{proof}
By theorem \ref{319}, for all $e,x\in {\mathcal N}_{\Phi}\cap
{\mathcal N}_{T_L}$ and $\eta\in D(_{\alpha}H_{\Phi},\mu^o)$, we
recall that:
$$(id\surl{\ _{\beta} \star_{\alpha}}_{\ \nu}
\omega_{J_{\Phi}\Lambda_{\Phi}(e),\eta} )(\Gamma(x))=(id*
\omega_{\Lambda_{\Phi}(x),J_{\Phi}e^*J_{\Phi}\eta})(U_H)$$ Then the
first two equalities are equivalent to formulas of the definition
and we get straightforward the equalities at the Hilbert level. The
last ones come from definition.
\end{proof}

\section{Uniqueness, modulus and scaling operator}
\label{croise} In this section, we obtain results about the modular
theory of the left-invariant operator-valued weight. We construct a
scaling operator and a modulus which link the left invariant
operator-valued weight $T_L$ and the right invariant operator-valued
weight $R\circ T_L\circ R$. We also prove that the modulus is a
co-character. We also establish uniqueness of the invariant
operator-valued weight.

\subsection{Definitions of modulus and scaling operators}

\begin{prop}\label{summer}
For all $t\in\mathbb R$, we have:
\begin{enumerate}
\item $\Gamma\circ\sigma_t^{\Phi}\tau_{-t}=(id\surl{\
_{\beta}\star_{\alpha}}_{\ N}\sigma_t^{\Phi}\tau_{-t})\circ\Gamma$
\item $R\circ T_L\circ
R\circ\sigma_t^{\Phi}\tau_{-t}=\beta\circ\gamma_t\sigma_{-t}^{\nu}\circ\beta^{-1}\circ
R\circ T_L\circ R$
\item $\Phi\circ R\circ\sigma_t^{\Phi}\tau_{-t}=\Phi\circ R$
\end{enumerate}
\end{prop}

\begin{proof}
For all $n\in N$ and $t\in\mathbb R$, we have:
$$\sigma_t^{\Phi}\tau_{-t}(\alpha(n))=\sigma_t^{\Phi}(\alpha(\sigma_{-t}^{\nu}(n)))=\alpha(n)$$
so that we can define $id\surl{\ _{\beta}\star_{\alpha}}_{\
N}\sigma_t^{\Phi}\tau_{-t}$. Then, the first statement comes
straightforward from definition property of $\tau$ and by
proposition \ref{easy}.

By right invariance of $T_R$, we deduce, for all $a\in {\mathcal
M}_{T_R}^+$:

$$\begin{aligned}
T_R\circ\sigma_t^{\Phi}\tau_{-t}(a)&=(\Psi\surl{\ _{\beta}
\star_{\alpha}}_{\ \nu}id)\Gamma(\sigma_t^{\Phi}\tau_{-t}(a))\\
&=\sigma_t^{\Phi}\tau_{-t}((\Psi\surl{\ _{\beta} \star_{\alpha}}_{\
\nu}id)\Gamma(a))=\sigma_t^{\Phi}\tau_{-t}\circ T_R(a)
\end{aligned}$$ Then, by hypothesis on $\tau$ and $T_L$, we get:
$$T_R\circ\sigma_t^{\Phi}\tau_{-t}=\sigma_t^{\Phi}\tau_{-t}\circ\beta\circ\beta^{-1}\circ T_R
=\sigma_t^{\Phi}\circ\beta\circ\sigma_{-t}^{\nu}\circ\beta^{-1}\circ
T_R=\beta\circ\gamma_t\sigma_{-t}^{\nu}\circ\beta^{-1}\circ T_R$$ To
conclude we just have to take $\nu\circ\beta^{-1}$ on the previous
relation and use invariance property of $\sigma^{\nu}$ and $\gamma$
w.r.t $\nu$.
\end{proof}

\begin{prop}
The one-parameter groups of automorphisms $\sigma^{\Phi}$ and $\tau$
(resp. $\sigma^{\Psi}$ and $\tau$) commute each other.
\end{prop}

\begin{proof}
We put $\kappa_t=\gamma_t\sigma_{-t}^{\nu}$. Since $\Psi$ is
$\kappa$-invariant, we have
$\sigma^{\Psi}_s\circ\sigma_t^{\Phi}\circ\tau_{-t}
=\sigma_t^{\Phi}\circ\tau_{-t}\circ\sigma_s^{\Psi}$, for all
$s,t\in\mathbb{R}$ so that:
$$\begin{aligned}
(id\surl{\ _{\beta} \star_{\alpha}}_{\
N}\kappa_t)\Gamma=\Gamma\circ\kappa_t
&=\Gamma\circ\sigma^{\Psi}_{-s}\circ\kappa_t\circ\sigma^{\Psi}_s
=(\sigma^{\Psi}_{-s}\surl{\ _{\beta} \star_{\alpha}}_{\
N}\tau_s)\circ\Gamma\circ\kappa_t\circ\sigma^{\Psi}_s\\
&=(\sigma^{\Psi}_{-s}\!\!\surl{\ _{\beta} \star_{\alpha}}_{\
N}\tau_s\circ\kappa_t)\circ\Gamma\circ\sigma^{\Psi}_s=(id\!\!\surl{\
_{\beta} \star_{\alpha}}_{\
N}\tau_s\circ\kappa_t\circ\tau_{-s})\circ\Gamma
\end{aligned}$$

So, for all $a\in M$, $\omega\in M_*^+$ and $k\in\mathbb{R}^+$ such
that $\omega\circ\beta\leq k\nu$, we get:
$$\sigma_t^{\Phi}\circ\tau_{-t}((\omega\surl{\ _{\beta} \star_{\alpha}}_{\ \nu}id)\Gamma(a))=
\tau_s\circ\sigma_t^{\Phi}\circ\tau_{-t}\circ\tau_{-s}((\omega\surl{\
_{\beta} \star_{\alpha}}_{\ \nu}id)\Gamma(a))$$ and by theorem
\ref{densevn1}, we easily obtain commutation between $\sigma^{\Phi}$
and $\tau$. By applying the co-involution $R$ to this commutation
relation, we end the proof.
\end{proof}

\begin{coro}\label{mod}
The one-parameter groups of automorphisms $\sigma^{\Phi}$ and
$\sigma^{\Psi}$ commute each other.
\end{coro}

\begin{proof}
By the previous proposition, we compute, for all $s,t\in\mathbb{R}$:
$$
\begin{aligned}
\Gamma\circ\sigma^{\Phi}_s\circ\sigma^{\Psi}_t=(\tau_s \surl{\
_{\beta} \star_{\alpha}}_{\ N} \sigma^{\Phi}_s)\circ
\Gamma\circ\sigma^{\Psi}_t &=(\tau_s\sigma_t^{\Psi} \surl{\ _{\beta}
\star_{\alpha}}_{\ N}
\sigma^{\Phi}_s\tau_{-t})\circ\Gamma\\
&=(\sigma_t^{\Psi}\tau_s \surl{\ _{\beta} \star_{\alpha}}_{\ N}
\tau_{-t}\sigma^{\Phi}_s)\circ\Gamma\\
&=(\sigma_t^{\Psi}\surl{\ _{\beta} \star_{\alpha}}_{\ N}
\tau_{-t})\circ\Gamma\circ\sigma^{\Phi}_s=\Gamma\circ\sigma^{\Psi}_t\circ\sigma^{\Phi}_s\\
\end{aligned}$$

Since $\Gamma$ is injective, we have done.
\end{proof}

By the previous proposition and by \cite{Vae} (proposition 2.5),
there exist a strictly positive operator $\delta$ affiliated with
$M$ and a strictly positive operator $\lambda$ affiliated to the
center of $M$ such that, for all $t\in\mathbb{R}$, we have
$[D\Phi\circ R :D\Phi]_t=\lambda^{\frac{1}{2}it^2}\delta^{it}$.
Modular groups of $\Phi$ and $\Phi\circ R$ are linked by
$\sigma_t^{\Phi\circ
R}(m)=\delta^{it}\sigma_t^{\Phi}(m)\delta^{-it}$ for all
$t\in\mathbb{R}$ and $m\in M$.

\begin{defi}
We call \textbf{scaling operator} the strictly positive operator
$\lambda$ affiliated to $Z(M)$ and \textbf{modulus} the strictly
positive operator $\delta$ affiliated to $M$ such that, for all
$t\in\mathbb{R}$, we have:
$$[D\Phi\circ R:D\Phi]_t=\lambda^{\frac{1}{2}it^2}\delta^{it}$$
\end{defi}

The following propositions give the compatibility of $\lambda$ and
$\delta$ w.r.t the structure of Hopf-bimodule.

\begin{lemm}\label{prep1}
For all $s,t\in\mathbb{R}$, we have $[D\Phi\circ\sigma_s^{\Phi\circ
R}:D\Phi]_t=\lambda^{ist}$.
\end{lemm}

\begin{proof}
The computation of the cocycle is straightforward:
$$
\begin{aligned}
\ [D\Phi\circ\sigma_s^{\Phi\circ R}:D\Phi]_t
&=[D\Phi\circ\sigma_s^{\Phi\circ R}:D\Phi\circ
R\circ\sigma_s^{\Phi\circ R}]_t
[D\Phi\circ R:D\Phi]_t\\
&=\sigma_{-s}^{\Phi\circ R}([D\Phi:D\Phi\circ R]_t)[D\Phi\circ R:D\Phi]_t\\
&=\delta^{-is}\sigma_{-s}^{\Phi}(\lambda^{-\frac{it^2}{2}}\delta^{-it})
\delta^{is}\lambda^{\frac{it^2}{2}}\delta^{it}\\
&=\delta^{-is}\lambda^{-\frac{it^2}{2}}\lambda^{ist}
\delta^{-it}\delta^{is}\lambda^{\frac{it^2}{2}}\delta^{it}=\lambda^{ist}
\end{aligned}$$
\end{proof}

\begin{prop}
We have $R(\lambda)=\lambda$, $R(\delta)=\delta^{-1}$ and
$\tau_t(\delta)=\delta$, $\tau_t(\lambda)=\lambda$ for all
$t\in\mathbb{R}$.
\end{prop}

\begin{proof}
Relations between $R$, $\lambda$ and $\delta$ come from uniqueness
of Radon-Nikodym cocycle decomposition. By proposition \ref{summer},
we have $\Phi\circ\tau_{-s}=\Phi\circ\sigma_s^{\Phi\circ R}$ for all
$s,t\in\mathbb{R}$, so:
$$\tau_s([D\Phi\circ R:D\Phi]_t)=[D\Phi\circ
R\circ\tau_{-s}:D\Phi\circ\tau_{-s}]_t=[D\Phi\circ\sigma_s^{\Phi\circ
R}\circ R:D\Phi\circ\sigma_s^{\Phi\circ R}]_t$$ Consequently, by the
previous lemma, we get:
$$
\begin{aligned}
&\ \quad\tau_s([D\Phi\circ R:D\Phi]_t)\\
&=[D\Phi\circ\sigma_s^{\Phi\circ R}\circ
R:D\Phi\circ R]_t[D\Phi\circ R:D\Phi]_t[D\Phi:D\Phi\circ\sigma_s^{\Phi\circ R}]_t\\
&=R([D\Phi\circ\sigma_s^{\Phi\circ R}:D\Phi]^*_{-t})
[D\Phi\circ R:D\Phi]_t[D\Phi\circ\sigma_s^{\Phi\circ R}:D\Phi]_t^*\\
&=R(\lambda^{ist})\lambda^{-\frac{it^2}{2}}\delta^{it}\lambda^{-ist}
=\lambda^{-\frac{it^2}{2}}\delta^{it}
\end{aligned}$$
\end{proof}

\subsection{First result of uniqueness for invariant operator-valued
weight}

Next, we want to precise where the scaling operator $\lambda$ sits.
We have to prove, first of all, a first result of uniqueness as far
as the invariant operator-valued weight is concerned.

Let $T_1$ and $T_2$ be two n.s.f left invariant operator-valued
weights from $M$ to $\alpha(N)$ such that $T_1\leq T_2$. For all
$i\in\{1,2\}$, we put $\Phi_i=\nu\circ\alpha^{-1}\circ T_i$ and
$\hat{\beta}_i(n)=J_{\Phi_i}\alpha(n^*)J_{\Phi_i}$.

We define, as we have done for $U_H$, an isometry $(U_2)_H$ by the
following formula:
$$(U_2)_H(v\surl{\ _{\alpha} \otimes_{\hat{\beta}_2}}_{\
\nu^o}\Lambda_{\Phi_2}(a))=\sum_{i\in I}\xi_i \surl{\ _{\beta}
\otimes_{\alpha}}_{\ \nu}\Lambda_{\Phi_2}((\omega_{v,\xi_i}\surl{\
_{\beta} \star_{\alpha}}_{\ \nu}id)(\Gamma(a)))$$ for all $v\in
D(H_{\beta},\nu^o)$ and $a\in{\mathcal N}_{\Phi_2}\cap {\mathcal
N}_{T_2}$. Then, we know that $(U_2)_H$ is unitary and
$\Gamma(m)=(U_2)_H(1\surl{\ _{\alpha} \otimes_{\hat{\beta}_2}}_{\ \
N^o}m)(U_2)_H^*$ for all $m\in M$.

Since $T_1\leq T_2$, there exists $F\in
\mathcal{L}(H_{\Phi_2},H_{\Phi_1})$ such that, for all $x\in
{\mathcal N}_{\Phi_2}\cap {\mathcal N}_{T_2}$, we have
$F\Lambda_{\Phi_2}(x)=\Lambda_{\Phi_1}(x)$. It is easy to verify
that, for all $n\in N$, we have $F\hat{\beta}_2(n)=
\hat{\beta}_1(n)F$. If we put $P=F^*F$, then $P$ belongs to $M'\cap
\hat{\beta}_2(N)'$ and $J_{\Phi_2}PJ_{\Phi_2}$ belongs to
$M\cap\alpha(N)'$.

\begin{lemm}
We have $\Gamma(J_{\Phi_2}PJ_{\Phi_2})=1\surl{\ _{\beta}
\otimes_{\alpha}}_{\ N}J_{\Phi_2}PJ_{\Phi_2}$.
\end{lemm}

\begin{proof}
We have, for all $v,w\in D(H_{\beta},\nu^o)$ and $a,b\in{\mathcal
N}_{\Phi_2}\cap {\mathcal N}_{T_2}$:

$$
\begin{aligned}
&\quad((1\surl{\ _{\beta} \otimes_{\alpha}}_{\ N}P)(U_2)_H(v\surl{\
_{\alpha} \otimes_{\hat{\beta}_2}}_{\
\nu^o}\Lambda_{\Phi_2}(a))|(U_2)_H(w\surl{\ _{\alpha}
\otimes_{\hat{\beta}_2}}_{\ \nu^o}\Lambda_{\Phi_2}(b)))\\
&=((U_1)_H(v\surl{\ _{\alpha} \otimes_{\hat{\beta}_1}}_{\
\nu^o}\Lambda_{\Phi_1}(a))|(U_1)_H(w\surl{\ _{\alpha}
\otimes_{\hat{\beta}_1}}_{\ \nu^o}\Lambda_{\Phi_1}(b)))
\end{aligned}$$
where $(U_1)_H$ is defined in the same way as $(U_2)_H$. The two
expressions are continuous in $v$ and $w$, so by density of
$D(H_{\beta},\nu^o)$ in $H$, we get, for all $v,w\in H$ and
$a,b\in{\mathcal N}_{\Phi_2}\cap {\mathcal N}_{T_2}$:

$$
\begin{aligned}
&\quad((1\surl{\ _{\beta} \otimes_{\alpha}}_{\ N}P)(U_2)_H(v\surl{\
_{\alpha} \otimes_{\hat{\beta}_2}}_{\
\nu^o}\Lambda_{\Phi_2}(a))|(U_2)_H(w\surl{\ _{\alpha}
\otimes_{\hat{\beta}_2}}_{\ \nu^o}\Lambda_{\Phi_2}(b)))\\
&=((U_1)_H(v\surl{\ _{\alpha} \otimes_{\hat{\beta}_1}}_{\
\nu^o}\Lambda_{\Phi_1}(a))|(U_1)_H(w\surl{\ _{\alpha}
\otimes_{\hat{\beta}_1}}_{\ \nu^o}\Lambda_{\Phi_1}(b)))\\
&=(v\surl{\ _{\alpha} \otimes_{\hat{\beta}_1}}_{\
\nu^o}\Lambda_{\Phi_1}(a)|w\surl{\ _{\alpha}
\otimes_{\hat{\beta}_1}}_{\ \nu^o}\Lambda_{\Phi_1}(b))\\
&=((1\surl{\ _{\alpha} \otimes_{\hat{\beta}_2}}_{\ N^o}P)(v\surl{\
_{\alpha} \otimes_{\hat{\beta}_2}}_{\
\nu^o}\Lambda_{\Phi_2}(a))|w\surl{\ _{\alpha}
\otimes_{\hat{\beta}_2}}_{\ \nu^o}\Lambda_{\Phi_2}(b))
\end{aligned}$$
so that $(U_2)_H^*(1\surl{\ _{\beta} \otimes_{\alpha}}_{\
N}P)(U_2)_H=1\surl{\ _{\alpha} \otimes_{\hat{\beta}_2}}_{\ N^o}P$.
In particular, if we take $H=H_{\Phi}$, then by \ref{PUJ} we get
$(U_2)_H(1\surl{\ _{\alpha} \otimes_{\hat{\beta}_2}}_{\
N^o}J_{\Phi_2}PJ_{\Phi_2})(U_2)_H^*=1\surl{\ _{\beta}
\otimes_{\alpha}}_{\ N}J_{\Phi_2}PJ_{\Phi_2}$. Finally, since
$J_{\Phi_2}PJ_{\Phi_2}\in M$, we have
$\Gamma(J_{\Phi_2}PJ_{\Phi_2})=1\surl{\ _{\beta}
\otimes_{\alpha}}_{\ N}J_{\Phi_2}PJ_{\Phi_2}$.
\end{proof}

\begin{prop}
If $T_1$ and $T_2$ are n.s.f left invariant weights from $M$ to
$\alpha(N)$ such that $T_1\leq T_2$, then there exists an injective
$p\in N$ such that $0\leq p\leq 1$ and, for all $t\in\mathbb{R}$:
$$[D\Phi_1:D\Phi_2]_t=\beta(p)^{it}$$
\end{prop}

\begin{proof}
By the previous lemma and proposition \ref{clef2}, there exists an
injective $p\in N$ such that $0\leq p\leq 1$ and, for all $x,y\in
{\mathcal N}_{\Phi_2}\cap {\mathcal N}_{T_2}$, we have
$(\Lambda_{\Phi_1}(x)|\Lambda_{\Phi_1}(y))=(J_{\Phi_2}\beta(p)J_{\Phi_2}\Lambda_{\Phi_2}(x)|
\Lambda_{\Phi_2}(y))$. By \cite{St} (proposition 3.13), we get that
$\beta(p)$ coincides with the analytic continuation in $-i$ of the
cocycle $[D\Phi_1:D\Phi_2]$. Then, we have, for all
$t\in\mathbb{R}$:
$$[D\Phi_1:D\Phi_2]_t=\beta(p)^{it}$$
\end{proof}

\begin{prop}\label{struc}
Let $T_1$ be a n.s.f left invariant operator-valued weight $\Phi_1$
is $\sigma^{\Phi}$-invariant. Then, there exists a strictly positive
operator $q$ which is affiliated to $N^{\gamma}$ such that
$\Phi_1=(\Phi)_{\beta(q)}$.
\end{prop}

\begin{proof}
We put $T_2=T_L+T_1$. Since $\Phi_1$ is $\sigma^{\Phi}$-invariant,
then the left invariant operator-valued weight $T_2$ is n.s.f.
Finally, since $T_1\leq T_2$ and $T_L\leq T_2$, there exists an
injective $p\in N$ between $0$ and $1$ such that
$\Phi_1=(\Phi_2)_{\beta(p)}$ and $\Phi=(\Phi_2)_{\beta(1-p)}$. By
\cite{St}, we have:
$$[D\Phi_1:D\Phi_2]_t=\beta(p)^{it} \text{ and } [D\Phi:D\Phi_2]_t=\beta(1-p)^{it}$$ Then, we have, for all
$t\in\mathbb{R}$:
$$[D\Phi_1:D\Phi]_t=[D\Phi_1:D\Phi_2]_t[D\Phi_2
:D\Phi]_t=\beta(\frac{p}{1-p})^{it}$$ that's why $q=\frac{p}{1-p}$
is the suitable element. Now, by \cite{St}, we have:
$$\beta(q)=\sigma_t^{\Phi}(\beta(q))=\beta(\gamma_t(q))$$ so that,
by injectivity of $\beta$, we get that $q$ is affiliated to
$N^{\gamma}$.
\end{proof}

\begin{lemm}
For all $t\in\mathbb{R}$, $\tau_{-t}\circ T_L\circ\tau_t$ is a n.s.f
left invariant operator-valued weight from $M$ to $\alpha(N)$.
Moreover, $\sigma_s^{\Phi\circ\tau_t}(\beta(n))=\beta(\gamma_s(n))$
for all $s,t\in\mathbb{R}$ and $n\in N$.
\end{lemm}

\begin{proof}
For all $t\in\mathbb{R}$, we have
$\nu\circ\alpha^{-1}\circ\tau_{-t}\circ
T_L\circ\tau_t=\Phi\circ\tau_t$. Then:
$$
\begin{aligned}
(id\surl{\ _{\beta} \star_{\alpha}}_{\
\nu}\nu\circ\alpha^{-1}\circ\tau_{-t}\circ
T_L\circ\tau_t)\circ\Gamma&=(id\surl{\ _{\beta} \star_{\alpha}}_{\
\nu}\Phi\circ\tau_t)\circ\Gamma\\
&=\tau_{-t}\circ(id\surl{\ _{\beta} \star_{\alpha}}_{\
\nu}\Phi)\circ\Gamma\circ\tau_t=\tau_{-t}\circ T_L\circ\tau_t
\end{aligned}$$

On the other hand, for all $s,t\in\mathbb{R}$ and $n\in N$, since
$\gamma$ and $\sigma^{\nu}$ commute, we have:
$$
\begin{aligned}
\sigma_s^{\Phi\circ\tau_t}(\beta(n))&=\tau_{-t}\circ\sigma_s^{\Phi}\circ\tau_t(\beta(n))
=\tau_{-t}\circ\sigma_s^{\Phi}(\beta(\sigma_t^{\nu}(n)))\\
&=\tau_{-t}(\beta(\gamma_s\sigma_t^{\nu}(n)))=\beta(\sigma_{-t}^{\nu}\gamma_s\sigma_t^{\nu}(n))=\beta(\gamma_s^{\nu}(n))
\end{aligned}$$
\end{proof}

\begin{prop}
There exists a strictly positive operator $q$ affiliated with $Z(N)$
such that the scaling operator $\lambda=\alpha(q)=\beta(q)$. In
particular, $\lambda$ is affiliated with
$Z(M)\cap\alpha(N)\cap\beta(N)$.
\end{prop}

\begin{proof}
By the previous lemma, $\tau_s\circ T_L\circ\tau_{-s}$ is left
invariant. Moreover, since $\sigma^{\Phi}$ and $\tau$ commute,
$\Phi\circ\tau_{-s}$ is $\sigma^{\Phi}$-invariant. That's why, we
are in conditions of proposition \ref{struc} so that we get a
strictly positive operator $q_s$ affiliated with $N^{\gamma}$ such
that $[D\Phi\circ\tau_{-s}:D\Phi]_t=\beta(q_s)^{it}$. On the other
hand, by lemma \ref{prep1}, we have $[D\Phi\circ\sigma_s^{\Phi\circ
R}:D\Phi]_t=\lambda^{ist}$. Since we have
$\Phi\circ\tau_{-s}=\Phi\circ\sigma_s^{\Phi\circ R}$, so we obtain
that $\lambda^{ist}=\beta(q_s)^{it}$ for all $s,t\in\mathbb{R}$. We
easily deduce that there exists a strictly positive operator $q$
affiliated with $Z(N)$ such that $\lambda=\beta(q)$. Finally, since
$R(\lambda)=\lambda$, we also have $\lambda=\alpha(q)$.
\end{proof}

\subsection{Properties of the modulus}
Now, we prove that the modulus $\delta$ is a co-character. This will
be a key-result for duality.

\begin{prop}
For all $n\in N$ and $t\in\mathbb R$, we have:
$$\delta^{it}\alpha(n)\delta^{-it}=\alpha(\gamma_t\sigma_t^{\nu}(n))\quad\text{
and
}\quad\delta^{it}\beta(n)\delta^{-it}=\beta(\gamma_t\sigma_t^{\nu}(n))$$
\end{prop}

\begin{proof}
By definition of $\gamma$, we have:
$$\beta(\sigma_{-t}^{\nu}(n))=\sigma_t^{\Psi}(\beta(n))=\delta^{it}\sigma_t^{\Phi}(\beta(n))\delta^{-it}
=\delta^{it}\beta(\gamma_t(n))\delta^{-it}$$ what gives the first
equality (we recall that $\gamma$ and $\sigma^{\nu}$ commute with
each other). Then, apply the co-involution to get the second one.
\end{proof}

Thanks to the commutation relations and by proposition
\ref{tenscomp}, we can define, for all $t\in \mathbb R$, a bounded
operator $\delta^{it}\surl{\ _{\beta} \otimes_{\alpha}}_{\ N}
\delta^{it}$ which naturally acts on elementary tensor products.

\begin{lemm}
There exists a strictly positive operator $P$ on $H_{\Phi}$
implementing $\tau$ such that, for all $x\in {\mathcal N}_{\Phi}$
and $t\in\mathbb{R}$, we have
$P^{it}\Lambda_{\Phi}(x)=\lambda^{\frac{t}{2}}\Lambda_{\Phi}(\tau_t(x))$.
\end{lemm}

\begin{proof}
Since $\Phi\circ R=\Phi_{\delta}$, by \cite{Vae} (5.3), we have:
$$\Lambda_{\Phi}(\sigma_t^{\Phi\circ
R}(x))=\delta^{it}J_{\Phi}
\lambda^{\frac{t}{2}}\delta^{it}J_{\Phi}\Delta^{it}_{\Phi}\Lambda_{\Phi}(x)$$
and since $\lambda$ is affiliated with $Z(M)$, we get
$||\Lambda_{\Phi}(\sigma_t^{\Phi\circ
R}(x))||=||\lambda^{\frac{t}{2}}\Lambda_{\Phi}(x)||$ for all $x\in
{\mathcal N}_{\Phi}\cap {\mathcal N}_{T_L}$ and $t\in\mathbb{R}$.
But, we know that $\Phi$ is $\sigma_t^{\Phi\circ
R}\circ\tau_t$-invariant, so
$||\Lambda_{\Phi}(x)||=||\lambda^{\frac{t}{2}}\Lambda_{\Phi}(\tau_t(x))||$.
Then, there exists $P_t$ on $H_{\Phi}$ such that:
$$P_t\Lambda_{\Phi}(x)=\lambda^{\frac{t}{2}}\Lambda_{\Phi}(\tau_t(x))$$
for all $x\in {\mathcal N}_{\Phi}\cap {\mathcal N}_{T_L}$ and
$t\in\mathbb{R}$. For all $s,t\in\mathbb{R}$, we verify that
$P_sP_t=P_{st}$ thanks to relation $\tau_t(\lambda)=\lambda$ and the
existence of $P$ follows. The fact that $P$ implements $\tau$ is
clear.
\end{proof}

\begin{lemm}
We have, for all $a,b\in {\mathcal N}_{\Phi}\cap {\mathcal N}_{T_L}$
and  $t\in\mathbb{R}$:
$$\omega_{J_{\Phi}\Lambda_{\Phi}(\lambda^{\frac{t}{2}}\tau_t(a))}
=\omega_{J_{\Phi}\Lambda_{\Phi}(a)}\circ\tau_{-t}\text{ and }
\omega_{J_{\Phi}\Lambda_{\Phi}(b)}\circ\sigma_t^{\Phi\circ R}=
\omega_{J_{\Phi}\Lambda_{\Phi}(\lambda^{\frac{t}{2}}\sigma_{-t}^{\Phi\circ
R}(b))}$$
\end{lemm}

\begin{proof}
Since $\tau$ is implemented by $P$, the first relation holds. By
\cite{Vae} (proposition 2.4), we know that $\Delta_{\Phi\circ
R}=J_{\Phi}\delta J_{\Phi}\delta\Delta_{\Phi}$ so that we can
compute, for all $x\in M$ and
$b\in\mathcal{N}_{\Phi}\cap\mathcal{N}_{T_L}$:
$$
\begin{aligned}
(\sigma_t^{\Phi\circ
R}(x)J_{\Phi}\Lambda_{\Phi}(b)|J_{\Phi}\Lambda_{\Phi}(b))
&=(x\Delta_{\Phi\circ R}^{-it}J_{\Phi}\Lambda_{\Phi}(b)|
\Delta_{\Phi\circ R}^{-it}J_{\Phi}\Lambda_{\Phi}(b))\\
&=(xJ_{\Phi}\delta^{it}\!\!J_{\Phi}\delta^{-it}\!\Delta_{\Phi}^{-it}\!\!J_{\Phi}\Lambda_{\Phi}(b)|
J_{\Phi}\delta^{it}\!J_{\Phi}\delta^{-it}\!\Delta_{\Phi}^{-it}\!J_{\Phi}\Lambda_{\Phi}(b))\\
&=(x\delta^{-it}J_{\Phi}\Lambda_{\Phi}(\sigma_{-t}^{\Phi}(b))|
\delta^{-it}J_{\Phi}\Lambda_{\Phi}(\sigma_{-t}^{\Phi}(b)))\\
&=(xJ_{\Phi}\Lambda_{\Phi}(\lambda^{\frac{t}{2}}\sigma_{-t}^{\Phi}(b)\delta^{it})|
J_{\Phi}\Lambda_{\Phi}(\lambda^{\frac{t}{2}}\sigma_{-t}^{\Phi}(b)\delta^{it}))\\
&=(xJ_{\Phi}\Lambda_{\Phi}(\lambda^{\frac{t}{2}}\delta^{-it}\sigma_{-t}^{\Phi}(b)\delta^{it})|
J_{\Phi}\Lambda_{\Phi}(\lambda^{\frac{t}{2}}\delta^{-it}\sigma_{-t}^{\Phi}(b)\delta^{it}))\\
&=(xJ_{\Phi}\Lambda_{\Phi}(\lambda^{\frac{t}{2}}\sigma_{-t}^{\Phi\circ
R}(b))|
J_{\Phi}\Lambda_{\Phi}(\lambda^{\frac{t}{2}}\sigma_{-t}^{\Phi\circ
R}(b)))
\end{aligned}$$

\end{proof}

\begin{prop}\label{fortau}
We have $\Gamma\circ\tau_t=(\sigma_t^{\Phi}\surl{\ _{\beta}
  \star_{\alpha}}_{\
  N}\sigma_{-t}^{\Phi\circ R})\circ\Gamma$ for all $t\in\mathbb{R}$.
\end{prop}

\begin{proof}
For all $a,b\in {\mathcal N}_{\Phi}\cap {\mathcal N}_{T_L}$ and
$t\in\mathbb{R}$, we compute:
$$
\begin{aligned}
&\ \quad (id\surl{\ _{\beta}\star_{\alpha}}_{\
\nu}\omega_{J_{\Phi}\Lambda_{\Phi}(b)})[(\sigma^{\Phi}_{-t}\surl{\
_{\beta}\star_{\alpha}}_{\ N}\sigma^{\Phi\circ
R}_t)\circ\Gamma\circ\tau_t(a^*a)]\\&=\sigma^{\Phi}_{-t}[(id\surl{\
_{\beta}\star_{\alpha}}_{\
\nu}\omega_{J_{\Phi}\Lambda_{\Phi}(b)}\circ\sigma_t^{\Phi\circ
R})(\Gamma\circ\tau_t(a^*a))]
\end{aligned}$$ By the previous lemma,
this last expression is equal to:
$$
\begin{aligned}
&\ \quad\sigma^{\Phi}_{-t}[(id\surl{\ _{\beta}\star_{\alpha}}_{\
\nu}\omega_{J_{\Phi}\Lambda_{\Phi}(\lambda^{\frac{t}{2}}\sigma_{-t}^{\Phi\circ
R}(b))})(\Gamma\circ\tau_t(a^*a))]\\
&=\sigma^{\Phi}_{-t}\circ R[(id\surl{\ _{\beta}\star_{\alpha}}_{\
\nu}\omega_{J_{\Phi}\Lambda_{\Phi}(\tau_t(a))}
)(\Gamma(\lambda^t\sigma_{-t}^{\Phi\circ R}(b^*b)))]\\
&=R\circ\sigma_t^{\Phi\circ R}[(id\surl{\ _{\beta}\star_{\alpha}}_{\
\nu}\omega_{J_{\Phi}\Lambda_{\Phi}(\lambda^{\frac{t}{2}}\tau_t(a))}
)(\Gamma\circ\sigma_{-t}^{\Phi\circ R}(b^*b))]
\end{aligned}$$ Again, by the previous lemma, this last expression
is equal to:
$$
\begin{aligned}
&\ \quad R\circ\sigma_t^{\Phi\circ R}[(id\surl{\
_{\beta}\star_{\alpha}}_{\
\nu}\omega_{J_{\Phi}\Lambda_{\Phi}(a)}\circ\tau_{-t}
(\Gamma\circ\sigma_{-t}^{\Phi\circ R}(b^*b))]\\
&=R[(id\surl{\ _{\beta}\star_{\alpha}}_{\
\nu}\omega_{J_{\Phi}\Lambda_{\Phi}(a)} )(\Gamma(b^*b))]=(id\surl{\
_{\beta}\star_{\alpha}}_{\ \nu}\omega_{J_{\Phi}\Lambda_{\Phi}(b)}
)(\Gamma(a^*a))
\end{aligned}$$ So, we conclude that $(\sigma^{\Phi}_{-t}\surl{\
_{\beta}\star_{\alpha}}_{\ N}\sigma^{\Phi\circ
R}_t)\circ\Gamma\circ\tau_t=\Gamma$ for all $t\in\mathbb{R}$.
\end{proof}

\begin{coro}\label{ouf}
For all $t\in\mathbb{R}$ and $m\in M$, we have:
$$(\delta^{it}\surl{\ _{\beta} \otimes_{\alpha}}_{\
N}\delta^{it})\Gamma(m)(\delta^{-it}\surl{\
_{\beta}\otimes_{\alpha}}_{\
N}\delta^{-it})=\Gamma(\delta^{it}m\delta^{-it})$$ In particular,
for all $s,t\in\mathbb{R}$, $\Gamma(\delta^{is})$ and
$\delta^{it}\surl{\ _{\beta}\otimes_{\alpha}}_{\ N}\delta^{it}$
commute each other.
\end{coro}

\begin{proof}
For all $t\in\mathbb{R}$, we have:
$$
\begin{aligned}
(\sigma_{-t}^{\Phi}\circ\sigma_t^{\Phi\circ R}\surl{\ _{\beta}
\star_{\alpha}}_{\ N}\sigma_{-t}^{\Phi}\circ\sigma_t^{\Phi\circ
R})\circ\Gamma&=(\sigma_{-t}^{\Phi}\surl{\
_{\beta}\star_{\alpha}}_{\
N}\sigma_{-t}^{\Phi}\circ\sigma_t^{\Phi\circ
R}\circ\tau_t)\circ\Gamma\circ\sigma_t^{\Phi\circ R}\\
&=(\sigma_{-t}^{\Phi}\circ\tau_t\surl{\ _{\beta} \star_{\alpha}}_{\
N}\sigma_t^{\Phi\circ R}\circ\tau_t)\circ
\Gamma\circ\sigma_{-t}^{\Phi}\circ\sigma_t^{\Phi\circ R}\\
&=(\sigma_{-t}^{\Phi}\surl{\ _{\beta}\star_{\alpha}}_{\
N}\sigma_t^{\Phi\circ
R})\circ\Gamma\circ\tau_t\sigma_{-t}^{\Phi}\circ\sigma_t^{\Phi\circ
R}\\
&=\Gamma\circ\sigma_{-t}^{\Phi}\circ\sigma_t^{\Phi\circ R}
\end{aligned}$$ We know that
$\sigma_{-t}^{\Phi}\sigma_t^{\Phi\circ
R}(m)=\delta^{it}m\delta^{-it}$ for all $m\in M$, that's why we get:
$$(\delta^{it}\surl{\ _{\beta} \otimes_{\alpha}}_{\
N}\delta^{it})\Gamma(m)(\delta^{-it}\surl{\
_{\beta}\otimes_{\alpha}}_{\
N}\delta^{-it})=\Gamma(\delta^{it}m\delta^{-it})$$ In particular,
for all $s\in\mathbb{R}$, we have:
$$(\delta^{it}\surl{\ _{\beta}
  \otimes_{\alpha}}_{\ N}\delta^{it})\Gamma(\delta^{is})(\delta^{-it}\surl{\ _{\beta}
  \otimes_{\alpha}}_{\
  N}\delta^{-it})=\Gamma(\delta^{it}\delta^{is}\delta^{-it})=\Gamma(\delta^{is})$$
\end{proof}

\begin{prop}\label{tech}
Let us denote by ${\mathcal T}_{\Phi,T_L}^{\Psi}$ made of elements
$a\in {\mathcal N}_{T_R}\cap {\mathcal N}_{\Phi} \cap {\mathcal
N}_{\Psi}$, analytic with respect to both $\Phi$ and $\Psi$ such
that, for all $z,z'\in {\mathcal C}$,
$\sigma_z^{\Psi}\circ\sigma_z^{\Phi}(a)$ belongs to ${\mathcal
N}_{T_R}\cap {\mathcal N}_{\Phi} \cap {\mathcal N}_{\Psi}$. This
linear space is weakly dense in $M$ and the set of
$\Lambda_{\Psi}(a)$ (resp. $\Lambda_{\Phi}(a)$), for all $a\in
{\mathcal T}_{\Phi,T_L}^{\Psi}$, is a linear dense subset in $H$.
Moreover, the subset $J_{\Psi}\Lambda_{\Psi}({\mathcal
T}_{\Phi,T_L}^{\Psi})$ is included in the domain of $\delta^z$, for
all $z\in {\mathcal C}$ and is an essential domain for $\delta^z$.
\end{prop}

\begin{proof}
Let us take $x\in {\mathcal T}_{\Phi,T_L}^+$ and let us write
$\lambda=\int_0^{\infty}\!t\ de_t$ and define
$f_p=\int_{\frac{1}{p}}^pde_t$. If we put:
$$x_{p,q}=f_p\sqrt{\frac{q}{\pi}}\int_{-\infty}^{+\infty}e^{-qt^2}\sigma_t^{\Psi}(x)\
dt$$ we obtain that $x_{p,q}$ belongs to ${\mathcal
T}_{\Phi,T_L}^+$, is analytical with respect to $\Psi$,
$\sigma_z^{\Psi}(x_{p,q})$ is weakly converging to $x$ and
$\Lambda_{\Phi}(x_{p,q})$ is weakly converging to
$\Lambda_{\Phi}(x)$. Moreover, $\Lambda_{T_L}(x_{p,q})$is weakly
converging to $\Lambda_{T_L}(x)$.

Since, for all $y\in M$ and $t\in\mathbb{R}$, we have:
$$\sigma_t^{\Psi}\sigma_{-t}^{\Phi}=\delta^{it}x\delta^{-it}$$ we
see that, for all such elements $x_{p,q}$ and $z\in\mathbb{C}$,
$\delta^{iz}x_{p,q}\delta^{-iz}$ is bounded and belongs to
${\mathcal T}_{\Phi,T_L}$. In particular,
$\delta^{-\frac{1}{2}}x_{p,q}\delta^{\frac{1}{2}}$ belongs to
${\mathcal M}_{T_L}\cap {\mathcal T}_{\Phi}$ and is analytic with
respect to both $\Phi$ and $\Psi$. Using then the operator $e_n$
introduced in \cite{Vae}, 1.1, which are analytic to both $\Phi$ and
$\Psi$ and converging to $1$ when $n$ goes to infinity, we get that
$e_nx_{p,q}$ belongs to ${\mathcal N}_{T_L}\cap {\mathcal
N}_{\Phi}$. On the other hand, since:
$$e_nx_{p,q}\delta^{\frac{1}{2}}=(e_n\delta^{\frac{1}{2}})\delta^{-\frac{1}{2}}x_{p,q}\delta^{\frac{1}{2}}$$
belongs to ${\mathcal N}_{\Phi}$, we see, by \cite{Vae}, 3.3, that
$e_nx_{p,q}$ belongs to ${\mathcal N}_{\Psi}$ and, therefore, to
${\mathcal T}_{\Phi,T_L}^{\Psi}$, from which we then get all the
results claimed.
\end{proof}

Let recall proposition 2.4 of \cite{E3}:

\begin{prop}
Let $a,b$ in ${\mathcal N}_{T_L}$. Then $T_L(a^*a)$ and $T_L(b^*b)$
are positive self-adjoint closed operators which verify:
$$\omega_{J_{\Psi}\Lambda_{\Psi}(a)}(T_L(b^*b))=\omega_{J_{\Psi}\Lambda_{\Psi}(b)}(T_L(a^*a))$$
\end{prop}

\begin{lemm}\label{inter}
Let $b\in {\mathcal N}_{T_R}\cap {\mathcal N}_{\Phi} \cap {\mathcal
N}_{\Psi}$ and $X$ positive affiliated to $M$ be such that
$\delta^{-\frac{1}{2}}X\delta^{-\frac{1}{2}}$ is bounded. Then the
element of the extended positive part $(id\surl{\ _{\beta}
\star_{\alpha}}_{\ \nu}\Psi)\Gamma(X)$ is such that:
$$\omega_{J_{\Psi}\Lambda_{\Psi}(b)}((id\surl{\ _{\beta}
\star_{\alpha}}_{\
\nu}\Psi)\Gamma(X))=\omega_{\delta^{\frac{1}{2}}J_{\Psi}\Lambda_{\Psi}(b)}
(T_L(\delta^{-\frac{1}{2}}X\delta^{-\frac{1}{2}}))$$ If $X$ is
bounded, such that $\delta^{-\frac{1}{2}}X\delta^{-\frac{1}{2}}$ is
bounded and in ${\mathcal M}_{T_L}^+$ then
$(\omega_{J_{\Psi}\Lambda_{\Psi}(b)}\surl{\ _{\beta}
\star_{\alpha}}_{\ \nu}id)\Gamma(X)$ belongs to ${\mathcal
M}_{T_L}^+\cap {\mathcal M}_{\Psi}^+$. If $Y$ is in ${\mathcal
M}_{T_L}^+$, we have:
$$\delta^{-\frac{1}{2}}T_L(Y)\delta^{\frac{1}{2}}=(id\surl{\ _{\beta}
\star_{\alpha}}_{\
\nu}\Psi)\Gamma(\delta^{\frac{1}{2}}X\delta^{\frac{1}{2}})$$
\end{lemm}

\begin{proof}
Let us assume that $a,b\in {\mathcal N}_{T_R}\cap {\mathcal
N}_{\Phi} \cap {\mathcal N}_{\Psi}$. By \cite{Vae},
$J_{\Phi}\Lambda_{\Phi}(a)$ is in the domain of
$\delta^{-\frac{1}{2}}$ and
$\delta^{-\frac{1}{2}}J_{\Phi}\Lambda_{\Phi}(a)=\lambda^{\frac{i}{4}}\delta^{-\frac{1}{2}}J_{\Psi}\Lambda_{\Psi}(a)$.
Then, we compute the following:
$$
\begin{aligned}
&\ \quad\omega_{J_{\Psi}\Lambda_{\Psi}(b)}((id\surl{\ _{\beta}
\star_{\alpha}}_{\ \nu}\Psi)\Gamma(a^*a))=\Phi\circ
R((\omega_{J_{\Psi}\Lambda_{\Psi}(b)}\surl{\ _{\beta}
\star_{\alpha}}_{\
\nu}id)\Gamma(a^*a))\\
&=\Phi((\omega_{J_{\Psi}\Lambda_{\Psi}(a)}\surl{\ _{\beta}
\star_{\alpha}}_{\ \nu}id)\Gamma(b^*b))
=\omega_{J_{\Psi}\Lambda_{\Psi}(a)}(T_L(b^*b))=\omega_{\delta^{-\frac{1}{2}}J_{\Phi}\Lambda_{\Phi}(a)}(T_L(b^*b))\\
&=\omega_{J_{\Phi}\Lambda_{\Phi}(a\delta^{-\frac{1}{2}})}(T_L(b^*b))=\omega_{J_{\Phi}\Lambda_{\Phi}(b)}(T_L(\delta^{-\frac{1}{2}}a^*a\delta^{-\frac{1}{2}}))\\
&=\omega_{\delta^{\frac{1}{2}}J_{\Psi}\Lambda_{\Psi}(b)}(T_L(\delta^{-\frac{1}{2}}a^*a\delta^{-\frac{1}{2}}))
\end{aligned}$$
If $X$ is positive such that
$\delta^{-\frac{1}{2}}X\delta^{-\frac{1}{2}}$ is bounded, we may
consider $X$ as the upper limit of elements of the type $a_i^*a_i$
where $a_i$ belongs to the dense left ideal ${\mathcal N}_{T_R}\cap
{\mathcal N}_{\Phi} \cap {\mathcal N}_{\Psi}$. Then every
$a_i\delta^{-\frac{1}{2}}$ is bounded and we get the first formula
by increasing limits. The proof of the second one is an easy
corollary of the first one because we are in the essential domain of
$\delta^{\frac{1}{2}}$.
\end{proof}

\begin{theo}
We have $\Gamma(\delta)=\delta\surl{\ _{\beta}\otimes_{\alpha}}_{\
N}\delta$.
\end{theo}

\begin{proof}
Applying $\Gamma$ to the second equality of the previous
proposition, we get for all $Y\in {\mathcal M}_{\Phi}^+$:
$$
\begin{aligned}
\Gamma(\delta^{\frac{1}{2}})(T_L(Y)\surl{\
_{\beta}\otimes_{\alpha}}_{\ N}1)\Gamma(\delta^{\frac{1}{2}})
&=\Gamma((id\surl{\ _{\beta}\star_{\alpha}}_{\
\nu}\Psi)\Gamma(\delta^{\frac{1}{2}}Y\delta^{\frac{1}{2}}))\\
&=(id\surl{\ _{\beta}\star_{\alpha}}_{\ \nu}id\surl{\
_{\beta}\star_{\alpha}}_{\ \nu}\Psi)(\Gamma\surl{\
_{\beta}\star_{\alpha}}_{\
\nu}id)\Gamma(\delta^{\frac{1}{2}}Y\delta^{\frac{1}{2}}))\\
&=(id\surl{\ _{\beta}\star_{\alpha}}_{\ \nu}id\surl{\
_{\beta}\star_{\alpha}}_{\ \nu}\Psi)(id\surl{\
_{\beta}\star_{\alpha}}_{\
\nu}\Gamma)\Gamma(\delta^{\frac{1}{2}}Y\delta^{\frac{1}{2}}))
\end{aligned}$$
Let now $b\in {\mathcal T}_{\Phi,T_L}^{\Psi}$ and define $Z$ by:
$$Z=(\omega_{J_{\Psi}\Lambda_{\Psi}(b)}\surl{\ _{\beta}\star_{\alpha}}_{\
\nu}id)\Gamma(\delta^{\frac{1}{2}}Y\delta^{\frac{1}{2}})$$ By
corollary \ref{ouf}, we have:
$$\delta^{-\frac{1}{2}}Z\delta^{-\frac{1}{2}}
=(\omega_{\delta^{\frac{1}{2}}J_{\Psi}\Lambda_{\Psi}(b)}\surl{\
_{\beta}\star_{\alpha}}_{\ \nu}id)\Gamma(Y)$$ which is bounded by
proposition \ref{tech}. By the previous proposition, we get for all
$b'\in{\mathcal T}_{T_R}\cap {\mathcal T}_{\Psi}$:
$$
\begin{aligned}
&\ \quad\omega_{J_{\Psi}\Lambda_{\Psi}(b)\surl{\
_{\beta}\otimes_{\alpha}}_{\
\nu}J_{\Psi}\Lambda_{\Psi}(b')}(\Gamma(\delta^{\frac{1}{2}})(T_L(Y)\surl{\
_{\beta}\otimes_{\alpha}}_{\ N}1)\Gamma(\delta^{\frac{1}{2}}))\\
&= \omega_{J_{\Psi}\Lambda_{\Psi}(b)\surl{\
_{\beta}\otimes_{\alpha}}_{\
\nu}J_{\Psi}\Lambda_{\Psi}(b')}((id\surl{\
_{\beta}\star_{\alpha}}_{\ \nu}id\surl{\ _{\beta}\star_{\alpha}}_{\
\nu}\Psi)(id\surl{\ _{\beta}\star_{\alpha}}_{\
\nu}\Gamma)\Gamma(\delta^{\frac{1}{2}}Y\delta^{\frac{1}{2}}))\\
&=\omega_{J_{\Psi}\Lambda_{\Psi}(b')}((id\surl{\
_{\beta}\star_{\alpha}}_{\
\nu}\Psi)\Gamma((\omega_{J_{\Psi}\Lambda_{\Psi}(b)}\surl{\
_{\beta}\otimes_{\alpha}}_{\
\nu}id)\Gamma(\delta^{\frac{1}{2}}Y\delta^{\frac{1}{2}})))\\
&=\omega_{\delta^{\frac{1}{2}}J_{\Psi}\Lambda_{\Psi}(b')}
(T_L(\delta^{-\frac{1}{2}}Z\delta^{-\frac{1}{2}}))=\omega_{\delta^{\frac{1}{2}}J_{\Psi}\Lambda_{\Psi}(b')}
(T_L((\omega_{\delta^{\frac{1}{2}}J_{\Psi}\Lambda_{\Psi}(b)}\surl{\
_{\beta}\star_{\alpha}}_{\ \nu}id)\Gamma(Y)))\\
&=\omega_{\delta^{\frac{1}{2}}J_{\Psi}\Lambda_{\Psi}(b)\surl{\
_{\beta}\otimes_{\alpha}}_{\
\nu}\delta^{\frac{1}{2}}J_{\Psi}\Lambda_{\Psi}(b')}(T_L(Y)\surl{\
_{\beta}\otimes_{\alpha}}_{\ N}1)
\end{aligned}$$
from which we infer, by increasing limits, that:
$$\omega_{J_{\Psi}\Lambda_{\Psi}(b)\surl{\
_{\beta}\otimes_{\alpha}}_{\
\nu}J_{\Psi}\Lambda_{\Psi}(b')}(\Gamma(\delta))=||\delta^{\frac{1}{2}}J_{\Psi}\Lambda_{\Psi}(b)\surl{\
_{\beta}\otimes_{\alpha}}_{\
\nu}\delta^{\frac{1}{2}}J_{\Psi}\Lambda_{\Psi}(b')||^2$$ which
finishes the proof by proposition \ref{tech}.
\end{proof}

\subsection{Uniqueness of invariant operator-valued weight}

\begin{theo}
If $T'$ a n.s.f left invariant operator-valued weight such that
$(\tau_t \surl{\ _{\beta} \star_{\alpha}}_{\ N}
\sigma^{\Phi'}_t)\circ \Gamma=\Gamma\circ\sigma^{\Phi'}_t$,
$\nu\circ\gamma'=\nu$ and $\gamma\circ\gamma'=\gamma'\circ\gamma$,
then there exists a strictly positive operator $h$ affiliated with
$Z(N)$ such that, for all $t\in\mathbb{R}$, we have:
$$\Phi'=\nu\circ\alpha^{-1}\circ T'=(\nu\circ\alpha^{-1}\circ
T_L)_{\beta(h)}\text{ and } [DT':DT_L]_t=\beta(h^{it})$$
\end{theo}

\begin{proof}
We put $\Phi'=\nu\circ\alpha^{-1}\circ T'$. We have for all
$s\in\mathbb{R}$:
$$\Gamma\circ\sigma_{-s}^{\Phi}\circ\sigma_s^{\Phi'}
=(\tau_{-s}\surl{\ _{\beta}\star_{\alpha}}_{\
N}\sigma_{-s}^{\Phi})\circ\Gamma\circ\sigma_s^{\Phi'}= (id\surl{\
_{\beta} \star_{\alpha}}_{\
N}\sigma_{-s}^{\Phi}\circ\sigma_s^{\Phi'})\circ\Gamma$$

By right invariance of $T_R$, we have for all $a\in {\mathcal
M}_{T_R}^+$:
$$
\begin{aligned}
T_R(\sigma_{-s}^{\Phi}\circ\sigma_s^{\Phi'}(a))&=(\Phi\circ R\surl{\
_{\beta} \star_{\alpha}}_{\
\nu}id)(\Gamma(\sigma_{-s}^{\Phi}\circ\sigma_s^{\Phi'}(a)))\\
&=\sigma_{-s}^{\Phi}\circ\sigma_s^{\Phi'}((\Phi\circ R\surl{\
_{\beta} \star_{\alpha}}_{\ \nu}id)\Gamma(a))
=\sigma_{-s}^{\Phi}\circ\sigma_s^{\Phi'}(T_R(a))
\end{aligned}$$ Since $\gamma$ and $\gamma'$ leave $\nu$ invariant,
we get that $\Phi\circ R$ is
$\sigma_{-s}^{\Phi}\circ\sigma_s^{\Phi'}$-invariant and, so
$\sigma_t^{\Phi\circ R}$ and
$\sigma_{-s}^{\Phi}\circ\sigma_s^{\Phi'}$ commute each other. But
$\sigma^{\Phi\circ R}$ and $\sigma^{\Phi}$ commute each other that's
why $\sigma^{\Phi\circ R}$ and $\sigma^{\Phi'}$ also commute each
other. For all $s,t\in\mathbb{R}$, we have:
$$\Gamma(\sigma_t^{\Phi'}(\delta^{is})\delta^{-is})=
(\tau_t\surl{\ _{\beta} \star_{\alpha}}_{\
N}\sigma_t^{\Phi'})(\Gamma(\delta^{is}))(\delta^{-is}\surl{\
_{\beta} \otimes_{\alpha}}_{\ N}\delta^{-is})=1\surl{\ _{\beta}
\otimes_{\alpha}}_{\ N}\sigma_t^{\Phi'}(\delta^{is})\delta^{-is}$$
Consequently $\sigma_t^{\Phi'}(\delta^{is})\delta^{-is}$ belongs to
$\beta(N)$. For all $n\in N$ and $s,t\in\mathbb{R}$, we have:
$$
\begin{aligned}
&\
\quad\sigma_t^{\Phi'}(\delta^{is})\beta(n)\sigma_t^{\Phi'}(\delta^{-is})
=\sigma_t^{\Phi'}((\delta^{is})\sigma_{-t}^{\Phi'}(\beta(n))\delta^{-is})
=\sigma_t^{\Phi'}((\delta^{is})\beta(\gamma'_{-t}(n))\delta^{-is})\\
&=\sigma_t^{\Phi'}(\beta(\gamma_s\sigma_s^{\nu}\gamma'_{-t}(n)))
=\beta(\gamma'_t\gamma_s\sigma_s^{\nu}\gamma'_{-t}(n)))
=\beta(\gamma_s\sigma_s^{\nu}(n)))=\delta^{is}\beta(n)\delta^{-is}
\end{aligned}$$

So $\sigma_t^{\Phi'}(\delta^{is})\delta^{-is}$ belongs to
$\beta(Z(N))$ and we easily get that there exists a strictly
positive operator $k$ affiliated with $Z(N)$ such that
$\sigma_t^{\Phi'}(\delta^{is})=\beta(k^{ist})\delta^{is}$. Then, we
have:
$$
\begin{aligned}
\sigma_s^{\Phi'}\circ\sigma_t^{\Phi}(m)
=\sigma_s^{\Phi'}(\delta^{-it}\sigma_t^{\Phi\circ
R}(m)\delta^{it})&=\beta(k^{-ist})\delta^{-it}\sigma_s^{\Phi'}\circ\sigma_t^{\Phi\circ
R}(m)\delta^{it}\beta(k^{ist})\\
&=\beta(k^{-ist})\sigma_t^{\Phi}\circ\sigma_s^{\Phi'}(m)\beta(k^{ist})
\end{aligned}$$ Take $m=\delta^{iu}$ to get $k$ is affiliated to $N^\gamma$. Apply
$\Phi$ to the previous formula and get:
$$\begin{aligned}
\Phi\circ\sigma_s^{\Phi'}\circ\sigma_s^{\Phi}(m^*m)
&=\Phi(\beta(k^{-ist})\sigma_t^{\Phi}\circ\sigma_s^{\Phi'}(m^*m)\beta(k^{ist}))\\
&=\Phi(\sigma_t^{\Phi}\circ\sigma_s^{\Phi'}(m^*m))=\Phi\circ\sigma_s^{\Phi'}(m^*m)
\end{aligned}$$
So, by \ref{struc} and left invariance $\sigma_{-s}^{\Phi'}\circ
T_L\circ\sigma_s^{\Phi'}$, there exists a strictly positive operator
$q_s$ affiliated with $Z(N)$ such that
$\Phi\circ\sigma_s^{\Phi'}=\Phi_{\beta(q_s)}$. By usual arguments,
we deduce that there exists a strictly positive $q$ affiliated to
$Z(N)$ such that $\Phi\circ\sigma_s^{\Phi'}=\Phi_{\beta(q^{-s})}$
and $[D\Phi\circ\sigma_s^{\Phi'}:D\Phi]_s=\beta(q^{-s})$. Then,
again by \ref{struc}, there exists a strictly positive operator $h$
affiliated to $Z(N)$ such that $\Phi=\Phi_{\beta(h)}$ avec
$[DT':DT_L]_t=\beta(h^{it})$.
\end{proof}

Also, we have a similar result for right invariant operator-valued
weight.

\begin{coro}
If $T_R$ a n.s.f right invariant operator-valued weight such that
$(\sigma^{\Psi'}_t\surl{\ _{\beta} \star_{\alpha}}_{\ N} \tau_{-t}
)\circ \Gamma=\Gamma\circ\sigma^{\Psi'}_t$, $\nu\circ\gamma'=\nu$
and $\gamma\circ\gamma'=\gamma'\circ\gamma$, then there exists a
strictly positive operator $h$ affiliated with $Z(N)$ such that:
$$T_R=(R\circ T_L\circ R)_{\alpha(h)}$$
\end{coro}

We state results of the section in the following theorems:

\begin{theo}\label{ensemble1}
Let $(N,M,\alpha,\beta,\Gamma,T_L,R,\tau,\nu)$ be a measured quantum
groupoid. If $T'$ a n.s.f left invariant operator-valued weight such
that $(\tau_t \surl{\ _{\beta} \star_{\alpha}}_{\ N}
\sigma^{\Phi'}_t)\circ \Gamma=\Gamma\circ\sigma^{\Phi'}_t$,
$\nu\circ\gamma'=\nu$ and $\gamma\circ\gamma'=\gamma'\circ\gamma$,
then there exists a strictly positive operator $h$ affiliated with
$Z(N)$ such that, for all $t\in\mathbb{R}$:
$$\nu\circ\alpha^{-1}\circ T'=(\nu\circ\alpha^{-1}\circ
T_L)_{\beta(h)}$$ We have a similar result for the right invariant
operator-valued weights.
\end{theo}

\begin{theo}\label{ensemble2}
Let $(N,M,\alpha,\beta,\Gamma,T_L,R,\tau,\nu)$ be a measured quantum
groupoid. Then there exists a strictly positive operator $\delta$
affiliated with $M$ called modulus and then there exists a strictly
positive operator $\lambda$ affiliated with
$Z(M)\cap\alpha(N)\cap\beta(N)$ called scaling operator such that
$[D\nu\circ\alpha^{-1}\circ T_L\circ R:D\nu\circ\alpha^{-1}\circ
T_L]_t=\lambda^{\frac{it^2}{2}}\delta^{it}$ for all
$t\in\mathbb{R}$.

Moreover, we have, for all $s,t\in\mathbb{R}$:
\begin{enumerate}[i)]
\item $\begin{aligned} &\ [D\nu\circ\alpha^{-1}\circ T_L\circ\tau_s
:D\nu\circ\alpha^{-1}\circ T_L]_t=\lambda^{-ist}\\
                       &\ [D\nu\circ\alpha^{-1}\circ T_L\circ R\circ\tau_s
:D\nu\circ\alpha^{-1}\circ T_L\circ R]_t=\lambda^{-ist}\\
                       &\ [D\nu\circ\alpha^{-1}\circ T_L\circ\sigma^{\nu\circ\alpha^{-1}\circ T_L\circ R}_s
:D\nu\circ\alpha^{-1}\circ T_L]_t=\lambda^{ist}\\
                       &\ [D\nu\circ\alpha^{-1}\circ T_L\circ R\circ\sigma^{\nu\circ\alpha^{-1}\circ T_L}_s:D\nu\circ\alpha^{-1}\circ T_L\circ R]_t=\lambda^{-ist}
       \end{aligned}$
\item $R(\lambda)=\lambda$, $R(\delta)=\delta^{-1}$,
$\tau_t(\delta)=\delta$ and $\tau_t(\lambda)=\lambda$ ;
\item $\delta$ is a group-like element i.e $\Gamma(\delta)
      =\delta\surl{\ _{\beta} \otimes_{\alpha}}_{\ N}\delta$.
\end{enumerate}
\end{theo}

\section{A density theorem}\label{dtheo}
In this section, we prove that there are sufficiently enough
operators which are both bounded under the left-invariant
operator-valued weight and the right-invariant operator-valued
weight. This allows, as a corollary, to found bounded elements for
both $\alpha$ and $\beta$ which will be useful for duality. This
chapter is mostly inspired by chapter $7$ of \cite{E3}.

\begin{lemm}
Let $y,z\in {\mathcal N}_{T_R}\cap {\mathcal N}_{\Psi}$ and $\xi\in
D((H_{\Psi})_{\beta},\mu^o)$, then we have:
$$\begin{aligned}
&\ \quad [(\omega_{\Lambda_{\Psi}(y),J_{\Psi}z^*J_{\Psi}\xi}*
id)(U'_H)]^*(\omega_{\Lambda_{\Psi}(y),J_{\Psi}z^*J_{\Psi}\xi}*
id)(U'_H)\\
&\leq ||R^{\beta,\nu^0}(\xi)||^2(\omega_{J_{\Psi}\Lambda_{\Psi}(z)}
\surl{\ _{\beta} \star_{\alpha}}_{\ \mu}id)(\Gamma(y^*y))
\end{aligned}$$ For all $y\in {\mathcal N}_{T_R}\cap {\mathcal
N}_{\Psi}$, $z\in M$ and $\xi\in D((H_{\Psi})_{\beta},\mu^o)$, then
we have:
$$
\begin{aligned}
&\ \quad R\left([(\omega_{\Lambda_{\Psi}(y),J_{\Psi}z^*J_{\Psi}\xi}*
id)(U'_H)]^*(\omega_{\Lambda_{\Psi}(y),J_{\Psi}z^*J_{\Psi}\xi}*
id)(U'_H)\right)\\
&\leq ||R^{\beta,\nu^0}(\xi)||^2(\omega_{J_{\Psi}\Lambda_{\Psi}(y)}
\surl{\ _{\beta} \star_{\alpha}}_{\ \mu}id)(\Gamma(z^*z))
\end{aligned}$$
\end{lemm}

\begin{proof}
The first inequality comes straightforward from theorem \ref{319}.
Then, apply $R$ to get for all $z\in {\mathcal N}_{T_R}\cap
{\mathcal N}_{\Psi}$:
$$
\begin{aligned}
&\ \quad R\left([(\omega_{\Lambda_{\Psi}(y),J_{\Psi}z^*J_{\Psi}\xi}*
id)(U'_H)]^*(\omega_{\Lambda_{\Psi}(y),J_{\Psi}z^*J_{\Psi}\xi}*
id)(U'_H)\right)\\
&\leq ||R^{\beta,\nu^0}(\xi)||^2R(\omega_{J_{\Psi}\Lambda_{\Psi}(z)}
\surl{\ _{\beta} \star_{\alpha}}_{\ \mu}id)(\Gamma(y^*y))\\
&=||R^{\beta,\nu^0}(\xi)||^2(\omega_{J_{\Psi}\Lambda_{\Psi}(y)}
\surl{\ _{\beta} \star_{\alpha}}_{\ \mu}id)(\Gamma(z^*z))
\end{aligned}$$

Let us assume now that $z\in M$. Using Kaplansky' s theorem, there
exist a family $z_i$ in ${\mathcal N}_{T_R}\cap {\mathcal
N}_{\Psi}$, weakly converging to $z$, with $||z_i||\leq ||z||$. Then
we infer that $R^{\beta,\nu^0}(J_{\Psi}z_i^*J_{\Psi}\xi)$ is weakly
converging to $R^{\beta,\nu^0}(J_{\Psi}z^*J_{\Psi}\xi)$ with:
$$||R^{\beta,\nu^0}(J_{\Psi}z^*J_{\Psi}\xi)||\leq
||R^{\beta,\nu^0}(J_{\Psi}z^*J_{\Psi}\xi)||$$ Therefore
$(\omega_{\Lambda_{\Psi}(y),J_{\Psi}z_i^*J_{\Psi}\xi}* id)(U'_H)$ is
weakly converging to
$(\omega_{\Lambda_{\Psi}(y),J_{\Psi}z^*J_{\Psi}\xi}* id)(U'_H)$
with:
$$||(\omega_{\Lambda_{\Psi}(y),J_{\Psi}z_i^*J_{\Psi}\xi}*
id)(U'_H)||\leq
||(\omega_{\Lambda_{\Psi}(y),J_{\Psi}z^*J_{\Psi}\xi}* id)(U'_H)||$$
which finishes the proof.
\end{proof}

\begin{prop}\label{dens1}
If $z\in {\mathcal N}_{T_L}$, $y\in {\mathcal N}_{T_R}\cap {\mathcal
N}_{\Psi}$ and $\xi\in D((H_{\Psi})_{\beta},\mu^o)$ then
$(\omega_{\Lambda_{\Psi}(y),J_{\Psi}z^*J_{\Psi}\xi}* id)(U'_H)$
belongs to ${\mathcal N}_{T_R}\cap {\mathcal N}_{\Psi}$.
\end{prop}

\begin{proof}
By the previous lemma and by right left-invariance of $\Phi$, we
have:
$$
\begin{aligned}
&\
\quad\Psi\left((\omega_{\Lambda_{\Psi}(y),J_{\Psi}z^*J_{\Psi}\xi}*
id)(U'_H)^*(\omega_{\Lambda_{\Psi}(y),J_{\Psi}z^*J_{\Psi}\xi}*
id)(U'_H)\right)\\&= \Phi\circ
R\left([(\omega_{\Lambda_{\Psi}(y),J_{\Psi}z^*J_{\Psi}\xi}*
id)(U'_H)]^*(\omega_{\Lambda_{\Psi}(y),J_{\Psi}z^*J_{\Psi}\xi}*
id)(U'_H)\right)\\
&\leq ||R^{\beta,\nu^0}(\xi)||^2\omega_{J_{\Psi}\Lambda_{\Psi}(y)}
(id\surl{\ _{\beta} \star_{\alpha}}_{\
\mu}\Phi)(\Gamma(z^*z))=||R^{\beta,\nu^0}(\xi)||^2\omega_{J_{\Psi}\Lambda_{\Psi}(y)}(T_L(z^*z))
\end{aligned}$$
Also, we have:
$$
\begin{aligned}
&\ \quad
T_R\left((\omega_{\Lambda_{\Psi}(y),J_{\Psi}z^*J_{\Psi}\xi}*
id)(U'_H)^*(\omega_{\Lambda_{\Psi}(y),J_{\Psi}z^*J_{\Psi}\xi}*
id)(U'_H)\right)\\&= R\circ T_L\circ
R\left([(\omega_{\Lambda_{\Psi}(y),J_{\Psi}z^*J_{\Psi}\xi}*
id)(U'_H)]^*(\omega_{\Lambda_{\Psi}(y),J_{\Psi}z^*J_{\Psi}\xi}*
id)(U'_H)\right)\\
&\leq
||R^{\beta,\nu^0}(\xi)||^2(\omega_{J_{\Psi}\Lambda_{\Psi}(y)}\surl{\
_{\beta} \star_{\alpha}}_{\ \mu}id)(id\surl{\ _{\beta}
\star_{\alpha}}_{\
\mu}T_L)(\Gamma(z^*z))\\
&=||R^{\beta,\nu^0}(\xi)||^2(\omega_{J_{\Psi}\Lambda_{\Psi}(y)}\surl{\
_{\beta} \star_{\alpha}}_{\ \mu}id)(T_L(z^*z)\surl{\ _{\beta}
\otimes_{\alpha}}_{\ \mu}1))\\
&\leq ||R^{\beta,\nu^0}(\xi)||^2||T_L(z^*z)||||T_R(y^*y)||1
\end{aligned}$$
\end{proof}

\begin{lemm}\label{dens2}
For all $y,z\in {\mathcal N}_{T_R}\cap {\mathcal N}_{\Psi}$ and
$\eta\in D((H_{\Psi})_{\beta},\mu^o)$, we have:
$$\begin{aligned}
&\ \quad R[(\omega_{y^*J_{\Psi}\eta,J_{\Psi}\Lambda_{\Psi}(z)}*
id)(U'_H)]R[(\omega_{y^*J_{\Psi}\eta,J_{\Psi}\Lambda_{\Psi}(z)}*
id)(U'_H)]^*\\
&\leq ||T_R(y^*y)||^2(\omega_{\eta}\surl{\ _{\beta}
\star_{\alpha}}_{\ \mu}id)(\Gamma(zz^*)) \end{aligned}$$
\end{lemm}

\begin{proof}
Let us compute:
$$
\begin{aligned}
&\ \quad R[(\omega_{y^*J_{\Psi}\eta,J_{\Psi}\Lambda_{\Psi}(z)}*
id)(U'_H)]R[(\omega_{y^*J_{\Psi}\eta,J_{\Psi}\Lambda_{\Psi}(z)}*
id)(U'_H)]^*\\
&= (\omega_{\Lambda_{\Psi}(z),J_{\Psi}y^*J_{\Psi}\eta}*
id)(U'_H)(\omega_{\Lambda_{\Psi}(z),J_{\Psi}y^*J_{\Psi}\eta}*
id)(U'_H)^*\\
&= (\omega_{J_{\Psi}\Lambda_{\Psi}(y),\eta}\surl{\ _{\beta}
\star_{\alpha}}_{\
\mu}id)(\Gamma(z))(\omega_{J_{\Psi}\Lambda_{\Psi}(y),\eta}\surl{\
_{\beta} \star_{\alpha}}_{\ \mu}id)(\Gamma(z))^*\\
&\leq ||T_R(y^*y)||^2(\omega_{\eta}\surl{\ _{\beta}
\star_{\alpha}}_{\ \mu}id)(\Gamma(zz^*))
\end{aligned}$$
\end{proof}

\begin{prop}
Let $y_1,z'\in {\mathcal N}_{T_R}\cap {\mathcal N}_{\Psi}$, $y_2\in
{\mathcal N}_{T_R}\cap {\mathcal N}_{\Psi}\cap {\mathcal N}_{\Phi}$,
$z\in R({\mathcal T}_{\Phi,T_L}^{\Psi})^*$ defined in proposition
\ref{tech} and $e_n$ the analytic elements associated to the
Radon-Nikodym derivative $\delta$ defined in \cite{Vae}. Then the
operators
$(\omega_{y_1^*\Lambda_{\Psi}(y_2),J_{\Psi}z^*e_n^*\Lambda_{\Psi}(z')}*
id)(U'_H)$ belong to ${\mathcal N}_{T_R}\cap {\mathcal N}_{\Psi}\cap
{\mathcal N}_{T_L}\cap {\mathcal N}_{\Phi}$.
\end{prop}

\begin{proof}
Let us write
$X=(\omega_{y_1^*\Lambda_{\Psi}(y_2),J_{\Psi}z^*e_n^*\Lambda_{\Psi}(z')}*
id)(U'_H)$. Since $y_1^*y_2$ belongs to ${\mathcal N}_{T_R}\cap
{\mathcal N}_{\Psi}$ and $z$ belongs to $R({\mathcal
N}_{T_R})^*={\mathcal N}_{T_L}$ and therefore $e_nz$ belongs to
${\mathcal N}_{T_R}$, we get, using proposition \ref{dens1}, that
$X$ belongs to ${\mathcal N}_{T_R}\cap {\mathcal N}_{\Psi}$. On the
other hand, since $y_1, y_2, z^*e_n^*z'$ belong to ${\mathcal
N}_{T_R}\cap {\mathcal N}_{\Psi}$, we can use lemma \ref{dens2} to
get that:
$$
\begin{aligned}
R(X)R(X)^*&\leq
||T_R(y_1^*y_1)||(\omega_{J_{\Psi}\Lambda_{\Psi}(y_2)}\surl{\
_{\beta} \star_{\alpha}}_{\ \mu}id)(\Gamma(z^*e_n^*z'z'^*e_nz))\\
&\leq
||T_R(y_1^*y_1)||||z'||^2(\omega_{J_{\Psi}\Lambda_{\Psi}(y_2)}\surl{\
_{\beta} \star_{\alpha}}_{\ \mu}id)(\Gamma(z^*e_n^*e_nz))
\end{aligned}$$
Let us apply $T_R$ to this inequality, we get that:
$$T_R(R(X)R(X)^*)\leq ||T_R(y_1^*y_1)||||z'||^2T_R(\omega_{J_{\Psi}\Lambda_{\Psi}(y_2)}\surl{\
_{\beta} \star_{\alpha}}_{\ \mu}id)(\Gamma(z^*e_n^*e_nz))$$ which is
equal, thanks to lemma \ref{inter}, to:
$$||T_R(y_1^*y_1)||||z'||^2\omega_{\delta^{-\frac{1}{2}}J_{\Psi}\Lambda_{\Psi}(y_2)}
(T_L(\delta^{-\frac{1}{2}}z^*e_n^*e_nz\delta^{-\frac{1}{2}}))$$ With
the hypothesis, we get that
$\delta^{\frac{1}{2}}z\delta^{-\frac{1}{2}}$ belongs to ${\mathcal
N}_{T_L}$ and therefore
$e_nz\delta^{-\frac{1}{2}}=(e_n\delta^{-\frac{1}{2}})\delta^{\frac{1}{2}}z\delta^{-\frac{1}{2}}$
belongs  also to ${\mathcal N}_{T_L}$. We also get that
$J_{\Psi}\Lambda_{\Psi}(y_2)$ belongs to the domain of
$\delta^{-\frac{1}{2}}$ which proves that $R(X)^*$ belongs to
${\mathcal N}_{T_R}$ and therefore $X$ belongs to ${\mathcal
N}_{T_L}$. We prove by similar computations that $X$ belongs to
${\mathcal N}_{\Phi}$.
\end{proof}

\begin{theo}
The left ideal ${\mathcal N}_{T_R}\cap {\mathcal N}_{\Psi}\cap
{\mathcal N}_{T_L}\cap {\mathcal N}_{\Phi}$ is dense in $M$ and
$\Lambda_{\Psi}({\mathcal N}_{T_R}\cap {\mathcal N}_{\Psi}\cap
{\mathcal N}_{T_L}\cap {\mathcal N}_{\Phi})$ is dense in $H$.
\end{theo}

\begin{proof}
Let $y$ be in ${\mathcal N}_{T_R}\cap {\mathcal N}_{\Psi}\cap
{\mathcal N}_{\Phi}$ and $z$ in ${\mathcal N}_{T_R}\cap {\mathcal
N}_{\Psi}$. Taking, by Kaplansky's theorem, a bounded family $e_i$
in ${\mathcal N}_{T_R}\cap {\mathcal N}_{\Psi}$ strongly converging
to $1$, we get that $R^{\hat{\alpha},\mu}(e_i^*\Lambda_{\Psi}(y))$
is weakly converging to $R^{\hat{\alpha},\mu}(\Lambda_{\Psi}(y))$.
Taking also a bounded family $f_k$ in $R({\mathcal
T}_{T_R,\Psi}^{\Phi})^*$ strongly converging to $1$, we get that
$R^{\beta,\mu^0}(J_{\Psi}f_k^*e_n^*\Lambda_{\Psi}(z))$ is weakly
converging, when $n,k$ go to infinity, to
$R^{\beta,\mu^0}(\Lambda_{\Psi}(z))$. Therefore, using the previous
proposition, we get that
$(\omega_{\Lambda_{\Psi}(y),J_{\Psi}\Lambda_{\Psi}(z)}*id)(U'_H)$
belongs to the weak closure of ${\mathcal N}_{T_R}\cap {\mathcal
N}_{\Psi}\cap {\mathcal N}_{T_L}\cap {\mathcal N}_{\Phi}$. By
proposition \ref{tech}, we get that, for any $x\in {\mathcal
T}_{T_R,\Psi}$, there exists $y_i$ in ${\mathcal N}_{T_R}\cap
{\mathcal N}_{\Psi}\cap {\mathcal N}_{\Phi}$ such that
$\Lambda_{T_R}(y_i)$ is weakly converging to $\Lambda_{T_R}(x)$ or
equivalently $R^{\hat{\alpha},\mu}(\Lambda_{\Psi}(y_i))$ is weakly
converging to $R^{\hat{\alpha},\mu}(\Lambda_{\Psi}(x))$. Therefore,
we get that
$(\omega_{\Lambda_{\Psi}(x),J_{\Psi}\Lambda_{\Psi}(z)}*id)(U'_H)$
belongs to the weak closure of ${\mathcal N}_{T_R}\cap {\mathcal
N}_{\Psi}\cap {\mathcal N}_{T_L}\cap {\mathcal N}_{\Phi}$. It
remains true for $x$ in ${\mathcal N}_{T_R}\cap {\mathcal
N}_{\Psi}\cap {\mathcal N}_{T_R}^*\cap {\mathcal N}_{\Psi}^*$ by
density. If now $x$ belongs to ${\mathcal N}_{T_R}\cap {\mathcal
N}_{\Psi}$, and $h_i$ is a bounded family in ${\mathcal N}_{T_R}\cap
{\mathcal N}_{\Psi}$, since
$\Lambda_{T_R}(h_i^*x)=h_i^*\Lambda_{T_R}(x)$ is weakly converging
to $\Lambda_{T_R}(x)$, we finally obtain that, for any $x,z$ in
${\mathcal N}_{T_R}\cap {\mathcal N}_{\Psi}$, the operator
$(\omega_{\Lambda_{\Psi}(y),J_{\Psi}\Lambda_{\Psi}(z)}*id)(U'_H)$
belongs to the weak closure of ${\mathcal N}_{T_R}\cap {\mathcal
N}_{\Psi}\cap {\mathcal N}_{T_L}\cap {\mathcal N}_{\Phi}$. By
density, for all $\xi\in D(\ _{\hat{\alpha}}H,\mu)$ and $\eta\in
D(H_{\beta},\mu^0)$, the operator $(\omega_{\xi,\eta}*id)(U'_H)$
belongs to the weak closure of ${\mathcal N}_{T_R}\cap {\mathcal
N}_{\Psi}\cap {\mathcal N}_{T_L}\cap {\mathcal N}_{\Phi}$. Which
proves the density of ${\mathcal N}_{T_R}\cap {\mathcal
N}_{\Psi}\cap {\mathcal N}_{T_L}\cap {\mathcal N}_{\Phi}$ in $M$ by
theorem \ref{densevn1}.

Let $g_n$ an increasing sequence of positive elements of ${\mathcal
M}_{T_R}\cap {\mathcal M}_{\Psi}\cap {\mathcal M}_{T_L}\cap
{\mathcal M}_{\Phi}$ strongly  converging to $1$. The operators:
$$h_n=\sqrt{\frac{1}{\pi}}\int_{-\infty}^{+\infty}\!e^{-t^2}\sigma_t^{\Psi}(g_n)\
dt$$ are in ${\mathcal M}_{T_R}\cap {\mathcal M}_{\Psi}$, analytic
with respect to $\Psi$, and, for any $z\in\mathbb{C}$,
$\sigma_z^{\Psi}(h_n)$ is a bounded sequence strongly converging to
$1$. Let now $\lambda=\int_0^{+\infty}\! t\ de_t$ be the scaling
operator. Let us write $h'_n=\left(\int_{\frac{1}{n}}^n
de_t\right)h_n$. These operators are in ${\mathcal N}_{T_R}\cap
{\mathcal N}_{\Psi}$, analytic with respect to $\Psi$, and, for any
$z\in\mathbb{C}$, $\sigma_z^{\Psi}(h'_n)$ is a bounded sequence
strongly converging to $1$. Moreover the operators $h'_n$ belong
also to ${\mathcal N}_{T_L}\cap {\mathcal N}_{\Phi}$ by lemma
\ref{prep1} and \cite{Vae}. Let now $x$ be in ${\mathcal N}_{\Psi}$.
We get that $xh'_n$ belongs to ${\mathcal N}_{T_R}\cap {\mathcal
N}_{\Psi}\cap {\mathcal N}_{T_L}\cap {\mathcal N}_{\Phi}$ and that:
$$\Lambda_{\Psi}(xh'_n)=J_{\Psi}\sigma_{-i/2}^{\Psi}(h'_n)J_{\Psi}\Lambda_{\Psi}(x)$$
is converging to $\Lambda_{\Psi}(x)$ which finishes the proof.
\end{proof}

\begin{theo}
Let ${\mathcal T}_{T_R,\Psi,T_L,\Phi}$ be the subset of elements $x$
in ${\mathcal N}_{T_R}\cap {\mathcal N}_{\Psi}\cap {\mathcal
N}_{T_L}\cap {\mathcal N}_{\Phi}$, analytic with respect to both
$\Phi$ and $\Psi$, and such that, for all $z,z'\in\mathbb{C}$,
$\sigma_z^{\Psi}\circ\sigma_{z'}^{\Phi}(x)$ belongs to ${\mathcal
N}_{T_R}\cap {\mathcal N}_{\Psi}\cap {\mathcal N}_{T_L}\cap
{\mathcal N}_{\Phi}$. Then ${\mathcal T}_{T_R,\Psi,T_L,\Phi}$ is
dense in $M$ and $\Lambda_{\Psi}({\mathcal T}_{T_R,\Psi,T_L,\Phi})$
is dense in $H$.
\end{theo}

\begin{proof}
Let $x$ be a positive operator in ${\mathcal M}_{T_R}\cap {\mathcal
M}_{\Psi}\cap {\mathcal M}_{T_L}\cap {\mathcal M}_{\Phi}$. Let now
$\lambda=\int_0^{\infty}\!t\ de_t$ be the scaling operator and let
us define:
$$x_n=\left(\int_{\frac{1}{n}}^nde_t\right)\frac{n}{\pi}\int_{-\infty}^{+\infty}\!\int_{-\infty}^{+\infty}\!
e^{-n(t^2+s^2)}\sigma_t^{\Psi}\sigma_s^{\Phi}(x)\ dsdt$$ It is not
so difficult to see that $x_n$ is analytic both with respect to
$\Phi$ and $\Psi$. By lemma \ref{prep1} and thanks to \cite{Vae} and
\cite{EN} 10.12, we see that the operators $\sigma_z^{\Psi}(x_n)$
and $\sigma_z^{\Phi}(x_n)$ are linear combinations of positive
elements in ${\mathcal M}_{T_R}\cap {\mathcal M}_{\Psi}\cap
{\mathcal M}_{T_L}\cap {\mathcal M}_{\Phi}$.
\end{proof}

\begin{coro}\label{formi}
There exist a dense linear subspace $E$ of ${\mathcal N}_{\Phi}$
such that $\Lambda_{\Phi}(E)$ is dense in $L^2(M,\Phi)=H$ and:
$$J_{\Phi}\Lambda_{\Phi}(E)\subset D(\ _{\alpha}H,\mu)\cap
D(H_{\beta},\nu^0)$$
\end{coro}

\begin{proof}
Let $E$ be the linear subspace spanned by the elements of the form
$e_nx$ where $e_n$ are the analytic elements associated to the
Radon-Nikodym derivative $\delta$, defined in \cite{Vae}, and $x$
belongs to ${\mathcal T}_{T_R,\Psi,T_L,\Phi}$. It is clear that $E$
is a subset of ${\mathcal N}_{\Phi}$, dense in $M$ and that
$\Lambda_{\Phi}(E)$ is dense in $H$. Since $E\subset {\mathcal
N}_{\Phi}\cap {\mathcal N}_{T_L}$, we have:
$$J_{\Phi}\Lambda_{\Phi}(E)\subset D(\ _{\alpha}H,\mu)$$ Using
\cite{Vae}, we get that:
$$J_{\Phi}\Lambda_{\Phi}(e_nx)=\delta^{-\frac{1}{2}}J_{\Psi}\Lambda_{\Psi}(e_nx)$$
Since
$e_nx\delta^{-\frac{1}{2}}=(e_n\delta^{-\frac{1}{2}})\delta^{\frac{1}{2}}x\delta^{-\frac{1}{2}}$
and, by the previous theorem, that
$\delta^{\frac{1}{2}}x\delta^{-\frac{1}{2}}$ is a bounded operator
in ${\mathcal N}_{T_R}$, so is $e_nx\delta^{-\frac{1}{2}}$ and
therefore, we have:
$$\delta^{-\frac{1}{2}}J_{\Psi}\Lambda_{\Psi}(e_nx)
=\lambda^{\frac{1}{4}}J_{\Psi}\Lambda_{\Psi}(e_nx\delta^{-\frac{1}{2}})\subset
J_{\Psi}\Lambda_{\Psi}({\mathcal N}_{\Psi}\cap {\mathcal N}_{T_R})$$
and we get that $J_{\Phi}\Lambda_{\Phi}(e_nx)$ belongs to
$D(H_{\beta},\mu^0)$. By linearity, we get the result.
\end{proof}

\section{Manageability of the fundamental unitary}\label{mtheo}

In this section, we prove that the fundamental unitary satisfies a
proposition similar to Woronowicz's manageability of \cite{W}.
Following \cite{E2} (definition 4.1), we define the notion of weakly
regular pseudo-multiplicative unitary which is interesting by itself
but it will be useful for us to get easily von Neumann algebra
structure on the dual structure.

\begin{defi}
We call \textbf{manageable operator} the strictly positive operator
$P$ on $H_{\Phi}$ such that
$P^{it}\Lambda_{\Phi}(x)=\lambda^{\frac{t}{2}}\Lambda_{\Phi}(\tau_t(x))$,
for all $x\in {\mathcal N}_{\Phi}$ and  $t\in\mathbb{R}$.
\end{defi}

\begin{prop}
For all $m\in M$, $n\in N$ and $t\in\mathbb{R}$, we have:
$$
\begin{aligned}
&P^{it}mP^{-it}=\tau_t(m)
&\quad P^{it}\alpha(n)P^{-it}=\alpha(\sigma^{\nu}_t(n))\\
&P^{it}\beta(n)P^{-it}=\beta(\sigma^{\nu}_t(n)) &\quad
P^{it}\hat{\beta}(n)P^{-it}=\hat{\beta}(\sigma^{\nu}_t(n))
\end{aligned}$$
\end{prop}

\begin{proof}
Straightforward.
\end{proof}

\noindent Then, we can define operators $P^{it}\surl{\
_{\beta}\otimes_{\alpha}}_{\ \nu}P^{it}$ on $H_{\Phi}\surl{\
_{\beta} \otimes_{\alpha}}_{\ \nu}H_{\Phi}$ and $P^{it}\surl{\
_{\alpha}\otimes_{\hat{\beta}}}_{\ \nu^o}P^{it}$ on $H_{\Phi}\surl{\
_{\alpha}\otimes_{\hat{\beta}}}_{\ \nu^o}H_{\Phi}$ for all
$t\in\mathbb{R}$.

\begin{theo}
The unitary $W$ satisfies a manageability relation. More exactly, we
have:
$$(\sigma_{\nu}W^*\sigma_{\nu}(q\surl{\ _{\hat{\beta}}
\otimes_{\alpha}}_{\ \nu}v)|p\surl{\ _{\alpha}\otimes_{\beta}}_{\
\nu^o}w)=(\sigma_{\nu^o}W\sigma_{\nu^o}(J_{\Phi}p\surl{\
_{\alpha}\otimes_{\beta}}_{\ \nu^o}P^{-1/2}v)| J_{\Phi}q\surl{\
_{\hat{\beta}} \otimes_{\alpha}}_{\ \nu}P^{1/2}w)$$ for all $v\in
{\mathcal D}(P^{-\frac{1}{2}})$, $w\in {\mathcal
D}(P^{\frac{1}{2}})$ and $p,q\in D(_{\alpha}H_{\Phi},\nu)\cap
D((H_{\Phi})_{\hat{\beta}},\nu^o)$. Moreover, for all
$t\in\mathbb{R}$, we have $W(P^{it}\surl{\
_{\beta}\otimes_{\alpha}}_{\ \nu}P^{it})=(P^{it}\surl{\
_{\alpha}\otimes_{\hat{\beta}}}_{\ \nu^o}P^{it})W$.  \label{mania}
\end{theo}

\begin{proof}
Let $p,q\in D(_{\alpha}H_{\Phi},\nu)\cap
D((H_{\Phi})_{\hat{\beta}},\nu^o)$. For all $v\in {\mathcal
D}(D^{1/2})$ and $w\in {\mathcal D}(D^{-1/2})$, we know that:
$$(I(id*\omega_{q,p})(W)Iv|w)=((id*\omega_{p,q})(W)P^{1/2}v|P^{-1/2}w)$$
for all $v\in {\mathcal D}(P^{1/2})$ and $w\in {\mathcal
D}(P^{-1/2})$. By \ref{PUJ}, we rewrite the formula:
$$(\sigma_{\nu}W^*\sigma_{\nu}(q\surl{\ _{\hat{\beta}}
\otimes_{\alpha}}_{\ \nu}v)|p\surl{\ _{\alpha}\otimes_{\beta}}_{\
\nu^o}w)=(\sigma_{\nu^o}W\sigma_{\nu^o}(J_{\Phi}p\surl{\
_{\alpha}\otimes_{\beta}}_{\ \nu^o}P^{-1/2}v)| J_{\Phi}q\surl{\
_{\hat{\beta}} \otimes_{\alpha}}_{\ \nu}P^{1/2}w)$$ Now, we have to
prove $W^*(P^{it}\surl{\ _{\alpha}\otimes_{\hat{\beta}}}_{\
\nu^o}P^{it})= (P^{it}\surl{\ _{\beta}\otimes_{\alpha}}_{\
\nu}P^{it})W^*$ for all $t\in\mathbb{R}$. First of all, because of
the commutation relation between $P$ and $\beta$,
$D((H_{\Phi})_{\beta},\nu^o)$ is $P^{it}$-invariant and if
$(\xi_i)_{i\in I}$ is a $(N^o,\nu^o)$-basis of $(H_{\Phi})_{\beta}$,
then $(P^{it}\xi_i)_{i\in I}$ is also. Let $v\in
D((H_{\Phi})_{\beta},\nu^o)$ and $a\in {\mathcal N}_{T_L}\cap
{\mathcal N}_{\Phi}$. We compute:

$$
\begin{aligned}
&\ \quad (P^{it}\surl{\ _{\beta}\otimes_{\alpha}}_{\ \nu}P^{it})
W^*(v\surl{\ _{\alpha}\otimes_{\hat{\beta}}}_{\
\nu^o}\Lambda_{\Phi}(a))
\\&=\sum_{i\in I}P^{it}\xi_i\surl{\
_{\beta} \otimes_{\alpha}}_{\
\nu}\lambda^{t/2}\Lambda_{\Phi}(\tau_t((\omega_{v,\xi_i}
\surl{\ _{\beta} \star_{\alpha}}_{\ \nu}id)(\Gamma(a))))\\
&=\sum_{i\in I}P^{it}\xi_i\surl{\ _{\beta} \otimes_{\alpha}}_{\
\nu}\Lambda_{\Phi}((\omega_{P^{it}v,P^{it}\xi_i}
\surl{\ _{\beta} \star_{\alpha}}_{\ \nu}id)(\Gamma(\lambda^{t/2}\tau_t(a))))\\
&=W^*(P^{it}v\surl{\ _{\alpha}\otimes_{\hat{\beta}}}_{\
\nu^o}\lambda^{t/2} \Lambda_{\Phi}(\tau_t(a)))=W^*(P^{it}\surl{\
_{\alpha}\otimes_{\hat{\beta}}}_{\ \nu^o}P^{it})(v\surl{\
_{\alpha}\otimes_{\hat{\beta}}}_{\ \nu^o}\Lambda_{\Phi}(a))
\end{aligned}$$

\end{proof}

\begin{defi}
A pseudo-multiplicative unitary ${\mathcal W}$ w.r.t
$\alpha,\beta,\hat{\beta}$ is said to be \textbf{weakly regular} if
the weakly closed linear span of
$(\lambda_v^{\alpha,\beta})^*\mathcal{W}\rho_w^{\hat{\beta},\alpha}$
where $v,w$ belongs to $D(_{\alpha}H,\nu)$ is equal to $\alpha(N)'$.
\end{defi}

\begin{prop}\label{reg}
The operator $\widehat{W}=\sigma_{\nu}W^*\sigma_{\nu}$ from
$H_{\Phi}\surl{\ _{\hat{\beta}}\otimes_{\alpha}}_{\ \nu}H_{\Phi}$
onto $H_{\Phi}\surl{\ _{\alpha}\otimes_{\beta}}_{\ \nu^o}H_{\Phi}$
is a pseudo-multiplicative unitary over $N$ w.r.t
$\alpha,\beta,\hat{\beta}$ which is weakly regular in the sense of
\cite{E2} (definition 4.1).
\end{prop}

\begin{proof}
By \cite{EV}, we know that $\widehat{W}$ is a pseudo-multiplicative
unitary. We also know that
$<(\lambda_v^{\alpha,\beta})^*\widehat{W}\rho_w^{\hat{\beta},\alpha}>
^{-\textsc{w}}\subset\alpha(N)'$. For all $v\in {\mathcal
D}(P^{-\frac{1}{2}})$, $w\in {\mathcal D}(P^{\frac{1}{2}})$ and
$p,q\in D(_{\alpha}H_{\Phi},\nu)\cap
D((H_{\Phi})_{\hat{\beta}},\nu^o)$, we have, by theorem \ref{mania}:
$$((\lambda_p^{\alpha,\beta})^*\widehat{W}\rho_v^{\hat{\beta},\alpha}q|w)
=(\sigma_{\nu^o}W\sigma_{\nu^o}(J_{\Phi}p\surl{\
_{\alpha}\otimes_{\beta}}_{\ \nu^o}P^{-1/2}v)| J_{\Phi}q\surl{\
_{\hat{\beta}} \otimes_{\alpha}}_{\ \nu}P^{1/2}w)$$ and on the other
hand:
$$
\begin{aligned}
(R^{\alpha,\nu}(v)R^{\alpha,\nu}(p)^*q|w)
&=(R^{\alpha,\nu}(v)J_{\nu}R^{\hat{\beta},\nu^o}(J_{\Phi}p)^*J_{\Phi}q|w)\\
&=(R^{\alpha,\nu}(v)J_{\nu}\Lambda_{\nu}(<J_{\Phi}q,J_{\Phi}p>_{\hat{\beta},\nu^o_L})|w)\\
&=(P^{-1/2}R^{\alpha,\nu}(v)J_{\nu}\Lambda_{\nu}(<J_{\Phi}q,J_{\Phi}p>_{\hat{\beta},\nu^o_L})|P^{1/2}w)\\
&=(R^{\alpha,\nu}(P^{-1/2}v)\Delta_{\nu}^{-1/2}J_{\nu}\Lambda_{\nu}(<J_{\Phi}q,J_{\Phi}p>_{\hat{\beta},\nu^o_L})|P^{1/2}w)\\
&=(R^{\alpha,\nu}(P^{-1/2}v)\Lambda_{\nu}(<J_{\Phi}p,J_{\Phi}q>_{\hat{\beta},\nu^o_L})|P^{1/2}w)\\
&=(\alpha(<J_{\Phi}p,J_{\Phi}q>_{\hat{\beta},\nu^o_L})P^{-1/2}v|P^{1/2}w)\\
&=(J_{\Phi}p\surl{\ _{\hat{\beta}} \otimes_{\alpha}}_{\
\nu}P^{-1/2}v| J_{\Phi}q\surl{\ _{\hat{\beta}} \otimes_{\alpha}}_{\
\nu}P^{1/2}w)
\end{aligned}$$
There exists $\Xi\in H_{\Phi}\surl{\
_{\hat{\beta}}\otimes_{\alpha}}_{\ \nu}H_{\Phi}$ such that
$\sigma_{\nu^o}W\sigma_{\nu^o}\Xi=J_{\Phi}p\surl{\ _{\hat{\beta}}
\otimes_{\alpha}}_{\ \nu}P^{-1/2}v$ since $W$ is onto. By
definition, there exists a net
$(\sum_{k=1}^{n(i)}J_{\Phi}p_k^i\surl{\ _{\alpha}\otimes_{\beta}}_{\
\nu^o}P^{-1/2}v_k^i)_{i\in I}$ which converges to $\Xi$. Then
$((\sum_{k=1}^{n(i)}(\lambda_{p_k^i}^{\alpha,\beta})^*
\widehat{W}\rho_{v_k^i}^{\hat{\beta},\alpha}q|w))_{i\in I}$
converges to:
$$\begin{aligned}
(\sigma_{\nu^o}W\sigma_{\nu^o}\Xi|J_{\Phi}q\surl{\ _{\hat{\beta}}
\otimes_{\alpha}}_{\ \nu}P^{1/2}w)&=(J_{\Phi}p\surl{\ _{\hat{\beta}}
\otimes_{\alpha}}_{\ \nu}P^{-1/2}v| J_{\Phi}q\surl{\ _{\hat{\beta}}
\otimes_{\alpha}}_{\
\nu}P^{1/2}w)\\
&=(R^{\alpha,\nu}(v)R^{\alpha,\nu}(p)^*q|w)
\end{aligned}$$

Then, we obtain $\alpha(N)'=
<R^{\alpha,\nu}(v)R^{\alpha,\nu}(p)^*>^{-\textsc{w}}\subset
<(\omega_{v,p}*id)(\widehat{W}\sigma_{\nu^o})>^{-\textsc{w}}$.
\end{proof}

\section{Duality}
In this section, a dual measured quantum groupoid is constructed
thanks to modulus and scaling operator. Then, we obtain a bi-duality
theorem which generalizes Pontryagin duality, locally compact
quantum groups duality and duality for groupoids. Finally, we get
Heisenberg's relations.

\subsection{Dual structure}

\begin{defi}
The weak closure of the linear span of $(\omega_{\xi,\eta}*id)(W)$,
where $\xi\in D((H_{\Phi})_{\beta},\nu^o)$ and $\eta\in
D(_{\alpha}H_{\Phi},\nu)$, is denoted by $\widehat{M}$. It's a von
Neumann algebra because weak regularity of $\hat{W}$ (prop.
\ref{reg}) and \cite{E2} (proposition 3.2).
\end{defi}

\begin{defi}
We put $\hat{\Gamma}$ the application from $\widehat{M}$ into
$\mathcal{L}(H_{\Phi}\surl{\ _{\hat{\beta}}\otimes_{\alpha}}_{\
\nu}H_{\Phi})$ such that, for all $x\in\widehat{M}$, we have:
$$\hat{\Gamma}(x)=\sigma_{\nu^o}W(x\surl{\
_{\beta}\otimes_{\alpha}}_{\ N}1)W^*\sigma_{\nu}$$
\end{defi}

\begin{prop}\label{deb1}
The 5-uple $(N,\widehat{M},\alpha,\hat{\beta},\hat{\Gamma})$ is a
Hopf-bimodule called dual Hopf-bimodule.
\end{prop}

\begin{proof}
The proposition comes from theorems 6.2 and 6.3 of \cite{EV} applied
to $\widehat{W}=\sigma_{\nu}W^*\sigma_{\nu}$.
\end{proof}

\begin{lemm}
Let call $M_*^{\alpha,\beta}$ the subspace of $M_*$ spanned by the
positive and normal forms such that there exists $k\in\mathbb{R}^+$
and both $\omega\circ\alpha$ and $\omega\circ\beta$ are dominated by
$k\nu$. Then, $M_*^{\alpha,\beta}$ is dense *-subalgebra of $M_*$
such that, for all $m\in M$, we have:
$$\omega\mu(m)=\mu((\omega\surl{\ _{\beta} \star_{\alpha}}_{\ \nu}
id)(\Gamma(m)))\quad\text{ and }\quad
\omega^*(m)=\overline{\omega\circ R (m^*)}$$
\end{lemm}

\begin{proof}
By definition $\omega\mu$ belongs to $M_*$. There exists $\xi\in
D(H_{\beta},\nu^o)$ such that $\omega=\omega_{\xi}$. For all $n\in
N$, we have:
$$
\begin{aligned}
&\ \quad\omega\mu(\alpha(n^*n))=\mu((\omega_{\xi}\surl{\ _{\beta}
\star_{\alpha}}_{\ \nu}id)(\Gamma(\alpha(n^*n))))\\
&=\mu((\lambda_{\xi}^{\beta,\nu^0})^*(\alpha(n^*n)\surl{\ _{\beta}
\otimes_{\alpha}}_{\ N}1)\lambda_{\xi}^{\beta,\nu^0})
=\mu((\lambda_{\alpha(n)\xi}^{\beta,\nu^0})^*\lambda_{\alpha(n)\xi}^{\beta,\nu^0})\\
&=\mu\circ\alpha(<\alpha(n)\xi,\alpha(n)\xi>_{\beta,\nu^0})\leq
k\nu(<\alpha(n)\xi,\alpha(n)\xi>_{\beta,\nu^0})=k||\alpha(n)\xi||^2\\
&=k\omega\circ\alpha(n^*n)\leq k^2\nu(n^*n)
\end{aligned}$$
Also, we can prove that $\omega\mu\circ\beta$ is dominated by
$k^2\nu$ so that $\omega\mu$ belongs to $M_*^{\alpha,\beta}$. Since
$R\circ\alpha=\beta$, $M_*^{\alpha,\beta}$ is *-stable. We have to
prove associativity of product and that
$(\omega\mu)^*=\mu^*\omega^*$. The first property comes from
co-associativity of co-product and the second one comes from
co-involution property. We only check the first one because the
second proof is very similar computation. Let $\omega,\mu,\chi\in
M_*^{\alpha,\beta}$ and $\xi,\xi',\xi''\in D(H_{\beta},\nu^o)$ the
corresponding vectors. Then, for all $m\in M$, it is easy to see
that:
$$
\begin{aligned}
(\omega\mu)\chi(x)&=((\Gamma\surl{\ _{\beta} \star_{\alpha}}_{\
\nu}id)(\Gamma(x))(\xi\surl{\ _{\beta} \otimes_{\alpha}}_{\
\nu}\xi'\surl{\ _{\beta} \otimes_{\alpha}}_{\ \nu}\xi'')|\xi\surl{\
_{\beta} \otimes_{\alpha}}_{\ \nu}\xi'\surl{\ _{\beta}
\otimes_{\alpha}}_{\
\nu}\xi'')\\
&=((id\surl{\ _{\beta} \star_{\alpha}}_{\
\nu}\Gamma)(\Gamma(x))(\xi\surl{\ _{\beta} \otimes_{\alpha}}_{\
\nu}\xi'\surl{\ _{\beta} \otimes_{\alpha}}_{\ \nu}\xi'')|\xi\surl{\
_{\beta} \otimes_{\alpha}}_{\ \nu}\xi'\surl{\ _{\beta}
\otimes_{\alpha}}_{\
\nu}\xi'')\\
&=(\Gamma((\omega_{\xi}\surl{\ _{\beta} \star_{\alpha}}_{\ \nu}
id)(\Gamma(x)))(\xi'\surl{\ _{\beta} \otimes_{\alpha}}_{\
\nu}\xi'') |\xi'\surl{\ _{\beta} \otimes_{\alpha}}_{\ \nu}\xi'')\\
&=\mu\chi((\omega_{\xi}\surl{\ _{\beta} \star_{\alpha}}_{\ \nu}
id)(\Gamma(x)))=\omega(\mu\chi)(x)
\end{aligned}$$ Density condition comes from corollary \ref{formi}
for example.
\end{proof}

\begin{coro}
The contractive application $\widehat{\pi}$ from
$M_*^{\alpha,\beta}$ to $\widehat{M}$ such that
$\widehat{\pi}(\omega)=(\omega*id)(W)$ is 1-1 and multiplicative.
\end{coro}

\begin{proof}
The application $\widehat{\pi}$ is injective because of theorem
\ref{densevn1}. We prove multiplicativity of $\widehat{\pi}$ for
positive linear forms because the general case comes then from
linearity. Let $\xi,\eta\in D(_{\alpha}H,\nu)\cap
D(H_{\beta},\nu^o)$, $\zeta_1\in D(_{\alpha}H,\nu)$ and $\zeta_2\in
D(H_{\hat{\beta}},\nu^o)$. By proposition \ref{calcule} of the first
part, we know that:
$$((\omega_{\xi}*id)(W)(\omega_{\eta}*id)(W)\zeta_1
|\zeta_2)$$ is equal to the scalar product of
$$(\sigma_{\nu^o}\!\! \surl{\ _{\alpha}
  \otimes_{\hat{\beta}}}_{\ \ N^o} 1)(1\!\! \surl{\ _{\alpha}
  \otimes_{\hat{\beta}}}_{\ \ N^o} W)\sigma_{2\nu}(1\!\! \surl{\ _{\beta}
  \otimes_{\alpha}}_{\ N} \sigma_{\nu^o})(1\!\! \surl{\ _{\beta}
  \otimes_{\alpha}}_{\ N} W)(\xi\surl{\
_{\beta} \otimes_{\alpha}}_{\ \nu}\eta\surl{\ _{\beta}
\otimes_{\alpha}}_{\ \nu}\zeta_1)$$ by $[\xi\surl{\ _{\beta}
\otimes_{\alpha}}_{\ \nu}\eta]\surl{\ _{\alpha}
\otimes_{\hat{\beta}}}_{\ \nu^o}\zeta_2$. Then, by
pseudo-multiplicativity of $W$, this equal to:
$$((W^*\surl{\ _{\alpha} \otimes_{\hat{\beta}}}_{\
N^o} 1)(1\surl{\ _{\alpha} \otimes_{\hat{\beta}}}_{\ N^o}
W)(W\surl{\ _{\beta} \otimes_{\alpha}}_{\ N} 1)(\xi\surl{\ _{\beta}
\otimes_{\alpha}}_{\ \nu}\eta\surl{\ _{\beta} \otimes_{\alpha}}_{\
\nu}\zeta_1)|[\xi\surl{\ _{\beta} \otimes_{\alpha}}_{\
\nu}\eta]\surl{\ _{\alpha} \otimes_{\hat{\beta}}}_{\
\nu^o}\zeta_2)$$
$$
\begin{aligned}
&=((1\surl{\ _{\alpha} \otimes_{\hat{\beta}}}_{\ N^o} W)
W(\xi\surl{\ _{\beta} \otimes_{\alpha}}_{\ \nu}\eta)\surl{\ _{\beta}
\otimes_{\alpha}}_{\ \nu}\zeta_1|W(\xi\surl{\ _{\beta}
\otimes_{\alpha}}_{\ \nu}\eta)\surl{\ _{\alpha}
\otimes_{\hat{\beta}}}_{\ \nu^o}\zeta_2)\\
&=((1\surl{\ _{\alpha}\otimes_{\hat{\beta}}}_{\
N^o}(id*\omega_{\zeta_1,\zeta_2})(W))W(\xi\surl{\ _{\beta}
\otimes_{\alpha}}_{\ \nu}\eta)|W(\xi\surl{\ _{\beta}
\otimes_{\alpha}}_{\ \nu}\eta))
\end{aligned}$$ Since $\Gamma$ is
implemented by $W$, this is equal to:
$$
\begin{aligned}
(\Gamma((id*\omega_{\zeta_1,\zeta_2})(W))(\xi\surl{\ _{\beta}
\otimes_{\alpha}}_{\ \nu}\eta)|(\xi\surl{\ _{\beta}
\otimes_{\alpha}}_{\
\nu}\eta))&=(\omega_{\xi}\omega_{\eta})((id*\omega_{\zeta_1,\zeta_2})(W))\\
&=(((\omega_{\xi}\omega_{\eta})*id)(W)\zeta_1|\zeta_2)
\end{aligned}$$ By density of $D(_{\alpha}H,\nu)$ and
$D(H_{\hat{\beta}},\nu^o)$ in $H$, we get that $\hat{\pi}$ is
multiplicative.
\end{proof}

To get a measured quantum groupoid from the dual Hopf-bimodule, we
have to exhibit, first of all, a co-involution. This is done and the
following proposition:

\begin{prop}\label{deb3}
There exists a unique *-anti-automorphism $\widehat{R}$ of
$\widehat{M}$ such that, for all $\omega\in M_*^{\alpha,\beta}$, we
have $\widehat{R}(\widehat{\pi}(\omega))=\widehat{\pi}(\omega\circ
R)$. Moreover $\widehat{R}(x)=J_{\Phi}x^*J_{\Phi}$ for all
$x\in\widehat{M}$ and $\widehat{R}$ is a co-involution.
\end{prop}

\begin{proof}
For all $\xi\in D(_{\alpha}(H_{\Phi}),\nu)$ and $\eta\in
D((H_{\Phi})_{\beta},\nu^o)$, we have:
$$
\begin{aligned}
&\ \quad (J_{\Phi}\hat{\pi}(\omega\circ
R)^*J_{\Phi}\xi|\eta)=(\hat{\pi}(\omega\circ
R)J_{\Phi}\eta|J_{\Phi}\xi)=((\omega\circ R*id)(W)J_{\Phi}\eta|J_{\Phi}\xi)\\
&=\omega\circ R((id*\omega_{J_{\Phi}\eta,J_{\Phi}\xi}(W))
=\omega((id*\omega_{\xi,\eta})(W))=(\hat{\pi}(\omega)\xi|\eta)
\end{aligned}$$ So, if we define $\widehat{R}$ by $\widehat{R}(x)=J_{\Phi}x^*J_{\Phi}$ for all
$x\in\widehat{M}$, we obtain a *-anti-automorphism of $\widehat{M}$
such that, for all $\omega\in M_*^{\alpha,\beta}$, we have
$\widehat{R}(\widehat{\pi}(\omega))=\widehat{\pi}(\omega\circ R)$.
Uniqueness comes from density of $\widehat{\pi}(M_*^{\alpha,\beta})$
in $\widehat{M}$. By definition, we have
$\widehat{R}\circ\alpha=\hat{\beta}$. So, we have to check
co-involution property to finish the proof. For all $\omega\in
M_*^{\alpha,\beta}$, we compute:
$$
\begin{aligned}
\hat{\Gamma}(\hat{\pi}(\omega))&=\hat{W}^*(1\surl{\ _{\alpha}
\otimes_{\beta}}_{\ N^o}(\omega
*id)(W))\hat{W}=\sigma_{\nu^o}W((\omega *id)(W)\surl{\ _{\beta}
\otimes_{\alpha}}_{\ N}1)W^*\sigma_{\nu}\\
&=\sigma_{\nu^o}(\omega*id*id)((1\surl{\ _{\alpha}
\otimes_{\hat{\beta}}}_{\ N^o}W)(W\surl{\ _{\beta}
\otimes_{\alpha}}_{\ N}1)(1\surl{\ _{\beta} \otimes_{\alpha}}_{\
N}W^*))\sigma_{\nu}
\end{aligned}$$
By pseudo-multiplicativity of $W$, this is equal to:
$$
\begin{aligned}
&\ \quad\sigma_{\nu^o}(\omega*id*id)((W \!\!\surl{\ _{\alpha}
  \otimes_{\hat{\beta}}}_{\ \ N^o} 1)(\sigma_{\nu^o}\!\! \surl{\ _{\alpha}
  \otimes_{\hat{\beta}}}_{\ \ N^o} 1)(1\!\! \surl{\ _{\alpha}
  \otimes_{\hat{\beta}}}_{\ \ N^o} W)\sigma_{2\nu}(1\!\! \surl{\ _{\beta}
  \otimes_{\alpha}}_{\ N} \sigma_{\nu^o}))\sigma_{\nu}\\
&=(\omega*id*id)((1\surl{\ _{\alpha}
  \otimes_{\hat{\beta}}}_{\ \ N^o}\sigma_{\nu^o})(W\surl{\ _{\alpha}
  \otimes_{\hat{\beta}}}_{\ \ N^o}1)(\sigma_{\nu^o}\surl{\ _{\alpha}
  \otimes_{\hat{\beta}}}_{\ \ N^o}1)(1\surl{\ _{\alpha}
  \otimes_{\hat{\beta}}}_{\ \ N^o}W)\sigma_{2\nu})\\
\end{aligned}$$
Then, we get:
$$\hat{\Gamma}\circ\hat{R}(\hat{\pi}(\omega))=\hat{\Gamma}(\hat{\pi}(\omega\circ
R))$$
$$=(\omega\circ R*id*id)((1\surl{\ _{\alpha}
  \otimes_{\hat{\beta}}}_{\ \ N^o}\sigma_{\nu^o})(W\surl{\ _{\alpha}
  \otimes_{\hat{\beta}}}_{\ \ N^o}1)(\sigma_{\nu^o}\surl{\ _{\alpha}
  \otimes_{\hat{\beta}}}_{\ \ N^o}1)(1\surl{\ _{\alpha}
  \otimes_{\hat{\beta}}}_{\ \ N^o}W)\sigma_{2\nu})$$
Now, by proposition \ref{PUJ}, we know that: $W=(I\surl{\
_{\beta}\otimes_{\alpha}}_{\ \ N}J_{\Phi})W^*(I\surl{\ _{\beta}
  \otimes_{\alpha}}_{\ \ N}J_{\Phi})$ so that:
$$(1\surl{\ _{\alpha}
  \otimes_{\hat{\beta}}}_{\ \ N^o}\sigma_{\nu^o})(W\surl{\ _{\alpha}
  \otimes_{\hat{\beta}}}_{\ \ N^o}1)(\sigma_{\nu^o}\surl{\ _{\alpha}
  \otimes_{\hat{\beta}}}_{\ \ N^o}1)(1\surl{\ _{\alpha}
  \otimes_{\hat{\beta}}}_{\ \ N^o}W)\sigma_{2\nu}$$
$$=(I\!\!\surl{\ _{\alpha}
  \otimes_{\beta}}_{\ \ N^o}\! J_{\Phi}\!\!\surl{\ _{\alpha}
  \otimes_{\hat{\beta}}}_{\ \ N^o}\! J_{\Phi})[(W \!\!\surl{\ _{\alpha}
  \otimes_{\hat{\beta}}}_{\ \ N^o}\!\! 1)(\sigma_{\nu^o}\!\! \surl{\ _{\alpha}
  \otimes_{\hat{\beta}}}_{\ \ N^o}\!\! 1)(1\!\! \surl{\ _{\alpha}
  \otimes_{\hat{\beta}}}_{\ \ N^o}\!\! W)\sigma_{2\nu}(1\!\! \surl{\ _{\beta}
  \otimes_{\alpha}}_{\ N}\!\! \sigma_{\nu^o})]^*(I\!\!\surl{\ _{\beta}
  \otimes_{\alpha}}_{\ \ N}\! J_{\Phi}\!\!\surl{\ _{\hat{\beta}}
  \otimes_{\alpha}}_{\ \ N}\! J_{\Phi})$$ Since $R$ is implemented
by $I$ and $\widehat{R}$ is implemented by $J_{\Phi}$, we have:

$$
\begin{aligned}
&\
\quad\hat{\Gamma}\circ\hat{R}(\hat{\pi}(\omega))\\
&=(\hat{R}\surl{\ _{\alpha}\star_{\beta}}_{\ \
N^o}\hat{R})((\omega*id)[(W \!\!\surl{\ _{\alpha}
  \otimes_{\hat{\beta}}}_{\ \ N^o} 1)(\sigma_{\nu^o}\!\! \surl{\ _{\alpha}
  \otimes_{\hat{\beta}}}_{\ \ N^o} 1)(1\!\! \surl{\ _{\alpha}
  \otimes_{\hat{\beta}}}_{\ \ N^o} W)\sigma_{2\nu}(1\!\! \surl{\ _{\beta}
  \otimes_{\alpha}}_{\ N} \sigma_{\nu^o})])\\
&=(\hat{R}\surl{\ _{\alpha}\star_{\hat{\beta}}}_{\ \
N^o}\hat{R})\circ\varsigma_N\circ\hat{\Gamma}(\hat{\pi}(\omega))
=\varsigma_{N^o}\circ(\hat{R}\surl{\
_{\hat{\beta}}\star_{\alpha}}_{\ \
N}\hat{R})\circ\hat{\Gamma}(\hat{\pi}(\omega))
\end{aligned}$$
A density argument enables us to conclude.
\end{proof}

Then, we have to construct a left-invariant operator-valued weight
$\widehat{T_L}$ from $\widehat{M}$ to $\alpha(N)$. We follow J.
Kustermans and S. Vaes' paper \cite{KV1}: we define in fact a GNS
construction $(H,\iota,\widehat{\Lambda})$ and we give a core for
$\widehat{\Lambda}$. Let introduce the space $\mathcal{I}$ of
$\omega\in M_*^{\alpha,\beta}$ such that there exists
$k\in\mathbb{R}^+$ and $|\omega(x^*)|\leq k||\Lambda_{\Phi}(x)||$
for all $x\in {\mathcal N}_{\Phi}\cap {\mathcal N}_{T_L}$. Then, by
Riesz' theorem, there exists $\xi(\omega)\in H$ such that:
$$\omega(x^*)=(\xi(\omega)|\Lambda_{\Phi}(x))$$

\begin{lemm}
The set $\{\xi(\omega)|\ \omega\in\mathcal{I}\}$ is dense in $H$.
\end{lemm}

\begin{proof}
Let $a,b\in E$ define in corollary \ref{formi}. Then
$\omega_{\Lambda_{\Phi}(a),\Lambda_{\Phi}(b)}$ belongs to
$M_*^{\alpha,\beta}$ and we have, for all $x\in {\mathcal
N}_{\Phi}\cap {\mathcal N}_{T_L}$:
$$\omega_{\Lambda_{\Phi}(a),\Lambda_{\Phi}(b)}(x^*)=\Phi(b^*x^*a)
=\Phi(x^*a\sigma_{-i}^{\Phi}(b^*))=(\Lambda_{\Phi}(a\sigma_{-i}^{\Phi}(b^*))|\Lambda_{\Phi}(x))$$
so that $\omega_{\Lambda_{\Phi}(a),\Lambda_{\Phi}(b)}$ belongs to
$\mathcal{I}$ and we have
$\xi(\omega_{\Lambda_{\Phi}(a),\Lambda_{\Phi}(b)})=\Lambda_{\Phi}(a\sigma_{-i}^{\Phi}(b^*))$
which is dense in $H$.
\end{proof}

In the following, for all form $\omega$, we denote by
$\overline{\omega}$ the form such that
$\overline{\omega}(x)=\overline{\omega(x^*)}$. Observe that
$\omega\in M_*^{\alpha,\beta}$ implies that $\overline{\omega}$
belongs also to $M_*^{\alpha,\beta}$.

\begin{prop}
The space $\mathcal{I}$ is a dense left ideal of
$M_*^{\alpha,\beta}$ such that, for all $\omega\in
M_*^{\alpha,\beta}$ and $\mu\in\mathcal{I}$, we have:
$$\xi(\omega\mu)=\hat{\pi}(\omega)\xi(\mu)$$
\end{prop}

\begin{proof}
If $\xi,\eta$ belong to $D(_{\alpha}H,\nu)\cap D(H_{\beta},\nu^0)$,
then $\omega_{\xi,\eta}$ belongs to $M_*^{\alpha,\beta}$. Moreover,
if $\eta$ belongs also to
$D(_{\text{id}}H_{\Phi},\Phi)=J_{\Phi}\Lambda_{\Phi}({\mathcal
N}_{\Phi})$, then we have:
$$|\omega_{\xi,\eta}(x^*)|=|(\xi|x\eta)|\leq||\xi||||x\eta||\leq k||\xi||||\Lambda_{\Phi}(x)||$$
so that, by corollary \ref{formi}, we can deduce that $\mathcal{I}$
is dense in $M_*^{\alpha,\beta}$ and therefore in $M_*$. Now, for
all $x\in {\mathcal N}_{\Phi}\cap {\mathcal N}_{T_L}$, we have:
$$
\begin{aligned}
\omega\mu(x^*)&=\mu((\omega\surl{\ _{\beta} \star_{\alpha}}_{\
\nu}id)\Gamma(x^*))=\mu(((\overline{\omega}\surl{\ _{\beta}
\star_{\alpha}}_{\
\nu}id)\Gamma(x))^*)\\&=(\xi(\mu)|\Lambda_{\Phi}((\overline{\omega}\surl{\
_{\beta} \star_{\alpha}}_{\ \nu}id)\Gamma(x)))
=(\xi(\mu)|(\overline{\omega}*id)(W^*)\Lambda_{\Phi}(x))\\
&= ((\omega *id)(W)\xi(\mu)|\Lambda_{\Phi}(x))
\end{aligned}$$ so that the proposition holds.
\end{proof}

\begin{defi}
For all $t\in\mathbb{R}$ and $\omega\in M_*$, we define elements of
$M_*$ such that, for all $x\in M$:
$$\tau^*_t(\omega)(w)=\omega\circ\tau_t(x),\ \quad
\delta^*_t(\omega)(x)=\omega(\delta^{it}x),\ \quad\text{ and }\
\quad\rho_t(\omega)(x)=\omega(\delta^{-it}\tau_{-t}(x))$$
\end{defi}

\begin{prop}\label{forP}
The applications $\tau^*,\delta^*$ and $\rho$ define strongly
continuous one-parameter groups of *-automorphisms of
$M_*^{\alpha,\beta}$. Moreover, they leave $\mathcal{I}$ stable and,
for all $t\in\mathbb{R}$ and $\omega\in\mathcal{I}$, we have:
$$\xi(\tau^*_t(\omega))=\lambda^{-\frac{t}{2}}P^{-it}\xi(\omega),\ \quad
\xi(\delta^*_t(\omega))=\lambda^{\frac{t}{2}}J_{\Phi}\delta^{-it}J_{\Phi}\xi(\omega),$$
$$\text{ and}\ \quad
\xi(\rho_t(\omega))=P^{it}J_{\Phi}\delta^{it}J_{\Phi}\xi(\omega)$$
\end{prop}

\begin{proof}
Since $\tau_t(\delta)=\delta$, it is easy to see that $\tau^*$ and
$\delta^*$ commute with each other and, for all $t\in\mathbb{R}$, we
have $\rho_t=\tau^*_{-t}\circ\delta^*_{-t}$ so that the last
statement comes from the two first one. Since $\tau$ is implemented
by $P$, $\tau^*$ defines a strongly continuous one-parameter
representation of $M_*$. It is the same for $\delta^*$. If $\omega$
belongs to $M_*^{\alpha,\beta}$, then there exists $k\in\mathbb{R}+$
such that, for all $t\in\mathbb{R}$, we have:
$$\tau^*_t(\omega)\circ\alpha=\omega\circ\tau_t\circ\alpha
=\omega\circ\alpha\circ\sigma_t^{\nu}\leq
k\nu\circ\sigma_t^{\nu}=k\nu$$ Moreover, there exists $\xi\in
D(_{\alpha}H,\nu)\cap D(H_{\beta},\nu^0)$ such that
$\omega=\omega_{\xi}$ and, for all $t\in\mathbb{R}$ and $n\in N$, we
have:
$$
\begin{aligned}
\delta_t(\omega)(\alpha(n^*n))&=(\delta^{it}\alpha(n^*n)\xi|\xi)
=(\alpha(n)\xi|\alpha(n)\delta^{-it}\xi)\\
&=(\alpha(n)\xi|\delta^{-it}\alpha(\gamma_t\sigma_t^{\nu}(n)\xi))
\end{aligned}$$
so that we get:
$$|\delta_t^*(\omega)(\alpha(n^*n))|\leq
k||\Lambda_{\nu}(n)||^2=k\nu(n^*n)$$ A similar proof with $\beta$
allows us to deduce that $\tau^*,\delta^*$ and $\rho$ belongs to
$M_*^{\alpha,\beta}$ as soon as $\omega$ belongs to
$M_*^{\alpha,\beta}$. It is also straightforward to check that
$\tau^*_t$ is a *-automorphism of $M_*^{\alpha,\beta}$ thanks to
$\Gamma\circ\tau_t=(\tau_t\surl{\ _{\beta} \star_{\alpha}}_{\
N}\tau_t)\circ\Gamma$ and the commutation between $\tau$ and $R$.
Also, it is also straightforward to check that $\delta^*_t$ is a
*-automorphism of $M_*^{\alpha,\beta}$ thanks to
$\Gamma(\delta)=\delta\surl{\ _{\beta} \otimes_{\alpha}}_{\
N}\delta$ and $R(\delta)=\delta^{-1}$. Finally, for all $x\in
{\mathcal N}_{\Phi}\cap {\mathcal N}_{T_L}$, we have, on one hand:
$$
\begin{aligned}
\tau_t(\omega)(x^*)&=\omega(\tau_t(x^*))=\omega\circ\tau_t(x^*)
=(\xi(\omega)|\Lambda_{\Phi}(\tau_t(x)))\\
&=(\xi(\omega)|\lambda^{\frac{-t}{2}}P^{it}\Lambda_{\Phi}(x))
=(\lambda^{\frac{-t}{2}}P^{-it}\xi(\omega)|\Lambda_{\Phi}(x))
\end{aligned}$$ and on the other hand:
$$
\begin{aligned}
\delta_t(\omega)(x^*)&=\omega((x\delta^{-it})^*)
=(\xi(\omega)|\Lambda_{\Phi}(x\delta^{-it}))\\
&=(\xi(\omega)|J_{\Phi}\delta^{-it}J_{\Phi}\lambda^{\frac{t}{2}}\Lambda_{\Phi}(x))
=(\lambda^{\frac{t}{2}}J_{\Phi}\delta^{it}J_{\Phi}\xi(\omega)|\Lambda_{\Phi}(x))
\end{aligned}$$ That finishes the proof.
\end{proof}

\begin{prop}\label{groupe}
There exists unique strongly continuous one-parameter groups
$\widehat{\tau},\widehat{\kappa}$ and $\widehat{\sigma}$ of
*-automorphisms of $\widehat{M}$ such that, for all $t\in\mathbb{R}$
and $\omega\in M_*^{\alpha,\beta}$, we have:
$$\widehat{\tau}_t(\widehat{\pi}(\omega))=\widehat{\pi}(\tau^*_{-t}(\omega)),\ \quad
\widehat{\kappa}_t(\widehat{\pi}(\omega))=\widehat{\pi}(\delta^*_{-t}(\omega))\
\quad\text{ and }\ \quad
\widehat{\sigma}_t(\widehat{\pi}(\omega))=\widehat{\pi}(\rho_t(\omega))
$$
Moreover, for all $t\in\mathbb{R}$ and $x\in\widehat{M}$, the
following properties hold:
\begin{itemize}
\item $\widehat{\tau}_t(x)=P^{it}xP^{-it}$, $\widehat{\kappa}_t(x)=J_{\Phi}\delta^{it}J_{\Phi}x
J_{\Phi}\delta^{-it}J_{\Phi}$\\ and
$\widehat{\sigma_t}(x)=P^{it}J_{\Phi}\delta^{it}J_{\Phi}x
J_{\Phi}\delta^{-it}J_{\Phi}P^{-it}$
\item $\widehat{\tau}$, $\widehat{\kappa}$ and $\widehat{\sigma}$ commute with each other. Also $\widehat{\tau}$ and $\widehat{R}$
do.
\item $\widehat{\kappa}\circ\alpha=\alpha$ and $\widehat{\tau}\circ\alpha=\alpha\circ\sigma_t^{\nu}=\widehat{\sigma}\circ\alpha$
\item $(\widehat{\tau}_t \surl{\ _{\beta} \star_{\alpha}}_{\ N}
\widehat{\tau}_t)\circ\widehat{\Gamma}=\widehat{\Gamma}\circ\widehat{\tau}_t$,
$(id\surl{\ _{\beta} \star_{\alpha}}_{\ N}
\widehat{\kappa}_t)\circ\widehat{\Gamma}=\widehat{\Gamma}\circ\widehat{\kappa}_t$
and $(\widehat{\tau}_t \surl{\ _{\beta} \star_{\alpha}}_{\ N}
\widehat{\sigma}_t)\circ\widehat{\Gamma}=\widehat{\Gamma}\circ\widehat{\sigma}_t$
\end{itemize}

\end{prop}

\begin{proof}
By definition, we have
$\widehat{\sigma}=\widehat{\tau}\circ\widehat{\kappa}=\widehat{\kappa}\circ\widehat{\tau}$
so that we just have to do the proof for $\widehat{\tau}$ and
$\widehat{\kappa}$. For all $\omega\in M^{\alpha,\beta}_*$ and
$t\in\mathbb{R}$, we compute the values of
$P^{it}\hat{\pi}(\omega)P^{-it}$ and
$J_{\Phi}\delta^{it}J_{\Phi}\hat{\pi}(\omega)J_{\Phi}\delta^{-it}J_{\Phi}$.
Let $\mu\in {\mathcal I}$. Since $\hat{\pi}(\omega)$ belongs to
$\beta(N)'$, we have on one hand:
$$
\begin{aligned}
P^{it}\hat{\pi}(\omega)P^{-it}\xi(\mu)
&=P^{it}\hat{\pi}(\omega)\lambda^{\frac{t}{2}}\xi(\tau^*_t(\mu))
=\lambda^{\frac{t}{2}}P^{it}\xi(\omega\tau^*_t(\mu))\\
&=\xi(\tau_{-t}^*(\omega)\mu)=\hat{\pi}(\tau_{-t}^*(\omega))\xi(\mu)
\end{aligned}$$ and on the other hand:
$$
\begin{aligned}
J_{\Phi}\delta^{it}J_{\Phi}\hat{\pi}(\omega)J_{\Phi}\delta^{-it}J_{\Phi}\xi(\mu)
&=J_{\Phi}\delta^{it}J_{\Phi}\hat{\pi}(\omega)\lambda^{\frac{-t}{2}}\xi(\delta^*_t(\mu))\\
&=\xi(\delta^*_{-t}(\omega)\mu)=\hat{\pi}(\delta^*_{-t}(\omega))\xi(\mu)
\end{aligned}$$ So, if we define $\widehat{\tau}_t$ by
$\widehat{\tau}_t(x)=P^{it}xP^{-it}$ and $\widehat{\kappa}_t$ by
$\widehat{\kappa}_t(x)=J_{\Phi}\delta^{it}J_{\Phi}xJ_{\Phi}\delta^{-it}J_{\Phi}$,
then we get strongly continuous *-automorphism of $\widehat{M}$
satisfying the first property. By definition, $\widehat{\tau}$ is
implemented by $P$ and $\widehat{R}$ by $J_{\Phi}$. Since $P$ and
$J_{\Phi}$ commute with each other, so $\widehat{\tau}$ and
$\widehat{R}$ do. Now, $\widehat{\tau}$ and $\tau$ coincide on
$\alpha(N)\subset M\cap\widehat{M}$ because they are both
i:mplemented by $P$. Also $\widehat{\tau}$ coincide with $id$ on
$M\cap\widehat{M}$ by definition. By the way, we can give a meaning
for formulas of the fourth point. Thanks to manageability of $W$, we
have, for all $t\in\mathbb{R}$ and $x\in\widehat{M}$:
$$
\begin{aligned}
\widehat{\Gamma}(\widehat{\tau}_t(x))
&=\sigma_{\nu}W(P^{it}xP^{-it}\surl{\
_{\beta}\otimes_{\alpha}}_{\ N}1)W^*\sigma_{\nu}\\
&=(P^{it}\surl{\ _{\hat{\beta}}\otimes_{\alpha}}_{\
N}P^{it})\sigma_{\nu}W(x\surl{\ _{\beta}\otimes_{\alpha}}_{\
N}1)W^*\sigma_{\nu}(P^{-it}\surl{\
_{\hat{\beta}}\otimes_{\alpha}}_{\
N}P^{-it})\\
&=(\widehat{\tau}_t\surl{\ _{\hat{\beta}}\star_{\alpha}}_{\
N}\widehat{\tau}_t)\widehat{\Gamma}(x)
\end{aligned}$$
Finally, since the left leg of $W$ leaves in $M$, we have:
$$
\begin{aligned}
\widehat{\Gamma}(\widehat{\kappa}_t(x))
&=\sigma_{\nu}W(J_{\Phi}\delta^{it}J_{\Phi}xJ_{\Phi}\delta^{-it}J_{\Phi}\surl{\
_{\beta}\otimes_{\alpha}}_{\ N}1)W^*\sigma_{\nu}\\
&=(1\surl{\ _{\hat{\beta}}\otimes_{\alpha}}_{\
N}J_{\Phi}\delta^{it}J_{\Phi})\sigma_{\nu}W(x\surl{\
_{\beta}\otimes_{\alpha}}_{\ N}1)W^*\sigma_{\nu}(1\surl{\
_{\hat{\beta}}\otimes_{\alpha}}_{\
N}J_{\Phi}\delta^{-it}J_{\Phi})\\
&=(id\surl{\ _{\hat{\beta}}\star_{\alpha}}_{\
N}\widehat{\kappa}_t)\widehat{\Gamma}(x)
\end{aligned}$$
\end{proof}

\begin{lemm}
We have $(\omega R*id)(W^*)=(\tau^*_{-i/2}(\omega)*id)(W)$ for all
$\omega\in {\mathcal D}(\tau^*_{-i/2})$.
\end{lemm}

\begin{proof}
We know that $(id*\mu)(W)$ belongs to ${\mathcal D}(S)$ and that
$S((id*\mu)(W)=(id*\mu)(W^*)$. So $(id*\mu)(W)$ belongs to $
{\mathcal D}(\tau_{-i/2})$ and
$\tau_{-i/2}((id*\mu)(W))=R((id*\mu)(W^*))$. By applying $\omega$ to
the previous equation, we easily get the result.
\end{proof}

Since $\Psi=\Phi\circ R$, there exists an anti-unitary ${\mathcal
J}$ from $H_{\Psi}$ onto $H_{\Phi}$ such that ${\mathcal
J}\Lambda_{\Psi}(x)=\Lambda_{\Phi}(R(x^*))$ for all $x\in {\mathcal
N}_{\Psi}\cap {\mathcal N}_{T_R}$.

\begin{prop}\label{antilineaire}
For all $\omega\in {\mathcal I}$ and $\mu\in {\mathcal
D}(\rho_{i/2})$, $\omega\mu$ belongs to ${\mathcal I}$ and we have:
$$\xi(\omega\mu)={\mathcal J}^*\hat{\pi}(\rho_{i/2}(\mu))^*{\mathcal
J}\xi(\omega)$$
\end{prop}

\begin{proof}
For all $n\in\mathbb{N}$, we put
$e_n=\frac{n}{\sqrt{\pi}}\int\!\exp(-n^2t^2)\delta^{it}dt$ so that
$e_n$is analytic with respect to $\sigma^{\Phi}$, ${\mathcal
N}_{\Phi}e_n\subset {\mathcal N}_{\Phi}$ and ${\mathcal
N}_{\Phi}\delta^{-\frac{1}{2}}e_n\subset {\mathcal N}_{\Psi}$. It is
sufficient to prove the proposition for all $\mu\in {\mathcal
D}(\tau_{-i/2}^*\delta_{i/2})$. Then, since $\delta$ is a
co-character, we can compute, for all $x\in {\mathcal N}_{\Phi}$:
$$
\begin{aligned}
\Lambda_{\Phi}((id\surl{\ _{\beta}
  \star_{\alpha}}_{\ \nu}\overline{\mu})\Gamma(xe_n))&=\Lambda_{\Psi}((id\surl{\
_{\beta} \star_{\alpha}}_{\
\nu}\overline{\mu})\Gamma(xe_n)\delta^{-\frac{1}{2}})\\
&=\Lambda_{\Psi}((id\surl{\ _{\beta} \star_{\alpha}}_{\
\nu}\overline{\mu})\Gamma(xe_n\delta^{-\frac{1}{2}})(1\surl{\
_{\beta}\otimes_{\alpha}}_{\ N}\delta^{-\frac{1}{2}}))
\end{aligned}$$ The computation goes on as follow:
$$
\begin{aligned}
&\ \quad\Lambda_{\Phi}((id\surl{\ _{\beta}
  \star_{\alpha}}_{\ \nu}\overline{\mu})\Gamma(xe_n))
=\Lambda_{\Psi}((id\surl{\ _{\beta} \star_{\alpha}}_{\
\nu}\overline{\delta_{-i/2}^*(\mu)})\Gamma(x\delta^{-\frac{1}{2}}e_n))\\
&={\mathcal J}^*\Lambda_{\Phi}(R((id\surl{\ _{\beta}
\star_{\alpha}}_{\
\nu}\delta_{-i/2}^*(\mu))\Gamma((x\delta^{-\frac{1}{2}}e_n)^*)))\\
&={\mathcal J}^*\Lambda_{\Phi}((\delta_{-i/2}^*(\mu)\circ R\surl{\
_{\beta} \star_{\alpha}}_{\
\nu}id)\Gamma(R(x\delta^{-\frac{1}{2}}e_n)^*))\\
&={\mathcal J}^*(\delta_{-i/2}^*(\mu)\circ
R*id)(W^*)\Lambda_{\Phi}(R(x\delta^{-\frac{1}{2}}e_n)^*))\\
&={\mathcal J}^*(\delta_{-i/2}^*(\mu)\circ
R*id)(W^*){\mathcal J}\Lambda_{\Psi}(x\delta^{-\frac{1}{2}}e_n))\\
&={\mathcal J}^*(\tau_{-i/2}^*\delta_{-i/2}^*(\mu)*id)(W){\mathcal
J}\Lambda_{\Phi}(xe_n))={\mathcal
J}^*(\rho_{i/2}(\mu)*id)(W){\mathcal J}\Lambda_{\Phi}(xe_n))
\end{aligned}$$
Now, we have:
$$
\begin{aligned}
(\omega\mu)((xe_n)^*)&=(\omega\surl{\ _{\beta} \star_{\alpha}}_{\
\nu}\mu)\Gamma((xe_n)^*)=\omega((id\surl{\ _{\beta}
\star_{\alpha}}_{\ \nu}\mu)\Gamma((xe_n)^*))\\
&=(\xi(\omega)|\Lambda_{\Phi}((id\surl{\ _{\beta} \star_{\alpha}}_{\
\nu}\overline{\mu})\Gamma(xe_n)))=(\xi(\omega)|{\mathcal
J}^*\hat{\pi}(\rho_{i/2}(\mu)){\mathcal
J}\Lambda_{\Phi}(xe_n))\\
&=({\mathcal J}^*\hat{\pi}(\rho_{i/2}(\mu))^*{\mathcal
J}\xi(\omega)|\Lambda_{\Phi}(xe_n))
\end{aligned}$$
Since $(xe_n)_{n\in\mathbb{N}}$ is converging to $x$ and
$(\Lambda_{\Phi}(xe_n))_{n\in\mathbb{N}}$ is converging to
$\Lambda_{\Phi}(x)$, we finally have:
$$(\omega\mu)(x^*)=({\mathcal J}^*\hat{\pi}(\rho_{i/2}(\mu))^*{\mathcal
J}\xi(\omega)|\Lambda_{\Phi}(x))$$ so that $\omega\mu\in {\mathcal
I}$ and $\xi(\omega\mu)={\mathcal
J}^*\hat{\pi}(\rho_{i/2}(\mu))^*{\mathcal J}\xi(\omega)$.
\end{proof}

\begin{coro}
There exists a unique closed densely defined operator
$\widehat{\Lambda}$ from ${\mathcal
D}(\hat{\Lambda})\subset\widehat{M}$ to $H_{\Phi}$ such that
$\widehat{\pi}(\mathcal{I})$ is a core for $\widehat{\Lambda}$ and
$\hat{\Lambda}(\hat{\pi}(\omega))=\xi(\omega)$ for all $\omega\in
{\mathcal I}$.
\end{coro}

\begin{proof}
Let $(\omega_n)_{n\in\mathbb{N}}$ be a sequence of ${\mathcal I}$
and let $w\in H_{\Phi}$ such that
$(\hat{\pi}(\omega_n))_{n\in\mathbb{N}}$ is converging to $0$ and
$(\xi_n)_{n\in\mathbb{N}}$ is converging to $w$. If $\mu$ belongs to
${\mathcal D}(\rho_{i/2})\cap {\mathcal I}$, then we have, by the
previous proposition, for all $n\in\mathbb{N}$:
$$\hat{\pi}(\omega_n)\xi(\mu)={\mathcal J}^*\hat{\pi}(\rho_{i/2}(\mu))^*{\mathcal
J}\xi(\omega_n)$$ Take the limit to get that $0={\mathcal
J}^*\hat{\pi}(\rho_{i/2}(\mu))^*{\mathcal J}w$. Since it is easy to
check that $\rho_{i/2}({\mathcal D}(\rho_{i/2})\cap {\mathcal I})$
is dense in ${\mathcal I}$ we get that $w=0$. So the formula of the
proposition defines a closable operator and its closure satisfy all
expected conditions.
\end{proof}

\begin{theo}\label{deb2}
There exists a unique normal semi-finite faithful weight
$\widehat{T_L}:\widehat{M}\rightarrow\alpha(N)$ such that the normal
semi-finite faithful weight
$\widehat{\Phi}=\nu\circ\alpha^{-1}\circ\widehat{T_L}$ admits
$(H,\iota,\widehat{\Lambda})$ as GNS construction. Moreover,
$\widehat{\sigma}$ is the modular group of $\widehat{\Phi}$, the
closure of $PJ_{\Phi}\delta J_{\Phi}$ ($P$ and $J_{\Phi}\delta
J_{\Phi}$ commute with each other) coincide with the modular
operator of $\widehat{\Phi}$ and
$\sigma_t^{\widehat{T_L}}(\hat{\beta}(n))=\hat{\beta}(\gamma_{-t}(n))$
for all $t\in\mathbb{R}$ and $n\in N$.
\end{theo}

\begin{proof}
Since $\widehat{\pi}$ is a multiplicative application and since
${\mathcal I}$ is a left ideal of $M_*^{\beta,\alpha}$, $xy$ belongs
to $\widehat{\pi}({\mathcal I})$ for all
$x\in\widehat{\pi}(M_*^{\beta,\alpha})$ and
$y\in\widehat{\pi}({\mathcal I})$ so, by definition, we have
$\widehat{\Lambda}(xy)=x\widehat{\Lambda}(y)$. Using the closeness
of $\widehat{\Lambda}$, we show that ${\mathcal
D}(\widehat{\Lambda})$ is a left ideal of $\widehat{M}$ and
$\widehat{\Lambda}(xy)=x\widehat{\Lambda}(y)$ for all
$x\in\widehat{M}$ and $y\in {\mathcal D}(\widehat{\Lambda})$.

By proposition \ref{groupe}, $\widehat{\sigma}_t(x)$ belongs to
${\mathcal D}(\widehat{\Lambda})$ for all
$x\in\widehat{\pi}({\mathcal I})$ and $t\in\mathbb{R}$ and
$\widehat{\Lambda}(\widehat{\sigma}_t(x))=P^{it}J_{\Phi}\delta^{it}J_{\Phi}\widehat{\Lambda}(x)$.
Using again the closeness of $\widehat{\Lambda}$, we get that
$\widehat{\sigma}_t(x)$ belongs to ${\mathcal D}(\widehat{\Lambda})$
for all $x\in {\mathcal D}(\widehat{\Lambda})$ and $t\in\mathbb{R}$
and we have:
$$\widehat{\Lambda}(\widehat{\sigma}_t(x))=P^{it}J_{\Phi}\delta^{it}J_{\Phi}\widehat{\Lambda}(x)$$

By proposition \ref{antilineaire}, for all $\omega\in {\mathcal
D}(\rho_{i/2})$ and $x\in\widehat{\pi}({\mathcal I})$,
$x\widehat{\pi}(\omega)$ belongs to ${\mathcal
D}(\widehat{\Lambda})$ and we have
$\widehat{\Lambda}(x\widehat{\pi}(\omega))={\mathcal
J}^*\widehat{\pi}(\rho_{i/2}(\omega))^*{\mathcal
J}\widehat{\Lambda}(x)={\mathcal
J}^*\widehat{\sigma}_{i/2}(\widehat{\pi}(\rho_{i/2})))^*{\mathcal
J}\widehat{\Lambda}(x)$. Since $\widehat{\pi}({\mathcal
D}(\rho_{i/2}))$ is dense in ${\mathcal D}(\widehat{\sigma}_{i/2})$
and $\widehat{\sigma}$-invariant, $\widehat{\pi}({\mathcal
D}(\rho_{i/2}))$ is a core for $\widehat{\sigma}$. The closeness of
$\widehat{\Lambda}$ allows us to conclude that $xy$ belongs to
${\mathcal D}(\widehat{\Lambda})$ for all $x\in {\mathcal
D}(\widehat{\Lambda})$ and $y\in {\mathcal
D}(\widehat{\sigma}_{i/2})$ and we have:
$$\widehat{\Lambda}(xy)={\mathcal J}^*\widehat{\sigma}_{i/2}(y)^*{\mathcal
J}\widehat{\Lambda}(x)$$

Therefore we know, by proposition 5.14 of \cite{Ku1}, that there
exists a normal semi-finite weight $\widehat{\Phi}$ on $\widehat{M}$
such that $(H,\iota,\widehat{\Lambda})$ is a GNS construction for
$\widehat{\Phi}$ and $\widehat{\sigma}$ is the modular group of
$\widehat{\Phi}$. Moreover, thanks to the previous equation, we
have:
$$\widehat{\Lambda}(xy)={\mathcal J}^*\widehat{\sigma}_{i/2}(y)^*{\mathcal
J}\widehat{\Lambda}(x)$$ for all $x\in {\mathcal
N}_{\widehat{\Phi}}$ and $y\in {\mathcal
D}(\widehat{\sigma}_{i/2})\cap {\mathcal N}_{\widehat{\Phi}}$. We
easily get faithfulness of $\widehat{\Phi}$ from this last relation.
We already know that $\alpha(N)\subseteq M\cap\widehat{M}$ and, by
proposition \ref{groupe} we have, for all $n\in N$:
$$\sigma_t^{\widehat{\Phi}}(\alpha(n))=\alpha(\sigma_t^{\nu}(n))
=\sigma_t^{\nu\circ\alpha^{-1}}(\alpha(n))$$ By Haagerup's existence
theorem, we get the normal semi-finite faithful weight
$\widehat{\Phi}$. Finally, we check the last property. For all $n\in
N$ and $t\in\mathbb{R}$, we have:
$$\begin{aligned}
\sigma_t^{\widehat{\Phi}}(\hat{\beta}(n))&=
P^{it}J_{\Phi}\delta^{it}\alpha(n^*)\delta^{-it}J_{\Phi}P^{-it}=
P^{it}J_{\Phi}\alpha(\gamma_{-t}\sigma^{\nu}_{-t}(n^*))J_{\Phi}P^{-it}\\
&=P^{it}\hat{\beta}(\gamma_{-t}\sigma^{\nu}_{-t}(n))P^{-it}=\hat{\beta}(\gamma_{-t}(n))
\end{aligned}$$ because $\gamma$ and $\sigma^{\nu}$ commute with each other.
\end{proof}

\begin{lemm}
For all $x\in {\mathcal N}_{\widehat{T_L}}\cap {\mathcal
N}_{\hat{\Phi}}$, $\widehat{\Lambda}(x)$ belongs to
$D(H_{\beta},\nu^0)$ and we have
$R^{\beta,\nu^o}(\hat{\Lambda}(x))=\Lambda_{\widehat{T_L}}(x)$.
\end{lemm}

\begin{proof}
By definition $J_{\widehat{\Phi}}$ and $\mathcal{J}$ implement the
same operator on $\alpha(N)\subset M\cap\widehat{M}$ so that
$J_{\widehat{\Phi}}\alpha(n^*)J_{\widehat{\Phi}}=\beta(n)$ for all
$n\in N$. Then the lemma is a consequence of proposition \ref{prem}.
\end{proof}

\begin{lemm}\label{dernier}
For all $\xi\in D((H_{\Phi})_{\hat{\beta}},\nu^o)$, all $\eta\in
D((H_{\Phi})_{\hat{\beta}},\nu^o)\cap D(_{\alpha}H_{\Phi},\nu^o)$
and all $x\in {\mathcal N}_{\hat{\Phi}}\cap {\mathcal
N}_{\widehat{T_L}}$, $(\omega_{\eta,\xi}\surl{\
_{\hat{\beta}}\otimes_{\alpha}}_{\ \nu}id)(\hat{\Gamma}(x))$ belongs
to ${\mathcal N}_{\hat{\Phi}}\cap {\mathcal N}_{\widehat{T_L}}$ and
we have:
$$\hat{\Lambda}((\omega_{\eta,\xi}\surl{\ _{\hat{\beta}}\otimes_{\alpha}}_{\
\nu}id)(\hat{\Gamma}(x)))=(id*\omega_{\eta,\xi})(W)\hat{\Lambda}(x)$$
\end{lemm}

\begin{proof}
Thanks to the pentagonal relation, we can compute for all
$\omega\in\mathcal{I}$:
$$
\begin{aligned}
&\ \quad (\omega_{\eta,\xi}\surl{\
_{\hat{\beta}}\otimes_{\alpha}}_{\
\nu}id)(\hat{\Gamma}(\hat{\pi}(\omega))) =(\omega_{\eta,\xi}\surl{\
_{\hat{\beta}}\otimes_{\alpha}}_{\
\nu}id)(\hat{\Gamma}((\omega*id)(W)))\\
&=(\omega_{\eta,\xi}\surl{\ _{\hat{\beta}}\otimes_{\alpha}}_{\
\nu}id)(\sigma_{\nu^o}W((\omega*id)(W)\surl{\
_{\beta}\otimes_{\alpha}}_{\ N}1)W^*\sigma_{\nu})\\
&=(\omega*id*\omega_{\eta,\xi})((1\surl{\
_{\alpha}\otimes_{\hat{\beta}}}_{\ N^o}W)(W\surl{\
_{\beta}\otimes_{\alpha}}_{\ N}1)(1\surl{\
_{\beta}\otimes_{\alpha}}_{\ N}W^*))\\
&=(\omega*\omega_{\eta,\xi}*id)((W\surl{\
_{\alpha}\otimes_{\hat{\beta}}}_{\ N^o}1)(\sigma_{\nu^o}\surl{\
_{\alpha}\otimes_{\hat{\beta}}}_{\ N^o}1)(1\surl{\
_{\alpha}\otimes_{\hat{\beta}}}_{\ N^o}W)\sigma_{2\nu}(1\surl{\
_{\beta}\otimes_{\alpha}}_{\ N}\sigma_{\nu^o}))\\
&=\hat{\pi}((id*\omega_{\eta,\xi})(W)\omega)
\end{aligned}$$
Then, by definition, $(\omega_{\eta,\xi}\surl{\
_{\hat{\beta}}\otimes_{\alpha}}_{\
\nu}id)(\hat{\Gamma}(\hat{\pi}(\omega)))$ belongs to ${\mathcal
N}_{\hat{\Phi}}\cap {\mathcal N}_{\widehat{T_L}}$ for all $\omega\in
{\mathcal I}$ and we have:
$$\hat{\Lambda}((\omega_{\eta,\xi}\surl{\
_{\hat{\beta}}\otimes_{\alpha}}_{\
\nu}id)(\hat{\Gamma}(\hat{\pi}(\omega))))
=(id*\omega_{\eta,\xi})(W)\hat{\Lambda}(\hat{\pi}(\omega))$$
Closeness of $\widehat{\Lambda}$ finishes the proof.
\end{proof}

\begin{prop}
The operator-valued weight $\widehat{T_L}$ is left invariant.
\end{prop}

\begin{proof}
Let $(\xi_i)_{i\in I}$ be a $(N^o,\nu^o)$-basis of
$(H_{\Phi})_{\hat{\beta}}$. For all $x\in {\mathcal
N}_{\widehat{\Phi}}\cap {\mathcal N}_{\widehat{T_L}}$ and $\eta\in
D((H_{\Phi})_{\hat{\beta}},\nu^o)\cap D(_{\alpha}H_{\Phi},\nu)$, we
have:
$$\begin{aligned}
&\ \quad\widehat{\Phi}((\omega_{\eta}\surl{\
_{\hat{\beta}}\star_{\alpha}}_{\
\nu}id)(\widehat{\Gamma}(x^*x)))=\sum_{i\in
I}\widehat{\Phi}((\omega_{\eta,\xi_i}\surl{\
_{\hat{\beta}}\star_{\alpha}}_{\
\nu}id)(\widehat{\Gamma}(x))^*(\omega_{\eta,\xi_i}\surl{\
_{\hat{\beta}}\star_{\alpha}}_{\ \nu}id)(\widehat{\Gamma}(x)))\\
&=\sum_{i\in I}||\widehat{\Lambda}((\omega_{\eta,\xi_i}\surl{\
_{\hat{\beta}}\star_{\alpha}}_{\ \nu}id)(\widehat{\Gamma}(x))||^2=\sum_{i\in I}||(id*\omega_{\eta,\xi_i})(W)\widehat{\Lambda}(x)||^2\\
&=((\rho^{\beta,\alpha}_{\eta})^*\rho^{\beta,\alpha}_{\eta}\widehat{\Lambda}(x)|\widehat{\Lambda}(x))=||\widehat{\Lambda}(x)\surl{\
_{\beta}\otimes_{\alpha}}_{\
\nu}\eta||^2\\
&=(\alpha(<\widehat{\Lambda}(x),\widehat{\Lambda}(x)>_{\beta,\nu^o}\eta|\eta)=(\widehat{T_L}(x^*x)\eta|\eta)
\end{aligned}$$

\end{proof}

To have a measured quantum groupoid, we need to check a relation
between the co-involution $\widehat{R}$ and $\Gamma$. By the way, it
will give a link between the two natural GNS constructions of
$\Phi_{\delta}=\Psi=\Phi\circ R$. We put
$S_{\widehat{\Phi}},J_{\widehat{\Phi}}$ and $
\Delta_{\widehat{\Phi}}$ to be the fundamental objects associated to
$\widehat{\Phi}$ by the Tomita's theory in the GNS construction
$(H,\iota,\widehat{\Lambda})$.

\begin{defi}
We put $\mathcal{I}^{\sharp}$ the subset of $\mathcal{I}$ consisting
of elements of the form
$\omega_{\Lambda_{\Phi}(a),\Lambda_{\Phi}(b)}$ where $a,b$ belong to
$E$.
\end{defi}

\begin{lemm}
We have that $\widehat{\pi}(\mathcal{I}^{\sharp})$ is a core for
$\widehat{\Lambda}$ and
$\widehat{\Lambda}(\widehat{\pi}(\mathcal{I}^{\sharp}))$ is a core
for $S_{\widehat{\Phi}}$ and for $\Delta_{\widehat{\Phi}}^z$ for all
$z\in\mathbb{C}$.
\end{lemm}

\begin{proof}
This lemma comes from standard arguments and by definition of
$\widehat{\Lambda}$.
\end{proof}

\begin{prop}
For all $x\in E$, we have $\Lambda_{\Phi}(x)$ belongs to
$\mathcal{D}(S_{\widehat{\Phi}}^*)$ and we have:
$$S_{\widehat{\Phi}}^*\Lambda_{\Phi}(x)=\Lambda_{\Phi}(S^{-1}(x)^*)$$
Moreover $\Lambda_{\Phi}(E)$ is a core for $S_{\widehat{\Phi}}^*$.
\end{prop}

\begin{proof}
Let $\omega\in\mathcal{I}^{\sharp}$. For all $\mu\in
M_*^{\alpha,\hat{\beta}}$, we have:
$$
\begin{aligned}
\mu(\widehat{\pi}(\omega)^*)&=\mu((\omega*id)(W)^*)=\overline{\omega}((id*\mu)(W^*))
=\overline{\omega}\circ
S((id*\mu)(W))\\&=\omega^*\circ\tau_{-\frac{i}{2}}((id*\mu)(W))
=\mu((\omega^*\circ\tau_{-\frac{i}{2}}*id)(W))=\mu(\widehat{\pi}(\omega^*\circ\tau_{-\frac{i}{2}}))
\end{aligned}$$
Then, we have:
$$
\begin{aligned}
&\
(S_{\widehat{\Phi}}\widehat{\Lambda}(\widehat{\pi}(\omega))|\Lambda_{\Phi})
=(\widehat{\Lambda}(\widehat{\pi}(\omega)^*)|\Lambda_{\Phi})
=(\widehat{\Lambda}(\widehat{\pi}(\omega^*\circ\tau_{-\frac{i}{2}})|\Lambda_{\Phi}))
=(\xi(\omega^*\circ\tau_{-\frac{i}{2}})|\Lambda_{\Phi})\\
&=\omega^*\circ\tau_{-\frac{i}{2}}(x^*)=\overline{\omega(S^{-1}(x))}
=(\Lambda_{\Phi}(S^{-1}(x)^*)|\xi(\omega))=(\Lambda_{\Phi}(S^{-1}(x)^*)|\widehat{\Lambda}(\widehat{\pi}(\omega)))
\end{aligned}$$ Thus the previous lemma and the fact that
$\xi(\omega_{\Lambda_{\Phi}(a),\Lambda_{\Phi}(b)})=\Lambda_{\Phi}(a\sigma_{-i}^{\Phi}(b^*))$
implies the proposition.
\end{proof}

\begin{prop}\label{hatj}
For all $x\in\mathcal{N}_{T_R}\cap\mathcal{N}_{\Psi}$, we have:
$$J_{\widehat{\Phi}}\Lambda_{\Phi_{\delta}}(x)=\Lambda_{\Phi}(R(x^*))$$
\end{prop}

\begin{proof}
Define the anti-unitary $\mathcal{J}$ of $H$ such that
$\mathcal{J}\Lambda_{\Phi_{\delta}}(x)=\Lambda_{\Phi}(R(x^*))$ for
all $x\in\mathcal{N}_{T_R}\cap\mathcal{N}_{\Psi}$. Let $a$ belongs
to $E$. For all $n\in\mathbb{N}$, we put
$e_n=\frac{n}{\sqrt{\pi}}\int\!\exp(-n^2t^2)\delta^{it}dt$ so that
$e_n$ is analytic with respect to $\sigma^{\Phi}$, ${\mathcal
N}_{\Phi}e_n\subset {\mathcal N}_{\Phi}$ and ${\mathcal
N}_{\Phi}\delta^{-\frac{1}{2}}e_n\subset {\mathcal N}_{\Psi}$. Since
$\tau_s(\delta)=\delta$, we see that $\tau_s(e_n)=e_n$ for all
$s\in\mathbb{R}$, hence $e_n\in\mathcal{D}(\tau_{\frac{i}{2}})$ and
$\tau_{\frac{i}{2}}(e_n)=e_n$. By assumption $a$ belongs to
$\mathcal{D}(\tau_{\frac{i}{2}})$ so that $ae_n$ belongs to
$\mathcal{D}(\tau_{\frac{i}{2}})$ and
$\tau_{\frac{i}{2}}(ae_n)=\tau_{\frac{i}{2}}(a)e_n$. Hence
$\tau_{\frac{i}{2}}(ae_n)\delta^{\frac{1}{2}}$ is a bounded operator
and its closure is equal to
$\tau_{\frac{i}{2}}(a)(\delta^{\frac{1}{2}}e_n)$. We recall that
$\kappa_t(x)$ is equal, by definition, to $\tau_t(m)\delta^{it}$ for
all $t\in\mathbb{R}$ and $m\in M$. Then $ae_n$ belongs to
$\mathcal{D}(\kappa_{\frac{i}{2}})$ and
$\kappa_{\frac{i}{2}}(ae_n)=\tau_{\frac{i}{2}}(a)(\delta^{\frac{1}{2}}e_n)$.
By assumption, $\tau_{\frac{i}{2}}$ belongs to
$\mathcal{N}_{\Psi}\cap\mathcal{N}_{T_R}$. So we see that
$\kappa_{\frac{i}{2}}(ae_n)\delta^{-\frac{1}{2}}$ is bounded and its
closure equals
$\tau_{\frac{i}{2}}(a)e_n\in\mathcal{N}_{\Psi}\cap\mathcal{N}_{T_R}$
implying that
$\kappa_{\frac{i}{2}}(ae_n)\in\mathcal{N}_{\Phi}\cap\mathcal{N}_{T_L}$
and:
$$\Lambda_{\Phi}(\kappa_{\frac{i}{2}}(ae_n))=
\Lambda_{\Phi_{\delta}}(\kappa_{\frac{i}{2}}(ae_n)\delta^{-\frac{1}{2}})
=\Lambda_{\Phi_{\delta}}(\tau_{\frac{i}{2}}(a)e_n)$$ By definition,
we have
$\Delta_{\widehat{\Phi}}^{it}\Lambda_{\Phi}(x)=\Lambda_{\Phi}(\kappa_t(x))$,
we easily get that $\Lambda_{\Phi}(ae_n)$ belongs to
$\mathcal{D}(\Delta_{\widehat{\Phi}}^{-\frac{1}{2}})$ and:
$$\Delta_{\widehat{\Phi}}^{-\frac{1}{2}}\Lambda_{\Phi}(ae_n)
=\Lambda_{\Phi}(\kappa_{\frac{i}{2}}(ae_n))
=\Lambda_{\Phi_{\delta}}(\tau_{\frac{i}{2}}(a)e_n)$$ By closedness
of $\Delta_{\widehat{\Phi}}^{-\frac{1}{2}}$, this implies that
$\Lambda_{\Phi}(a)$ belongs to
$\mathcal{D}(\Delta_{\widehat{\Phi}}^{-\frac{1}{2}})$ and:
$$\Delta_{\widehat{\Phi}}^{-\frac{1}{2}}\Lambda_{\Phi}(a)
=\Lambda_{\Phi_{\delta}}(\tau_{\frac{i}{2}}(a))$$ Consequently, we
have:
$$\mathcal{J}\Delta_{\widehat{\Phi}}^{-\frac{1}{2}}\Lambda_{\Phi}(a)
=\mathcal{J}\Lambda_{\Phi_{\delta}}(\tau_{\frac{i}{2}}(a))=\Lambda_{\Phi}(S^{-1}(a)^*)
=S_{\widehat{\Phi}}^*\Lambda_{\Phi}(a)
=J_{\widehat{\Phi}}\Delta_{\widehat{\Phi}}^{-\frac{1}{2}}\Lambda_{\Phi}(a)$$
Since $\Lambda_{\Phi}(E)$ is a core for
$\Delta_{\widehat{\Phi}}^{-\frac{1}{2}}=J_{\widehat{\Phi}}S^*_{\widehat{\Phi}}$,
we have done.
\end{proof}

Finally, we have to recognize what is $\widehat{W}$.

\begin{prop}\label{transfo}
The unitary $\sigma_{\nu}W^*\sigma_{\nu}$ is the fundamental unitary
associated with the dual Hopf-bimodule structure.
\end{prop}

\begin{proof}
The fundamental unitary associated with the dual quantum groupoid is
denoted by $\widehat{W}$. By definition of $\widehat{W}$ and lemma
\ref{dernier}, we have for all $\xi\in D(_{\alpha}H_{\Phi},\nu)\cap
D((H_{\Phi})_{\hat{\beta}},\nu^o)$, $\eta\in
D((H_{\Phi})_{\hat{\beta}},\nu^o)$ and $x\in {\mathcal
N}_{\hat{\Phi}}\cap {\mathcal N}_{\widehat{T_L}}$:
$$\begin{aligned}
(\omega_{\xi,\eta}*id)(\widehat{W}^*)\widehat{\Lambda}(x)&=\widehat{\Lambda}((\omega_{\xi,\eta}\surl{\
_{\hat{\beta}}\star_{\alpha}}_{\ \nu}id)(\widehat{\Gamma}(x)))\\
&=(id*\omega_{\xi,\eta})(W)\widehat{\Lambda}(x)
=(\omega_{\xi,\eta}*id)(\sigma_{\nu^o}W\sigma_{\nu^o})\widehat{\Lambda}(x)
\end{aligned}$$ from which we easily deduce that
$\widehat{W}=\sigma_{\nu}W^*\sigma_{\nu}$.
\end{proof}

\begin{theo}\label{caradual}
$(N,\widehat{M},\alpha,\hat{\beta},\widehat{\Gamma},\widehat{T_L},\widehat{R},\widehat{\tau},\nu)$
is a measured quantum groupoid called \textbf{dual quantum groupoid}
of $(N,M,\alpha,\beta,\Gamma,T_L,R,\tau,\nu)$. Fundamental objects
of the dual quantum groupoid
$(N,\widehat{M},\alpha,\hat{\beta},\widehat{\Gamma},\widehat{R},\widehat{T_L},\widehat{\tau},\nu)$
are given, for all $x\in\widehat{M}$ and $t\in\mathbb{R}$, by:
\begin{enumerate}[i)]
\item $\widehat{W}=\sigma_{\nu}W^*\sigma_{\nu}$ is the fundamental
unitary,
\item $\widehat{R}(x)=J_{\Phi}x^*J_{\Phi}$ is the unitary antipode and $\widehat{\tau}_t(x)=P^{it}xP^{-it}$ is the scaling group,
\item $\widehat{\lambda}=\lambda^{-1}$ is the scaling operator and
the closure of $P^{-1}J_{\Phi}\delta
J_{\Phi}\delta^{-1}\Delta_{\Phi}^{-1}$ is the modulus
$\widehat{\delta}$,
\item $\widehat{P}=P$ is the manipulation operator,
\item in the GNS construction
$(H,\iota,\widehat{\Lambda})$, the modular operator
$\Delta_{\widehat{\Phi}}$ is the closure of
$PJ_{\Phi}\delta^{-1}J_{\Phi}$ and the modular conjugation satisfies
$J_{\widehat{\Phi}}\Lambda_{\Phi_{\delta}}(x)=\Lambda_{\Phi}(R(x^*))$
for all $x\in\mathcal{N}_{T_R}\cap\mathcal{N}_{\Psi}$.
\end{enumerate}
\end{theo}

\begin{proof}
By proposition \ref{deb1},
$(N,\widehat{M},\alpha,\hat{\beta},\widehat{\Gamma})$ is a
Hopf-bimodule. By theorem \ref{deb2}, it admits a normal semi-finite
faithful left-invariant operator-valued weight $\widehat{T_L}$. By
proposition \ref{deb3}, $\widehat{R}$ is a co-involution for this
structure and, by definition, we have
$\widehat{R}(\widehat{\pi}(\omega))=\widehat{\pi}(\omega\circ R)$
for all $\omega\in M_*^{\alpha,\beta}$. Since
$J_{\widehat{\Phi}}=\mathcal{J}$ implement $R$ on $M$ and since
$=\widehat{W}=\sigma_{\nu}W^*\sigma_{\nu}$, we get
$R((id*\omega_{J_{\Phi}v,w})(\widehat{W}))
=(id*\omega_{J_{\Phi}w,v})(\widehat{W})$. By proposition
\ref{groupe}, $\widehat{\tau}$ is a scaling group. We just have to
check that the one-parameter group of automorphisms
$\widehat{\gamma}$ of $N$ leaves $\nu$ invariant. However, we have
already noticed, in theorem \ref{deb2}, that we have
$\widehat{\gamma}_t=\gamma_{-t}$ for all $t\in\mathbb{R}$. By
hypothesis over $\gamma$, we have done.

By proposition \ref{forP} and by definition of $\widehat{\tau}$,
$\widehat{\pi}(\mathcal{I})$ is stable under $\widehat{\tau}_t$
$t\in\mathbb{R}$ and we have, for all $\omega\in\mathcal{I}$:
$$\widehat{\Lambda}(\widehat{\tau}_t(\widehat{\pi}(\omega)))
=\widehat{\Lambda}(\widehat{\pi}(\omega\circ\tau_{-t}))
=\xi(\omega\circ\tau_{-t})=\lambda^{\frac{t}{2}}P^{it}\widehat{\Lambda}(\widehat{\pi}(\omega))$$
Now, by closeness of $\widehat{\Lambda}$, we get that
$P^{it}\widehat{\Lambda}(x)=\lambda^{-\frac{t}{2}}\widehat{\Lambda}(\widehat{\tau}_t(x))$
for all
$x\in\mathcal{N}_{\widehat{T_L}}\cap\mathcal{N}_{\widehat{\Phi}}$
and $t\in\mathbb{R}$. From this and from lemma \ref{prep1}, we get
that:
$$\lambda^{-ist}=[D\widehat{\Phi}\circ\widehat{\tau}_{-s}:D\widehat{\Phi}]_t
=[D\widehat{\Phi}\circ\widehat{\sigma}_s^{\widehat{\Phi}\circ\widehat{R}}:D\widehat{\Phi}]_t
=\widehat{\lambda}^{ist}$$ $$\text{ and }\ \quad
P^{it}\widehat{\Lambda}(x)=\widehat{\lambda}^{\frac{t}{2}}\widehat{\Lambda}(\widehat{\tau}_t(x))
=\widehat{P}^{it}\widehat{\Lambda}(x)$$
\end{proof}

The whole picture is not completely drawn yet because the value of
$\widehat{\delta}$ is missing. For this, we need the bi-duality
theorem. The expression will finally be given in \ref{hatdelta}.

\subsection{Bi-duality theorem}
In this section, we compute fundamentals objects of the dual
structure. Also, we can construct the bi-dual quantum groupoid that
is the dual quantum groupoid of the dual quantum groupoid and we
establish a bi-duality theorem.

\begin{theo}
The measured quantum groupoid
$(N,M,\alpha,\beta,\Gamma,T_L,R,\tau,\nu)$ and its bi-dual
$(N,\widehat{\widehat{M}},\alpha,\hat{\hat{\beta}},
\widehat{\widehat{\Gamma}},\widehat{\widehat{T_L}},\widehat{\widehat{R}},
\widehat{\widehat{\tau}},\nu)$ coincide. Moreover, we have
$\widehat{\widehat{\Lambda}}=\Lambda_{\Phi}$.
\end{theo}

\begin{proof}
We know that $J_{\widehat{\Phi}}=\mathcal{J}$. Then, on
$\alpha(N)\subset M\cap\widehat{M}$, we have:
$$\hat{\hat{\beta}}(n)=J_{\widehat{\Phi}}\alpha(n)^*J_{\widehat{\Phi}}
=\mathcal{J}\alpha(n)^*\mathcal{J}=R(\alpha(n))=\beta(n)$$ By
proposition \ref{transfo}, we have:
$$\widehat{\widehat{W}}=\sigma_{\nu}\widehat{W}^*\sigma_{\nu}=W$$
so that we deduce that the Hopf-bimodule and its bi-dual coincide.
We denote by $\widehat{\widehat{\pi}}(\omega)=(\omega
*id)(\widehat{W})=(id*\omega)(W^*)$ for all $\omega\in
M_*^{\alpha,\hat{\beta}}$. By definition of $\widehat{R}$ and
$\widehat{\widehat{R}}$, we have for all $\xi,\eta\in
D(_{\alpha}H,\nu)$:
$$\begin{aligned}
\widehat{\widehat{R}}((id*\omega_{J_{\Phi}\xi,\eta})(W^*))&=
\widehat{\widehat{R}}(\widehat{\widehat{\pi}}(\omega_{J_{\Phi}\xi,\eta}))=
\widehat{\widehat{\pi}}(\omega_{J_{\Phi}\xi,\eta}\circ\widehat{R})\\
&=
\widehat{\widehat{\pi}}(\omega_{J_{\Phi}\eta,\xi})=(id*\omega_{J_{\Phi}\eta,\xi})(W^*)
\end{aligned}$$ so that $\widehat{\widehat{R}}=R$. Let
$\omega\in\widehat{{\mathcal I}}$. On note
$a=\widehat{\widehat{\pi}}(\omega)$. Then, for all $\Theta\in
{\mathcal I}$, we have:
$$
\begin{aligned}
\omega(\hat{\pi}(\Theta)^*)&=\omega((\Theta*id)(W)^*)=\omega((\overline{\Theta}*id)(W^*))=
\overline{\Theta}((id*\omega)(W^*))\\
&=\overline{\Theta(a^*)}=\overline{(\xi(\Theta)|\Lambda_{\Phi}(a))}
=(\Lambda_{\Phi}(a)|\hat{\Lambda}(\hat{\pi}(\Theta)))
\end{aligned}$$
Since $\widehat{\pi}({\mathcal I})$ is a core for
$\widehat{\Lambda}$, this implies
$\omega(x^*)=(\Lambda_{\Phi}(a)|\widehat{\Lambda}(x))$ for all $x\in
{\mathcal N}_{\widehat{\Phi}}$. By definition of
$\widehat{\widehat{\Lambda}}$, we get
$\widehat{\widehat{\Lambda}}(\widehat{\widehat{\pi}}(\omega))=\Lambda_{\Phi}(a)
=\Lambda_{\Phi}(\widehat{\widehat{\pi}}(\omega))$. Since
$\widehat{\widehat{\pi}}(\widehat{{\mathcal I}})$ is a core for
$\widehat{\widehat{\Lambda}}$ and by closeness of $\Lambda_{\Phi}$
we have $\widehat{\widehat{\Lambda}}(y)=\Lambda_{\Phi}(y)$ for all
$y\in {\mathcal N}_{\widehat{\Phi}}$. In particular
$\widehat{\widehat{T_L}}=T_L$. Finally, we have to compute
$\widehat{\widehat{\tau}}$. For example, we can use proposition
\ref{fortau}, to get for all $t\in\mathbb{R}$:
$$\Gamma\circ\widehat{\widehat{\tau}}_t=\widehat{\widehat{\Gamma}}\circ\widehat{\widehat{\tau}}_t=
(\sigma_t^{\widehat{\widehat{\Phi}}}\surl{\ _{\beta}
\star_{\alpha}}_{\ N}\sigma_{-t}^{\widehat{\widehat{\Phi}}\circ
\widehat{\widehat{R}}})\circ\widehat{\widehat{\Gamma}}=
(\sigma_t^{\Phi}\surl{\ _{\beta}\star_{\alpha}}_{\
N}\sigma_{-t}^{\Phi\circ R})\circ\Gamma=\Gamma\circ\tau$$ and we can
conclude by injectivity of $\Gamma$.
\end{proof}

\begin{prop}\label{hatdelta}
For all $t\in\mathbb{R}$, we have:
$$\widehat{\delta}^{it}=P^{-it}J_{\Phi}\delta^{-it}J_{\Phi}\delta^{-it}\Delta_{\Phi}^{-it}$$
\end{prop}

\begin{proof}
By theorem \ref{caradual}, we know that
$\Delta_{\widehat{\Phi}}^{it}=P^{it}J_{\Phi}\delta^{it}J_{\Phi}$ so
that we get, thanks to the bi-duality theorem that:
$$\widehat{\delta}^{it}=\widehat{P}^{-it}J_{\widehat{\Phi}}\Delta^{it}J_{\widehat{\Phi}}
=P^{-it}J_{\widehat{\Phi}}\Delta^{it}J_{\widehat{\Phi}}$$ From the
previous proposition, it is easy to check on
$\Lambda_{\Phi_{\delta}}(x)$ that $J_{\widehat{\Phi}}\Delta
J_{\widehat{\Phi}}$ coincide with the modular operator of $\Psi$ in
the GNS construction $(H,\iota,\Lambda_{\Phi_{\delta}})$. Now, by
proposition 2.5 of \cite{Vae}, this last modular operator is equal
to the closure of $J_{\Phi}\delta^{-1}J_{\Phi}\delta\Delta_{\Phi}$
so that we get the result.
\end{proof}

\begin{rema}
From this last expression of $\widehat{\delta}$, we can directly
verify the following properties which should be satisfied by
duality, for all $x\in\widehat{M}$ and $s,t\in\mathbb{R}$:
$$\sigma^{\widehat{\Phi}}_s(\widehat{\delta}^{it})=\widehat{\lambda}^{ist}\widehat{\delta}^{it},\
\quad\sigma^{\widehat{\Phi}\circ\widehat{R}}_t(x)=\widehat{\delta}^{it}\sigma^{\widehat{\Phi}}_t(x)\widehat{\delta}^{-it}\
\quad\text{ and }\ \quad
\widehat{\Gamma}(\delta^{it})=\delta^{it}\surl{\
_{\hat{\beta}}\otimes_{\alpha}}_{\ \nu}\delta^{it}$$
\end{rema}

\begin{theo}
The following properties and their dual hold:
\begin{itemize}
\item
$\tau_t(m)=\Delta_{\widehat{\Phi}}^{it}m\Delta_{\widehat{\Phi}}^{-it}$\quad
and\quad $R(m)=J_{\widehat{\Phi}}m^*J_{\widehat{\Phi}}$\quad for all
$t\in\mathbb{R}$ and $m\in M$
\item $W(\Delta_{\widehat{\Phi}}\surl{\
_{\beta}\otimes_{\alpha}}_{\
\nu}\Delta_{\Phi})=(\Delta_{\widehat{\Phi}}\surl{\
_{\alpha}\otimes_{\hat{\beta}}}_{\ \nu^0}\Delta_{\Phi})W$\\ and\quad
$W(J_{\widehat{\Phi}}\surl{\ _{\alpha}\otimes_{\hat{\beta}}}_{\
\nu^0}J_{\Phi})=(J_{\widehat{\Phi}}\surl{\
_{\beta}\otimes_{\alpha}}_{\ \nu}J_{\Phi})W^*$
\item
$\Delta_{\widehat{\Phi}}^{it}\Delta_{\Phi}^{is}
=\lambda^{ist}\Delta_{\Phi}^{is}\Delta_{\widehat{\Phi}}^{it}$,\quad
$\Delta_{\Phi}^{it}\delta^{is}=\lambda^{ist}\delta^{is}\Delta_{\Phi}^{it}$\quad
and\quad
$\Delta_{\widehat{\Phi}}^{it}\delta^{is}=\delta^{is}\Delta_{\widehat{\Phi}}^{it}$
\item
$J_{\widehat{\Phi}}J_{\Phi}=\lambda^{\frac{i}{4}}J_{\Phi}J_{\widehat{\Phi}}$,\quad
$J_{\Phi}PJ_{\Phi}=P^{-1}$\quad and\quad $J_{\widehat{\Phi}}\delta
J_{\widehat{\Phi}}=\delta^{-1}$
\item $P^{is}\Delta_{\Phi}^{it}=\Delta_{\Phi}^{it}P^{is}$\quad
and\quad $P^{is}\delta_{\Phi}^{it}=\delta_{\Phi}^{it}P^{is}$
\end{itemize}
\end{theo}

\begin{proof}
Since $\delta$ is affiliated to $M$, $J_{\Phi}\delta J_{\Phi}$ is
affiliated to $M'$ so that, for all $t\in\mathbb{R}$ and $m\in M$,
we have:
$$\Delta_{\widehat{\Phi}}^{it}m\Delta_{\widehat{\Phi}}^{-it}=
P^{it}J_{\Phi}\delta^{it}J_{\Phi}mJ_{\Phi}\delta^{-it}J_{\Phi}P^{-it}
=P^{it}mP^{-it}=\tau_t(m)$$ We have already noticed that $R$ is
implemented by $J_{\widehat{\Phi}}$ by definition of
$\widehat{\Phi}$ but we can recover this point thanks to the
bi-duality theorem and the fact that, by definition, $\widehat{R}$
is implemented by $J_{\Phi}$. Now, since we have
$R((id*\omega_{\xi,J_{\Phi}\eta})(W))=(id*\omega_{\eta,J_{\Phi}\xi})(W)$
for all $\xi,\eta\in D(_{\alpha}H,\nu)$, we easily get the second
equality of the second point from the first point. Also, we know
that $\tau_t((id*\omega_{\xi,J_{\Phi}\eta})(W))
=(id*\omega_{\Delta_{\Phi}^{it}\xi,\Delta_{\Phi}^{it}J_{\Phi}\eta})(W)$
for all $t\in\mathbb{R}$ from which and from the first point we get
the first equality of the second point. Since $\tau$ and $\sigma$
commute each other, it is easy to check on $\Lambda_{\Phi}(x)$ the
first equality of the last point. Since $\tau(\delta)=\delta$, we
get the last equality of the last point. The last equality of the
third point comes from the fact that $\tau$ is implemented by
$\Delta_{\widehat{\Phi}}$ and that $\tau(\delta)=\delta$. By
proposition 5.2 of \cite{Vae}, we have
$\sigma_t^{\Phi}(\delta^{is})=\lambda^{ist}\delta^{it}$ so that we
get the second equality of the third point. Then, for all
$s,t\in\mathbb{R}$, we have:
$$\begin{aligned}
\Delta_{\widehat{\Phi}}^{it}\Delta_{\Phi}^{is}&=
P^{it}J_{\Phi}\delta^{it}J_{\Phi}\Delta_{\Phi}^{is}
=P^{it}J_{\Phi}\delta^{it}\Delta_{\Phi}^{is}J_{\Phi}\\
&=P^{it}J_{\Phi}\lambda^{-ist}\Delta_{\Phi}^{is}\delta^{it}J_{\Phi}
=\lambda^{ist}\Delta_{\Phi}^{is}P^{it}J_{\Phi}\delta^{it}J_{\Phi}
=\lambda^{ist}\Delta_{\Phi}^{is}\Delta_{\widehat{\Phi}}^{it}
\end{aligned}$$ As far as the fourth point is concerned, the last equality
comes from the fact that $R$ is implemented by $J_{\widehat{\Phi}}$
and $R(\delta)=\delta^{-1}$. The second one can be directly checked
on $\Lambda_{\Phi}(x)$. Let us prove the first equality. Let $x$
belongs to
$\mathcal{N}_{\Psi}\cap\mathcal{D}(\sigma^{\Psi}_{\frac{i}{2}})$.
Then, it is easy to see that $R(x^*)$ belongs to
$\mathcal{N}_{\Phi}\cap\mathcal{D}(\sigma^{\Phi}_{\frac{i}{2}})$.
Remembering that the modular conjugation of $\Psi=\Phi_{\delta}$
associated with the GNS construction
$(H,\iota,\Lambda_{\Phi_{\delta}})$ is equal to
$\lambda^{\frac{i}{4}}$ by proposition 2.5 of \cite{Vae}, we get:
$$
\begin{aligned}
J_{\widehat{\Phi}}J_{\Phi}\Lambda_{\Phi_{\delta}}(x)
&=\lambda^{\frac{i}{4}}J_{\widehat{\Phi}}\lambda^{\frac{i}{4}}J_{\Phi}\Lambda_{\Phi_{\delta}}(x)
=\lambda^{\frac{i}{4}}J_{\widehat{\Phi}}\Lambda_{\Phi_{\delta}}(\sigma^{\Psi}_{-\frac{i}{2}}(x^*))
=\lambda^{\frac{i}{4}}\Lambda_{\Phi}(R\circ\sigma^{\Psi}_{\frac{i}{2}}(x))\\
&=\lambda^{\frac{i}{4}}\Lambda_{\Phi}(\sigma^{\Psi}_{\frac{-i}{2}}(R(x^*)^*)
=\lambda^{\frac{i}{4}}J_{\Phi}\Lambda_{\Phi}(R(x^*))
=\lambda^{\frac{i}{4}}J_{\Phi}J_{\widehat{\Phi}}\Lambda_{\Phi_{\delta}}(x)
\end{aligned}$$
\end{proof}

\subsection{Heisenberg's relations}
We recall that $\alpha(N)\cup\beta(N)\subset
M\subset\widehat{\beta}(N)'$ and
$\alpha(N)\cup\widehat{\beta}(N)\subset \widehat{M}\subset\beta(N)'$
in ${\mathcal L}(H)$.

\begin{prop}
For all $x\in M'$ and $y\in\widehat{M}'$, we have:
$$W(x\surl{\
_{\beta}\otimes_{\alpha}}_{\ N^0} y)=(x\surl{\
_{\alpha}\otimes_{\hat{\beta}}}_{\ N^0} y)W$$
\end{prop}

\begin{proof}
Straightforward by proposition \ref{appartenance} and by definition
of $\widehat{M}$.
\end{proof}

\begin{prop}
The following equalities hold:
$$
\begin{array}{rllcrll}
i)&\ M\cap\widehat{M}&=\alpha(N)&\quad\qquad & ii)&\ M'\cap\widehat{M}&=\widehat{\beta}(N)\\
iii)&\ M\cap\widehat{M}'&=\beta(N) &\quad\qquad & iv)&\
M'\cap\widehat{M}'&=J_{\Phi}\beta(N)J_{\Phi}
\end{array}$$
\end{prop}

\begin{proof}
We start to prove i). We already know that
$M\cap\widehat{M}\supset\beta(N)$. In the other way, let $m\in M
\cap\widehat{M}$. Then, we have by the previous proposition and the
unitarity of $W$:
$$\Gamma(m)=W^*(1\surl{\
_{\alpha}\otimes_{\hat{\beta}}}_{\ N^0}m)W=W^*W(1\surl{\
_{\beta}\otimes_{\alpha}}_{\ N}m)=1\surl{\
_{\beta}\otimes_{\alpha}}_{\ N}m$$ so that $m$ belongs to $\beta(N)$
by proposition \ref{clef2}. Apply $R$ to get iii) and then apply
$\widehat{R}$ to get iv). Finally apply $\widehat{R}$ to i) to get
ii).
\end{proof}

\part{Examples}

In this part, we present a variety of measured quantum groupoids.
First of all, we are interested in the so-called adapted measured
quantum groupoids. These are a class of measured quantum groupoids
with much less complicated axioms because we are able to construct
the antipode. The axiomatic is inspired by J. Kustermans and S.
Vaes' locally quantum groups with a weak condition on the basis.
That is what we develop first. We also characterize adapted measured
quantum groupoids and their dual among measured quantum groupoids.
Then, we give different examples of adapted measured quantum
groupoids and, in particular, the case of groupoids and quantum
groups. In a second time, we investigate inclusions of von Neumann
algebras of depth 2 which can be seen as measured quantum groupoids
but they are not in general of adapted measured quantum groupoids'
type. Finally, we explain how to produce new examples from well
known measured quantum groupoids thanks to simple operations.

We want to lay stress on a fact: historically speaking, the notion
of adapted measured quantum groupoid was the first one we introduce.
The main interest of the structure is the rather quite simple
axioms. So it is easier to find examples (see sections
\ref{groupoids}, \ref{finite} \ref{quantumgroups},
\ref{compactcase}). But we discovered examples of quantum space
quantum groupoid (section \ref{qsqg}) and pairs quantum groupoid
(section \ref{pqg}) duals of which are not adapted measured quantum
groupoid anymore that is we have not a dual structure within
category of adapted measured quantum groupoid. Moreover this
category do not cover all inclusions of von Neumann algebras
(section \ref{ivnasur}). That's why we introduce a larger category
the now so-called measured quantum groupoid which answer all the
problems.

\section{Adapted measured quantum groupoids}\label{fornu}

In this section, we introduce a new natural hypothesis which gives a
link between the right (resp. left) invariant operator-valued weight
and the (resp. anti-) representation of the basis.

\subsection{Definitions} \label{pop}
\begin{defi}
We say that a n.s.f operator-valued weight $T_L$ from $M$ to
$\alpha(N)$ is {\bf $\beta$-adapted} if there exists a n.s.f weight
$\nu_L$ on $N$ such that:
$$\sigma_t^{T_L}(\beta(n))=\beta(\sigma_{-t}^{\nu_L}(n))$$ for all
$n\in N$ and $t\in\mathbb{R}$. We also say that $T_L$ is
$\beta$-adapted w.r.t $\nu_L$.

We say that a n.s.f operator-valued weight $T_R$ from $M$ to
$\beta(N)$ is {\bf $\alpha$-adapted} if there exists a n.s.f weight
$\nu_R$ on $N$ such that:
$$\sigma_t^{T_R}(\alpha(n))=\alpha(\sigma_t^{\nu_R}(n))$$ for all
$n\in N$ and $t\in\mathbb{R}$. We also say that $T_R$ is
$\alpha$-adapted w.r.t $\nu_R$.
\end{defi}

\begin{defi}
A Hopf bimodule $(N,M,\alpha,\beta,\Gamma)$ with left (resp. right)
invariant n.s.f operator-valued weight $T_L$ (resp. $T_R$) from $M$
to $\alpha(N)$ (resp. $\beta(N)$) is said to be a {\bf adapted
measured quantum groupoid} if there exists a n.s.f weight $\nu$ on
$N$ such that $T_L$ is $\beta$-adapted w.r.t $\nu$ and $T_R$ is
$\alpha$-adapted w.r.t $\nu$. Then, we denote by
$(N,M,\alpha,\beta,\Gamma,\nu,T_L,T_R)$ the adapted measured quantum
groupoid and we say that $\nu$ is {\bf quasi-invariant}.
\end{defi}

\begin{rema}
If a n.s.f operator-valued weight $T_L$ from $M$ to $\alpha(N)$ is
$\beta$-adapted w.r.t $\nu$ and if $R$ is a co-involution of $M$,
then the n.s.f operator-valued weight $R\circ T_L\circ R$ from $M$
to $\beta(N)$ is $\alpha$-adapted w.r.t the same weight $\nu$.
\end{rema}

\begin{lemm}\label{timpo}
If $\mu$ is a n.s.f weight on $N$ and if an operator-valued weight
$T_L$ is $\beta$-adapted w.r.t $\nu$, then there exists an
operator-valued weight $S^{\mu}$ from $M$ to $\beta(N)$, which is
$\alpha$-adapted w.r.t $\mu$ such that $\mu\circ\alpha^{-1}\circ
T_L=\nu\circ\beta^{-1}\circ S^{\mu}$. Also, if $\chi$ is a n.s.f
weight on $N$ and if an operator-valued weight $T_R$ is
$\alpha$-adapted w.r.t $\nu$, then there exists an operator-valued
weight $S_{\chi}$ from $M$ to $\alpha(N)$ normal, which is
$\beta$-adapted w.r.t $\chi$ such that $\chi\circ\beta^{-1}\circ
T_R=\nu\circ\beta^{-1}\circ S_{\chi}$.
\end{lemm}

\begin{proof}
For all $n\in N$ and $t\in\mathbb{R}$, we have
$\sigma_t^{\mu\circ\alpha^{-1}\circ
T_L}(\beta(n))=\sigma_t^{\nu\circ\beta^{-1}}(\beta(n))$. By
Haagerup's theorem, we obtain the existence of $S^{\mu}$ which is
clearly adapted. The second part of the lemma is very similar.
\end{proof}

Let $(N,M,\alpha,\beta,\Gamma,\nu,T_L,T_R)$ be a adapted measured
quantum groupoid. Then the opposite adapted measured quantum
groupoid is
$(N^o,M,\beta,\alpha,\varsigma_N\circ\Gamma,\nu^o,T_R,T_L)$. We put:
$$\Phi=\nu\circ\alpha^{-1}\circ T_L \quad\text{ and }\quad \Psi=\nu\circ\beta^{-1}\circ
T_R$$ We also put $S^{\nu}=S_L$ and $S_{\nu}=S_R$. By \ref{prem} and
\ref{evi}, we have: $$\Lambda_{\Phi}({\mathcal
T}_{\Phi,S_L})\subseteq J_{\Phi}\Lambda_{\Phi}({\mathcal
N}_{\Phi}\cap {\mathcal N}_{S_L})\subseteq
D((H_{\Phi})_{\beta},\nu^o)$$ and we have
$R^{\beta,\nu^o}(J_{\Phi}\Lambda_{\Phi}(a))=J_{\Phi}\Lambda_{S_L}(a)J_{\nu}$
for all $a\in {\mathcal N}_{\Phi}\cap {\mathcal N}_{S_L}$.

\subsection{Antipode}\label{fornu}
Then we construct a closed antipode with polar decomposition which
leads to a co-involution and a one-parameter group of automorphisms
of $M$ called scaling group.

\subsubsection{The operator $G$}

We construct now an closed unbounded operator on $H_{\Phi}$ with
polar decomposition which gives needed elements to construct the
antipode. We have the following lemmas:

\begin{lemm}\label{jrel}
For all $\lambda \in \mathbb C$, $x\in {\mathcal
D}(\sigma_{i\lambda}^{\nu})$ and $\xi,\xi' \in
\Lambda_{\Phi}({\mathcal T}_{\Phi,T_L})$, we have:

\begin{equation}\label{rela}
\begin{aligned}
\alpha(x)\Delta_{\Phi}^{\lambda}
 &\subseteq \Delta_{\Phi}^{\lambda}
\alpha(\sigma_{i\lambda}^{\nu}(x)) \\
R^{\alpha,\nu}(\Delta_{\Phi}^{\lambda}\xi)\Delta_{\nu}^{\lambda}
&\subseteq \Delta_{\Phi}^{\lambda}R^{\alpha,\nu}(\xi) \\
\text{and }
\sigma_{i\lambda}^{\nu}(<\Delta_{\Phi}^{\lambda}\xi,\xi'>_{\alpha,\nu})&=<\xi,\Delta_{\Phi}^{\overline{\lambda}}\xi'>_{\alpha,\nu}
\end{aligned}
\end{equation}

and:

\begin{equation}\label{rela2}
\begin{aligned}
\hat{\beta}(x)\Delta_{\Phi}^{\lambda}
 &\subseteq \Delta_{\Phi}^{\lambda}
\hat{\beta}(\sigma_{i\lambda}^{\nu}(x)) \\
R^{\hat{\beta},\nu^o}(\Delta_{\Phi}^{\lambda}\xi)\Delta_{\nu}^{\lambda}
&\subseteq \Delta_{\Phi}^{\lambda}R^{\hat{\beta},\nu^o}(\xi) \\
\text{and }
\sigma_{i\lambda}^{\nu}(<\Delta_{\Phi}^{\lambda}\xi,\xi'>_{\hat{\beta},\nu^o})&=<\xi,\Delta_{\Phi}^{\overline{\lambda}}\xi'>_{\hat{\beta},\nu^o}.
\end{aligned}
\end{equation}
\end{lemm}

\begin{proof}
Straightforward.
\end{proof}

Then, by \cite{S3} and proposition \ref{tenscomp}, we can define a
closed operator $\Delta_{\Phi}^{\lambda} \surl{\ _{\alpha}
\otimes_{\hat{\beta}}}_{\ \ \nu^o} \Delta_{\Phi}^{\lambda}$ which
naturally acts on elementary tensor products for all $\lambda \in
\mathbb C$. Moreover, for all $n\in N$, we have
$J_{\Phi}\alpha(n)=\hat{\beta}(n^*)J_{\Phi}$, so that we can define
a unitary anti-linear operator:
$$J_{\Phi}\surl{\ _{\alpha} \otimes_{\hat{\beta}}}_{\ \ \nu^o} J_{\Phi}:
H_{\Phi}\surl{\ _{\alpha} \otimes_{\hat{\beta}}}_{\ \ \nu^o}H_{\Phi}
\rightarrow H_{\Phi}\surl{\ _{\hat{\beta}} \otimes_{\alpha}}_{\ \nu}
H_{\Phi}$$ such that the adjoint is $J_{\Phi} \surl{\ _{\hat{\beta}}
\otimes_{\alpha}}_{\ \nu} J_{\Phi}$. Also, by composition, it is
possible to define a natural closed anti-linear operator:
$$S_{\Phi} \surl{\ _{\alpha}
\otimes_{\hat{\beta}}}_{\ \ \nu^o} S_{\Phi}: H_{\Phi}\surl{\
_{\alpha} \otimes_{\hat{\beta}}}_{\ \ \nu^o} H_{\Phi}\rightarrow
H_{\Phi}\surl{\ _{\hat{\beta}} \otimes_{\alpha}}_{\ \nu}H_{\Phi}$$
In the same way, if $F_{\Phi}=S_{\Phi}^*$, then it is possible to
define a natural closed anti-linear operator: $F_{\Phi}\surl{\
_{\hat{\beta}}\otimes_{\alpha}}_{\ \nu}F_{\Phi} : H_{\Phi}\surl{\
_{\hat{\beta}}\otimes_{\alpha}}_{\ \nu} H_{\Phi}\rightarrow
H_{\Phi}\surl{\ _{\alpha}\otimes_{\hat{\beta}}}_{\ \ \nu^o}H_{\Phi}$
and we have:
$$(S_{\Phi} \surl{\ _{\alpha} \otimes_{\hat{\beta}}}_{\ \ \nu^o}
S_{\Phi})^*=F_{\Phi}\surl{\ _{\hat{\beta}}\otimes_{\alpha}}_{\
\nu}F_{\Phi}$$

\begin{lemm}\label{chia}
For all $c\in ({\mathcal N}_{\Phi}\cap {\mathcal
N}_{T_L})^*({\mathcal N}_{\Psi}\cap {\mathcal N}_{T_R})$, $e\in
{\mathcal N}_{\Phi}\cap {\mathcal N}_{T_L}$ and all net $(e_k)_{k\in
K}$ of elements of ${\mathcal N}_{\Psi}\cap {\mathcal N}_{T_R}$
weakly converging to $1$, then
$(\lambda_{J_{\Psi}\Lambda_{\Psi}(e_k)}^{\beta,\alpha})^*(1\surl{\
_{\beta} \otimes_{\alpha}}_{\
N}J_{\Phi}eJ_{\Phi})U_{H_{\Psi}}\rho_{\Lambda_{\Phi}(c^*)}^{\alpha,\hat{\beta}}
$ converges to
$(\lambda_{\Lambda_{\Psi}(c)}^{\hat{\alpha},\beta})^*U'^*_{H_{\Phi}}
\rho_{J_{\Phi}\Lambda_{\Phi}(e)}^{\beta,\alpha}$ in the weak
topology.
\end{lemm}

\begin{proof}
By \ref{raccourci}, we have, for all $k\in K$:

$$
\begin{aligned}
&\ \quad
(\lambda_{J_{\Psi}\Lambda_{\Psi}(e_k)}^{\beta,\alpha})^*(1\surl{\
_{\beta} \otimes_{\alpha}}_{\
\nu}J_{\Phi}eJ_{\Phi})U_{H_{\Psi}}\rho_{\Lambda_{\Phi}(c^*)}^{\alpha,\hat{\beta}}
\\
&=(\lambda_{J_{\Psi}\Lambda_{\Psi}(e_k)}^{\beta,\alpha})^*\Gamma(c^*)
\rho_{J_{\Phi}\Lambda_{\Phi}(e)}^{\beta,\alpha}=\left(
\Gamma(c)\lambda_{J_{\Psi}\Lambda_{\Psi}(e_k)}^{\beta,\alpha}\right)^*
\rho_{J_{\Phi}\Lambda_{\Phi}(e)}^{\beta,\alpha}\\
&=\left( (J_{\Psi}e_kJ_{\Psi}\surl{\ _{\beta} \otimes_{\alpha}}_{\
N}1)U'_{H_{\Phi}}\lambda_{\Lambda_{\Psi}(c)}^{\hat{\alpha},\beta}\right)^*
\rho_{J_{\Phi}\Lambda_{\Phi}(e)}^{\beta,\alpha}\\
&=(\lambda_{\Lambda_{\Psi}(c)}^{\hat{\alpha},\beta})^*U'^*_{H_{\Phi}}(J_{\Psi}e_k^*J_{\Psi}\surl{\
_{\beta} \otimes_{\alpha}}_{\ N}1)
\rho_{J_{\Phi}\Lambda_{\Phi}(e)}^{\beta,\alpha}=(\lambda_{\Lambda_{\Psi}(c)}^{\hat{\alpha},\beta})^*U'^*_{H_{\Phi}}
\rho_{J_{\Phi}\Lambda_{\Phi}(e)}^{\beta,\alpha}J_{\Psi}e_k^*J_{\Psi}
\end{aligned}$$
This computation implies the lemma.
\end{proof}

\begin{lemm}\label{dchia}
If $c\in ({\mathcal N}_{\Phi}\cap {\mathcal N}_{T_L})^*({\mathcal
N}_{\Psi}\cap {\mathcal N}_{T_R})$, $e\in {\mathcal N}_{\Phi}\cap
{\mathcal N}_{T_L}$, $\eta\in H_{\Psi}$, $v\in H_{\Phi}$ and a net
$(e_k)_{k\in K}$ of ${\mathcal N}_{\Psi}\cap {\mathcal N}_{T_R}$
converges weakly to $1$, then the net:
$$((U_{H_{\Psi}}(\eta\surl{\ _{\alpha} \otimes_{\hat{\beta}}}_{\ \
\nu^o} \Lambda_{\Phi}(c^*))|J_{\Psi}\Lambda_{\Psi}(e_k)\surl{\
_{\beta} \otimes_{\alpha}}_{\ \nu} J_{\Phi}e^*J_{\Phi}v))_{k\in K}$$
converges to $(
\eta|(\rho_{J_{\Phi}\Lambda_{\Phi}(e)})^*U'_{H_{\Phi}}(\Lambda_{\Psi}(c)\surl{\
_{\hat{\alpha}} \otimes_{\beta}}_{\ \ \nu^o}v))$.
\end{lemm}

\begin{proof}
It's a re-formulation of the previous lemma.
\end{proof}

\begin{prop}
Let $(\eta_i)_{i\in I}$ be a $(N,\nu)$-basis of $\ _{\alpha}H$,
$\Xi\in H_{\Psi}\surl{\ _{\beta} \otimes_{\alpha}}_{\ \nu} H$, $u\in
D(_{\alpha}H,\nu)$, $c\in ({\mathcal N}_{\Phi}\cap {\mathcal
N}_{T_L})^*({\mathcal N}_{\Psi}\cap {\mathcal N}_{T_R})$, $h\in
{\mathcal N}_{\Phi}\cap {\mathcal N}_{T_L}$ and $e$ be an element of
${\mathcal N}_{\Phi}\cap {\mathcal N}_{T_L}\cap {\mathcal
N}_{\Phi}^*\cap {\mathcal N}_{T_L}^*$. Then, we have:

$$\lim_k\sum_{i\in I} (\eta_i\surl{\ _{\alpha}
\otimes_{\hat{\beta}}}_{\ \ \nu^o}
h^*(\lambda^{\beta,\alpha}_{J_{\Phi}\Lambda_{\Phi}(e_k)})^*U_{H_{\Psi}}(
(\rho_{\eta_i}^{\beta,\alpha})^*\Xi\surl{\ _{\alpha}
\otimes_{\hat{\beta}}}_{\ \ \nu^o}\Lambda_{\Phi}(c^*))| u\surl{\
_{\alpha} \otimes_{\hat{\beta}}}_{\ \
\nu^o}J_{\Phi}\Lambda_{\Phi}(e^*))$$ exists and is equal to
$((\rho^{\beta,\alpha}_u)^*\Xi|
(\rho^{\beta,\alpha}_{J_{\Phi}\Lambda_{\Phi}(e)})^*
U'_{H_{\Psi}}(\Lambda_{\Psi}(c)\surl{\ _{\hat{\alpha}}
\otimes_{\beta}}_{\ \ \nu^o}\Lambda_{\Phi}(h)))$.
\end{prop}

\begin{proof}
By \ref{base} and \ref{comm}, we can compute,  for all $i\in I$ and
$k\in K$:

$$\begin{aligned}
&\quad\ (\eta_i\surl{\ _{\alpha} \otimes_{\hat{\beta}}}_{\ \ \nu^o}
h^*(\lambda^{\beta,\alpha}_{J_{\Phi}\Lambda_{\Phi}(e_k)})^*U_{H_{\Psi}}(
(\rho_{\eta_i}^{\beta,\alpha})^*\Xi\surl{\ _{\alpha}
\otimes_{\hat{\beta}}}_{\ \ \nu^o}\Lambda_{\Phi}(c^*))| u\surl{\
_{\alpha} \otimes_{\hat{\beta}}}_{\ \
\nu^o}J_{\Phi}\Lambda_{\Phi}(e^*))\\
&=\!\!(\hat{\beta}(<\eta_i,u>_{\alpha,\nu})
(\lambda^{\beta,\alpha}_{J_{\Phi}\Lambda_{\Phi}(e_k)})^*U_{H_{\Psi}}(
(\rho_{\eta_i}^{\beta,\alpha})^*\Xi\surl{\ _{\alpha}
\otimes_{\hat{\beta}}}_{\ \
\nu^o}\Lambda_{\Phi}(c^*))|gJ_{\Phi}\Lambda_{\Phi}(e^*))\\
&=\!\!((\lambda^{\beta,\alpha}_{J_{\Phi}\Lambda_{\Phi}(e_k)}\!)^*\!(1\!\!\surl{\
_{\beta} \otimes_{\alpha}}_{\
\nu}\!\!\hat{\beta}\!(<\eta_i,u>_{\alpha,\nu}\!)U_{H_{\Psi}}\!\!(
(\rho_{\eta_i}^{\beta,\alpha})^*\Xi\!\!\surl{\ _{\alpha}
\otimes_{\hat{\beta}}}_{\ \
\nu^o}\!\Lambda_{\Phi}(c^*)\!)|J_{\Phi}e^*J_{\Phi}\!\Lambda_{\Phi}(h)))\\
&=\!\!((\lambda^{\beta,\alpha}_{J_{\Phi}\Lambda_{\Phi}(e_k)})^*U_{H_{\Psi}}(
\beta(<\eta_i,u>_{\alpha,\nu})(\rho_{\eta_i}^{\beta,\alpha})^*
\Xi\surl{\ _{\alpha} \otimes_{\hat{\beta}}}_{\ \
\nu^o}\Lambda_{\Phi}(c^*))|J_{\Phi}e^*J_{\Phi}\Lambda_{\Phi}(h))
\end{aligned}$$
Take the sum over $i$ to obtain:
$$\sum_{i\in I}(\eta_i\surl{\ _{\alpha} \otimes_{\hat{\beta}}}_{\ \
\nu^o}
h^*(\lambda^{\beta,\alpha}_{J_{\Phi}\Lambda_{\Phi}(e_k)})^*U_{H_{\Psi}}(
(\rho_{\eta_i}^{\beta,\alpha})^*\Xi\surl{\ _{\alpha}
\otimes_{\hat{\beta}}}_{\ \ \nu^o}\Lambda_{\Phi}(c^*))| u\surl{\
_{\alpha} \otimes_{\hat{\beta}}}_{\ \
\nu^o}J_{\Phi}\Lambda_{\Phi}(e^*))$$
$$=(U_{H_{\Psi}}((\rho_u^{\beta,\alpha})^* \Xi\surl{\
_{\alpha} \otimes_{\hat{\beta}}}_{\ \
\nu^o}\Lambda_{\Phi}(c^*))|J_{\Phi}\Lambda_{\Phi}(e_k)\surl{\
_{\beta} \otimes_{\alpha}}_{\
\nu}J_{\Phi}e^*J_{\Phi}\Lambda_{\Phi}(h))$$ so that lemma
\ref{dchia} implies:
$$\lim_k\sum_{i\in I}(\eta_i\surl{\ _{\alpha} \otimes_{\hat{\beta}}}_{\ \
\nu^o}
h^*(\lambda^{\beta,\alpha}_{J_{\Phi}\Lambda_{\Phi}(e_k)})^*U_{H_{\Psi}}(
(\rho_{\eta_i}^{\beta,\alpha})^*\Xi\surl{\ _{\alpha}
\otimes_{\hat{\beta}}}_{\ \ \nu^o}\Lambda_{\Phi}(c^*))| u\surl{\
_{\alpha} \otimes_{\hat{\beta}}}_{\ \
\nu^o}J_{\Phi}\Lambda_{\Phi}(e^*))$$
$$=((\rho^{\beta,\alpha}_u)^*\Xi|
(\rho^{\beta,\alpha}_{J_{\Phi}\Lambda_{\Phi}(e)})^*
U'_{H_{\Phi}}(\Lambda_{\Psi}(c)\surl{\ _{\hat{\alpha}}
\otimes_{\beta}}_{\ \ \nu^o}\Lambda_{\Phi}(h)))\\
$$
\end{proof}

\begin{prop}\label{cx}
For all $a,c\in ({\mathcal N}_{\Phi}\cap {\mathcal
N}_{T_L})^*({\mathcal N}_{\Psi}\cap {\mathcal N}_{T_R})$, $b,d\in
{\mathcal T}_{\Psi,T_R}$ and $g,h\in {\mathcal T}_{\Phi,S_L}$, the
following vector:
$$U_{H_{\Phi}}^*\Gamma(g^*)(\Lambda_{\Phi}(h)\surl{\ _{\beta}
\otimes_{\alpha}}_{\ \nu}(\lambda^{\beta,\alpha}_{\Lambda_{\Psi}
(\sigma_{-i}^{\Psi}(b^*))})^*U_{H_{\Psi}}(\Lambda_{\Psi}(a) \surl{\
_{\alpha} \otimes_{\hat{\beta}}}_{\ \
\nu^o}\Lambda_{\Phi}((cd)^*)))$$ belongs to ${\mathcal
D}(S_{\Phi}\surl{\ _{\alpha} \otimes_{\hat{\beta}}}_{\ \ \nu^o}
S_{\Phi})$ and the value of $\sigma_{\nu}(S_{\Phi}\surl{\ _{\alpha}
\otimes_{\hat{\beta}}}_{\ \ \nu^o} S_{\Phi})$ on this vector is
equal to:
$$U_{H_{\Phi}}^*\Gamma(h^*)(\Lambda_{\Phi}(g)\surl{\ _{\beta}
\otimes_{\alpha}}_{\
\nu}(\lambda^{\beta,\alpha}_{\Lambda_{\Psi}(\sigma_{-i}^{\Psi}(d^*))})^*
U_{H_{\Psi}}(\Lambda_{\Psi}(c) \surl{\ _{\alpha}
\otimes_{\hat{\beta}}}_{\ \ \nu^o}\Lambda_{\Phi}((ab)^*)))$$
\end{prop}

\begin{proof}
For the proof, let denote by
$\Xi_1=U'_{H_{\Phi}}(\Lambda_{\Psi}(ab)\surl{\ _{\hat{\alpha}}
\otimes_{\beta}}_{\ \ \nu^o}\Lambda_{\Phi}(h))$ and by
$\Xi_2=U'_{H_{\Phi}}(\Lambda_{\Psi}(cd)\surl{\ _{\hat{\alpha}}
\otimes_{\beta}}_{\ \ \nu^o}\Lambda_{\Phi}(g))$. Then, for all
$e,f\in {\mathcal N}_{T_L}\cap {\mathcal N}_{\Phi}\cap{\mathcal
N}_{T_L}^*\cap {\mathcal N}_{\Phi}^*$, the scalar product of
$F_{\Phi}J_{\Phi}\Lambda_{\Phi}(e^*)\!\!\surl{\ _{\alpha}
\otimes_{\hat{\beta}}}_{\ \
\nu^o}\!F_{\Phi}J_{\Phi}\Lambda_{\Phi}(f)$ by:
$$U_{H_{\Phi}}^*\Gamma(g^*)(\Lambda_{\Phi}(h)\!\!\surl{\
_{\beta} \otimes_{\alpha}}_{\
\nu}(\lambda^{\beta,\alpha}_{\Lambda_{\Psi}
(\sigma_{-i}^{\Psi}(b^*))})^* U_{H_{\Psi}}(\Lambda_{\Psi}(a)
\!\!\surl{\ _{\alpha} \otimes_{\hat{\beta}}}_{\ \
\nu^o}\!\Lambda_{\Phi}((cd)^*)))$$ is equal to the scalar product of
$J_{\Phi}\Lambda_{\Phi}(e)\surl{\ _{\alpha}
\otimes_{\hat{\beta}}}_{\ \ \nu^o}J_{\Phi}\Lambda_{\Phi}(f^*)$ by:
$$U_{H_{\Phi}}^*\Gamma(g^*)(\Lambda_{\Phi}(h)\surl{\
_{\beta} \otimes_{\alpha}}_{\
\nu}(\lambda^{\beta,\alpha}_{\Lambda_{\Psi}
(\sigma_{-i}^{\Psi}(b^*))})^* U_{H_{\Psi}}(\Lambda_{\Psi}(a) \surl{\
_{\alpha} \otimes_{\hat{\beta}}}_{\ \
\nu^o}\Lambda_{\Phi}((cd)^*)))$$ By \ref{precau}, this scalar
product is equal to the limit over $k$ of the sum over $i$ of:
$$(J_{\Phi}\Lambda_{\Phi}(e)\!\!\surl{\ _{\alpha}
\otimes_{\hat{\beta}}}_{\ \
\nu^o}\!J_{\Phi}\Lambda_{\Phi}(f^*)|\eta_i\!\!\surl{\ _{\alpha}
\otimes_{\hat{\beta}}}_{\ \ \nu^o}\!\!
g^*(\lambda^{\beta,\alpha}_{J_{\Psi}\Lambda_{\Psi}(e_k)})^*U_{H_{\Psi}}\!(
(\rho_{\eta_i}^{\beta,\alpha})^*\Xi_1\!\!\surl{\ _{\alpha}
\otimes_{\hat{\beta}}}_{\ \ \nu^o}\!\!\Lambda_{\Phi}((cd)^*)))$$ By
the previous proposition applied with $\Xi=\Xi_1$, we get the
symmetric expression:
$$((\rho_{J_{\Phi}\Lambda_{\Phi}(f)}^{\beta,\alpha})^*\Xi_2|
(\rho_{J_{\Phi}\Lambda_{\Phi}(e)}^{\beta,\alpha})^*\Xi_1)$$ so that,
again by the previous proposition applied, this time, with
$\Xi=\Xi_2$ we obtain the limit over $k$ of the sum over $i$ of:
$$(\eta_i\!\!\surl{\ _{\alpha} \otimes_{\hat{\beta}}}_{\ \ \nu^o}\!
h^*(\lambda^{\beta,\alpha}_{J_{\Psi}\Lambda_{\Psi}(e_k)})^*U_{H_{\Psi}}(
(\rho_{\eta_i}^{\beta,\alpha})^*\Xi_2\!\!\surl{\ _{\alpha}
\otimes_{\hat{\beta}}}_{\ \ \nu^o}\!\Lambda_{\Phi}((ab)^*))|
J_{\Phi}\Lambda_{\Phi}(f)\!\!\surl{\ _{\alpha}
\otimes_{\hat{\beta}}}_{\ \ \nu^o}\! J_{\Phi}\Lambda_{\Phi}(e^*))$$

This last expression is equal to the scalar product of:
$$U_{H_{\Phi}}^*\Gamma(h^*)(\Lambda_{\Phi}(g)\surl{\ _{\beta}
\otimes_{\alpha}}_{\
\nu}(\lambda^{\beta,\alpha}_{\Lambda_{\Psi}(\sigma_{-i}^{\Psi}(d^*))})^*
U_{H_{\Psi}}(\Lambda_{\Psi}(c) \surl{\ _{\alpha}
\otimes_{\hat{\beta}}}_{\ \ \nu^o}\Lambda_{\Phi}((ab)^*)))$$ by
$J_{\Phi}\Lambda_{\Phi}(f)\surl{\ _{\alpha}
\otimes_{\hat{\beta}}}_{\ \ \nu^o}J_{\Phi}\Lambda_{\Phi}(e^*)$ and
to the scalar product of:
$$\sigma_{\nu^o}U_{H_{\Phi}}^*\Gamma(h^*)(\Lambda_{\Phi}(g)\surl{\
_{\beta} \otimes_{\alpha}}_{\
\nu}(\lambda^{\beta,\alpha}_{\Lambda_{\Psi}(\sigma_{-i}^{\Psi}(d^*))})^*
U_{H_{\Psi}}(\Lambda_{\Psi}(c) \surl{\ _{\alpha}
\otimes_{\hat{\beta}}}_{\ \ \nu^o}\Lambda_{\Phi}((ab)^*)))$$ by
$J_{\Phi}\Lambda_{\Phi}(e^*)\surl{\ _{\hat{\beta}}
\otimes_{\alpha}}_{\ \ \nu}J_{\Phi}\Lambda_{\Phi}(f)$. Since the
linear span of $J_{\Phi}\Lambda_{\Phi}(e^*)\surl{\ _{\hat{\beta}}
\otimes_{\alpha}}_{\ \ \nu}J_{\Phi}\Lambda_{\Phi}(f)$ where $e,f\in
{\mathcal N}_{T_L}\cap {\mathcal N}_{\Phi}\cap{\mathcal
N}_{T_L}^*\cap {\mathcal N}_{\Phi}^*$ is a core of $F_{\Phi}\surl{\
_{\hat{\beta}}\otimes_{\alpha}}_{\ \nu}F_{\Phi}$, we get that:
$$U_{H_{\Phi}}^*\Gamma(g^*)(\Lambda_{\Phi}(h)\surl{\ _{\beta}
\otimes_{\alpha}}_{\ \nu}(\lambda^{\beta,\alpha}_{\Lambda_{\Psi}
(\sigma_{-i}^{\Psi}(b^*))})^*U_{H_{\Psi}}^*(\Lambda_{\Psi}(a)
\surl{\ _{\alpha} \otimes_{\hat{\beta}}}_{\ \
\nu^o}\Lambda_{\Phi}((cd)^*)))$$ belongs to ${\mathcal
D}(S_{\Phi}\surl{\ _{\alpha} \otimes_{\hat{\beta}}}_{\ \ \nu^o}
S_{\Phi})$ and the value of $S_{\Phi}\surl{\ _{\alpha}
\otimes_{\hat{\beta}}}_{\ \ \nu^o} S_{\Phi}$ on this vector is:
$$\sigma_{\nu^o}U_{H_{\Phi}}^*\Gamma(h^*)(\Lambda_{\Phi}(g)\surl{\ _{\beta}
\otimes_{\alpha}}_{\
\nu}(\lambda^{\beta,\alpha}_{\Lambda_{\Psi}(\sigma_{-i}^{\Psi}(d^*))})^*
U_{H_{\Psi}}(\Lambda_{\Psi}(c) \surl{\ _{\alpha}
\otimes_{\hat{\beta}}}_{\ \ \nu^o}\Lambda_{\Phi}((ab)^*)))$$
\end{proof}

\begin{prop}
There exists a closed densely defined anti-linear operator $G$ on
$H_{\Phi}$ such that the linear span of:
$$(\lambda^{\beta,\alpha}_{\Lambda_{\Psi}(\sigma_{-i}^{\Psi}(b^*))})^*
U_{H_{\Psi}}(\Lambda_{\Psi}(a) \surl{\ _{\alpha}
\otimes_{\hat{\beta}}}_{\ \ \nu^o}\Lambda_{\Phi}((cd)^*))$$ with $
a,c\in ({\mathcal N}_{\Phi}\cap {\mathcal N}_{T_L})^*({\mathcal
N}_{\Psi}\cap {\mathcal N}_{T_R}), b,d\in {\mathcal T}_{\Psi,T_R} $,
is a core of $G$ and we have:

$$G(\lambda^{\beta,\alpha}_{\Lambda_{\Psi}(\sigma_{-i}^{\Psi}(b^*))})^*
U_{H_{\Psi}}(\Lambda_{\Psi}(a) \surl{\ _{\alpha}
\otimes_{\hat{\beta}}}_{\ \ \nu^o}\Lambda_{\Phi}((cd)^*))$$
$$=(\lambda^{\beta,\alpha}_{\Lambda_{\Psi}(\sigma_{-i}^{\Psi}(d^*))})^*
U_{H_{\Psi}}(\Lambda_{\Psi}(c) \surl{\ _{\alpha}
\otimes_{\hat{\beta}}}_{\ \ \nu^o}\Lambda_{\Phi}((ab)^*))$$
Moreover, $G{\mathcal D}(G)={\mathcal D}(G)$ and $G^2=
id_{|{\mathcal D}(G)}$. \label{ferme}
\end{prop}

\begin{proof}
For all $n\in \mathbb N$, let $k_n \in\mathbb N$, $a(n,l),c(n,l)\in
({\mathcal N}_{\Phi}\cap {\mathcal N}_{T_L})^*({\mathcal
N}_{\Psi}\cap {\mathcal N}_{T_R})$ and $b(n,l),d(n,l)\in {\mathcal
T}_{\Psi,T_R}$ and let $w\in H_{\Phi}$ such that:
$$v_n=\sum_{l=1}^{k_n} (\lambda^{\beta,\alpha}_{\Lambda_{\Psi}
(\sigma_{-i}^{\Psi}(b(n,l)^*))})^* U_{H_{\Psi}}
(\Lambda_{\Psi}(a(n,l)) \surl{\ _{\alpha} \otimes_{\hat{\beta}}}_{\
\ \nu^o}\Lambda_{\Phi}((c(n,l)d(n,l))^*))\rightarrow 0$$
$$w_n=\sum_{l=1}^{k_n} (\lambda^{\beta,\alpha}_{\Lambda_{\Psi}
(\sigma_{-i}^{\Psi}(d(n,l)^*))})^*U_{H_{\Psi}}(\Lambda_{\Psi}(c(n,l))
\surl{\ _{\alpha} \otimes_{\hat{\beta}}}_{\ \
\nu^o}\Lambda_{\Phi}((a(n,l)b(n,l))^*))\rightarrow w$$ We have
$U_{H_{\Phi}}^*\Gamma(g^*)(\Lambda_{\Phi}(h)\surl{\ _{\beta}
\otimes_{\alpha}}_{\ \nu}v_n)\in {\mathcal D}(S_{\Phi}\surl{\
_{\alpha} \otimes_{\hat{\beta}}}_{\ \ \nu^o} S_{\Phi})$ for all
$g,h\in {\mathcal T}_{\Phi,S_L}$ and $n\in\mathbb{N}$ by the
previous proposition. Moreover, we have:
$$\sigma_{\nu}(S_{\Phi}\surl{\ _{\alpha}
\otimes_{\hat{\beta}}}_{\ \ \nu^o} S_{\Phi})
U_{H_{\Phi}}^*\Gamma(g^*)(\Lambda_{\Phi}(h)\surl{\ _{\beta}
\otimes_{\alpha}}_{\
\nu}v_n)=U_{H_{\Phi}}^*\Gamma(h^*)(\Lambda_{\Phi}(g)\surl{\ _{\beta}
\otimes_{\alpha}}_{\ \nu}w_n)$$

Since $\Lambda_{\Phi}(g)$ and $\Lambda_{\Phi}(h)$ belongs to
$D((H_{\Phi})_{\beta},\nu^o)$, we obtain:
$$\sigma_{\nu}(S_{\Phi}\surl{\ _{\alpha}
\otimes_{\hat{\beta}}}_{\ \ \nu^o} S_{\Phi})
U_{H_{\Phi}}^*\Gamma(g^*)\lambda_{\Lambda_{\Phi}(h)}^{\beta,\alpha}v_n
=U_{H_{\Phi}}^*\Gamma(h^*)\lambda_{\Lambda_{\Phi}(g)}^{\beta,\alpha}w_n$$
The closure of $S_{\Phi}\surl{\ _{\alpha} \otimes_{\hat{\beta}}}_{\
\ \nu^o} S_{\Phi}$ implies that
$U_{H_{\Phi}}^*\Gamma(h^*)\lambda_{\Lambda_{\Phi}(g)}^{\beta,\alpha}w=0$.
So, apply $U_{H_{\Phi}}$, to get
$\Gamma(h^*)\lambda_{\Lambda_{\Phi}(g)}^{\beta,\alpha}w=0$. Now,
${\mathcal T}_{\Phi,S_L}$ is dense in $M$ that's why
$\lambda_{\Lambda_{\Phi}(g)}^{\beta,\alpha}w=0$ for all $g\in
{\mathcal T}_{\Phi,S_L}$. Then, by \ref{evi}, we have:
$$||\lambda_{\Lambda_{\Phi}(g)}^{\beta,\alpha}w||^2
=(\alpha(<\Lambda_{\Phi}(g),\Lambda_{\Phi}(g)>_{\beta,\nu^o})w|w)\\
=(S_L(\sigma_{i/2}^{\Phi}(g)\sigma_{-i/2}^{\Phi}(g^*))w|w)$$ By
density of ${\mathcal T}_{\Phi,S_L}$, we obtain $||w||^2=0$ i.e
$w=0$. Consequently, the formula given in the proposition for $G$
gives rise to a closable densely defined well-defined operator on
$H_{\Phi}$. So the required operator is the closure of the previous
one.
\end{proof}

Thanks to polar decomposition of the closed operator $G$, we can
give the following definitions:

\begin{defi}\label{Ndefi}
We denote by $D$ the strictly positive operator $G^*G$ on $H_{\Phi}$
(that means positive, self-adjoint and injective) and by $I$ the
anti-unitary operator on $H_{\Phi}$ such that $G=ID^{1/2}$.
\end{defi}

Since $G$ is involutive, we have $I=I^*$, $I^2=1$ and $IDI=D^{-1}$.

\subsubsection{A fundamental commutation relation}
In this section, we establish a commutation relation between $G$ and
the elements $(\omega_{v,w}*id)(U'_{H_{\Phi}})$. We recall that
$W'=U'_{H_{\Psi}}$. We begin by two lemmas borrowed from \cite{E2}.

\begin{lemm}
Let $\xi_i$ be a $(N^o,\nu^o)$-basis of $(H_{\Psi})_{\beta}$. For
all $w'\in D(_{\hat{\alpha}}H_{\Psi},\nu)$ and $w\in H_{\Psi}$, we
have: $$W'(w'\surl{\ _{\hat{\alpha}} \otimes_{\beta}}_{\ \ \nu^o}
w)=\sum_{i}\xi_i \surl{\ _{\beta} \otimes_{\alpha}}_{\ \nu}
(\omega_{w',\xi_i} * id)(W')w$$ If we put
$\delta_i=(\omega_{w',\xi_i}
* id)(W')w$, then
$\alpha(<\xi_i,\xi_i>_{\beta,\nu^o})\delta_i=\delta_i$. Moreover, if
$w\in D(_{\hat{\alpha}}(H_{\Psi}),\nu)$, then $\delta_i\in
D(_{\hat{\alpha}}(H_{\Psi}),\nu)$.

For all $v,v'\in D((H_{\Psi})_{\beta},\nu^o)$ and $i\in I$, there
exists \mbox{$\zeta_i\in D((H_{\Psi})_{\beta},\nu^o)$} such that
$\alpha(<\xi_i,\xi_i>_{\beta,\nu^o})\zeta_i=\zeta_i$ and:
$$W'(v'\surl{\ _{\hat{\alpha}} \otimes_{\beta}}_{\ \ \nu^o}
v)=\sum_{i}\xi_i \surl{\ _{\beta} \otimes_{\alpha}}_{\ \nu}
\zeta_i$$
\end{lemm}

\begin{proof}
Lemma 3.4 of \cite{E2}.
\end{proof}

\begin{rema}
If $v,v'\in \Lambda_{\Psi}({\mathcal T}_{\Psi,T_R}) \subseteq
D(_{\hat{\alpha}}H, \nu) \cap D(H_{\beta},\nu^o)$, then, with
notations of the previous lemma, we have $\zeta_i\in
D(_{\hat{\alpha}}H, \nu) \cap D(H_{\beta},\nu^o)$.
\end{rema}

\begin{lemm}\label{alg}
Let $v,v'\in D(H_{\beta},\nu^o)$ and $w,w'\in D(_{\hat{\alpha}}H,
\nu)$. With notations of the previous lemma, we have:
$$(\omega_{v,w}*id)(U_H'{}^{\hspace{-.15cm}*})(\omega_{v',w'}*
id)(U_H'{}^{\hspace{-.15cm}*})=\sum_i(\omega_{\zeta_i,\delta_i} *
id)(U_H'{}^{\hspace{-.15cm}*})$$ in the norm convergence (and also
in the weak convergence).
\end{lemm}

\begin{proof}
Proposition 3.6 of \cite{E2}.
\end{proof}

\begin{lemm}
Let $a,c$ belonging to $({\mathcal N}_{\Phi}\cap {\mathcal
N}_{T_L})^*({\mathcal N}_{\Psi}\cap {\mathcal N}_{T_R})$. For all
$b,d,a',b',c',d'\in {\mathcal T}_{\Psi,T_R}$, the value of
$(\lambda_{\Lambda_{\Psi}(\sigma_{-i}^{\Psi}(b'^*))}
^{\beta,\alpha})^*U_{H_{\Psi}}$ on the sum over $i$ of:
$$\Lambda_{\Psi}((\omega_{\Lambda_{\Psi}(ab),\xi_i}\!*id)
(W')a')\!\!\surl{\ _{\alpha} \otimes_{\hat{\beta}}}_{\ \
\nu^o}\!\Lambda_{\Phi}(\!(c'd')^*
(\omega_{\xi_i,\Lambda_{\Psi}(cd)}*id)(W'^*))$$ is equal to:
$$(\omega_{\Lambda_{\Psi}(a'b'),\Lambda_{\Psi}(c'd')}*id)(U_{H_{\Phi}}'{}^{\hspace{-.3cm}*}\,)
(\lambda_{\Lambda_{\Psi}(\sigma_{-i}^{\Psi}(b^*))}^{\beta,\alpha})^*
U_{H_{\Psi}}(\Lambda_{\Psi}(a)\surl{\ _{\alpha}
\otimes_{\hat{\beta}}}_{\ \ \nu^o}\Lambda_{\Phi}((cd)^*))$$
\end{lemm}

\begin{proof}
First, let's suppose that $a\in {\mathcal T}_{\Psi,T_R}$. By
\ref{rap} and \ref{corres}, we have:
$$
\begin{aligned}
&\ \quad
(\omega_{\Lambda_{\Psi}(a'b'),\Lambda_{\Psi}(c'd')}*id)(U_{H_{\Phi}}'{}^{\hspace{-.3cm}*}\,)
(\lambda_{\Lambda_{\Psi}(\sigma_{-i}^{\Psi}(b^*))}^{\beta,\alpha})^*
U_{H_{\Psi}}(\Lambda_{\Psi}(a)\surl{\ _{\alpha}
\otimes_{\hat{\beta}}}_{\ \ \nu^o}\Lambda_{\Phi}((cd)^*))\\
&=(\omega_{\Lambda_{\Psi}(a'b'),\Lambda_{\Psi}(c'd')}*id)(U_{H_{\Phi}}'{}^{\hspace{-.3cm}*}\,)
\Lambda_{\Phi}((\omega_{\Lambda_{\Psi}(a),
\Lambda_{\Psi}(\sigma_{-i}^{\Psi}(b^*))}\surl{\ _{\beta}
\otimes_{\alpha}}_{\ \nu}id)(\Gamma((cd)^*)))\\
&=(\omega_{\Lambda_{\Psi}(a'b'),\Lambda_{\Psi}(c'd')}*id)(U_{H_{\Phi}}'{}^{\hspace{-.3cm}*}\,)
\Lambda_{\Phi}((\omega_{\Lambda_{\Psi}(ab),\Lambda_{\Psi}(cd)}*id)(U_{H_{\Phi}}'{}^{\hspace{-.3cm}*}\,))
\end{aligned}$$
By \ref{alg} and the closure of $\Lambda_{\Phi}$, this expression is
equal to the sum over $i\in I$ of:
$$\Lambda_{\Phi}((\omega_{(\omega_{\Lambda_{\Psi}(ab),\xi_i}*id)
(W')\Lambda_{\Psi}(a'b'),(\omega_{\Lambda_{\Psi}(cd),\xi_i}*id)(W')
\Lambda_{\Psi}(c'd')}*id)(U_{H_{\Phi}}'{}^{\hspace{-.3cm}*}\,))$$
Again, \ref{rap} and \ref{corres}, we obtain the sum over $i\in I$
of the value of $(\lambda_{\Lambda_{\Psi}(\sigma_{-i}^{\Psi}(b'^*))}
^{\beta,\alpha})^*U_{H_{\Psi}}$ on:
$$\Lambda_{\Psi}((\omega_{\Lambda_{\Psi}(ab),\xi_i}*id)
(W')a')\!\!\surl{\ _{\alpha} \otimes_{\hat{\beta}}}_{\ \
\nu^o}\Lambda_{\Phi}((c'd')^*(\omega_{\xi_i,\Lambda_{\Psi}(cd)}*id)(W'^*))$$
A density argument finishes the proof.
\end{proof}

\begin{prop}
If $v,w \in \Lambda_{\Psi}({\mathcal T}_{\Psi,T_R}^2) \subseteq
D(_{\hat{\alpha}}(H_{\Psi}), \nu) \cap D((H_{\Psi})_{\beta},\nu^o)$,
then we have:
\begin{align} (\omega_{v,w}* id)(U_{H_{\Phi}}'{}^{\hspace{-.3cm}*}\,)G & \subseteq
G(\omega_{w,v}* id)(U_{H_{\Phi}}'{}^{\hspace{-.3cm}*}\,) \label{inclu1}\\
\text{and} \ (\omega_{v,w}* id)(U'_{H_{\Phi}})G^* & \subseteq
G^*(\omega_{v,w}* id)(U'_{H_{\Phi}})\label{inclu2}
\end{align}

\end{prop}

\begin{proof}
Let $a,c\in ({\mathcal N}_{\Phi}\cap {\mathcal N}_{T_L})^*({\mathcal
N}_{\Psi}\cap {\mathcal N}_{T_R})$ and $b,d,a',b',c',d'\in {\mathcal
T}_{\Psi,T_R}$. By definition of $G$, we have:
$$(\lambda_{\Lambda_{\Psi}(\sigma_{-i}^{\Psi}(d^*))}^{\beta,\alpha})^*
U_{H_{\Psi}}(\Lambda_{\Psi}(c)\surl{\ _{\alpha}
\otimes_{\hat{\beta}}}_{\ \ \nu^o}\Lambda_{\Phi}((ab)^*))\in
{\mathcal D}(G)$$ and:
$$(\omega_{\Lambda_{\Psi}(a'b'),\Lambda_{\Psi}(c'd')}*id)(U_{H_{\Phi}}'{}^{\hspace{-.3cm}*}\,)
G(\lambda_{\Lambda_{\Psi}(\sigma_{-i}^{\Psi}(d^*))}^{\beta,\alpha})^*
U_{H_{\Psi}}(\Lambda_{\Psi}(c)\surl{\ _{\alpha}
\otimes_{\hat{\beta}}}_{\ \ \nu^o}\Lambda_{\Phi}((ab)^*))$$
$$=(\omega_{\Lambda_{\Psi}(a'b'),\Lambda_{\Psi}(c'd')}*id)(U_{H_{\Phi}}'{}^{\hspace{-.3cm}*}\,)
(\lambda_{\Lambda_{\Psi}(\sigma_{-i}^{\Psi}(b^*))}^{\beta,\alpha})^*
U_{H_{\Psi}}(\Lambda_{\Psi}(a)\surl{\ _{\alpha}
\otimes_{\hat{\beta}}}_{\ \ \nu^o}\Lambda_{\Phi}((cd)^*))$$ By the
previous lemma, this is the sum over $i\in I$ of
$G(\lambda_{\Lambda_{\Psi}(\sigma_{-i}^{\Psi}(d'^*))}
^{\beta,\alpha})^*U_{H_{\Psi}}$ on:
$$\Lambda_{\Psi}((\omega_{\Lambda_{\Psi}(cd),\xi_i}*id)
(W')c')\!\!\surl{\ _{\alpha} \otimes_{\hat{\beta}}}_{\ \
\nu^o}\!\!\Lambda_{\Phi}((a'b')^*(\omega_{\xi_i,\Lambda_{\Psi}(ab)}*id)(W'^*))$$
Now, $G$ is a closed operator, so that the sum over $i\in I$ of
$(\lambda_{\Lambda_{\Psi}(\sigma_{-i}^{\Psi}(d'^*))}
^{\beta,\alpha})^*U_{H_{\Psi}}$ on:
$$\Lambda_{\Psi}((\omega_{\Lambda_{\Psi}(cd),\xi_i}*id)
(W')c')\!\!\surl{\ _{\alpha} \otimes_{\hat{\beta}}}_{\ \
\nu^o}\Lambda_{\Phi}((a'b')^*
(\omega_{\xi_i,\Lambda_{\Psi}(ab)}*id)(W'^*))$$ belongs to $
{\mathcal D}(G)$ and by the previous lemma, we obtain:
$$
\begin{aligned}
&\
\quad(\omega_{\Lambda_{\Psi}(a'b'),\Lambda_{\Psi}(c'd')}*id)(W'^*)
G(\lambda_{\Lambda_{\Psi}(\sigma_{-i}^{\Psi}(d^*))}^{\beta,\alpha})^*
U_{H_{\Psi}}(\Lambda_{\Psi}(c)\surl{\ _{\alpha}
\otimes_{\hat{\beta}}}_{\ \ \nu^o}\Lambda_{\Phi}((ab)^*))\\
&=G(\omega_{\Lambda_{\Psi}(c'd'),\Lambda_{\Psi}(a'b')}*id)(U_{H_{\Phi}}'{}^{\hspace{-.3cm}*}\,)
(\lambda_{\Lambda_{\Psi}(\sigma_{-i}^{\Psi}(d^*))}^{\beta,\alpha})^*
U_{H_{\Psi}}(\Lambda_{\Psi}(c)\surl{\ _{\alpha}
\otimes_{\hat{\beta}}}_{\ \ \nu^o}\Lambda_{\Phi}((ab)^*))
\end{aligned}$$
Now the linear span:
$$(\lambda^{\beta,\alpha}_{\Lambda_{\Psi}(\sigma_{-i}^{\Psi}(b^*))})^*
U_{H_{\Psi}}(\Lambda_{\Psi}(a) \surl{\ _{\alpha}
\otimes_{\hat{\beta}}}_{\ \ \nu^o}\Lambda_{\Phi}((cd)^*))$$ with
$a,c\in ({\mathcal N}_{\Phi}\cap {\mathcal N}_{T_L})^*({\mathcal
N}_{\Psi}\cap {\mathcal N}_{T_R}), b,d\in {\mathcal T}_{\Psi,T_R}
\}$, is a core for $G$ that's why the first inclusion holds. The
second one is the adjoint of the first one.
\end{proof}

\begin{coro}\label{coN}
For all $v,w \in \Lambda_{\Psi}({\mathcal T}_{\Psi,T_R}^2)$, we
have:
$$(\omega_{v,w}*id)(U'_{H_{\Phi}})D \subseteq
D(\omega_{\Delta_{\Psi}^{-1}v,\Delta_{\Psi}w}* id)(U'_{H_{\Phi}})$$
where $D=G^*G$ is defined in \ref{Ndefi}.
\end{coro}

\begin{proof}
We have:
$$
\begin{aligned}
(\omega_{w,v}* id)(U'_{H_{\Phi}})G &=
(\omega_{S_{\Psi}w,\Delta_{\Psi}S_{\Psi}v} *
id)(U_{H_{\Phi}}'{}^{\hspace{-.3cm}*}\hspace{.1cm})G &&\quad
\text{by lemma
\ref{switch}} \\
& \subseteq G(\omega_{\Delta_{\Psi}S_{\Psi}v,S_{\Psi}w} \star
id)(U_{H_{\Phi}}'{}^{\hspace{-.3cm}*}\hspace{.1cm}) && \quad
\text{by inclusion\  \eqref{inclu1}} \\
&= G(\omega_{\Delta_{\Psi}^{-1}v,\Delta_{\Psi}w} *
id)(U_{H_{\Phi}}'{}^{\hspace{-.3cm}*}\hspace{.1cm}) && \quad
\text{by lemma \ref{switch}}
\end{aligned}$$

In the same way, we can finish the proof:

$$
\begin{aligned}
(\omega_{v,w}* id)(U'_{H_{\Phi}})D &= (\omega_{v,w}*
id)(U'_{H_{\Phi}})G^* &&
\quad \text{by definition \ref{Ndefi}} \\
& \subseteq G^*(\omega_{w,v}* id)(U'_{H_{\Phi}})G &&\quad \text{by
inclusion \
\eqref{inclu2}} \\
& \subseteq G^*G(\omega_{\Delta_{\Psi}^{-1}v,\Delta_{\Psi}w}*
id)(U'_{H_{\Phi}}) &&\\
&= D(\omega_{\Delta_{\Psi}^{-1}v,\Delta_{\Psi}w} *
id)(U'_{H_{\Phi}}) &&\quad \text{by definition \ref{Ndefi}}.
\end{aligned}$$

\end{proof}

\subsubsection{Scaling group}
In this section, we give a sense and we prove the following
commutation relation $U'_{H_{\Phi}}(\Delta_{\Psi} \surl{\
_{\hat{\alpha}} \otimes_{\beta}}_{\ \ \nu^o} D)=(\Delta_{\Psi}
\surl{\ _{\beta} \otimes_{\alpha}}_{\ \nu} D)U'_{H_{\Phi}}$ so as to
construct the scaling group $\tau$.

\begin{lemm}
For all $\lambda \in \mathbb C$ and $x$ analytic w.r.t $\nu$, we
have:
$$\alpha(x)D^{\lambda} \subseteq
D^{\lambda}\alpha(\sigma_{-i\lambda}^{\nu}(x))\text{ and }
\beta(x)D^{\lambda} \subseteq
D^{\lambda}\beta(\sigma_{-i\lambda}^{\nu}(x))$$\label{commuNa}
\end{lemm}

\begin{proof}
For all $a,c \in ({\mathcal N}_{\Phi}\cap {\mathcal
N}_{T_L})^*({\mathcal N}_{\Psi}\cap {\mathcal N}_{T_R})$, $b,d\in
{\mathcal T}_{\Psi,T_R}$ and $x$ analytic w.r.t $\nu$, we have by
\ref{base} and \ref{comm}:
$$
\begin{aligned}
&\ \quad
\beta(x)G(\lambda^{\beta,\alpha}_{\Lambda_{\Psi}(\sigma_{-i}^{\Psi}(b^*))})^*
U_{\Psi}(\Lambda_{\Psi}(a) \surl{\ _{\alpha}
\otimes_{\hat{\beta}}}_{\ \ \nu^o}\Lambda_{\Phi}((cd)^*))\\
&=\beta(x)(\lambda^{\beta,\alpha}_{\Lambda_{\Psi}(\sigma_{-i}^{\Psi}(d^*))})^*
U_{\Psi}(\Lambda_{\Psi}(c) \surl{\ _{\alpha}
\otimes_{\hat{\beta}}}_{\ \ \nu^o}\Lambda_{\Phi}((ab)^*))\\
&=(\lambda^{\beta,\alpha}_{\Lambda_{\Psi}(\sigma_{-i}^{\Psi}(d^*))})^*(1\surl{\
_{\beta} \otimes_{\alpha}}_{\
\nu}\beta(x))U_{\Psi}(\Lambda_{\Psi}(c) \surl{\ _{\alpha}
\otimes_{\hat{\beta}}}_{\ \ \nu^o}\Lambda_{\Phi}((ab)^*))\\
&=(\lambda^{\beta,\alpha}_{\Lambda_{\Psi}(\sigma_{-i}^{\Psi}(d^*))})^*U_{\Psi}(\Lambda_{\Psi}(c)
\surl{\ _{\alpha} \otimes_{\hat{\beta}}}_{\ \
\nu^o}\Lambda_{\Phi}(\beta(x)b^*a^*))\\
&=G(\lambda^{\beta,\alpha}_{\Lambda_{\Psi}(\beta(\sigma_i^{\nu}(x))\sigma_{-i}^{\Psi}(b^*))})^*
U_{\Psi}(\Lambda_{\Psi}(a) \surl{\ _{\alpha}
\otimes_{\hat{\beta}}}_{\ \ \nu^o}\Lambda_{\Phi}((cd)^*))\\
&=G\alpha(\sigma_{-i/2}^{\nu}(x^*))(\lambda^{\beta,\alpha}_{\Lambda_{\Psi}(\sigma_{-i}^{\Psi}(b^*))})^*
U_{\Psi}(\Lambda_{\Psi}(a) \surl{\ _{\alpha}
\otimes_{\hat{\beta}}}_{\ \ \nu^o}\Lambda_{\Phi}((cd)^*))
\end{aligned}$$

Now, the linear span of:
$$(\lambda^{\beta,\alpha}_{\Lambda_{\Psi}(\sigma_{-i}^{\Psi}(b^*))})^*
U_{\Psi}(\Lambda_{\Psi}(a) \surl{\ _{\alpha}
\otimes_{\hat{\beta}}}_{\ \ \nu^o}\Lambda_{\Phi}((cd)^*))$$ where
$a,c \in ({\mathcal N}_{\Phi}\cap {\mathcal N}_{T_L})^*({\mathcal
N}_{\Psi}\cap {\mathcal N}_{T_R}), b,d\in {\mathcal T}_{\Psi,T_R}
\}$, is a core for $G$, so that we have:
$$\beta(x)G \subseteq G\alpha(\sigma_{-i/2}^{\nu}(x^*))$$
Take adjoint to obtain $\alpha(x)G^* \subseteq
G^*\beta(\sigma_{i/2}^{\nu}(x^*))$. So, we conclude by:
$$\alpha(x)D=\alpha(x)G^*G \subseteq
G^*\beta(\sigma_{i/2}^{\nu}(x^*))G \subseteq
D\alpha(\sigma_{-i}^{\nu}(x))$$ The second part of the lemma can be
proved in a very similar way.
\end{proof}

We now state two lemmas analogous to relations (\ref{rela}) and
(\ref{rela2}) for $\Psi$ and we justify the existence of natural
operators:

\begin{lemm}
For all $\lambda\in\mathbb C$, $x\in {\mathcal
D}(\sigma_{-i\lambda}^{\nu})$ and $\xi,\xi' \in
\Lambda_{\Psi}({\mathcal T}_{\Psi,T_R})$, we have:

\begin{equation}\label{relat1}
\begin{aligned}
\beta(x)\Delta_{\Psi}^{\lambda}
 &\subseteq \Delta_{\Psi}^{\lambda}
\beta(\sigma_{-i\lambda}^{\nu}(x)) \\
R^{\beta,\nu^o}(\Delta_{\Psi}^{\lambda}\xi)\Delta_{\nu}^{-\lambda}
&\subseteq \Delta_{\Psi}^{\lambda}R^{\beta,\nu^o}(\xi) \\
\text{and }
\sigma_{-i\lambda}^{\nu}(<\Delta_{\Psi}^{\lambda}\xi,\xi'>_{\beta,\nu^o})&
=<\xi,\Delta_{\Psi}^{\overline{\lambda}}\xi'>_{\beta,\nu^o}
\end{aligned}
\end{equation}

and:

\begin{equation}\label{relat2}
\begin{aligned}
\hat{\alpha}(x)\Delta_{\Psi}^{\lambda}
 &\subseteq \Delta_{\Psi}^{\lambda}
\hat{\alpha}(\sigma_{-i\lambda}^{\nu}(x)) \\
R^{\hat{\alpha},\nu^o}(\Delta_{\Psi}^{\lambda}\xi)\Delta_{\nu}^{-\lambda}
&\subseteq \Delta_{\Psi}^{\lambda}R^{\hat{\alpha},\nu^o}(\xi) \\
\text{and }
\sigma_{-i\lambda}^{\nu}(<\Delta_{\Psi}^{\lambda}\xi,\xi'>_{\hat{\alpha},\nu^o})
&=<\xi,\Delta_{\Psi}^{\overline{\lambda}}\xi'>_{\hat{\alpha},\nu^o}
\end{aligned}
\end{equation}

\end{lemm}

\begin{proof}
It is sufficient to apply \ref{jrel} to the opposite adapted
measured quantum groupoid for example.
\end{proof}

Then, we can define, for all $\lambda \in \mathbb C$, a closed
linear operator $\Delta_{\Psi}^{\lambda} \surl{\ _{\beta}
\otimes_{\alpha}}_{\ \nu} D^{\lambda}$ which naturally acts on
elementary tensor products of $H_{\Psi} \surl{\ _{\beta}
\otimes_{\alpha}}_{\ \nu} H_{\Phi}$. With relations \eqref{relat2}
in hand, we also get a closed linear operator
$\Delta_{\Psi}^{\lambda} \surl{\ _{\hat{\alpha}} \otimes_{\beta}}_{\
\ \nu^o} D^{\lambda}$ on $H_{\Psi} \surl{\ _{\hat{\alpha}}
\otimes_{\beta}}_{\ \ \nu^o} H_{\Phi}$.

\begin{prop}
The following relation holds:
\begin{equation}\label{fond}
U'_{H_{\Phi}}(\Delta_{\Psi} \surl{\ _{\hat{\alpha}}
\otimes_{\beta}}_{\ \ \nu^o} D)=(\Delta_{\Psi} \surl{\ _{\beta}
\otimes_{\alpha}}_{\ \nu} D)U'_{H_{\Phi}}
\end{equation}
\end{prop}

\begin{proof}
By \ref{coN}, we have, for all $v,w \in \Lambda_{\Psi}({\mathcal
T}_{\Psi,T_R})$ and $v',w' \in {\mathcal D}(D)$:

$$
\begin{aligned}
(U'_{H_{\Phi}}(v\surl{\ _{\hat{\alpha}} \otimes_{\beta}}_{\ \ \nu^o}
v')|\Delta_{\Psi}w\surl{\ _{\beta} \otimes_{\alpha}}_{\ \nu} Dw')
&=((\omega_{v,\Delta_{\Psi}w}* id)(U'_{H_{\Phi}})v'|Dw') \\
&=(D(\omega_{\Delta_{\Psi}^{-1}(\Delta_{\Psi}v),\Delta_{\Psi}w}* id)(U'_{H_{\Phi}})v'|w') \\
&= ((\omega_{\Delta_{\Psi}v,w}*id)(U'_{H_{\Phi}})Dv'|w')\\
&=(U'_{H_{\Phi}}(\Delta_{\Psi}v\surl{\ _{\hat{\alpha}}
\otimes_{\beta}}_{\ \ \nu^o} Dv')|w \surl{\ _{\beta}
\otimes_{\alpha}}_{\ \nu} w')
\end{aligned}$$

By definition, we know that $\Lambda_{\Psi}({\mathcal
T}_{\Psi,T_R})\odot {\mathcal D}(D)$ is a core for
$\Delta_{\Psi}\surl{\ _{\beta} \otimes_{\alpha}}_{\ \nu}D$ so, for
all $u \in {\mathcal D}(\Delta_{\Psi}\surl{\ _{\beta}
\otimes_{\alpha}}_{\ \nu}D)$, we have:
$$(U'_{H_{\Phi}}(v\surl{\ _{\hat{\alpha}}
\otimes_{\beta}}_{\ \ \nu^o} v')|(\Delta_{\Psi}\surl{\ _{\beta}
\otimes_{\alpha}}_{\ \nu} D)u)=(U'_{H_{\Phi}}(\Delta_{\Psi}v\surl{\
_{\hat{\alpha}} \otimes_{\beta}}_{\ \ \nu^o} Dv')|u)$$ Since
$\Delta_{\Psi}\surl{\ _{\beta} \otimes_{\alpha}}_{\ \nu}D$ is
self-adjoint, we get:
$$(\Delta_{\Psi}\surl{\ _{\beta} \otimes_{\alpha}}_{\ \nu}D)U'_{H_{\Phi}}(v\surl{\
_{\hat{\alpha}} \otimes_{\beta}}_{\ \ \nu^o}
v')=U'_{H_{\Phi}}(\Delta_{\Psi}v\surl{\ _{\hat{\alpha}}
\otimes_{\beta}}_{\ \ \nu^o}Dv')$$ Finally, since
$\Lambda_{\Psi}({\mathcal T}_{\Psi,T_R})\odot {\mathcal D}(D)$ is a
core for $\Delta_{\Psi}\surl{\ _{\hat{\alpha}} \otimes_{\beta}}_{\ \
\nu^o}D$ and by closeness of $\Delta_{\Psi}\surl{\ _{\beta}
\otimes_{\alpha}}_{\ \nu}D$, we deduce that:
$$U'_{H_{\Phi}}(\Delta_{\Psi}\surl{\ _{\hat{\alpha}}
\otimes_{\beta}}_{\ \ \nu^o}D) \subseteq (\Delta_{\Psi}\surl{\
_{\beta} \otimes_{\alpha}}_{\ \nu} D)U'_{H_{\Phi}}$$ Because of
unitarity of $U'_{H_{\Phi}}$, we get that $(\Delta_{\Psi}\surl{\
_{\hat{\alpha}} \otimes_{\beta}}_{\ \ \nu^o}
D)U_{H_{\Phi}}'{}^{\hspace{-.3cm}*} \subseteq
U_{H_{\Phi}}'{}^{\hspace{-.3cm}*}(\Delta_{\Psi}\surl{\ _{\beta}
\otimes_{\alpha}}_{\ \nu}D)$ and by taking the adjoint, we get the
reverse inclusion:
$$(\Delta_{\Psi}\surl{\ _{\beta} \otimes_{\alpha}}_{\ \nu}D)U'_{H_{\Phi}}
\subseteq U'_{H_{\Phi}}(\Delta_{\Psi}\surl{\ _{\hat{\alpha}}
\otimes_{\beta}}_{\ \ \nu^o}D)$$
\end{proof}

We know begin the construction of the scaling group $\tau$ strictly
speaking. We also prove a theorem which state that $A(U'_H)=M$ and
generalize proposition 1.5 of \cite{KV2}.

\begin{defi}
We denote by $M_R$ the weakly closed linear span of: $$\{
(\omega\surl{\ _{\beta} \star_{\alpha}}_{\ \nu} id)(\Gamma(x))\ | \
x\in M,\ \omega \in M^+_* \text{ s.t } \exists k\in {\mathbb
R}^+\!\!,\ \omega\circ\beta \leq k\nu\}$$ Also, we denote by $M_L$
the weakly closed linear span of:
$$\{ (id\surl{\
_{\beta} \star_{\alpha}}_{\ \nu} \omega)(\Gamma(x))\ | \ x\!\in\!
M,\ \omega\! \in\! M^+_* \text{ s.t } \exists k\in {\mathbb
R}^+\!\!, \ \omega\circ\alpha \leq k\nu\}$$
\end{defi}

By \ref{corres} and \ref{situa}, $M_R$ is equal to the von Neumann
subalgebra $A(U'_H)$ of $M$. Also, $M_L$ is a von Neumann subalgebra
of $M$. Moreover, we know $\alpha(N)\subseteq M_R$ and $\beta(N)
\subseteq M_L$, so that $M_L \surl{\ _{\beta} \star_{\alpha}}_{\
\nu} M_R$ makes sense. Also, we have, for all $m\in M$:

\begin{equation}\label{appart}
\Gamma(m)\in M_L \surl{\ _{\beta} \star_{\alpha}}_{\ N} M_R
\end{equation}

\begin{lemm}\label{tau}
There exists a unique strongly continuous one-parameter group $\tau$
of automorphisms of $M_R$ such that $\tau_t(x)=D^{-it}xD^{it}$ for
all $t \in \mathbb R$ and $x\in M_R$.
\end{lemm}

\begin{proof}
By commutation relation \eqref{fond}, for all $t\in \mathbb R$ and
$v,w \in \Lambda_{\Psi}({\mathcal T}_{\Psi,T_R})$, we get that:
$$D^{-it}(\omega_{v,w} * id)(U'_{H_{\Phi}})D^{it}
=(\omega_{\Delta_{\Psi}^{-it}v,\Delta_{\Psi}^{it}w}*
id)(U'_{H_{\Phi}})$$ Consequently, we obtain $D^{-it}M_RD^{it}=M_R$
which is the only point to show.
\end{proof}

\begin{lemm}\label{mieux}
We have $\tau_t(\alpha(n))=\alpha(\sigma_t^{\nu}(n))$ for all $n\in
N$ and $t\in\mathbb{R}$.
\end{lemm}

\begin{proof}
Straightforward by lemma \ref{commuNa}.
\end{proof}

\begin{prop}\label{egal}
We have $(\sigma_t^{\Psi} \surl{\ _{\beta} \star_{\alpha}}_{\ N}
\tau_{-t})\circ\Gamma=\Gamma \circ \sigma_t^{\Psi}$ for all $t\in
\mathbb R$.
\end{prop}

\begin{proof}
By proposition \ref{tenscomp} and thanks to the previous lemma, it
is possible to define a normal *-automorphism $\sigma_t^{\Psi}
\surl{\ _{\beta} \star_{\alpha}}_{\ N} \tau_{-t}$ of $M\surl{\
_{\beta} \star_{\alpha}}_{\ N} M_R$. By relation \eqref{appart}, the
formula makes sense ($\tau$ is just defined on $M_R$). By relation
\eqref{fond}, we can compute for all $m\in M$ and $t\in\mathbb R$:
$$
\begin{aligned}
(\sigma_t^{\Psi} \surl{\ _{\beta} \star_{\alpha}}_{\ \nu}
\tau_{-t})\circ\Gamma(m)&=(\Delta_{\Psi}^{it} \surl{\ _{\beta}
\otimes_{\alpha}}_{\ \nu}D^{it})\Gamma(m)(\Delta_{\Psi}^{-it}
\surl{\ _{\beta} \otimes_{\alpha}}_{\ \nu}
D^{-it})\\
&=(\Delta_{\Psi}^{it} \surl{\ _{\beta} \otimes_{\alpha}}_{\
\nu}D^{it})U'_{H_{\Phi}}(m\surl{\ _{\hat{\alpha}}
\otimes_{\beta}}_{\ \ \nu^o} \!\!
1)U_{H_{\Phi}}'{}^{\hspace{-.3cm}*} (\Delta_{\Psi}^{-it} \surl{\
_{\beta} \otimes_{\alpha}}_{\ \nu}D^{-it})\\
&=U'_{H_{\Phi}}(\Delta_{\Psi}^{it} \surl{\ _{\hat{\alpha}}
\otimes_{\beta}}_{\ \nu^o}D^{it})(m\surl{\ _{\hat{\alpha}}
\otimes_{\beta}}_{\ \ \nu^o} \!\! 1)(\Delta_{\Psi}^{-it} \surl{\
_{\hat{\alpha}} \otimes_{\beta}}_{\ \nu^o}
D^{-it})U_{H_{\Phi}}'{}^{\hspace{-.3cm}*}\\
&=U'_{H_{\Phi}}(\sigma_t^{\Psi}(m)\surl{\ _{\hat{\alpha}}
\otimes_{\beta}}_{\ \nu^o}
1)U_{H_{\Phi}}'{}^{\hspace{-.3cm}*}=\Gamma(\sigma_t^{\Psi}(m))
\end{aligned}$$

\end{proof}

We are now able to prove that we can re-construct $M$ thanks to the
fundamental unitary.

\begin{theo}\label{densevn}
If $<F>^{-\textsc{w}}$ is the weakly closed linear span of $F$ in
$M$, then following vector spaces:
$$\begin{aligned}
M_R&=<(\omega\surl{\ _{\beta} \star_{\alpha}}_{\ \nu}
id)(\Gamma(m))\ | \ m\in M,\omega\in M^+_*, k\in {\mathbb
R}^+\text{ s.t } \omega\circ\beta \leq k\nu>^{-\textsc{w}}\\
A(U'_H)&=<(\omega_{v,w}* id)(U'_H) \ |\ v \in
D(_{\hat{\alpha}}(H_{\Psi}),\mu), w \in
D((H_{\Psi})_{\beta},\mu^o)>^{-\textsc{w}}\\
M_L&=<(id\surl{\ _{\beta}\star_{\alpha}}_{\ \nu} \omega)(\Gamma(m))\
| \ m\in M, \omega\in M^+_*, k\in {\mathbb
R}^+ \ \text{s.t } \omega\circ\alpha \leq k\nu>^{-\textsc{w}}\\
A(U_H)&=<(id*\omega_{v,w})(U_H) \ |\ v \in
D((H_{\Psi}),\mu^o)_{\hat{\beta}}, w \in
D(_{\alpha}(H_{\Psi}),\mu)>^{-\textsc{w}}
\end{aligned}$$
are equal to the whole von Neumann algebra $M$.
\end{theo}

\begin{proof}
We have already noticed that $M_R=A(U'_H)$ and $M_L=A(U_H)$. Then,
we get inspired by \cite{KV2}. By \ref{mieux}, we have
$\tau_t(\alpha(n))=\alpha(\sigma_t^{\nu}(n))$ so:
$$M_L=<(id\surl{\ _{\beta} \star_{\alpha}}_{\ \nu}
\omega\circ\tau_t)(\Gamma(m))\ | \ m\in M,\omega \in (M_R)^+_* ,
k\in {\mathbb R}^+ \ \text{s.t } \omega\circ\alpha \leq
k\nu>^{-\textsc{w}}$$ By \ref{egal}, we have
$\sigma_t^{\Psi}((id\surl{\ _{\beta} \star_{\alpha}}_{\ \nu}
\omega)\Gamma(m))=(id\surl{\ _{\beta} \star_{\alpha}}_{\ \nu}
\omega\circ\tau_t)\Gamma(\sigma_t^{\Psi}(m))$ that's why
$\sigma_t^{\Psi}(M_L)=M_L$ for all $t\in\mathbb R$. On the other
hand, by \ref{semi}, restriction of $\Psi$ to $M_L$ is semi-finite.
By Takesaki's theorem (\cite{St}, theorem 10.1), there exists a
unique normal and faithful conditional expectation $E$ from $M$ to
$M_L$ such that $\Psi(m)=\Psi(E(m))$ for all $m\in M^+$. Moreover,
if $P$ is the orthogonal projection on the closure of
$\Lambda_{\Psi}({\mathcal N}_{\Psi} \cap M_L)$ then $E(m)P=PmP$.

So the range of $P$ contains $\Lambda_{\Psi}((id \surl{\ _{\beta}
\star_{\alpha}}_{\ \nu} \omega )\Gamma(x))$ for all $\omega$ and
$x\in {\mathcal N}_{\Psi}$. By right version of \ref{dense} implies
that $P=1$ so that $E$ is the identity and $M=M_L$. If we apply the
previous result to the opposite adapted measured quantum groupoid,
then we get that $M=M_R$.
\end{proof}

\begin{coro}
There exists a unique strongly continuous one-parameter group $\tau$
of automorphisms of $M$ such that, for all $t \in \mathbb R$, $m\in
M$ and $n\in N$:
$$\tau_t(m)=D^{-it}mD^{it},\ \tau_t(\alpha(n))=\alpha(\sigma_t^{\nu}(n)) \text{
and } \tau_t(\beta(n))=\beta(\sigma_t^{\nu}(n))$$
\end{coro}

\begin{proof}
Straightforward from the previous theorem and \ref{tau}. First
property comes from \ref{mieux} and the second one from
\ref{commuNa}.
\end{proof}

\begin{defi}
The one-parameter group $\tau$ is called {\bf scaling group}.
\end{defi}

Let us notice that it is possible to define normal *-automorphisms
$\tau_t \surl{\ _{\beta} \star_{\alpha}}_{\ N} \tau_t$ and $\tau_t
\surl{\ _{\beta} \star_{\alpha}}_{\ N} \sigma_t^{\Phi}$ of $M
\surl{\ _{\beta} \otimes_{\alpha}}_{\ N}M$ for all $t\in\mathbb{R}$,
thanks to the previous commutation relations and recalls about
tensor products.

\begin{prop}\label{clef20}
We have $\Gamma \circ \tau_t=(\tau_t \surl{\ _{\beta}
\star_{\alpha}}_{\ N}\tau_t)\circ\Gamma$ for all $t\in \mathbb R$.
\end{prop}

\begin{proof}
By \ref{egal} and co-product relation, we have for all $t\in\mathbb
R$:
$$\begin{aligned}
(id \surl{\ _{\beta} \star_{\alpha}}_{\ \nu} \Gamma)(\sigma_t^{\Psi}
\surl{\ _{\beta} \star_{\alpha}}_{\ \nu} \tau_{-t})\circ\Gamma &=(id
\surl{\ _{\beta} \star_{\alpha}}_{\
\nu} \Gamma)\Gamma \circ \sigma_t^{\Psi}\\
&=(\Gamma \surl{\ _{\beta} \star_{\alpha}}_{\ \nu} id)\Gamma \circ
\sigma_t^{\Psi}=(\Gamma \circ\sigma_t^{\Psi}\surl{\ _{\beta}
\star_{\alpha}}_{\
\nu} \tau_{-t})\Gamma\\
&=(\sigma_t^{\Psi}\surl{\ _{\beta} \star_{\alpha}}_{\ \nu}
\tau_{-t}\surl{\ _{\beta} \star_{\alpha}}_{\ \nu}
\tau_{-t})(\Gamma\surl{\ _{\beta} \star_{\alpha}}_{\ \nu}
id)\Gamma\\
&=(\sigma_t^{\Psi}\surl{\ _{\beta}\star_{\alpha}}_{\ \nu}
[(\tau_{-t}\surl{\ _{\beta}\star_{\alpha}}_{\ \nu}
\tau_{-t})\circ\Gamma])\circ\Gamma
\end{aligned}$$
Consequently, for all $m\in M$, $\omega\in M_*^+$,
$k\in\mathbb{R}^+$ such that $\omega\circ\beta\leq k\nu$, we have:

$$
\begin{aligned}
\Gamma\circ\tau_{-t}\circ((\omega\circ\sigma_t^{\Psi})\surl{\
_{\beta} \star_{\alpha}}_{\ \nu} id)\Gamma&=(\omega\surl{\ _{\beta}
\star_{\alpha}}_{\ \nu} id\surl{\ _{\beta} \star_{\alpha}}_{\ \nu}
id)(\sigma_t^{\Psi}\!\!\surl{\ _{\beta}
\star_{\alpha}}_{\ \nu}\!(\Gamma\circ\tau_{-t}))\circ\Gamma \\
&=(\omega\surl{\ _{\beta}\star_{\alpha}}_{\ \nu} id\surl{\ _{\beta}
\star_{\alpha}}_{\ \nu} id)(\sigma_t^{\Psi}\!\!\surl{\ _{\beta}
\star_{\alpha}}_{\ \nu}\![(\tau_{-t}\surl{\
_{\beta}\star_{\alpha}}_{\ \nu}\tau_{-t})\circ\Gamma])\\
&=[(\tau_{-t}\surl{\ _{\beta}\star_{\alpha}}_{\ \nu}
\tau_{-t})\circ\Gamma]\circ((\omega\circ\sigma_t^{\Psi})\surl{\
_{\beta} \star_{\alpha}}_{\ \nu} id)\Gamma
\end{aligned}$$

The theorem \ref{densevn} allows us to conclude.
\end{proof}

\begin{prop}\label{clef2}
For all $x\in M\cap\alpha(N)'$, we have $\Gamma(x)=1 \!\surl{\
_{\beta} \otimes_{\alpha}}_{\ N} x \Leftrightarrow x\in\beta(N)$.
Also, for all $x\in M\cap\beta(N)'$, we have $\Gamma(x)=x\! \surl{\
_{\beta} \otimes_{\alpha}}_{\ N} 1 \Leftrightarrow x\in\alpha(N)$.
\end{prop}

\begin{proof}
Let $x\in M\cap\alpha(N)'$ such that $\Gamma(x)=1\surl{\ _{\beta}
\otimes_{\alpha}}_{\ N} x$. For all $n \in \mathbb N$, we define in
the strong topology:
$$x_n=\frac{n}{\sqrt{\pi}}\int\!
exp(-n^2t^2)\sigma_t^{\Psi}(x)\ dt \quad \text{analytic w.r.t }
\sigma^{\Psi},$$ and:
$$y_n=\frac{n}{\sqrt{\pi}}\int\!
exp(-n^2t^2)\tau_{-t}(x)\ dt \quad \text{belongs to } \alpha(N)'.$$
By \ref{egal}, we have $\Gamma(x_n)=1\!\surl{\ _{\beta}
\otimes_{\alpha}}_{\ N}y_n$. If $d \in ({\mathcal M}_{\Psi} \cap
{\mathcal M}_{T_R})^+$, then, for all $n\in\mathbb N$, we have $dx_n
\in {\mathcal M}_{\Psi} \cap {\mathcal M}_{T_R}$. Let $\omega \in
M^+_*$ and $k\in {\mathbb R}^+$ such that $\omega\circ\alpha\leq
k\nu$. By right invariance, we get:

$$
\begin{aligned}
\omega\circ T_R(dx_n)&=\omega((\Psi \surl{\ _{\beta}
\star_{\alpha}}_{\ \nu} id)(\Gamma(dx_n)))\\
&=\Psi((id \surl{\ _{\beta} \star_{\alpha}}_{\ \nu}
\omega)(\Gamma(dx_n)))=\Psi((id \surl{\ _{\beta}
\star_{\alpha}}_{\ \nu} (y_n\omega))(\Gamma(d)))\\
&=\omega((\Psi \surl{\ _{\beta} \star_{\alpha}}_{\ \nu}
id)(\Gamma(d))y_n)=\omega(T_R(d)y_n)
\end{aligned}$$

Take the limit over $n\in\mathbb{N}$ to obtain $T_R(dx)=T_R(d)x$ for
all $d\in {\mathcal M}_{\Psi} \cap {\mathcal M}_{T_R}$ and, by
semi-finiteness of $T_R$, we conclude that $x$ belongs to
$\beta(N)$. Reverse inclusion comes from axioms. If we apply this
result to the opposite adapted measured quantum groupoid, then we
get the second point.
\end{proof}

\subsubsection{The antipode and its polar decomposition}
We now approach definition of the antipode.

\begin{lemm}\label{ense}
We have $(\omega_{v,w} * id)(U'_{H_{\Phi}})D^{\lambda} \subset
D^{\lambda}(\omega_{\Delta_{\Psi}^{-\lambda}v,\Delta_{\Psi}^{\lambda}w}
* id)(U'_{H_{\Phi}})$ for all $\lambda \in \mathbb C$ and $v,w\in
\Lambda_{\Psi}({\mathcal T}_{\Psi,T_R})$.
\end{lemm}

\begin{proof}
Straightforward from relation \eqref{fond}.
\end{proof}

\begin{prop}\label{proi}
If $I$ is the unitary part of the polar decomposition of $G$, then,
for all $v,w\in D((H_{\Psi})_{\beta},\nu^o)$, we have:
$$I(\omega_{J_{\Psi}w,v}*id)(U_{H_{\Phi}}'{}^{\hspace{-.3cm}*}\hspace{.1cm})I
=(\omega_{J_{\Psi}v,w}*id)(U'_{H_{\Phi}})$$
\end{prop}

\begin{proof}
We have $(\omega_{v,w} * id)(U'_{H_{\Phi}})D^{1/2}\subseteq
D^{1/2}(\omega_{\Delta_{\Psi}^{-1/2}v,\Delta_{\Psi}^{1/2}w}\star
id)(U'_{H_{\Phi}})$ for all $v,w\in\Lambda_{\Psi}({\mathcal
T}_{\Psi,T_R})$ by the previous lemma. On the other hand, by
inclusion \eqref{inclu2}, we have:
$$
(\omega_{v,w} * id)(U'_{H_{\Phi}})D^{1/2}= (\omega_{v,w} \star
id)(U'_{H_{\Phi}})G^*I\subseteq D^{1/2}I(\omega_{w,v} *
id)(U'_{H_{\Phi}})I$$ So $I(\omega_{w,v}*id)(U'_{H_{\Phi}})I=
(\omega_{\Delta_{\Psi}^{-1/2}v,\Delta_{\Psi}^{1/2}w} *
id)(U'_{H_{\Phi}})$ and, by \ref{switch}, we have:
$$I(\omega_{w,v}*id)(U_{H_{\Phi}}'{}^{\hspace{-.3cm}*}\hspace{.1cm})I
=(\omega_{\Delta_{\Psi}^{1/2}w,\Delta_{\Psi}^{-1/2}v} * id)(U_{H_{\Phi}}'{}^{\hspace{-.3cm}*}\,) \\
=(\omega_{J_{\Psi}v,J_{\Psi}w}* id)(U'_{H_{\Phi}})$$
\end{proof}

\begin{coro}
There exists a *-anti-automorphism $R$ of $M$ defined by
$R(m)=Im^*I$ such that $R^2=id$. (We recall that $I$ denotes the
unitary part of the polar decomposition of $G$).
\end{coro}

\begin{proof}
Straightforward from the previous proposition and theorem
\ref{densevn}.
\end{proof}

\begin{defi}
The unique *-anti-automorphism $R$ of $M$ such that $R(m)=Im^*I$,
where $I$ denotes the unitary part of the polar decomposition of
$G$, is called \textbf{unitary antipode}.
\end{defi}

\begin{defi}\label{antipode}
The application $S=R\tau_{-i/2}$ is called {\bf antipode}.
\end{defi}

The next proposition states elementary properties of the antipode.
Straightforward proofs are omitted.

\begin{prop}\label{numero}The antipode $S$ satisfies:
\begin{center}
\begin{minipage}{11cm}
\begin{enumerate}[i)]
\item for all $t\in\mathbb R$, we have $\tau_t\circ R=R\circ\tau_t$ and $\tau_t\circ S=S\circ\tau_t$
\item $SR=RS$ and $S^2=\tau_{-i}$
\item $S$ is densely defined and has dense range
\item $S$ is injective and $S^{-1}=R\tau_{i/2}=\tau_{i/2}R$
\item for all $x\in {\mathcal D}(S)$, $S(x^*)\in {\mathcal D}(S)$
and $S(S(x)^*)^*=x$ \label{test}
\end{enumerate}
\end{minipage}
\end{center}
\end{prop}

\subsubsection{Characterization of the antipode}

In \ref{antipode}, we define the antipode by giving its polar
decomposition. However, we have to verify that $S$ is what it should
be.

\paragraph{Usual characterization of the antipode.}

\begin{prop}\label{espoir}
For all $v,w\in \Lambda_{\Psi}({\mathcal T}_{\Psi,T_R})$,
$(\omega_{w,v} * id)(U'_{H_{\Phi}})$ belongs to ${\mathcal D}(S)$
and we have: $$S((\omega_{w,v}*id)(U'_{H_{\Phi}}))=(\omega_{w,v} *
id)(U_{H_{\Phi}}'{}^{\hspace{-.3cm}*}\hspace{.1cm})$$ Moreover, the
linear span of $(\omega_{v,w} * id)(U'_{H_{\Phi}})$, where $v,w\in
\Lambda_{\Psi}({\mathcal T}_{\Psi,T_R})$, is a core for $S$.
\end{prop}

\begin{proof}
By \ref{ense}, we have $(\omega_{w,v}\star id)(U'_{H_{\Phi}}) \in
{\mathcal D}(\tau_{-i/2})={\mathcal D}(S)$ and:
$$
\begin{aligned}
S((\omega_{w,v} * id)(U'_{H_{\Phi}})) &=
R((\omega_{\Delta_{\Psi}^{-1/2}w,\Delta_{\Psi}^{1/2}v}*id)(U'_{H_{\Phi}})) \\
&= (\omega_{S_{\Psi}v,\Delta_{\Psi}S_{\Psi}w}*
id)(U'_{H_{\Phi}}) &&\text{by proposition \ref{proi},} \\
&= (\omega_{w,v} * id)(U_{H_{\Phi}}'{}^{\hspace{-.3cm}*}\,)
&&\text{by lemma \ref{switch}}.
\end{aligned}$$

The involved subspace of $M$ is included in ${\mathcal
D}(\tau_{-i/2})$ by \ref{ense}, weakly dense in $M$ by theorem
\ref{densevn} and $\tau$-invariant by \ref{tau} which finishes the
proof.
\end{proof}

\begin{coro}\label{gammas}
For $a,b,c,d\in {\mathcal T}_{\Psi,T_R}$,
$(\omega_{\Lambda_{\Psi}(a),\Lambda_{\Psi}(b)} \surl{\ _{\beta}
\star_{\alpha}}_{\ \nu} id)(\Gamma(cd))$ belongs to ${\mathcal
D}(S)$ and we have:
$$S((\omega_{\Lambda_{\Psi}(a),\Lambda_{\Psi}(b)} \surl{\ _{\beta}
\star_{\alpha}}_{\
\nu}id)(\Gamma(cd)))=(\omega_{\Lambda_{\Psi}(c),\Lambda_{\Psi}(\sigma_{-i}^{\Psi}(d^*))}
\surl{\ _{\beta} \star_{\alpha}}_{\ \nu}
id)(\Gamma(\sigma_{i}^{\Psi}(a)b^*))$$
\end{coro}

\begin{proof}
By \ref{corres}, we know that:
$$(\omega_{\Lambda_{\Psi}(a),\Lambda_{\Psi}(b)}
\surl{\ _{\beta} \star_{\alpha}}_{\ \nu}
id)(\Gamma(cd))=(\omega_{\Lambda_{\Psi}(cd),\Lambda_{\Psi}
(b\sigma_{-i}^{\Psi}(a^*))}*id)(U'_{H_{\Phi}})$$ which belongs to
${\mathcal D}(S)$. Then, by \ref{corres} and \ref{switch}, we have:
$$
\begin{aligned}
S((\omega_{\Lambda_{\Psi}(a),\Lambda_{\Psi}(b)} \surl{\ _{\beta}
\star_{\alpha}}_{\ \nu} id)(\Gamma(cd)))&=
S((\omega_{\Lambda_{\Psi}(cd),\Lambda_{\Psi}
(b\sigma_{-i}^{\Psi}(a^*))}*id)(W'))\\
&=(\omega_{\Lambda_{\Psi}(cd),\Lambda_{\Psi}
(b\sigma_{-i}^{\Psi}(a^*))}*id)(W'^*)\\
&=(\omega_{\Lambda_{\Psi}(\sigma_{i}^{\Psi}(a)b^*),\Lambda_{\Psi}
(\sigma_{-i}^{\Psi}(d^*c^*))}*id)(W')\\
&=(\omega_{\Lambda_{\Psi}(c),
\Lambda_{\Psi}(\sigma_{-i}^{\Psi}(d^*))} \surl{\ _{\beta}
\star_{\alpha}}_{\ \nu} id)(\Gamma(\sigma_{i}^{\Psi}(a)b^*))
\end{aligned}$$
\end{proof}

\paragraph{The co-involution $R$}

In this section, we give a new expression of $R$ and we show that it
is a co-involution of the adapted measured quantum groupoid.

\begin{prop}\label{defR}
For all $a,b\in {\mathcal N}_{\Psi}\cap {\mathcal N}_{T_R}$, we
have:
$$R((\omega_{J_{\Psi}\Lambda_{\Psi}(a)}\surl{\ _{\beta}\star_{\alpha}}_{\ \
\nu}id)(\Gamma(b^*b)))=(\omega_{J_{\Psi}\Lambda_{\Psi}(b)}\surl{\
_{\beta} \star_{\alpha}}_{\ \ \nu}id)(\Gamma(a^*a))$$
\end{prop}

\begin{proof}
The proposition comes from the following computation:
$$
\begin{aligned}
&\ \quad
R((\omega_{J_{\Psi}\Lambda_{\Psi}(a),J_{\Psi}\Lambda_{\Psi}(a)}\surl{\
_{\beta} \star_{\alpha}}_{\ \nu}id)(\Gamma(b^*b)))&\\
&=R((\omega_{\Lambda_{\Psi}(b^*b),J_{\Psi}\Lambda_{\Psi}(a^*a)}*id)(U'_{H_{\Phi}})
&\text{by corollary \ref{corres},}\\
&=(\omega_{\Lambda_{\Psi}(a^*a),J_{\Psi}\Lambda_{\Psi}(b^*b)}*id)(U'_{H_{\Phi}})
&\text{by definition of } R,\\
&=(\omega_{J_{\Psi}\Lambda_{\Psi}(b),J_{\Psi}\Lambda_{\Psi}(b)}\surl{\
_{\beta} \star_{\alpha}}_{\ \nu}id)(\Gamma(a^*a))
&\text{by corollary \ref{corres}.}\\
\end{aligned}$$

\end{proof}

\begin{rema}
We notice that $R$ is $T_L$-independent.
\end{rema}

\begin{prop}
We have $I\alpha(n^*)=\beta(n)I$ for all $n\in N$ and
$R\circ\alpha=\beta$.
\end{prop}

\begin{proof}
By \ref{commuNa}, we have, for all $x \in {\mathcal T}_{\Psi,T_R}$:
$$\beta(x)GD^{-1/2}\subseteq G\alpha(\sigma_{-i/2}((x^*))
\subseteq GD^{-1/2}\alpha(x^*)\subseteq I\alpha(x^*)$$ and, on the
other hand, $\beta(x)GD^{-1/2}\subseteq \beta(x)I$ so that
$I\alpha(x^*)=\beta(x)I$. The result holds by normality of $\alpha$
and $\beta$.
\end{proof}

By \cite{S2}, there exists a unitary and anti-linear operator $I
\surl{\ _{\beta} \otimes_{\alpha}}_{\ \nu}I$ from $H\!\!\! \surl{\
_{\beta} \otimes_{\alpha}}_{\ \nu}\!\! H$ onto $H\!\!\!\surl{\
_{\alpha} \otimes_{\beta}}_{\ \ \nu^o}\!\! H$, the adjoint of which
is $I\surl{\ _{\alpha} \otimes_{\beta}}_{\ \ \nu^o}I$. Also, there
exists an anti-isomorphism $R\surl{\ _{\beta} \star_{\alpha}}_{\ N}
R$ from $M\surl{\ _{\beta} \star_{\alpha}}_{\ N}M$ onto $M\surl{\
_{\alpha} \star_{\beta}}_{\ N^o}M$ and, by definition of $R$, we
have, for all $X\in M\surl{\ _{\beta} \star_{\alpha}}_{\ N}M$:
$$(R\surl{\ _{\beta} \star_{\alpha}}_{\ N} R)(X)=(I \surl{\
_{\beta} \otimes_{\alpha}}_{\ \nu}I)X^*(I\surl{\ _{\alpha}
\otimes_{\beta}}_{\ \ \nu^o} I)$$ We underline the fact that, if
$\omega\in M_*^+$, then $\omega\circ R \in M_*^+$ and, if there
exists $k\in\mathbb{R}^+$ such that $\omega\circ\alpha\leq k\nu$,
then $\omega\circ R\circ\beta\leq k\nu$. Also, if $\theta\in M_*^+$
and $k'\in{\mathbb R}^+$ are such that $\theta\circ\beta\leq k'\nu$,
then $\theta\circ R\circ\alpha\leq k\nu$. Then, we have $\omega
R\surl{\ _{\beta} \star_{\alpha}}_{\ \nu} \theta R=(\omega \surl{\
_{\alpha} \star_{\beta}}_{\ \ \nu^o}\theta)\circ(R\surl{\ _{\beta}
\star_{\alpha}}_{\ \nu} R)$.

\begin{lemm}
For all $a,x\in {\mathcal N}_{T_R}\cap {\mathcal N}_{\Psi}$,
$\omega\in M^+_*$ and $k\in\mathbb{R}^+$ such that
$\omega\circ\alpha\leq k\nu$, we have: $$\omega\circ
R((\omega_{J_{\Psi}\Lambda_{\Psi}(a)}\surl{\ _{\beta}
\star_{\alpha}}_{\ \nu}id)(\Gamma(x)))=(\Lambda_{\Psi}((id\surl{\
_{\beta} \star_{\alpha}}_{\
\nu}\omega)(\Gamma(a^*a)))|J_{\Psi}\Lambda_{\Psi}(x))$$
\end{lemm}

\begin{proof}
Let $b\in {\mathcal N}_{T_R}\cap {\mathcal N}_{\Psi}$. By
\ref{defR}, we can compute:
$$\begin{aligned}
\omega\circ R((\omega_{J_{\Psi}\Lambda_{\Psi}(a)}\surl{\ _{\beta}
\star_{\alpha}}_{\
\nu}id)(\Gamma(b^*b)))&=\omega((\omega_{J_{\Psi}\Lambda_{\Psi}(b)}\surl{\
_{\beta} \star_{\alpha}}_{\ \nu}id)(\Gamma(a^*a)))\\
&=((id\surl{\ _{\beta} \star_{\alpha}}_{\
\nu}\omega)(\Gamma(a^*a))J_{\Psi}\Lambda_{\Psi}(b)|J_{\Psi}\Lambda_{\Psi}(b))\\
&=(J_{\Psi}bJ_{\Psi}\Lambda_{\Psi}((id\!\!\surl{\ _{\beta}
\star_{\alpha}}_{\ \nu}\omega)(\Gamma(a^*a)))|J_{\Psi}\Lambda_{\Psi}(b))\\
&=(\Lambda_{\Psi}((id\surl{\ _{\beta} \star_{\alpha}}_{\
\nu}\omega)(\Gamma(a^*a)))|J_{\Psi}\Lambda_{\Psi}(b^*b))
\end{aligned}$$

Linearity and normality of the expressions imply the lemma.
\end{proof}

\begin{prop}
We have $\varsigma_{N^o}\circ(R\surl{\ _{\beta} \star_{\alpha}}_{\
N} R)\circ\Gamma=\Gamma \circ R$.
\end{prop}

\begin{proof}
Let $a,b\in {\mathcal N}_{T_R}\cap {\mathcal N}_{\Psi}$,
$\omega,\theta\in M^+_*$ and $k,k'\in\mathbb{R}^+$ such that
$\omega\circ\alpha\leq k\nu$ and $\theta\circ\beta\leq k'\nu$. Then,
we can compute by \ref{defR} and the previous lemma:
$$
\begin{aligned}
&\ \quad(\theta\surl{\ _{\beta} \star_{\alpha}}_{\
\nu}\omega)(\Gamma\circ R((\omega_{J_{\Psi}\Lambda_{\Psi}(a)}\surl{\
_{\beta}
\star_{\alpha}}_{\ \nu}id)(\Gamma(b^*b))))\\
&=(\theta\surl{\ _{\beta} \star_{\alpha}}_{\
\nu}\omega)(\Gamma((\omega_{J_{\Psi}\Lambda_{\Psi}(b)}\surl{\
_{\beta}\star_{\alpha}}_{\ \nu}id)(\Gamma(a^*a))))\\
&=(\omega_{J_{\Psi}\Lambda_{\Psi}(b)}\surl{\
_{\beta}\star_{\alpha}}_{\ \nu}\theta\surl{\ _{\beta}
\star_{\alpha}}_{\ \nu}\omega)(id\surl{\
_{\beta}\star_{\alpha}}_{\ N}\Gamma)(\Gamma(a^*a))\\
&=(\omega_{J_{\Psi}\Lambda_{\Psi}(b)}\surl{\
_{\beta}\star_{\alpha}}_{\ \nu}\theta\surl{\ _{\beta}
\star_{\alpha}}_{\ \nu}\omega)(\Gamma\surl{\
_{\beta}\star_{\alpha}}_{\ N}id)(\Gamma(a^*a))\\
&=(\omega_{J_{\Psi}\Lambda_{\Psi}(b)}\surl{\
_{\beta}\star_{\alpha}}_{\ \nu}\theta)[\Gamma((id\surl{\ _{\beta}
\star_{\alpha}}_{\ \nu}\omega)(\Gamma(a^*a)))]\\
&=(\Lambda_{\Psi}((id\surl{\ _{\beta} \star_{\alpha}}_{\
\nu}\theta\circ R)(\Gamma(b^*b)))|J_{\Psi}\Lambda_{\Psi}((id\surl{\
_{\beta} \star_{\alpha}}_{\ \nu}\omega)(\Gamma(a^*a))))
\end{aligned}$$
Observe the symmetry of the last expression and use it to proceed
towards the computation:
$$
\begin{aligned}
&\ \quad (\Lambda_{\Psi}((id\surl{\ _{\beta} \star_{\alpha}}_{\
\nu}\omega)(\Gamma(a^*a)))|J_{\Psi}\Lambda_{\Psi}((id\surl{\
_{\beta} \star_{\alpha}}_{\ \nu}\theta\circ
R)(\Gamma(b^*b))))\\
&=(\omega_{J_{\Psi}\Lambda_{\Psi}(a)}\surl{\
_{\beta}\star_{\alpha}}_{\ \nu}\omega\circ R)[\Gamma((id\surl{\
_{\beta} \star_{\alpha}}_{\ \nu}\theta\circ R)(\Gamma(b^*b)))]\\
&=(\omega_{J_{\Psi}\Lambda_{\Psi}(a)}\surl{\
_{\beta}\star_{\alpha}}_{\ \nu}\omega\circ R\surl{\ _{\beta}
\star_{\alpha}}_{\ \nu}\theta\circ R)(\Gamma\surl{\
_{\beta}\star_{\alpha}}_{\ N}id)(\Gamma(b^*b))\\
&=(\omega_{J_{\Psi}\Lambda_{\Psi}(a)}\surl{\
_{\beta}\star_{\alpha}}_{\ \nu}\omega\circ R\surl{\ _{\beta}
\star_{\alpha}}_{\ \nu}\theta\circ R)(id\surl{\
_{\beta}\star_{\alpha}}_{\ N}\Gamma)(\Gamma(b^*b))\\
&=(\omega\circ R\surl{\ _{\beta} \star_{\alpha}}_{\ \nu}\theta\circ
R)(\Gamma((\omega_{J_{\Psi}\Lambda_{\Psi}(a)}\surl{\ _{\beta}
\star_{\alpha}}_{\ \nu}id)(\Gamma(b^*b))))\\
&=(\omega\surl{\ _{\alpha} \star_{\beta}}_{\ \nu^o}\theta)(R\surl{\
_{\beta} \star_{\alpha}}_{\
N}R)(\Gamma((\omega_{J_{\Psi}\Lambda_{\Psi}(a)}\surl{\ _{\beta}
\star_{\alpha}}_{\ \nu}id)(\Gamma(b^*b))))\\
&=(\theta\surl{\ _{\beta} \star_{\alpha}}_{\
\nu}\omega)\varsigma_{N^o}(R\surl{\ _{\beta} \star_{\alpha}}_{\
N}R)(\Gamma((\omega_{J_{\Psi}\Lambda_{\Psi}(a)}\surl{\ _{\beta}
\star_{\alpha}}_{\ \nu}id)(\Gamma(b^*b))))
\end{aligned}$$

Theorem \ref{densevn} easily implies the result.
\end{proof}

\paragraph{Left strong invariance w.r.t the antipode.}

In this section, $T'$ denotes a left invariant n.s.f weight from $M$
to $\alpha(N)$. We put $\Phi'=\nu\circ\alpha^{-1}\circ T'$,
$J_{\Phi'}$ the anti-linear operator and $\Delta_{\Phi'}$ the
modular operator which come from Tomita's theory of $\Phi'$,
$\sigma^{\Phi'}$ its modular group and $V=(U_{T'})_{H_{\Phi}}^*$ i.e
the fundamental unitary associated with $T'$. The next proposition
is the left strong invariance w.r.t $S$.

\begin{prop}\label{invafort}
Elements $(id*\omega_{v,w})(V)$ belong to the domain of $S$ for all
$v,w\in\Lambda_{\Phi'}({\mathcal T}_{\Phi',T'})$ and we have
$S((id*\omega_{v,w})(V))=(id*\omega_{v,w})(V^*)$.
\end{prop}

\begin{proof}
By \ref{corres}, we have $(id*\omega)(V)=(\omega\circ
R*id)(U'_{H_{\Phi}})$ for all $\omega$. If
$\overline{\omega}(x)=\overline{\omega(x^*)}$, then, by
\ref{espoir}, we have:
$$\begin{aligned}
S((id*\omega)(V))=S((\omega\circ R*id)(U'_{H_{\Phi}}))
&=(\omega\circ R*id)(U_{H_{\Phi}}'{}^{\hspace{-.3cm}*}\,)\\
&=[(\overline{\omega}\circ R*id)(U'_{H_{\Phi}})]^*\\
&=[(id*\overline{\omega})(V)]^*=(id*\omega)(V^*)
\end{aligned}$$

\end{proof}

\begin{lemm}
For all $v\in {\mathcal D}(D^{1/2})$ and $w\in {\mathcal
D}(D^{1/2})$, we have:
$$(\omega_{v,w}*id)(V)^*=(\omega_{ID^{-1/2}v,ID^{1/2}w}*id)(V)$$
\end{lemm}

\begin{proof}
We have $(id*\omega_{w',v'})(V)\in {\mathcal D}(S)={\mathcal
D}(\tau_{-i/2})$ for all $v',w'$ belonging to
$\Lambda_{\Phi'}({\mathcal T}_{\Phi',T'})$ by \ref{invafort} and,
since $\tau$ is implemented by $D^{-1}$, we have:

$$
\begin{aligned}
(id*\omega_{w',v'})(V)D^{1/2} &\subseteq
D^{1/2}\tau_{-i/2}((id*\omega_{w',v'})(V)) \\
&=D^{1/2}R(S((id*\omega_{w',v'})(V))) \\
&=D^{1/2}I[(id*\omega_{w',v'})(V^*)]^*I \\
&=D^{1/2}I(id*\omega_{v',w'})(V)I.
\end{aligned}$$

Then, for all $v\in {\mathcal D}(D^{1/2})$ and $w\in {\mathcal
D}(D^{1/2})$, we have:

$$
\begin{aligned}
((\omega_{ID^{-1/2}v,ID^{1/2}w} * id)(V)w'|v') &=
((id*\omega_{w',v'})(V)
D^{1/2}Iv|D^{-1/2}Iw) \\
&=(D^{1/2}I(id*\omega_{v',w'})(V)v|D^{-1/2}Iw) \\
&=(w|(id*\omega_{v',w'})v) \\
&=((\omega_{v,w} * id)(V)^*w',v')
\end{aligned}$$

Then, the proposition holds.
\end{proof}

\begin{prop}\label{clef}
The following relations are satisfied:
\begin{center}
\begin{minipage}{10cm}
\begin{enumerate}[i)]
\item $(I \surl{\ _{\alpha} \otimes_{\epsilon}}_{\ \ N^o} J_{\Phi'})V
=V^*(I \surl{\ _{\beta} \otimes_{\alpha}}_{\ N} J_{\Phi'})$;
\item $(D^{-1}\surl{\ _{\alpha} \otimes_{\epsilon}}_{\ \ \nu^o} \Delta_{\Phi'})V
=V(D^{-1}\surl{\ _{\beta} \otimes_{\alpha}}_{\ \nu}
\Delta_{\Phi'})$;
\item $(\tau_t \surl{\ _{\beta} \star_{\alpha}}_{\ N}
\sigma^{\Phi'}_t)\circ \Gamma=\Gamma\circ\sigma^{\Phi'}_t$ for all
$t\in\mathbb R$ .
\end{enumerate}
\end{minipage}
\end{center}
where $\epsilon(n)=J_{\Phi'}\alpha(n^*)J_{\Phi'}$ for all $n\in N$.
\end{prop}

\begin{proof}
We denote by $S_{\Phi'}$ the operator of Tomita's theory associated
with $\Phi'$ and defined as the closed operator on $H_{\Phi'}$ such
that $\Lambda_{\Phi'}({\mathcal N}_{\Phi'} \cap {\mathcal
N}_{\Phi'}^*)$ is a core for $S_{\Phi'}$ and
$S_{\Phi'}\Lambda_{\Phi'}(x)=\Lambda_{\Phi'}(x^*)$ for all $x\in
{\mathcal N}_{\Phi'} \cap {\mathcal N}_{\Phi'}^*$. Then, by
definition, we have $\Delta_{\Phi'}=S_{\Phi'}^*S_{\Phi'}$ and
$S_{\Phi'}=J_{\Phi'}\Delta_{\Phi'}^{1/2}$. Moreover, for all $m\in
M$ and $t\in\mathbb R$, we have
$\sigma^{\Phi'}_t(m)=\Delta_{\Phi'}^{it}m\Delta_{\Phi'}^{-it}$.

First of all, we verify these relations make sense. We have to prove
some commutation relations. We can write for all $n\in {\mathcal
T}_{\nu}$ and $y\in {\mathcal N}_{\Phi'} \cap {\mathcal
N}_{\Phi'}^*$:
$$S_{\Phi'}\alpha(n)\Lambda_{\Phi'}(y)=S_{\Phi'}\Lambda_{\Phi'}(\alpha(n)y)$$
$$=\Lambda_{\Phi'}(y^*\alpha(n^*))=
\hat{\alpha}(\sigma_{-i/2}^{\nu}(n))S_{\Phi'}\Lambda_{\Phi'}(y)$$ so
$\hat{\alpha}(\sigma_{-i/2}^{\nu}(n))S_{\Phi'} \subseteq
S_{\Phi'}\alpha(n)$ and by adjoint $\alpha(n)S_{\Phi'}^* \subseteq
S_{\Phi'}^*\hat{\alpha}(\sigma_{i/2}^{\nu}(n)$. Then:
$$\alpha(n)\Delta_{\Phi'}=\alpha(n)S_{\Phi'}^*S_{\Phi'} \subseteq
S_{\Phi'}^*\hat{\alpha}(\sigma_{i/2}^{\nu}(n)S_{\Phi'} \subseteq
\Delta_{\Phi'}\alpha(\sigma_i^{\nu}(n))$$ Since $\beta(n)D^{-1}
\subseteq D^{-1}\beta(\sigma_i^{\nu}(n))$, the second relation makes
sense. On an other hand, we know that $I\beta(n)=\alpha(n^*)I$ and
${\mathcal J}\alpha(n)=\epsilon(n^*)J_{\Phi'}$ to terms of the first
relation. Finally, for all $t\in\mathbb R$, we have:
$$\tau_t \circ\beta=\beta\circ\sigma_t^{\nu} \quad\text{ and }
\sigma^{\Phi'}_t(\alpha(n))=\Delta_{\Phi'}^{it}\alpha(n)\Delta_{\Phi'}^{-it}=
\alpha(\sigma_t^{\nu}(n))$$ which finishes verifications.

Let $v,w \in \Lambda_{\Phi}({\mathcal T}_{\Phi,S_L})$. By
\ref{prem}, we know that $(\omega_{v,w} \surl{\ _{\beta}
\star_{\alpha}}_{\ \nu} id)(\Gamma(y))$ belongs to ${\mathcal
N}_{T'}\cap {\mathcal N}_{\Phi'} \cap {\mathcal N}_{T'}^*\cap
{\mathcal N}_{\Phi'}^*$ for all $y \in {\mathcal N}_{T'}\cap
{\mathcal N}_{\Phi'} \cap {\mathcal N}_{T'}^*\cap {\mathcal
N}_{\Phi'}^*$. By \ref{lienGV}, we can write $(\omega_{v,w} \star
id)(V^*)\Lambda_{\Phi'}(y)=\Lambda_{\Phi'}((\omega_{v,w} \surl{\
_{\beta} \star_{\alpha}}_{\ \nu} id)(\Gamma(y)))$ so that
$(\omega_{v,w} * id)(V^*)\Lambda_{\Phi'}(y)$ belongs to ${\mathcal
D}(S_{\Phi'})$. Then, we compute:
$$
\begin{aligned}
S_{\Phi'}(\omega_{v,w}*id)(V^*)\Lambda_{\Phi'}(y)&=
S_{\Phi'}\Lambda_{\Phi'}((\omega_{v,w} \surl{\ _{\beta}
\star_{\alpha}}_{\ \nu}
id)(\Gamma(y)))\\
&=\Lambda_{\Phi'}((\omega_{w,v} \surl{\ _{\beta} \star_{\alpha}}_{\
\nu}
id)(\Gamma(y^*)))\\
&=(\omega_{w,v}*id)(V^*)\Lambda_{\Phi'}(y^*)\\
&=(\omega_{w,v}*id)(V^*)S_{\Phi'} \Lambda_{\Phi'}(y)
\end{aligned}$$ Since $\Lambda_{\Phi'}({\mathcal
N}_{T'}\cap {\mathcal N}_{\Phi'} \cap {\mathcal N}_{T'}^*\cap
{\mathcal N}_{\Phi'}^*)$ is a core for $S_{\Phi'}$, this implies:
\begin{equation}\label{parallele1}
(\omega_{w,v} * id)(V^*)S_{\Phi'} \subseteq S_{\Phi'}(\omega_{v,w}
* id)(V^*)
\end{equation}
Take adjoint so as to get:
\begin{equation}\label{parallele2}
(\omega_{w,v} * id)(V)S_{\Phi'}^* \subseteq S_{\Phi'}^*(\omega_{v,w}
* id)(V)
\end{equation}
Then, we deduce by the previous lemma:

$$
\begin{aligned}
(\omega_{v,w}* id)(V)\Delta_{\Phi'} &= (\omega_{v,w}*
id)(V)S_{\Phi'}^*S_{\Phi'}
\\
&\subseteq S_{\Phi'}^*(\omega_{v,w} * id)(V)S_{\Phi'} \\
&=S_{\Phi'}^*[(\omega_{ID^{-1/2}w,ID^{1/2}v} * id)(V)]^*S_{\Phi'}
\end{aligned}$$\\
Then by inclusion \eqref{parallele1} and the previous lemma, we
have:

$$
\begin{aligned}
(\omega_{v,w}* id)(V)\Delta_{\Phi'} &\subseteq
S_{\Phi'}^*S_{\Phi'}[(\omega_{ID^{1/2}v,ID^{-1/2}w} * id)(V)]^*\\
&=\Delta_{\Phi'} (\omega_{D^{1/2}IID^{1/2}v,D^{-1/2}IID^{-1/2}w}*id)(V)\\
&=\Delta_{\Phi'} (\omega_{Dv,N^{-1}w} * id)(V)
\end{aligned}$$
Consequently, like relation \eqref{fond}, we easily deduce that:
$$(D^{-1}\surl{\ _{\alpha} \otimes_{\epsilon}}_{\ \ \nu^o}
\Delta_{\Phi'})V =V(D^{-1}\surl{\ _{\beta} \otimes_{\alpha}}_{\ \nu}
\Delta_{\Phi'})$$

Let's prove the first relation. By inclusion \eqref{parallele1}, for
all $v\in {\mathcal D}(N^{-1/2})$ and $w\in {\mathcal D}(D^{1/2})$,
we have:

\begin{equation}\label{comenc}
\begin{aligned}
J_{\Phi'}(\omega_{w,v}* id)(V^*)J_{\Phi'}\Delta_{\Phi'}^{1/2}
&=J_{\Phi'}(\omega_{w,v}* id)(V^*)S_{\Phi'} \\
&\subseteq J_{\Phi'}S_{\Phi'}(\omega_{v,w}* id)(V^*)\\
&=\Delta_{\Phi'}^{1/2}(\omega_{v,w}* id)(V^*)
\end{aligned}
\end{equation}

For all $p,q\in {\mathcal D}(\Delta_{\Phi'}^{1/2})$, we have by ii):

$$
\begin{aligned}
((\omega_{v,w}* id)(V^*)p,\Delta_{\Phi'}^{1/2}q) &= (V^*(v\surl{\
_{\alpha} \otimes_{\epsilon}}_{\ \ \nu^o} p)|w\surl{\ _{\beta}
\otimes_{\alpha}}_{\ \nu} \Delta_{\Phi'}^{1/2}q) \\
&=(V^*(v \surl{\ _{\alpha} \otimes_{\epsilon}}_{\ \
\nu^o}p)|D^{-1/2}(D^{1/2}w) \surl{\
_{\beta} \otimes_{\alpha}}_{\ \nu} \Delta_{\Phi'}^{1/2}q) \\
&=((D^{-1/2} \surl{\ _{\beta} \otimes_{\alpha}}_{\ \nu}
\Delta_{\Phi'}^{1/2})V^*(v\surl{\ _{\alpha} \otimes_{\epsilon}}_{\ \
\nu^o} p)|D^{1/2}w\surl{\ _{\beta}
\otimes_{\alpha}}_{\ \nu} q) \\
&=(V^*(D^{-1/2}v\surl{\ _{\alpha} \otimes_{\epsilon}}_{\ \ \nu^o}
\Delta_{\Phi'}^{1/2}p)|D^{1/2}w\surl{\ _{\beta} \otimes_{\alpha}}_{\
\nu}
q) \\
&=((\omega_{D^{-1/2}v,D^{1/2}w} * id)(V^*)\Delta_{\Phi'}^{1/2}p|q).
\end{aligned}$$

Since $\Delta_{\Phi'}^{1/2}$ is self-adjoint, we get:
$$(\omega_{D^{-1/2}v,D^{1/2}w} * id)(V^*)\Delta_{\Phi'}^{1/2} \subseteq
\Delta_{\Phi'}^{1/2}(\omega_{v,w}* id)(V^*)$$ Also, by the previous
lemma, we have:
$$
\begin{aligned}
(\omega_{D^{-1/2}v,D^{1/2}w} *
id)(V^*)&=(\omega_{D^{1/2}w,D^{-1/2}v} * id) (V)^*\\
&=(\omega_{Iw,Iv}* id)(V)
\end{aligned}$$ That's why
$(\omega_{Iw,Iv}* id)(V)\Delta_{\Phi'}^{1/2}\subseteq
\Delta_{\Phi'}^{1/2}(\omega_{v,w}* id)(V^*)$. Since
$\Delta_{\Phi'}^{1/2}$ has dense range, this last inclusion and
\eqref{comenc} imply that:
$$(\omega_{Iw,Iv}*id)(V)=J_{\Phi'}(\omega_{v,w}*
id)(V^*)J_{\Phi'}$$ Then, we can compute:
$$
\begin{aligned}
&\quad((I\surl{\ _{\beta} \otimes_{\alpha}}_{\ \nu}J_{\Phi'})V^*(I
\surl{\ _{\beta} \otimes_{\alpha}}_{\ \nu}J_{\Phi'})(v\surl{\
_{\beta} \otimes_{\alpha}}_{\ \nu} q)|w\surl{\ _{\alpha}
\otimes_{\epsilon}}_{\ \ \nu^o} q)\\
&=(V(Iw \surl{\ _{\beta} \otimes_{\alpha}}_{\
\nu}J_{\Phi'}q)|Iv\surl{\ _{\alpha} \otimes_{\epsilon}}_{\ \
\nu^o}J_{\Phi'}p)\\
&=((\omega_{Iw,Iv}\star id)(V)J_{\Phi'}q|J_{\Phi'}p)
=(J_{\Phi'}(\omega_{w,v}\star id)(V^*)q|J_{\Phi'}p)\\
&=((\omega_{v,w}* id)(V)p|q)=(V(v\surl{\ _{\beta}
\otimes_{\alpha}}_{\ \nu} q)|w\surl{\ _{\alpha}
\otimes_{\epsilon}}_{\ \ \nu^o} q)
\end{aligned}$$
so that the first relation is proved. We end the proof by the last
equality. We know that $\Gamma$ is implemented by $V$,
$\sigma^{\Phi'}$ by $\Delta_{\Phi'}$ and $\tau$  by $D$ so that the
relation comes from $(D^{-1}\surl{\ _{\alpha} \otimes_{\epsilon}}_{\
\nu} \Delta_{\Phi'})V =V(D^{-1}\surl{\ _{\beta} \otimes_{\alpha}}_{\
\nu} \Delta_{\Phi'})$ like \ref{egal}.
\end{proof}

If we take $T'=T_L$ then $V=W^*$, $J_{\Phi'}=J_{\Phi}$ and
$\Delta_{\Phi'}=\Delta_{\Phi}$ so that we have the following
propositions:

\begin{prop}
For all $v,w\in \Lambda_{\Phi}({\mathcal T}_{\Phi,S_L})$,
$(id*\omega_{v,w})(W)$ belongs to ${\mathcal D}(S)$ and:
$$S((id*\omega_{v,w})(W))=(id*\omega_{v,w})(W^*)$$
\end{prop}

\begin{prop}
We have $(\omega_{v,w}*id)(W^*)^*=(\omega_{ID^{-1/2}v,ID^{1/2}w}
*id)(W^*)$ for all $v\in {\mathcal D}(D^{1/2})$ and $w\in {\mathcal
D}(D^{1/2})$.
\end{prop}

\begin{prop}\label{besoin}
The following relations are satisfied:
\begin{center}
\begin{minipage}{10cm}
\begin{enumerate}[i)]
\item $(I \surl{\ _{\alpha} \otimes_{\hat{\beta}}}_{\ \ N^o} J_{\Phi})W^*
=W(I \surl{\ _{\beta} \otimes_{\alpha}}_{\ N} J_{\Phi})$ ;
\item $(D^{-1}\surl{\ _{\beta} \otimes_{\alpha}}_{\ \nu}
\Delta_{\Phi})W^* =W^*(D^{-1}\surl{\ _{\beta} \otimes_{\alpha}}_{\
\nu} \Delta_{\Phi})$ ;
\item $(\tau_t \surl{\ _{\beta} \star_{\alpha}}_{\ N}
\sigma^{\Phi}_t)\circ \Gamma=\Gamma\circ\sigma^{\Phi}_t$ for all
$t\in\mathbb R$.
\end{enumerate}
\end{minipage}
\end{center}
\end{prop}

We summarize the results of this section in the three following
theorems:

\begin{theo}
Let $(N,M,\alpha,\beta,\Gamma,\nu,T_L,T_R)$ be a adapted measured
quantum groupoid and $W$ the pseudo-multiplicative unitary
associated with. Then the closed linear span of
$(id*\omega_{v,w})(W)$ for all $v\in D(_{\alpha}H_{\Phi},\nu)$ and
$w\in D((H_{\Phi})_{\hat{\beta}},\nu^o)$ is equal to the whole von
Neumann algebra $M$.
\end{theo}

\begin{theo}\label{invforte}
Let $(N,M,\alpha,\beta,\Gamma,\nu,T_L,T_R)$ be a adapted measured
quantum groupoid and $W$ the pseudo-multiplicative associated with.
If we put $\Phi=\nu\circ\alpha^{-1}\circ T_L$, then there exists an
unbounded antipode $S$ which satisfies:
\begin{enumerate}[i)]
\item for all $x\in {\mathcal D}(S)$, $S(x)^*\in {\mathcal D}(S)$
and $S(S(x)^*)^*=x$
\item for all $v,w\in \Lambda_{\Phi}({\mathcal T}_{\Phi,S_L})$, $(id*\omega_{v,w})(W)$
belongs to ${\mathcal D}(S)$ and:
$$S((id*\omega_{v,w})(W))=(id*\omega_{v,w})(W^*)$$
\end{enumerate}
$S$ has the following polar decomposition $S=R\tau_{i/2}$, where $R$
is a co-involution of $M$ satisfying $R^2=id$, $R\circ\alpha=\beta$
and $\varsigma_{N^o}\circ(R\surl{\ _{\beta} \star_{\alpha}}_{\
N}R)\circ\Gamma=\Gamma\circ R$, and where $\tau$, the so-called
scaling group, is a one-parameter group of automorphisms such that
$\tau_t\circ\alpha=\alpha\circ\sigma_t^{\nu}$,
$\tau_t\circ\beta=\beta\circ\sigma_t^{\nu}$ satisfying
$\Gamma\circ\tau_t=(\tau_t\surl{\ _{\beta} \star_{\alpha}}_{\
N}\tau_t)\circ\Gamma$ for all $t\in\mathbb{R}$. $S,R$ and $\tau$ are
independent of $T_L$ and of $T_R$.

Moreover, $R\circ T_L\circ R$ is a n.s.f operator-valued weight
which is right invariant and $\alpha$-adapted w.r.t $\nu$.
\end{theo}

\begin{theo}\label{ex1}
Let $(N,M,\alpha,\beta,\Gamma,\nu,T_L,T_R)$ be a adapted measured
quantum groupoid. If $R$ is the co-involution and $\tau$ the scaling
group, then $(N,M,\alpha,\beta,\Gamma,T_L,R,\tau,\nu)$ becomes a
measured quantum groupoid.
\end{theo}

\begin{proof}
By hypothesis, we know that $\gamma_t=\sigma_{-t}^{\nu}$ for all
$t\in\mathbb R$ so that $\gamma$ leaves $\nu$ invariant. By theorem
\ref{invforte} and proposition \ref{defR}, we can construct a
co-involution $R$ and a scaling group $\tau$ such that
$(N,M,\alpha,\beta,\Gamma,T_L,R,\tau,\nu)$ becomes a measured
quantum groupoid.
\end{proof}

\subsection{Uniqueness, modulus and scaling operator}

By the general theory of measured quantum groupoids, theorems
\ref{ensemble1} and \ref{ensemble2} can be applied and we get the
following two theorems in the adapted measured quantum groupoids
case:

\begin{theo}
Let $(N,M,\alpha,\beta,\Gamma,\nu,T_L,T_R)$ be a adapted measured
quantum groupoid. If $T'$ is a left invariant operator-valued weight
which is $\beta$-adapted w.r.t $\nu$, then there exists a strictly
positive operator $h$ affiliated with $Z(N)$ such that, for all
$t\in\mathbb{R}$:
$$\nu\circ\alpha^{-1}\circ T'=(\nu\circ\alpha^{-1}\circ
T_L)_{\beta(h)}$$ We have a similar result for the right invariant
operator-valued weights.
\end{theo}

\begin{theo}
Let $(N,M,\alpha,\beta,\Gamma,\nu,T_L,R\circ T_L\circ R)$ be a
adapted measured quantum groupoid. Then there exists a strictly
positive operator $\delta$ affiliated with $M$ called modulus and
then there exists a strictly positive operator $\lambda$ affiliated
with $Z(M)\cap\alpha(N)\cap\beta(N)$ called scaling operator such
that $[D\nu\circ\alpha^{-1}\circ T_L\circ
R:D\nu\circ\alpha^{-1}\circ
T_L]_t=\lambda^{\frac{it^2}{2}}\delta^{it}$ for all
$t\in\mathbb{R}$.

Moreover, we have, for all $s,t\in\mathbb{R}$:
\begin{enumerate}[i)]
\item $\begin{aligned} &\ [D\nu\circ\alpha^{-1}\circ T_L\circ\tau_s
:D\nu\circ\alpha^{-1}\circ T_L]_t=\lambda^{-ist}\\
                       &\ [D\nu\circ\alpha^{-1}\circ T_L\circ R\circ\tau_s
:D\nu\circ\alpha^{-1}\circ T_L\circ R]_t=\lambda^{-ist}\\
                       &\ [D\nu\circ\alpha^{-1}\circ T_L\circ\sigma^{\nu\circ\alpha^{-1}\circ T_L\circ R}_s
:D\nu\circ\alpha^{-1}\circ T_L]_t=\lambda^{ist}\\
                       &\ [D\nu\circ\alpha^{-1}\circ T_L\circ R\circ\sigma^{\nu\circ\alpha^{-1}\circ T_L}_s:D\nu\circ\alpha^{-1}\circ T_L\circ R]_t=\lambda^{-ist}
       \end{aligned}$
\item $R(\lambda)=\lambda$, $R(\delta)=\delta^{-1}$,
$\tau_t(\delta)=\delta$ and $\tau_t(\lambda)=\lambda$ ;
\item $\delta$ is a group-like element i.e $\Gamma(\delta)
      =\delta\surl{\ _{\beta} \otimes_{\alpha}}_{\ N}\delta$.
\end{enumerate}
\end{theo}

Nevertheless, in the setting of adapted measured quantum groupoids,
we can improve the previous results. We want to precise where
$\delta$ sits and the dependence of fundamental elements with
respect to the quasi-invariant weight.

\begin{prop}\label{deltadot}
The scaling operator does not depend on the quasi-invariant weight
but just on the modular group associated with. If $\dot{\delta}$ is
the class of $\delta$ for the equivalent relation
$\delta_1\sim\delta_2$ if, and only if there exists a strictly
positive operator $h$ affiliated to $Z(N)$ such that
$\delta_2^{it}=\beta(h^{it})\delta_1^{it}\alpha(h^{-it})$, then
$\dot{\delta}$ does not depend on the quasi-invariant weight but
just on the modular group associated with.
\end{prop}

\begin{proof}
If $\nu'$ is a n.s.f weight on $N$ such that
$\sigma^{\nu'}=\sigma^{\nu}$, then there exists a strictly positive
$h$ affiliated to $Z(N)$ such that $\nu'=\nu_h$. We just have to
compute:
$$
\begin{aligned}
&\ \quad [D\nu'\circ\alpha^{-1}\circ T_L\circ
R:D\nu'\circ\alpha^{-1}\circ
T_L]_t\\
&=[D\nu_h\circ\alpha^{-1}\circ T_L\circ R:D\Phi\circ
R]_t[D\Phi\circ R:D\Phi]_t[D\Phi:D\nu_h\circ\alpha^{-1}\circ T_L]_t\\
&=\beta([D\nu_h:D\nu]^*_{-t})\lambda^{\frac{1}{2}it^2}\delta^{it}\alpha([D\nu:D\nu_h]_t)=\lambda^{\frac{1}{2}it^2}\beta(h^{it})\delta^{it}\alpha(h^{-it})
\end{aligned}$$
\end{proof}

\begin{prop}
The modulus $\delta$ is affiliated with
$M\cap\alpha(N)'\cap\beta(N)'$.
\end{prop}

\begin{proof}
Since $\Phi=\nu\circ\beta^{-1}\circ S_L$, with the notation of
section \ref{fornu}, we have:
$$\lambda^{\frac{it^2}{2}}\delta^{it}=[D\Phi\circ
R:D\Phi]_t=[DR\circ T_L\circ R:DS_L]_t$$ which belongs to
$M\cap\beta(N)'$. Since $\lambda$ is affiliated with $Z(M)$, we get
that $\delta$ is affiliated with $M\cap\beta(N)'$. Finally, since
$R(\delta)=\delta$, we obtain that $\delta$ is affiliated with
$M\cap\alpha(N)'\cap\beta(N)'$.
\end{proof}

Let $\nu'$ be a n.s.f weight on $N$ such that there exist strictly
positive operator $h$ and $k$ affiliated with $N$ strongly commuting
and $[D\nu':D\nu]_t=k^{\frac{it^2}{2}}h^{it}$ for all
$t\in\mathbb{R}$. By \cite{Vae} (proposition 5.1), it is equivalent
to $\sigma_t^{\nu}(h^{is})=k^{ist}h^{is}$ for all $s,t\in\mathbb{R}$
and $\nu'=\nu_h$ in the sense of \cite{Vae}. This hypothesis is
satisfied, in particular, if $\sigma^{\nu}$ and $\sigma^{\nu'}$
commute each other. In this cas, $k$ is affiliated with $Z(N)$.

\begin{prop}
There exists a n.s.f operator-valued weight $T_L'$ from $M$ to
$\alpha(N)$ which is $\beta$-adapted w.r.t $\nu'$ such that, for all
$t\in\mathbb{R}$, we have:
$$[DT_L':DT_L]_t=\beta(k^{\frac{-it^2}{2}}h^{it})$$
\end{prop}

\begin{proof}
By \ref{timpo}, there exists a n.s.f operator-valued weight $S_L$
from $M$ to $\beta(N)$ such that $\nu\circ\alpha^{-1}\circ
T_L=\nu\circ\beta^{-1}\circ S_L$ so that $S_L$ is $\alpha$-adapted
w.r.t $\nu$. Then, again by \ref{timpo}, there exists a n.s.f
operator-valued weight $T_L'$ from $M$ to $\alpha(N)$ such that
$\nu'\circ\beta^{-1}\circ S=\nu\circ\alpha^{-1}\circ T_L'$ so that
$T_L'$ is $\beta$-adapted w.r.t $\nu'$. Then, we compute the
Radon-Nikodym cocycle for all $t\in\mathbb{R}$:
$$
\begin{aligned}
\ [DT_L':DT_L]_t &=[D\nu\circ\alpha^{-1}\circ
T_L':D\nu\circ\alpha^{-1}\circ
T_L]_t\\
&=[D\nu'\circ\beta^{-1}\circ
S:D\nu\circ\beta^{-1}\circ S]_t\\
&=\beta([D\nu':D\nu]_{-t}^*)=\beta(k^{\frac{-it^2}{2}}h^{it})
\end{aligned}$$

\end{proof}

\begin{coro}
We have:
$$\nu\circ\alpha^{-1}\circ
T_L'=(\nu\circ\alpha^{-1}\circ T_L)_{\beta(h)}\quad\text{ and
}\quad\nu'\circ\alpha^{-1}\circ T_L'=(\nu\circ\alpha^{-1}\circ
T_L)_{\alpha(h)\beta(h)}$$
\end{coro}

\begin{proof}
Come from \cite{Vae} (proposition 5.1) and the following equality,
for all $t\in\mathbb{R}$, $[D\nu'\circ\alpha^{-1}\circ
T_L':D\nu\circ\alpha^{-1}\circ
T_L]_t=\alpha(k^{\frac{it^2}{2}})\beta(k^{\frac{-it^2}{2}})\alpha(h^{it})\beta(h^{it})$.
\end{proof}

\begin{prop}
$T_L'$ is left invariant.
\end{prop}

\begin{proof}
Let $a\in {\mathcal M}_{T_L'}^+$. By left invariance of $T_L$, we
have:
$$
\begin{aligned}
(id\surl{\ _{\beta} \star_{\alpha}}_{\
\nu'}\nu'\circ\alpha^{-1}\circ T_L')(\Gamma(a))&=(id\surl{\ _{\beta}
\star_{\alpha}}_{\ \nu}\nu\circ\alpha^{-1}\circ
T_L')(\Gamma(a))\\
&=(id\surl{\ _{\beta} \star_{\alpha}}_{\
\nu}(\nu\circ\alpha^{-1}\circ
T_L)_{\beta(h)})(\Gamma(a))\\
&=(id\surl{\ _{\beta} \star_{\alpha}}_{\
\nu}\nu\circ\alpha^{-1}\circ
T_L)(\Gamma(\beta(h^{1/2})a\beta(h^{1/2})))\\
&=T_L(\beta(h^{1/2})a\beta(h^{1/2}))=T'(a)
\end{aligned}$$

\end{proof}

We state the right version of these results:

\begin{prop}
There exists a n.s.f right invariant operator-valued weight $T_R'$
which is $\alpha$-adapted w.r.t $\nu'$ such that, for all
$t\in\mathbb{R}$, we have:
$$[DT_R':DT_R]_t=\alpha(k^{\frac{it^2}{2}}h^{it})$$
Moreover, we have:
$$\nu\circ\beta^{-1}\circ
T_R'=(\nu\circ\beta^{-1}\circ T_R)_{\alpha(h)}\quad\text{ and
}\quad\nu'\circ\beta^{-1}\circ T_R'=(\nu\circ\beta^{-1}\circ
T_R)_{\alpha(h)\beta(h)}$$
\end{prop}

\begin{lemm}
The application $I_{\nu}^{\nu'}$ defined by the following formula:
$$I_{\nu}^{\nu'}(\xi\surl{\ _{\beta} \otimes_{\alpha}}_{\
\nu}\eta)=\beta(k^{-i/8})\xi\surl{\ _{\beta} \otimes_{\alpha}}_{\
\nu'}\alpha(h^{1/2})\eta$$ for all $\xi\in H$ and  $\eta\in
D(_{\alpha}H,\nu)\cap {\mathcal D}(\alpha(h^{1/2}))$, is an
isomorphism of $\beta(N)'-\alpha(N)'^o$-bimodules from $H\surl{\
_{\beta} \otimes_{\alpha}}_{\ \nu}H$ onto $H\surl{\ _{\beta}
\otimes_{\alpha}}_{\ \nu'}H$.
\end{lemm}

\begin{proof}
For all $x\in {\mathcal N}_{\nu'}$, we have:
$$\alpha(x)\alpha(h^{1/2})\eta=\alpha(xh^{1/2})\eta
=R^{\alpha,\nu}(\eta)\Lambda_{\nu}(xh^{1/2})=R^{\alpha,\nu}(\eta)\Lambda_{\nu'}(x)$$
so that $\alpha(h^{1/2})\eta\in D(_{\alpha}H,\nu)$ and
$R^{\alpha,\nu'}(\alpha(h^{1/2})\eta)=R^{\alpha,\nu}(\eta)$. Also,
we recall that $J_{\nu'}=J_{\nu}k^{-i/8}J_{\nu}k^{i/8}J_{\nu}$ by
\cite{Vae} (proposition 2.5). Then, we have:
$$
\begin{aligned}
&\quad\ (\beta(k^{-i/8})\xi_1\surl{\ _{\beta} \otimes_{\alpha}}_{\
\nu}\alpha(h^{1/2})\eta_1|\beta(k^{-i/8})\xi_2\surl{\ _{\beta}
\otimes_{\alpha}}_{\ \nu}\alpha(h^{1/2})\eta_2)\\
&=(\beta(J_{\nu'}<\alpha(h^{1/2})\eta_1,\alpha(h^{1/2})\eta_2>_{\alpha,\nu'}^*J_{\nu'})
\beta(k^{-i/8})\xi_1|\beta(k^{-i/8})\xi_2)\\
&=(\beta(k^{-i/8}J_{\nu}k^{-i/8}J_{\nu}k^{i/8}J_{\nu}<\eta_1,\eta_2>_{\alpha,\nu}^*
J_{\nu}k^{-i/8}J_{\nu}k^{i/8}J_{\nu}k^{i/8}) \xi_1|\xi_2)\\
&=(\beta(J_{\nu}<\eta_1,\eta_2>_{\alpha,\nu}^*J_{\nu})
\xi_1|\xi_2)=(\xi_1\surl{\ _{\beta} \otimes_{\alpha}}_{\
\nu}\eta_1|\xi_2\surl{\ _{\beta} \otimes_{\alpha}}_{\ \nu}\eta_2)
\end{aligned}$$
\end{proof}

\begin{rema}
For all $\xi\in D(H_{\beta},\nu^o)$ and $\eta\in D(_{\alpha}H,\nu)$,
we have:
$$
\begin{aligned}
I_{\nu}^{\nu'}(\xi\surl{\ _{\beta} \otimes_{\alpha}}_{\
\nu}\eta)&=\beta(k^{-i/8})\xi\surl{\ _{\beta} \otimes_{\alpha}}_{\
\nu'}\alpha(h^{1/2})\eta=\xi\surl{\ _{\beta} \otimes_{\alpha}}_{\
\nu'}\alpha(k^{-i/8}h^{1/2})\eta\\
&=\beta(\sigma_{i/2}^{\nu}(h^{1/2}))\xi\surl{\
_{\beta}\otimes_{\alpha}}_{\
\nu'}\alpha(k^{-i/8})\eta=\beta(\sigma_{i/2}^{\nu}(k^{-i/8}h^{1/2}))\xi\surl{\ _{\beta} \otimes_{\alpha}}_{\ \nu'}\eta\\
&=\beta(k^{i/8})\xi\surl{\ _{\beta} \otimes_{\alpha}}_{\
\nu'}\alpha(\sigma_{-i/2}^{\nu}(h^{1/2}))\eta=\beta(k^{i/8}h^{1/2})\xi\surl{\
_{\beta} \otimes_{\alpha}}_{\ \nu'}\eta
\end{aligned}$$
\end{rema}

\begin{prop}
Let $(N,M,\alpha,\beta,\Gamma,\nu,T_L,T_R)$ be a adapted measured
quantum groupoid. There exists a adapted measured quantum groupoid
$(N,M,\alpha,\beta,\Gamma,\nu',T_L',T_R')$ fundamental objects of
which, $R'$, $\tau'$, $\lambda'$, $\delta'$ and $P'$, are expressed,
for all $t\in\mathbb{R}$, in the following way:
\begin{center}
\begin{minipage}{14cm}
\begin{enumerate}[i)]
\item $R'=R$, $\lambda'=\lambda$ and $\delta'=\delta$
\item
$\tau_t'=Ad_{\alpha(k^{\frac{-it^2}{2}}h^{-it})\beta(k^{\frac{it^2}{2}}h^{it})}\circ\tau_t
=Ad_{\alpha([D\nu':D\nu]_t^*)\beta([D\nu':D\nu]_t)}\circ\tau_t$
\item $P'^{it}=\alpha(k^{\frac{it^2}{2}}h^{it})\beta(k^{\frac{-it^2}{2}}h^{-it})
J_{\Phi}\alpha(k^{\frac{it^2}{2}}h^{it})\beta(k^{\frac{-it^2}{2}}h^{-it})J_{\Phi}P^{it}$
\end{enumerate}
\end{minipage}
\end{center}
\end{prop}

\begin{proof}
The existence of $(N,M,\alpha,\beta,\Gamma,\nu',T_L',T_R')$ has been
already proved. We put $\Phi'=\nu'\circ\alpha^{-1}\circ T_L'$ and
$\Psi'=\nu'\circ\beta^{-1}\circ T_R'$. Let $x,y\in {\mathcal
N}_{T_R'}\cap {\mathcal N}_{\Psi'}$. By \cite{Vae} (proposition
2.5), we have:
$$
\begin{aligned}
J_{\Psi'}\Lambda_{\Psi'}(x)&=J_{\Psi}\alpha(k^{-i/8})\beta(k^{i/8})
J_{\Psi}\alpha(k^{i/8})\beta(k^{-i/8})J_{\Psi}\Lambda_{\Psi}(x\alpha(h^{1/2})\beta(h^{1/2}))\\
\omega_{J_{\Psi'}\Lambda_{\Psi'}(x)}
&=\omega_{\alpha(k^{i/8})\beta(k^{-i/8})J_{\Psi}\Lambda_{\Psi}(x\alpha(h^{1/2})\beta(h^{1/2}))}
\end{aligned}$$
Then, we easily verify
$$\lambda^{\beta,\alpha,\nu'}
_{\alpha(k^{i/8})\beta(k^{-i/8})J_{\Psi}\Lambda_{\Psi}(x\alpha(h^{1/2})\beta(h^{1/2}))}=
I_{\nu}^{\nu'}\lambda^{\beta,\alpha,\nu}_
{\alpha(k^{i/8})J_{\Psi}\Lambda_{\Psi}(x\alpha(h^{1/2}))}$$ We
compute:
$$
\begin{aligned}
&\ \quad(\omega_{J_{\Psi'}\Lambda_{\Psi'}(x)}\surl{\ _{\beta}
\star_{\alpha}}_{\ \nu'}id)(\Gamma(y^*y))\\
&=(\omega_{\alpha(k^{i/8})\beta(k^{-i/8})J_{\Psi}\Lambda_{\Psi}(x\alpha(h^{1/2})\beta(h^{1/2}))}\surl{\
_{\beta} \star_{\alpha}}_{\ \nu'}id)(\Gamma(y^*y))\\
&=(\omega_{\alpha(k^{i/8})J_{\Psi}\Lambda_{\Psi}(x\alpha(h^{1/2}))}\surl{\
_{\beta} \star_{\alpha}}_{\
\nu}id)(\Gamma(y^*y))\\
&=(\omega_{J_{\Psi}\Lambda_{\Psi}(x\alpha(k^{-i/8}h^{1/2}))}\surl{\
_{\beta} \star_{\alpha}}_{\ \nu}id)(\Gamma(y^*y))
\end{aligned}$$
Apply $R$ to get:
$$
\begin{aligned}
&\ \quad R[(\omega_{J_{\Psi'}\Lambda_{\Psi'}(x)}\surl{\ _{\beta}
\star_{\alpha}}_{\ \nu'}id)(\Gamma(y^*y))]\\
&=(\omega_{J_{\Psi}\Lambda_{\Psi}(y)}\surl{\
_{\beta} \star_{\alpha}}_{\ \nu}\!id)(\Gamma(\alpha(k^{i/8}h^{1/2})x^*x\alpha(k^{-i/8}h^{1/2})))\\
&=(\omega_{\alpha(k^{-i/8}h^{1/2})J_{\Psi}\Lambda_{\Psi}(y)}\surl{\
_{\beta} \star_{\alpha}}_{\ \nu}id)(\Gamma(x^*x))\\
&=(\omega_{\alpha(k^{i/8})J_{\Psi}\Lambda_{\Psi}(y\alpha(h^{1/2}))}\surl{\
_{\beta} \star_{\alpha}}_{\ \nu}id)(\Gamma(x^*x))\\
&=(\omega_{J_{\Psi'}\Lambda_{\Psi'}(y)}\surl{\ _{\beta}
\star_{\alpha}}_{\
\nu'}id)(\Gamma(x^*x))=R'[(\omega_{J_{\Psi'}\Lambda_{\Psi'}(x)}\surl{\
_{\beta}
\star_{\alpha}}_{\ \nu'}id)(\Gamma(y^*y))]\\
\end{aligned}$$
so that $R=R'$. For all $a\in M$, $\xi\in D(H_{\beta},\nu'^o)$ and
$t\in\mathbb{R}$, we have:
$$
\begin{aligned}
&\quad\ \tau_t((\omega_{\xi}\surl{\ _{\beta} \star_{\alpha}}_{\
\nu'}id)(\Gamma(a)))\\
&=\tau_t(\alpha(k^{-i/8}h^{-1/2})(\omega_{\xi}\surl{\ _{\beta}
\star_{\alpha}}_{\
\nu}id)(\Gamma(a))\alpha(k^{i/8}h^{-1/2}))\\
&=\alpha(\sigma_t^{\nu}(k^{-i/8}h^{-1/2}))\tau_t((\omega_{\xi}\surl{\
_{\beta} \star_{\alpha}}_{\
\nu}id)(\Gamma(a)))\alpha(\sigma_t^{\nu}(k^{i/8}h^{-1/2}))\\
&=\alpha(k^{-t/2-i/8}h^{-1/2})(\omega_{\Delta_{\Psi}^{-it}\xi}\surl{\
_{\beta} \star_{\alpha}}_{\
\nu}id)(\Gamma(\sigma_{-t}^{\Psi}(a)))\alpha(k^{-t/2+i/8}h^{-1/2})
\end{aligned}$$
By \cite{Vae} (proposition 2.4 and corollaire 2.6), we know that:
$$
\begin{aligned}
&\ \quad(\omega_{\Delta_{\Psi}^{-it}\xi}\surl{\ _{\beta}
\star_{\alpha}}_{\ \nu}id)(\Gamma(\sigma_{-t}^{\Psi}(a)))\\
&=(\omega_{\alpha(k^{\frac{-it^2}{2}}h^{it})
\beta(k^{\frac{it^2}{2}}h^{it})\Delta_{\Psi'}^{-it}\xi}\surl{\
_{\beta} \star_{\alpha}}_{\
\nu}id)(\Gamma(Ad_{\alpha(k^{\frac{it^2}{2}}h^{-it})\beta(k^{\frac{-it^2}{2}}h^{-it})}\circ
\sigma_{-t}^{\Psi'}(a)))
\end{aligned}$$ so that:
$$
\begin{aligned}
&\quad\ \tau_t((\omega_{\xi}\surl{\ _{\beta} \star_{\alpha}}_{\
\nu'}id)(\Gamma(a)))\\
&=Ad_{\alpha(k^{-t/2+i/8}h^{-1/2})\beta(k^{\frac{it^2}{2}}h^{it})}\circ
(\omega_{\beta(k^{\frac{it^2}{2}}h^{it})\Delta_{\Psi'}^{-it}\xi}\surl{\
_{\beta} \star_{\alpha}}_{\
\nu}id)(\Gamma(\sigma_{-t}^{\Psi'}(a)))\\
&=\alpha(k^{\frac{-it^2}{2}}h^{-it})\beta(k^{\frac{it^2}{2}}h^{it})
(\omega_{\Delta_{\Psi'}^{-it}\xi}\surl{\ _{\beta} \star_{\alpha}}_{\
\nu'}id)(\Gamma(\sigma_{-t}^{\Psi'}(a)))
\alpha(k^{\frac{it^2}{2}}h^{it})\beta(k^{\frac{-it^2}{2}}h^{-it})\\
&=\alpha(k^{\frac{-it^2}{2}}h^{-it})\beta(k^{\frac{it^2}{2}}h^{it})
\tau'_t((\omega_{\xi}\surl{\ _{\beta} \star_{\alpha}}_{\
\nu'}id)(\Gamma(a)))\alpha(k^{\frac{it^2}{2}}h^{it})\beta(k^{\frac{-it^2}{2}}h^{-it})
\end{aligned}$$
Consequently, we have:
$$\tau_t'(z)=\alpha(k^{\frac{it^2}{2}}h^{it})\beta(k^{\frac{-it^2}{2}}h^{-it})
\tau_t(z)\alpha(k^{\frac{-it^2}{2}}h^{-it})\beta(k^{\frac{it^2}{2}}h^{it})$$
for all $z\in M$ and  $t\in\mathbb{R}$. Now, we compute the
Radon-Nikodym cocycle:
$$
\begin{aligned}
&\quad\ [D\nu'\circ\alpha^{-1}\circ T'\circ
R:D\nu'\circ\alpha^{-1}\circ
T']_t\\
&=[D\nu'\alpha^{-1}T'R:D\nu\alpha^{-1}TR]_t[D\nu\alpha^{-1}
TR:D\nu\alpha^{-1}T]_t[D\nu\alpha^{-1}T:D\nu'\alpha^{-1}T']_t\\
&=\alpha([D\nu':D\nu]_t)\beta([D\nu':D\nu]_{-t}^*)\lambda^{\frac{it^2}{2}}\delta^{it}
\alpha([D\nu:D\nu']_t)\beta([D\nu:D\nu']_{-t}^*)
\end{aligned}$$
which is equal to $\lambda^{\frac{it^2}{2}}\delta^{it}$. Finally, we
express the manageable operator $P'$ in terms of $P$. We have, for
all $x\in {\mathcal N}_{T_L'}\cap {\mathcal N}_{\Phi'}$ and
$t\in\mathbb{R}$:
$$
\begin{aligned}
&\quad\ P'^{it}\Lambda_{\Phi'}(x)=\lambda'^{t/2}\Lambda_{\Phi'}(\tau'_t(x))\\
&=\lambda^{t/2}\Lambda_{\Phi}(\alpha(k^{\frac{it^2}{2}}h^{it})\beta(k^{\frac{-it^2}{2}}h^{-it})
\tau_t(x)\alpha(k^{\frac{-it^2}{2}}h^{-it})\beta(k^{\frac{it^2}{2}}h^{it})
\alpha(h^{1/2})\beta(h^{1/2}))
\end{aligned}$$
which is equal to the value of:
$$\lambda^{t/2}\alpha(k^{\frac{it^2}{2}}h^{it})\beta(k^{\frac{-it^2}{2}}h^{-it})
J_{\Phi}\alpha(k^{\frac{it^2}{2}}h^{it})\beta(k^{\frac{-it^2}{2}}h^{-it})
\alpha(k^{t/2})\beta(k^{t/2})J_{\Phi}$$ on
$\Lambda_{\Phi}(\tau_t(x)\alpha(h^{1/2})\beta(h^{1/2}))$ and the
value of:
$$\lambda^{t/2}\alpha(k^{\frac{it^2}{2}}h^{it})\beta(k^{\frac{-it^2}{2}}h^{-it})
J_{\Phi}\alpha(k^{\frac{it^2}{2}}h^{it})\beta(k^{\frac{-it^2}{2}}h^{-it})J_{\Phi}$$
on $\Lambda_{\Phi}(\tau_t(x\alpha(h^{1/2})\beta(h^{1/2})))$ which
is:
$$\alpha(k^{\frac{it^2}{2}}h^{it})\beta(k^{\frac{-it^2}{2}}h^{-it})
J_{\Phi}\alpha(k^{\frac{it^2}{2}}h^{it})\beta(k^{\frac{-it^2}{2}}h^{-it})J_{\Phi}
P^{it}\Lambda_{\Phi'}(x)$$
\end{proof}

Thanks to these formulas, we verify for example that
$\tau'_t(\alpha(n))=\alpha(\sigma_t^{\nu'}(n))$,
$\tau'_t(\beta(n))=\beta(\sigma_t^{\nu'}(n))$ and $\tau'$ is
implemented by $P'$.

\begin{prop}
Let $(N,M,\alpha,\beta,\Gamma,\nu,T_L,T_R)$ be adapted measured
quantum groupoid and let $\tilde{T_L}$ be an other n.s.f left
invariant operator-valued weight which is $\beta$-adapted w.r.t
$\nu$. Then fundamental objects $\tilde{R}$, $\tilde{\tau}$,
$\tilde{\lambda}$, $\tilde{\delta}$ and $\tilde{P}$ of the adapted
measured quantum groupoid
$(N,M,\alpha,\beta,\Gamma,\nu,\tilde{T_L},T_R)$ can be expressed in
the following way:
\begin{center}
\begin{minipage}{14cm}
\begin{enumerate}[i)]
\item $\tilde{R}=R$, $\tilde{\tau}=\tau$, $\tilde{\lambda}=\lambda$
and $\tilde{P}=P$
\item $\tilde{\delta}=\delta\alpha(h)\beta(h^{-1})$
where $h$ is affiliated with $Z(N)$ s.t.
$\tilde{T_L}=(T_L)_{\beta(h)}$
\end{enumerate}
\end{minipage}
\end{center}
\end{prop}

\begin{proof}
By uniqueness theorem, there exists a strictly positive operator $h$
affiliated with $Z(N)$ such that
$\nu\circ\alpha^{-1}\circ\tilde{T_L}=(\nu\circ\alpha^{-1}\circ
T_L)_{\beta(h)}$ and, for all $t\in\mathbb{R}$, we have
$[D\tilde{T_L}:DT_L]_t=\beta(h^{it})$. We have already noticed that
$R$ and $\tau$ are independent w.r.t left invariant operator-valued
weight and $\beta$-adapted w.r.t $\nu$. We compute then
Radon-Nykodim cocycle:
$$
\begin{aligned}
&\ \quad
[D\nu\beta^{-1}R\tilde{T_L}R:D\nu\alpha^{-1}\tilde{T_L}]_t\\
&=[D\nu\beta^{-1}R\tilde{T_L}R:D\nu\beta^{-1}
RT_LR]_t[D\Psi:D\Phi]_t
[D\nu\alpha^{-1}T_L:D\nu\alpha^{-1}\tilde{T_L}]_t\\
&=R([D\tilde{T_L}:DT_L]^*_{-t})[D\Psi:D\Phi]_t
[DT_L:D\tilde{T_L}]_t\\
&=\alpha(h^{it})\lambda^{\frac{it^2}{2}}\delta^{it}\beta(h^{-it})
=\lambda^{\frac{it^2}{2}}\delta^{it}\alpha(h^{it})\beta(h^{-it})\\
\end{aligned}$$
Then, it remains to compute $\tilde{P}$. If, we put
$\tilde{\Phi}=\nu\circ\alpha^{-1}\circ\tilde{T_L}$, we have, for all
$t\in\mathbb{R}$ and  $x\in {\mathcal N}_{\tilde{T_L}}\cap {\mathcal
N}_{\tilde{\Phi}}$:
$$
\begin{aligned}
\tilde{P}^{it}\Lambda_{\tilde{\Phi}}(x)
=\tilde{\lambda}^{t/2}\Lambda_{\tilde{\Phi}}(\tilde{\tau}_t(x))
&=\lambda^{t/2}\Lambda_{\Phi}(\tau_t(x)\beta(h^{1/2}))
=\lambda^{t/2}\Lambda_{\Phi}(\tau_t(x\beta(h^{1/2}))\\
&=P^{it}\Lambda_{\Phi}(x\beta(h^{1/2}))\!=\!P^{it}\!\Lambda_{\tilde{\Phi}}(x)
\end{aligned}$$

\end{proof}

\begin{theo}
Let $(N,M,\alpha,\beta,\Gamma,\nu,T_L,T_R)$ and
$(N,M,\alpha,\beta,\Gamma,\nu',T_L',T_R')$ be adapted measured
quantum groupoids such that there exist strictly positive operators
$h$ and $k$ affiliated with $N$ which strongly commute and
$[D\nu':D\nu]_t=k^{\frac{it^2}{2}}h^{it}$ for all $t\in\mathbb{R}$.
For all $t\in\mathbb{R}$, fundamental objects of the two structures
are linked by:
\begin{center}
\begin{minipage}{14cm}
\begin{enumerate}[i)]
\item $R'=R$
\item
$\tau_t'=Ad_{\alpha(k^{\frac{-it^2}{2}}h^{-it})\beta(k^{\frac{it^2}{2}}h^{it})}\circ\tau_t
=Ad_{\alpha([D\nu':D\nu]_t^*)\beta([D\nu':D\nu]_t)}\circ\tau_t$
\item $\lambda'=\lambda$
\item $\dot{\delta'}=\dot{\delta}$ where $\dot{\delta}$ and $\dot{\delta}'$
have been defined in proposition \ref{deltadot}
\item $P'^{it}=\alpha(k^{\frac{it^2}{2}}h^{it})\beta(k^{\frac{-it^2}{2}}h^{-it})
J_{\Phi}\alpha(k^{\frac{it^2}{2}}h^{it})\beta(k^{\frac{-it^2}{2}}h^{-it})J_{\Phi}P^{it}$
\end{enumerate}
\end{minipage}
\end{center}
\end{theo}

\begin{proof}
We successively apply the two previous propositions.
\end{proof}

We summarize results concerning the change of quasi-invariant weight
in the following theorem:

\begin{theo}
Let $(N,M,\alpha,\beta,\Gamma,\nu,T_L,R\circ T_L\circ R)$ be a
adapted measured quantum groupoid. If $\nu'$ is a n.s.f weight on
$N$ and $h$, $k$ are strictly positive operators, affiliated with
$N$, strongly commuting and satisfying
$[D\nu':D\nu]_t=k^{\frac{it^2}{2}}h^{it}$ for all $t\in\mathbb{R}$,
then there exists a n.s.f left invariant operator-valued weight
$\tilde{T_L}$ which is $\beta$-adapted w.r.t $\nu'$. Moreover, if
$(N,M,\alpha,\beta,\Gamma,\nu',T_L',T_R')$ is an other adapted
measured quantum groupoid, then, for all $t\in\mathbb{R}$,
fundamental objects are linked by:
\begin{center}
\begin{minipage}{14cm}
\begin{enumerate}[i)]
\item $R'=R$
\item
$\tau_t'=Ad_{\alpha(k^{\frac{-it^2}{2}}h^{-it})\beta(k^{\frac{it^2}{2}}h^{it})}\circ\tau_t
=Ad_{\alpha([D\nu':D\nu]_t^*)\beta([D\nu':D\nu]_t)}\circ\tau_t$
\item $\lambda'=\lambda$
\item $\dot{\delta'}=\dot{\delta}$ where $\dot{\delta}$ and $\dot{\delta}'$
have been defined in proposition \ref{deltadot}
\item $P'^{it}=\alpha(k^{\frac{it^2}{2}}h^{it})\beta(k^{\frac{-it^2}{2}}h^{-it})
J_{\Phi}\alpha(k^{\frac{it^2}{2}}h^{it})\beta(k^{\frac{-it^2}{2}}h^{-it})J_{\Phi}P^{it}$
\end{enumerate}
\end{minipage}
\end{center}
\end{theo}

\subsection{Characterization}
In theorem \ref{ex1}, we explain how a adapted measured quantum
groupoid can be seen as a generalized quantum groupoid. But it is
easy to characterize them among measured quantum groupoids.

\begin{theo}
A measured quantum groupoid is a adapted measured quantum groupoid
if, and only if $\gamma=\sigma^{\nu}$ if, and only if $\delta$ is
affiliated with $M\cap\alpha(N)'\cap\beta(N)'$.
\end{theo}

\begin{proof}
Straightforward.
\end{proof}

In general, we have not a duality within adapted measured quantum
groupoid category that is the dual structure coming from measured
quantum groupoid is not a adapted measured quantum groupoid anymore.
We can be even more precise by characterizing dual objects of
adapted measured quantum groupoids.

\begin{theo}
A measured quantum groupoid is the dual of a adapted measured
quantum groupoid if, and only if $\gamma_t=\sigma_{-t}^{\nu}$ for al
$t\in\mathbb{R}$.
\end{theo}

\begin{proof}
Let us denote by $M$ a measured quantum groupoid and by
$\widehat{M}$ its dual. By the bi-duality theorem and the previous
theorem, $M$ is the dual of a adapted measured quantum groupoid if,
and only if $\widehat{M}$ is a adapted measured quantum groupoid if,
and only if $\gamma_{-t}=\hat{\gamma}_t=\sigma_t^{\nu}$ for all
$t\in\mathbb{R}$.
\end{proof}

Also, we can deduce a precise result concerning duality within
adapted measured quantum groupoids:

\begin{theo}
For all adapted measured quantum groupoid
$(N,M,\alpha,\beta,\Gamma,\nu,T_L,T_R)$, the dual measured quantum
groupoid is a adapted measured quantum groupoid if, and only if the
basis $N$ is semi-finite.
\end{theo}

\begin{proof}
$(N,\widehat{M},\alpha,\hat{\beta},\hat{\Gamma})$ equipped with
$\widehat{T_L}$ et $\widehat{R}\circ\widehat{T_L}\circ\widehat{R}$
is a adapted measured quantum groupoid if, and only if there exists
a nsf weight $\widehat{\nu}$ on $N$ such that, for all
$t\in\mathbb{R}$, we have
$\sigma_t^{\nu}=\sigma_{-t}^{\widehat{\nu}}$. In this case,
$\hat{\nu}$ is $\sigma^{\nu}$ invariant, so there exists a strictly
positive operator $h$ affiliated to the centralizer of $\nu$ such
that $[D\hat{\nu}:D\nu]_t=h^{it}$. Then, for all $x\in N$, we have
$\sigma_{-t}^{\nu}(x)=h^{it}\sigma_t^{\nu}(x)h^{-it}$ and
$\sigma_{-2t}^{\nu}(x)=h^{it}xh^{-it}$. Then $\sigma^{\nu}_t$ is
inner for all $t\in\mathbb{R}$ and $N$ is semi-finite by theorem
3.14 of \cite{T}. Conversely, if $N$ is semi-finite, there exists a
nsf trace $tr$ on $N$ and a strictly positive operator $h$such that
$\nu=tr(h.)$. So $\hat{\nu}=tr(h^{-1}.)$ satisfies conditions.
\end{proof}


\section{Groupoids}\label{groupoids}

\begin{defi}
A \textbf{groupoid} $G$ is a small category in which each morphism
$\gamma:x\rightarrow y$ is an isomorphism the inverse of which is
$\gamma^{-1}$. Let $G^{\{0\}}$ the set of objects of $G$ that we
identify with $\{\gamma\in G | \gamma\circ\gamma=\gamma\}$. For all
$\gamma\in G$, $\gamma:x\rightarrow y$, we denote
$x=\gamma^{-1}\gamma=s(\gamma)$ we call source object and
$y=\gamma\gamma^{-1}=r(\gamma)$ we call range object. If $G^{\{2\}}$
is the set of pairs $(\gamma_1,\gamma_2)$ of $G$ such that
$s(\gamma_1)=r(\gamma_2)$, then composition of morphisms makes sense
in $G^{\{2\}}$.
\end{defi}

In \cite{R1}, J. Renault defines the structure of locally compact
groupoid $G$ with a Haar system $\{\lambda^u,u\in G^{\{0\}}\}$ and a
quasi-invariant measure $\mu$ on $G^{\{0\}}$. We refer to \cite{R1}
for definitions and notations. We put $\nu=\mu\circ\lambda$. We
refer to \cite{C3} and \cite{ADR} for discussions about transversal
measures.

If $G$ is $\sigma$-compact, J.M Vallin constructs in \cite{V1} two
co-involutive Hopf bimodules on the same basis
$N=L^{\infty}(G^{\{0\}},\mu)$, following T. Yamanouchi's works in
\cite{Y}. The underlying von Neumann algebras are
$L^{\infty}(G,\nu)$ which acts by multiplication on $H=L^2(G,\nu)$
and $\mathcal{L}(G)$ generated by the left regular representation
$L$ of $G$.

We define a (resp. anti-) representation $\alpha$ (resp. $\beta$)
from $N$ in $L^{\infty}(G,\nu)$ such that, for all $f\in N$:
$$\alpha(f)=f\circ r\quad\text{ and }\quad\beta(f)=f\circ s$$

For all $i,j\in\{\alpha,\beta\}$, we define $G^{\{2\}}_{i,j}\subset
G\times G$ and a measure $\nu^2_{i,j}$ such that:
$$H\surl{\ _i\otimes_j}_{\ N}H\text{ is identified with }L^2(G^{\{2\}}_{i,j},\nu^2_{i,j})$$
For example, $G^{\{2\}}_{\beta,\alpha}$ is equal to $G^{\{2\}}$ and
$\nu^2_{\beta,\alpha}$ to $\nu^2$. Then, we construct a unitary
$W_G$ from $H\surl{\ _{\alpha}\otimes_{\alpha}}_{\ \mu}H$ onto
$H\surl{\ _{\beta}\otimes_{\alpha}}_{\ \mu}H$, defined for all
$\xi\in L^2(G^{\{2\}}_{\alpha,\alpha},\nu^2_{\alpha,\alpha})$ by:
$$W_G\xi(s,t)=\xi(s,st)$$ for $\nu^2$-almost all $(s,t)$ in $G^{\{2\}}$.

This leads to define co-products $\Gamma_G$ and $\widehat{\Gamma_G}$
by formulas:
$$\Gamma_G(f)=W_G(1\surl{\ _{\alpha}
  \otimes_{\alpha}}_{\ N}f)W_G^*\quad\text{ and }\quad \widehat{\Gamma_G}(k)=W_G^*(k\surl{\ _{\beta}
  \otimes_{\alpha}}_{\ N}1)W_G$$
for all $f\in L^{\infty}(G,\nu)$ and $k\in\mathcal{L}(G)$, this
explicitly gives:
$$\Gamma_G(f)(s,t)=f(st)$$
for all $f\in L^{\infty}(G,\nu)$ and $\nu^2$-almost all $(s,t)$ in
$G^{\{2\}}$,
$$\widehat{\Gamma_G}(L(h))\xi(x,y)=\int_Gh(s)\xi(s^{-1}x,s^{-1}y)d\lambda^{r(x)}(s)$$
for all $\xi\in
L^2(G^{\{2\}}_{\alpha,\alpha},\nu^2_{\alpha,\alpha})$, $h$ a
continuous function with compact support on $G$ and
$\nu^2_{\alpha,\alpha}$-almost all $(x,y)$ in
$G^{\{2\}}_{\alpha,\alpha}$. Moreover, we define two co-involutions
$j_G$ and $\widehat{j_G}$ by:
$$j_G(f)(x)=f(x^{-1})$$ for all $f\in L^{\infty}(G,\nu)$ and
almost all $x$,
$$\widehat{j_G}(g)=Jg^*J$$ for all $g\in\mathcal{L}(G)$ and where
$J$ is the involution $J\xi=\overline{\xi}$ for all $\xi\in L^2(G)$.
Finally, we define two n.s.f left invariant operator-valued weights
$P_G$ and $\widehat{P_G}$:
$$P_G(f)(y)=\int_Gf(x)d\lambda^{r(y)}(x)\quad\text{ and
}\quad\widehat{P_G}(L(f))=\alpha(f_{|G^{\{0\}}})$$ for all
continuous with compact support $f$ on $G$ $\nu$-almost all $y$ in
$G$.

\begin{theo}
Let $G$ be a $\sigma$-compact, locally compact groupoid with a Haar
system and a quasi-invariant measure $\mu$on units. Then:
$$
(L^{\infty}(G^{\{0\}},\mu),L^{\infty}(G,\nu),\alpha,\beta,\Gamma_G,\mu,P_G,j_GP_Gj_G)$$
is a commutative adapted measured quantum groupoid and:
$$(L^{\infty}(G^{\{0\}},\mu),{\mathcal
L}(G),\alpha,\alpha,\widehat{\Gamma_G},\mu,\widehat{P_G},
\widehat{j_G}\widehat{P_G}\widehat{j_G})
$$ is a symmetric adapted measured quantum groupoid. The unitary $V_G=W_G^*$
is the fundamental unitary of the commutative structure.
\end{theo}

\begin{proof}
By \cite{V1} (th. 3.2.7 and 3.3.7),
$(L^{\infty}(G^{\{0\}},\mu),L^{\infty}(G,\nu),\alpha,\beta,\Gamma_G)$
and $(L^{\infty}(G^{\{0\}},\mu),{\mathcal
L}(G),\alpha,\alpha,\widehat{\Gamma_G})$ are co-involutive Hopf
bimodules with left invariant operator-valued weights; to get right
invariants operator-valued weights, we consider $j_GP_Gj_G$ and
$\widehat{j_G}\widehat{P_G}\widehat{j_G}$.

Since $L^{\infty}(G,\nu)$ is commutative, $P_G$ is adapted w.r.t
$\mu$ by \cite{V1} (theorem 3.3.4),
$\sigma_t^{\mu\circ\alpha^{-1}\circ\widehat{P_G}}$ fixes point by
point $\alpha(N)$ so that $\widehat{P_G}$ is adapted w.r.t $\mu$.

Finally, for all $e,f,g$ continuous functions with compact support
and almost all $(s,t)$ in $G^{\{2\}}$, we have, by \ref{raccourci}:
$$
\begin{aligned}
(1\surl{\ _{\beta}\otimes_{\alpha}}_{\ N}JeJ)W_G(f\surl{\
_{\alpha}\otimes_{\alpha}}_{\ \mu}g)(s,t)&=\overline{e(t)}f(s)g(st)=\Gamma_G(g)(f\surl{\ _{\beta}\otimes_{\alpha}}_{\ \mu}\overline{e})(s,t)\\
&=(1\surl{\ _{\beta}\otimes_{\alpha}}_{\ N}JeJ)U_H(f\surl{\
_{\alpha}\otimes_{\alpha}}_{\ \mu}g)(s,t)
\end{aligned}$$
so that we get $U_H=W_G$.
\end{proof}

\begin{rema}
In the commutative structure, modular function
$\frac{d\nu^{-1}}{d\nu}$ and modulus coincide and the scaling
operator is trivial.
\end{rema}

We have a similar result for adapted measured quantum groupoids in
the sense of Hahn (\cite{Ha1} and \cite{Ha2}):

\begin{theo}
From all measured groupoid $G$, we construct a commutative adapted
measured quantum groupoid
$(L^{\infty}(G^{\{0\}},\mu),L^{\infty}(G,\nu),\alpha,\beta,\Gamma_G,\mu,P_G,j_GP_Gj_G)$
and a symmetric one $(L^{\infty}(G^{\{0\}},\mu),{\mathcal
L}(G),\alpha,\alpha,\widehat{\Gamma_G},\mu,\widehat{P_G},\widehat{j_G}\widehat{P_G}\widehat{j_G})$.
Objects are defined in a similar way as in the locally compact case.
The unitary $V_G$ is the fundamental unitary of the commutative
structure.
\end{theo}

\begin{proof}
Results come from \cite{Y} for the symmetric case. It is sufficient
to apply in this case, technics of \cite{V1} for the commutative
case and invariant operator-valued weights.
\end{proof}

\begin{conj}
If $(N,M,\alpha,\beta,\Gamma,\mu,T_L,T_R)$ is a adapted measured
quantum groupoid such that $M$ is commutative, then there exists a
locally compact groupoid $G$ such that:
$$(N,M,\alpha,\beta,\Gamma,\mu,T_L,T_R)\simeq
(L^{\infty}(G^{\{0\}},\mu),L^{\infty}(G,\nu),\alpha,\beta,\Gamma_G,\mu,P_G,j_G\circ
P_G\circ j_G)$$
\end{conj}

\section{Finite quantum groupoids}\label{finite}

\begin{defi}(Weak Hopf C*-algebras \cite{BSz})
We call \textbf{weak Hopf C*-algebra } or finite quantum groupoid
all $(M,\Gamma,\kappa,\varepsilon)$ where $M$ is a finite
dimensional C*-algebra with a co-product $\Gamma:M\rightarrow
M\otimes M$, a co-unit $\varepsilon$ and an antipode
$\kappa:M\rightarrow M$ such that, for all $x,y\in M$:
\begin{center}
\begin{minipage}{12cm}
\begin{enumerate}[i)]
\item $\Gamma$ is a *-homomorphism (not necessary unital);
\item Unit and co-unit satisfy the following relation:
$$(\varepsilon\otimes\varepsilon)((x\otimes
1)\Gamma(1)(1\otimes y))=\varepsilon(xy)$$
\item $\kappa$ is an anti-homomorphism of algebra and co-algebra such that:
\begin{itemize}
\item $(\kappa\circ *)^2=\iota$
\item $(m(\kappa\otimes id)\otimes id)(\Gamma\otimes id)\Gamma
(x)=(1\otimes x)\Gamma(1)$.
\end{itemize}
where $m$ denote the product on $M$.
\end{enumerate}
\end{minipage}
\end{center}

\end{defi}

We recall some results \cite{NV1}, \cite{NV2} and \cite{BNSz}. If
$(M,\Gamma,\kappa,\varepsilon)$ is a weak Hopf C*-algebra. We call
co-unit range (resp. source) the application
$\varepsilon_t=m(id\otimes \kappa)\Gamma$ (resp.
$\varepsilon_s=m(\kappa\otimes id)\Gamma$). We have
$\kappa\circ\varepsilon_t=\varepsilon_s\circ\kappa$. There exists a
unique faithful positive linear form $h$, called normalized Haar
measure of $(M,\Gamma,\kappa,\varepsilon)$ which is
$\kappa$-invariant, such that $(id\otimes h)(\Gamma(1))=1$ and, for
all $x,y\in M$, we have:
$$(id\otimes h)((1\otimes
y)\Gamma(x))=\kappa((i\otimes h)(\Gamma(y)(1\otimes x)))$$ Moreover,
$E^s_h=(h\otimes id)\Gamma$ (resp. $E^t_h=(id\otimes h)\Gamma$) is a
Haar conditional expectation to the source (resp. range) Cartan
subalgebra $\varepsilon_s(M)$ (resp. range $\varepsilon_t(M)$) such
that $h\circ E^s_h=h$ (resp. $h\circ E^t_h=h$). Range and source
Cartan subalgebras commute.

By \cite{MC}, \cite{Ni}, and \cite{V5}, we can always assume that
$\kappa^2_{|\varepsilon_t(M)}=id$ thanks to a deformation.
\textbf{In the following, we assume that the condition holds.}

Since $h\circ\kappa=h$ and
$\kappa\varepsilon_t=\varepsilon_s\kappa$, we have
$h\circ\varepsilon_t=h\circ\varepsilon_s$.

\begin{theo}
Let $(M,\Gamma,\kappa,\varepsilon)$ be a weak Hopf C*-algebra, $h$
its normalized Haar measure, $E^s_h$ (resp. $E^t_h$) its source
(resp. range) Haar conditional expectation and $\varepsilon_t(M)$
its range Cartan subalgebra. We put $N=\varepsilon_t(M)$,
$\alpha=id_{|N}$, $\beta=\kappa_{|N}$, $\tilde{\Gamma}$ the
co-product $\Gamma$ viewed as an operator which takes value in:
$$M\surl{\ _{\beta} \star_{\alpha}}_{\ N}M\simeq (M\otimes
M)_{\Gamma(1)}$$ and $\mu=h\circ\alpha=h\circ\beta$. Then
$(N,M,\alpha,\beta,\tilde{\Gamma},\mu,E^t_h,E^s_h)$ is a adapted
measured quantum groupoid.
\end{theo}

\begin{proof}
$\alpha$ is a representation from $N$ in $M$ and, since
$\kappa^2_{|\varepsilon_t(M)}=id$, $\beta$ is a anti-representation
from $N$ in $M$. They commute each other because Cartan subalgebras
commute and $\kappa\circ\varepsilon_t=\varepsilon_s\circ\kappa$. For
all $n\in N$, there exists $m\in M$ such that $n=\varepsilon_t(m)$.
So, we have:
$$\tilde{\Gamma}(\alpha(n))=\tilde{\Gamma}(\varepsilon_t(m))=\Gamma(1)(\varepsilon_t(m)\otimes 1)\Gamma(1)
=\alpha(n)\surl{\ _{\beta}\otimes_{\alpha}}_{\ N}1$$ Also, we have
$\tilde{\Gamma}(\beta(n))=1\surl{\ _{\beta} \otimes_{\alpha}}_{\
N}\beta(n)$ and $\tilde{\Gamma}$ is a co-product. Then
$(N,M,\alpha,\beta,\Gamma)$ is Hopf bimodule. Moreover, for all
$n\in N$ and $t\in\mathbb{R}$, we have:
$$
\begin{aligned} \sigma_t^{E_h^t}(\beta(n))=\sigma_t^{h\circ
E_h^t}(\beta(n))=\sigma_t^{h\circ
E_h^s}(\beta(n))&=\sigma_t^{h_{|\beta(N)}}(\beta(n))\\
&=\beta(\sigma_{-t}^{h_{|\beta(N)}\circ\beta}(n))
=\beta(\sigma_{-t}^{\mu}(n))
\end{aligned}$$ and $E_h^t$ is
$\beta$-adapted w.r.t $\mu$. Since $E_h^s=\kappa\circ
E_h^t\circ\kappa$, then $E_h^s$ is $\alpha$-adapted w.r.t $\mu$.
\end{proof}

\begin{theo}
Let $(N,M,\alpha,\beta,\Gamma,\nu,T_L,T_R)$ be a adapted measured
quantum groupoid such that $M$ is finite dimensional. Then, there
exist $\tilde{\Gamma}$, $\kappa$ and $\varepsilon$ such that
$(M,\tilde{\Gamma},\kappa,\varepsilon)$ is a weak Hopf C*-algebra.
\end{theo}

\begin{proof}
By \ref{iden}, we identify via $I_{\beta,\alpha}^{\nu}$,
$L^2(M)\surl{\ _{\beta}\otimes_{\alpha}}_{\ N}L^2(M)$ with a
subspace of $L^2(M)\otimes L^2(M)$. We put
$\tilde{\Gamma}(x)=I_{\beta,\alpha}^{\nu}\Gamma(x)(I_{\beta,\alpha}^{\nu})^*$.
By \cite{V3} (definition 2.2.3), the fundamental
pseudo-multiplicative unitary becomes a multiplicative partial
isometry on $L^2(M)\otimes L^2(M)$ of basis
$(N,\alpha,\hat{\beta},\beta)$ by
$I=I_{\alpha,\hat{\beta}}^{\nu}W(I_{\beta,\alpha}^{\nu})^*$. $I$ is
regular in the sense of \cite{V3} (definition 2.6.3) by \ref{reg}.
Moreover, if we put $H=L^2(M)$, then $Tr_H(R(m))=Tr_H(m)$ for all
$m\in M$ because $R$ is implemented by an anti-unitary, so
$Tr_H\circ\beta=Tr_H\circ\alpha=Tr_H\circ\hat{\beta}$ and we
conclude by \cite{V3} (proposition 3.1.3).
\end{proof}

\begin{rema}
With notations of section \ref{iden}, $\kappa$ and $S$ are linked
by:
$$\kappa(x)=\alpha(n_o^{1/2}d^{1/2})\beta(n_o^{-1/2}d^{-1/2})S(x)
\alpha(n_o^{-1/2}d^{-1/2})\beta(n_o^{1/2}d^{1/2})$$
\end{rema}

\section{Quantum groups}\label{quantumgroups}

\begin{theo}
adapted measured quantum groupoids, basis $N$ on which is equal to
$\mathbb{C}$ are exactly locally compact quantum groups (in the von
Neumann setting) introduced by J. Kustermans and S. Vaes in
\cite{KV2}.
\end{theo}

\begin{proof}
In this case, the notion of relative tensor product is just usual
tensor product of Hilbert spaces, the notion of fibered product is
just tensor product of von Neumann algebras and the notion of
operator-valued weight is just weight.
\end{proof}

\section{Compact case}\label{compactcase}

In this section, we show that pseudo-multiplicative unitaries of
compact type in the sense of \cite{E2} correspond exactly to adapted
measured quantum groupoids with a Haar conditional expectation.

\begin{defi}
Let $W$ be a pseudo-multiplicative unitary over $N$ w.r.t
$\alpha,\beta,\hat{\beta}$. Let $\nu$ be a n.s.f weight on $N$. We
say that $W$ is of \textbf{compact type} w.r.t $\nu$ if there exists
$\xi\in H$ such that:
\begin{center}
\begin{minipage}{13cm}
\begin{enumerate}[i)]
\item $\xi$ belongs to $D(H_{\hat{\beta}},\nu^o)\cap D(_{\alpha}H,\nu)\cap D(H_{\beta},\nu^o)$;
\item
$<\xi,\xi>_{\hat{\beta},\nu^o}=<\xi,\xi>_{\alpha,\nu}=<\xi,\xi>_{\beta,\nu^o}=1$
\item we have $W(\xi\surl{\ _{\hat{\beta}}\otimes_{\alpha}}_{\ \nu}\eta)
=\xi\surl{\ _{\alpha}\otimes_{\beta}}_{\ \nu^o}\eta$ for all
$\eta\in H$.
\end{enumerate}
\end{minipage}
\end{center}
In this case, $\xi$ is said to be \textbf{fixed and bi-normalized}.
We also say that $W$ is of \textbf{discrete type} w.r.t $\nu$ if
$\hat{W}$ is of compact type.
\end{defi}

By \cite{E2} (proposition 5.11), we recall that, if $W$ is of
compact type w.r.t $\nu$ and $\xi$ is a fixed and bi-normalized
vector, then $\nu$ shall be a faithful, normal, positive form on $N$
which is equal to
$\omega_{\xi}\circ\alpha=\omega_{\xi}\circ\beta=\omega_{\xi}\circ\hat{\beta}$
and it is called \textbf{canonical form}.

\begin{prop}\label{compact}
Let $(N,M,\alpha,\beta,\Gamma)$ be a Hopf bimodule. Assume there
exist:
\begin{center}
\begin{minipage}{13cm}
\begin{enumerate}[i)]
\item a n.f left invariant conditional expectation from $E$
to $\alpha(N)$;
\item a n.f right invariant conditional expectation from $F$
to $\beta(N)$;
\item a n.f state $\nu$ on $N$ such that $\nu\circ\alpha^{-1}
\circ E=\nu\circ\beta^{-1}\circ F$.
\end{enumerate}
\end{minipage}
\end{center}

Then $(N,M,\alpha,\beta,\Gamma,\nu,E,F)$ is a adapted measured
quantum groupoid. Moreover, if $R,\tau,\lambda$ and $\delta$ are
fondamental objects of the structure, then we have $F=R\circ E\circ
R$ and $\lambda=\delta=1$. Finally,
$\Lambda_{\nu\circ\alpha^{-1}\circ E}(1)$ is co-fixed and
bi-normalized, and the fundamental pseudo-multiplicative unitary $W$
is weakly regular and of discrete type in sense of \cite{E2}
(paragraphe 5).
\end{prop}

\begin{proof}
For all $t\in\mathbb{R}$ and $n\in N$, we have:
$$\sigma_t^E(\beta(n))=\sigma_t^{\nu\circ\alpha^{-1}
\circ E}(\beta(n))=\sigma_t^{\nu\circ\beta^{-1} \circ
F}(\beta(n))=\beta(\sigma_{-t}^{\nu}(n))$$ Also, we have:
$$\sigma_t^F(\alpha(n))=\sigma_t^{\nu\circ\beta^{-1}
\circ F}(\alpha(n))=\sigma_t^{\nu\circ\alpha^{-1} \circ
E}(\alpha(n))=\alpha(\sigma_t^{\nu}(n))$$ so that
$(N,M,\alpha,\beta,\Gamma,\nu,E,F)$ is a adapted measured quantum
groupoid. By definition, we have:
$$[D\nu\circ\alpha^{-1}
\circ E\circ R:D\nu\circ\alpha^{-1} \circ
E]_t=\lambda^{\frac{it^2}{2}}\delta^{it}$$ On the other hand, since
$\nu\circ\alpha^{-1}\circ E=\nu\circ\beta^{-1}\circ F$ and by
uniqueness, there exists a strictly positive element $h$ affiliated
with $Z(N)$:
$$
\begin{aligned}
\ [D\nu\circ\alpha^{-1}\circ E\circ R:D\nu\circ\alpha^{-1}\circ
E]_t&=[DR\circ E\circ R:DF]_t=\alpha(h^{it})
\end{aligned}$$
We deduce that $\lambda=1$ and $\delta=\alpha(h)$, so
$\alpha(h^{-1})=\delta^{-1}=R(\delta)=\beta(h)$ and by \cite{E1}
(5.2), we get $h=1$.

We put $\Phi=\nu\circ\alpha^{-1}\circ E$. If $(\xi_i)_{i\in I}$ is a
$(N^o,\nu^o)$-basis of $(H_{\Phi})_{\beta}$ then, for all $v\in
D(H_{\beta},\nu^o)$:
$$
\begin{aligned}
U_H(v \surl{\ _{\alpha} \otimes_{\hat{\beta}}}_{\ \ \nu^o}
\Lambda_{\Phi}(1))&=\sum_{i
    \in I} \xi_{i} \surl{\ _{\beta} \otimes_{\alpha}}_{\ \nu}
  \Lambda_{\Phi} ((\omega_{v,\xi_i} \surl{\
      _{\beta} \star_{\alpha}}_{\ \nu} id)(\Gamma(1)))\\
&=\sum_{i \in I} \xi_{i} \surl{\ _{\beta} \otimes_{\alpha}}_{\ \nu}
\alpha(<v,\xi_i>_{\beta,\nu^o})\Lambda_{\Phi}(1)=v\surl{\ _{\beta}
\otimes_{\alpha}}_{\ \nu}\Lambda_{\Phi}(1)
\end{aligned}$$
It is easy to see that $\Lambda_{\Phi}(1)$ belongs to
$D((H_{\Phi})_{\hat{\beta}},\nu^o)\cap D(_{\alpha}H_{\Phi},\nu)$ and
satisfies
$<\Lambda_{\Phi}(1),\Lambda_{\Phi}(1)>_{\hat{\beta},\nu^o}=
<\Lambda_{\Phi}(1),\Lambda_{\Phi}(1)>_{\alpha,\nu}=1$ so that, by
continuity, we get $U_H(v \surl{\ _{\alpha}
\otimes_{\hat{\beta}}}_{\ \ \nu^o} \Lambda_{\Phi}(1))=v\surl{\
_{\beta} \otimes_{\alpha}}_{\ \nu}\Lambda_{\Phi}(1)$ for all $v\in
H$ i.e $\Lambda_{\Phi}(1)$ is co-fixed and bi-normalized. Since
$\nu\circ\alpha^{-1}\circ E=\Phi=\nu\circ\beta^{-1}\circ F$, we have
by \ref{prem}, for all $n\in \mathcal{N}_{\nu}$:
$$\beta(n^*)\Lambda_{\Phi}(1)=\beta(n^*)J_{\Phi}\Lambda_{\Phi}(1)
=J_{\Phi}\Lambda_F(1)\Lambda_{\nu}(n)$$ so that $\Lambda_{\Phi}(1)$
is $\beta$-bounded w.r.t $\nu^o$ and
$R^{\beta,\nu^o}(\Lambda_{\Phi}(1))=J_{\Phi}\Lambda_F(1)J_{\nu}$.
Consequently, $\Lambda_{\Phi}(1)$ is bi-normalized and $W$ is of
discrete type.
\end{proof}

\begin{coro}
Let $W$ be a weakly regular pseudo-multiplicative unitary over $N$
w.r.t $\alpha,\beta,\hat{\beta}$ of compact type w.r.t the canonical
form $\nu$. If $\xi$ a fixed and bi-normalized vector, we put:
\begin{center}
\begin{minipage}{13cm}
\begin{enumerate}[i)]
\item ${\mathcal A}$ the von Neumann algebra generated by right leg of $W$;
\item $\Gamma(x)=\sigma_{\nu^o}W(x\surl{\ _{\alpha}\otimes_{\beta}}_{\ \nu^o}1)
W^*\sigma_{\nu}$ for all $x\in {\mathcal A}$ ;
\item $E=(\omega_{\xi}\surl{\ _{\beta}\star_{\alpha}}_{\ \nu}id)\circ\Gamma$
and $F=(id\surl{\ _{\beta}\star_{\alpha}}_{\
\nu}\omega_{\xi})\circ\Gamma$.
\end{enumerate}
\end{minipage}
\end{center}
Then $(N,{\mathcal A},\alpha,\beta,\Gamma,\nu,E,F)$ is a adapted
measured quantum groupoid such that $E$ and $F$ are n.f conditional
expectations. Moreover, if $R,\tau,\lambda$ and $\delta$ are the
fundamental objects of the structure, we have $F=R\circ E\circ R$,
$\lambda=\delta=1$ and the fundamental unitary is $\hat{W}$.
\end{coro}

\begin{proof}
By \cite{EV} (6.3), we know that $(N,{\mathcal
A},\alpha,\beta,\Gamma)$ is a Hopf bimodule. By \cite{E2} (theorem
6.6), $E$ is a n.f left invariant conditional expectation from
${\mathcal A}$ to $\alpha(N)$. By \cite{E2} (propositions 6.2 and
6.3), $F$ is a n.f right invariant conditional expectation from
${\mathcal A}$ to $\beta(N)$. Moreover, we clearly have
$\omega_{\xi}\circ E=\omega_{\xi}\circ F$ so that
$\nu\circ\alpha^{-1}\circ E=\nu\circ\beta^{-1}\circ F$. We are in
conditions of the previous proposition an we get that $(N,{\mathcal
A},\alpha,\beta,\Gamma,\nu,E,F)$ is a adapted measured quantum
groupoid, $F=R\circ E\circ R$ and $\lambda=\delta=1$. Finally, by
\cite{E2} (corollaire 7.7), $\hat{W}$ is the fundamental unitary.
(More exactly, it is $\sigma_{\nu^o}W_s^*\sigma_{\nu}$ where $W_s$
is the standard form of $W$ in th sense of \cite{E2} (paragraph 7)).
\end{proof}

The converse is also true and so we characterize the compact case:

\begin{coro}
Let $(N,M,\alpha,\beta,\Gamma)$ be a Hopf bimodule. We assume there
exist:
\begin{center}
\begin{minipage}{13cm}
\begin{enumerate}[i)]
\item a co-involution $R$;
\item a n.f left invariant conditional expectation from $E$
to $\alpha(N)$.
\end{enumerate}
\end{minipage}
\end{center}
Then there exists a n.f state $\nu$ on $N$ such that
$(N,M,\alpha,\beta,\Gamma,\nu,E,R\circ E\circ R)$ is a adapted
measured quantum groupoid with trivial modulus and scaling operator
and the fundamental unitary of which is of discrete type w.r.t
$\nu$.
\end{coro}

\begin{proof}
We put $F=R\circ E\circ R$ which is a n.f right invariant
conditional expectation from $M$ to $\beta(N)$. We also put:
$$\tilde{E}=E_{|\beta(N)}:\beta(N)\rightarrow\alpha(Z(N))\text{ and }
\tilde{F}=F_{|\alpha(N)}:\alpha(N)\rightarrow\beta(Z(N))$$ We have,
for all $m\in M$:
$$
\begin{aligned}
\tilde{F}E(m)\surl{\ _{\beta} \otimes_{\alpha}}_{\ N}1&=(F\surl{\
_{\beta} \star_{\alpha}}_{\ N}id)(E(m)\surl{\ _{\beta}
\otimes_{\alpha}}_{\ N}1)\\
&=(F\surl{\ _{\beta} \star_{\alpha}}_{\ N}id)(id\surl{\ _{\beta}
\star_{\alpha}}_{\ N}E)\Gamma(m)\\
&=(id\surl{\ _{\beta} \star_{\alpha}}_{\ N}E)(F\surl{\ _{\beta}
\star_{\alpha}}_{\ N}id)\Gamma(m)\\
&=(id\surl{\ _{\beta} \star_{\alpha}}_{\ N}E)(1\surl{\ _{\beta}
\otimes_{\alpha}}_{\ N}F(m))=1\surl{\ _{\beta} \otimes_{\alpha}}_{\
N}\tilde{E}F(m)
\end{aligned}$$
so, if $\tilde{F}E(m)=\beta(n)$ for some $n\in Z(N)$, then
$\tilde{E}F(m)=\alpha(n)$. Moreover, we have:
$$
\begin{aligned}
\tilde{E}F(m)\surl{\ _{\beta} \otimes_{\alpha}}_{\ N}1&=EF(m)\surl{\
_{\beta} \otimes_{\alpha}}_{\ N}1=(id\surl{\
_{\beta} \star_{\alpha}}_{\ N}E)\Gamma(F(m))\\
&=(id\surl{\ _{\beta} \star_{\alpha}}_{\ N}E)(1\surl{\ _{\beta}
\otimes_{\alpha}}_{\ N}F(m))=1\surl{\ _{\beta} \otimes_{\alpha}}_{\
N}\tilde{E}F(m)
\end{aligned}$$
so that $\alpha(n)=\beta(n)$. Consequently
$\tilde{E}F(m)=\tilde{F}E(m)$ and $EF=FE$ is a n.f conditional
expectation from $M$ to: $$\tilde{N}=\alpha(\{n\in
Z(N),\alpha(n)=\beta(n)\})=\beta(\{n\in Z(N),\alpha(n)=\beta(n)\})$$
Also, we have $R_{|\tilde{N}}=id$. So, if $\omega$ is a n.f state on
$\tilde{N}$, we have
$\omega\circ\tilde{E}\circ\beta=\omega\circ\tilde{F}\circ\alpha$ and
$\nu=\omega\circ\tilde{E}\circ\beta=\omega\circ\tilde{F}\circ\alpha$
satisfies hypothesis of \ref{compact}: then, corollary holds.
\end{proof}

\begin{coro}
Let $(N,M,\alpha,\beta,\Gamma,\nu,T_L,T_R)$ be a adapted measured
quantum groupoid such that $T_L$ is a conditional expectation. Then
there exists a n.f state $\nu'$ on $N$ such that
$\sigma^{\nu'}=\sigma^{\nu}$ and the fundamental unitary is of
discrete type w.r.t $\nu'$.
\end{coro}

\begin{proof}
Let $R$ be the co-involution. By the previous corollary, there
exists a n.f state $\nu'$ on $N$ such that
$(N,M,\alpha,\beta,\Gamma,\nu',T_L,R\circ T_L\circ R)$ is a adapted
measured quantum groupoid. Since $T_L$ is $\beta$-adapted w.r.t
$\nu$ and $\nu'$, we have $\sigma^{\nu'}=\sigma^{\nu}$. We easily
verify that the fundamental unitary of the first structure coincides
with that of the last one which is of discrete type w.r.t $\nu'$ by
the previous corollary.
\end{proof}

\section{Quantum space quantum groupoid}\label{qsqg}
\subsection{Definition} Let $M$ be a von Neumann algebra. $M$ acts on
$H=L^2(M)=L^2_{\nu}(M)$ where $\nu$ is a n.s.f weight on $M$. We
denote by $M'$, (resp. $Z(M)'$) the commutant of $M$ (resp. $Z(M)$)
in $\mathcal{L}(L^2(M))$. Let $tr$ be a n.s.f trace on $Z(M)$.
$M'\surl{\star}_{Z(M)}M=M'\surl{\otimes}_{Z(M)}M$ acts on
$L^2(M)\surl{\otimes}_{tr}L^2(M)$. There exists a n.s.f
operator-valued weight $T$ from $M$ to $Z(M)$ such that $\nu=tr\circ
T$.

Let $\alpha$ (resp. $\beta$) be the (resp. anti-) representation of
$M$ to $M'\surl{\otimes}_{Z(M)}M$ such that
$\alpha(m)=1\surl{\otimes}_{Z(M)}m$ (resp.
$\beta(m)=j(m)\surl{\otimes}_{Z(M)}1$) where
$j(x)=J_{\nu}x^*J_{\nu}$ for all $x\in \mathcal{L}(L^2_{\nu}(M))$.

\begin{prop}\label{renv2}
The following formula:
$$
\begin{aligned}
I: [L^2(M)\surl{\otimes}_{tr}L^2(M)]\surl{\ _{\beta}
  \otimes_{\alpha}}_{\ \nu}[L^2(M)\surl{\otimes}_{tr}L^2(M)]
  &\rightarrow L^2(M)\surl{\otimes}_{tr}L^2(M)
  \surl{\otimes}_{tr}L^2(M)\\
[\Lambda_{\nu}(y)\surl{\otimes}_{tr}\eta]\surl{\ _{\beta}
  \otimes_{\alpha}}_{\ \nu}\Xi&\mapsto \alpha(y)\Xi\surl{\otimes}_{tr}\eta
\end{aligned}$$
for all $\eta\in L^2(M),\Xi\in L^2(M)\surl{\otimes}_{tr}L^2(M)$ and
$y\in M$, defines a canonical isomorphism such that we have
$I([m\surl{\otimes}_{Z(M)}z]\surl{\ _{\beta}\otimes_{\alpha}}_{\
\nu}Z)=(\alpha(M)Z\surl{\otimes}_{Z(M)}z)I$, for all $m\in M$, $z\in
Z(M)'$ and $Z\in \mathcal{L}(L^2(M))\surl{\star}_{Z(M)}M'$.
\end{prop}

\begin{proof}
Straightforward.
\end{proof}
We identify $(M'\surl{\otimes}_{Z(M)}M)\surl{\
_{\beta}\star_{\alpha}}_{\ M}(M'\surl{\otimes}_{Z(M)}M)$ with
$M'\surl{\otimes}_{Z(M)}Z(M)\surl{\otimes}_{Z(M)}M$ and so with
$M'\surl{\otimes}_{Z(M)}M$. We define a normal *-homomorphism
$\Gamma$ by:
$$
\begin{aligned}
M'\surl{\otimes}_{Z(M)}M &\rightarrow (M'\surl{\otimes}_{Z(M)}M)
\surl{\ _{\beta}\star_{\alpha}}_{\ \nu} (M'\surl{\otimes}_{Z(M)}M)\\
n\surl{\otimes}_{Z(M)}m&\mapsto
I^*(n\surl{\otimes}_{Z(M)}1\surl{\otimes}_{Z(M)}m)I=[1\surl{\otimes}_{Z(M)}m]
\surl{\ _{\beta}\otimes_{\alpha}}_{\ \nu}[n\surl{\otimes}_{Z(M)}1]
\end{aligned}$$
$\Gamma$ is, in fact, the identity trough the previous isomorphism.

\begin{theo}\label{space}
If we put $T_R=id\surl{\star}_{Z(M)}T$ and
$R=\varsigma_{Z(M)}\circ(j\surl{\otimes}_{Z(M)}j)$, then
$(M,M'\surl{\otimes}_{Z(M)}M,\alpha,\beta,\Gamma,\nu,R\circ T_R\circ
R,T_R)$ becomes a adapted measured quantum groupoid w.r.t $\nu$
called \textbf{quantum space quantum groupoid}.
\end{theo}

\begin{proof}
By definition, $\Gamma$ is a morphism of Hopf bimodule. We have to
prove co-product relation. For all $m\in M$ and $n\in M'$, we have:
$$
\begin{aligned}
(\Gamma\surl{\ _{\beta}\star_{\alpha}}_{\
\nu}id)\circ\Gamma(n\surl{\otimes}_{Z(M)}m )
&=[1\surl{\otimes}_{Z(M)}m]\surl{\ _{\beta}\otimes_{\alpha}}_{\
\nu}[1\surl{\otimes}_{Z(M)}1]
\surl{\ _{\beta}\otimes_{\alpha}}_{\ \nu}[n\surl{\otimes}_{Z(M)}1]\\
&=(id\surl{\ _{\beta}\star_{\alpha}}_{\
\nu}\Gamma)\circ\Gamma(n\surl{\otimes}_{Z(M)}m )
\end{aligned}$$
Now, we show that $T_R$ is right invariant and $\alpha$-adapted
w.r.t $\nu$. So, for all $m\in M,n\in M'$ and $\xi\in
D(_{\alpha}(L^2(M)\surl{\otimes}_{tr} L^2(M)),\nu^o)$, we put
$\Psi=\nu\circ\beta^{-1}\circ T_R$ and we compute:
$$
\begin{aligned}
\omega_{\xi}((\Psi\surl{\ _{\beta}\star_{\alpha}}_{\
\nu}id)\Gamma(n\surl{\otimes}_{Z(M)} m))&=\Psi((id\surl{\
_{\beta}\star_{\alpha}}_{\
\nu}\omega_{\xi})([1\surl{\otimes}_{Z(M)}m]\surl{\ _{\beta}
\otimes_{\alpha}}_{\ \nu}[n\surl{\otimes}_{Z(M)}1]))\\
&=\Psi([1\surl{\otimes}_{Z(M)}m]\beta(<[n\surl{\otimes}_{Z(M)}1]\xi,\xi>_{\alpha,\nu}))\\
&=\nu(<[n\surl{\otimes}_{Z(M)}T(m)]\xi,\xi>_{\alpha,\nu})\\
&=\omega_{\xi}(n\surl{\otimes}_{Z(M)}T(m))=\omega_{\xi}(T_R(n\surl{\otimes}_{Z(M)}m))
\end{aligned}$$ Finally, we have
for all $t\in\mathbb{R}$:
$$
\begin{aligned}
\sigma_t^{T_R}&=\sigma_t^{\nu'\surl{\star}_{Z(M)}\nu}\
_{|(M'\surl{\otimes}_{Z(M)} M)\cap
\beta(M)'}=\sigma_t^{\nu'\surl{\star}_{Z(M)}\nu}\
_{|(M'\surl{\otimes}_{Z(M)}
M)\cap (M\surl{\star}_{Z(M)}\mathcal{L}(L^2(M)))}\\
&=\sigma_t^{\nu'\surl{\star}_{Z(M)}\nu}\
_{|Z(M)\surl{\otimes}_{Z(M)}
M}=(id\surl{\otimes}_{Z(M)}\sigma_t^{\nu})\
_{|1\surl{\otimes}_{Z(M)}M}=1\surl{\otimes}_{Z(M)}\sigma_t^{\nu}
\end{aligned}$$ so that
$\sigma_t^{T_R}\circ\alpha(m)=1\surl{\otimes}_{Z(M)}\sigma_t^{\nu}(m)
=\alpha(\sigma_t^{\nu}(m))$ for all $t\in\mathbb{R}$ and $m\in M$.
Since it is easy to see that $R$ is a co-involution, we have done.
\end{proof}

\subsection{Fundamental elements} By \ref{raccourci}, we can compute
the pseudo-multiplicative unitary. Let first notice that
$\Phi=\nu'\surl{\star}_{Z(M)}\nu =\Psi$ so that $\lambda=\delta=1$
and:
$$\alpha=1\surl{\otimes}_{Z(M)}id, \hat{\alpha}=id\surl{\otimes}_{Z(M)}1,
\beta=j\surl{\otimes}_{Z(M)}1\text{ and }
\hat{\beta}=1\surl{\otimes}_{Z(M)}j$$ For example, we have
$D((H\surl{\otimes}_{tr} H)_{\hat{\beta},\nu^o})\supset
H\surl{\otimes}_{tr}
D(H_j,\nu^o)=H\surl{\otimes}_{tr}\Lambda_{\nu}({\mathcal N}_{\nu})$
and for all $\eta\in H$ and $y\in {\mathcal N}_{\nu}$, we have
$R^{\hat{\beta},\nu^o}(\eta\surl{\otimes}_{tr}\Lambda_{\nu}(y))=\lambda^{tr}_{\eta}
R^{j,\nu^o}(\Lambda_{\nu}(y))=\lambda^{tr}_{\eta}y$.

\begin{lemm}\label{simex2}
We have, for all $\eta\in H$ and $e\in {\mathcal N}_{\nu}$:
$$I\rho^{\beta,\alpha}_{\eta\surl{\otimes}_{tr}
J_{\nu}\Lambda_{\nu}(e)}=\lambda^{tr}_{\eta}J_{\nu}eJ_{\nu}
\surl{\otimes}_{Z(M)}1\text{ and }
I\lambda_{\Lambda_{\nu}(y)\otimes\eta}^{\beta,\alpha}=\rho_{\eta}^{tr}(1\surl{\otimes}_{Z(M)}
y)$$
\end{lemm}

\begin{proof}
Straightforward.
\end{proof}

\begin{prop}
We have, for all $\Xi\in H\surl{\otimes}_{tr}H,\eta\in H$ and $m\in
{\mathcal N}_{\nu}$:
$$W^*(\Xi\surl{\ _{\alpha} \otimes_{\hat{\beta}}}_{\
\nu^o}(\eta\surl{\otimes}_{tr}\Lambda_{\nu}(m)))
=I^*(\eta\surl{\otimes}_{tr}(1\surl{\otimes}_{Z(M)}m)\Xi)$$
\end{prop}

\begin{proof}
For all $m,e\in {\mathcal N}_{\nu}$ and $m',e'\in {\mathcal
N}_{\nu'}$, we have by the previous lemma:
$$
\begin{aligned}
I\Gamma(m'\surl{\otimes}_{Z(M)}m)\rho^{\beta,\alpha}_{J_{\nu'}
\Lambda_{\nu'}(e')\surl{\otimes}_{tr}
J_{\nu}\Lambda_{\nu}(e)}&=(m'\surl{\otimes}_{Z(M)}1\surl{\otimes}_{Z(M)}m)I
\rho^{\beta,\alpha}_{J_{\nu'}\Lambda_{\nu'}(e')\surl{\otimes}_{tr}
J_{\nu}\Lambda_{\nu}(e)}\\
&=(m'\otimes 1\otimes
m)\lambda^{tr}_{J_{\nu'}\Lambda_{\nu'}(e')}J_{\nu}eJ_{\nu}
\surl{\otimes}_{Z(M)}1\\
&=\lambda^{tr}_{J_{\nu'}e'J_{\nu'}\Lambda_{\nu'}(m')}J_{\nu}eJ_{\nu}\surl{\otimes}_{Z(M)}m
\end{aligned}$$
On the other hand, we have by \ref{renv2}:
$$
\begin{aligned}
&\ \quad I([1\surl{\otimes}_{Z(M)}1]\surl{\ _{\beta}
\otimes_{\alpha}}_{\ \nu} [J_{\nu'}e'J_{\nu'}\surl{\otimes}_{Z(M)}
J_{\nu}eJ_{\nu}])W^*
\rho^{\alpha,\hat{\beta}}_{\Lambda_{\nu'}(m')\surl{\otimes}_{tr}\Lambda_{\nu'}(m')}\\
&=(J_{\nu'}e'J_{\nu'}\surl{\otimes}_{Z(M)}J_{\nu}eJ_{\nu}\surl{\otimes}_{Z(M)}
1)IW^*\rho^{\alpha,\hat{\beta}}_{\Lambda_{\nu'}(m')\surl{\otimes}_{tr}\Lambda_{\nu'}(m')}
\end{aligned}$$
Then, by \ref{raccourci} and taking the limit over $e$ and $e'$
which go to $1$, we get for all $\Xi\in H\surl{\otimes}_{tr}H$:
$$W^*(\Xi\surl{\ _{\alpha} \otimes_{\hat{\beta}}}_{\
\nu^o}(\Lambda_{\nu'}(m')\surl{\otimes}_{tr}\Lambda_{\nu}(m)))
=I^*(\Lambda_{\nu'}(m')\surl{\otimes}_{tr}(1\surl{\otimes}_{Z(M)}m)\Xi)$$
Now, if $\Xi\in D(_{\alpha}(H\surl{\otimes}_{tr}H),\nu)$, by
continuity and density of $\Lambda_{\nu'}({\mathcal N}_{\nu'})$ we
have for all $\Xi\in D(_{\alpha}(H\surl{\otimes}_{tr}H),\nu)$:
$$W^*(\Xi\surl{\ _{\alpha} \otimes_{\hat{\beta}}}_{\
\nu^o}(\eta\surl{\otimes}_{tr}\Lambda_{\nu}(m)))
=I^*(\eta\surl{\otimes}_{tr}(1\surl{\otimes}_{Z(M)}m)\Xi)$$ Since
$\eta\surl{\otimes}_{tr}\Lambda_{\nu}(m)\in D((H\surl{\otimes}_{tr}
H)_{\hat{\beta},\nu^o})$, the relation holds by continuity for all
$\Xi\in H\surl{\otimes}_{tr}H$.
\end{proof}

\begin{rema}
If $\sigma_{tr}$ is the flip of $L^2(M)\surl{\otimes}_{tr}L^2(M)$,
then $\sigma_{tr}\circ\hat{\beta}=\beta\circ\sigma_{tr}$ and if
$I'=(1\surl{\otimes}_{Z(M)}\sigma_{tr})I(\sigma_{tr}\surl{\
_{\hat{\beta}}\otimes_{\alpha}}_{\ \nu} [1\surl{\otimes}_{Z(M)}
1])\sigma_{\nu^o}$, then $I'$ is the identification:

$$
\begin{aligned}
I': [L^2(M)\surl{\otimes}_{tr}L^2(M)]\surl{\ _{\alpha}
\otimes_{\hat{\beta}}}_{\ \nu^o}[L^2(M)\surl{\otimes}_{tr}
L^2(M)]&\rightarrow L^2(M)\surl{\otimes}_{tr}L^2(M)
\surl{\otimes}_{tr}L^2(M)\\
\Xi\surl{\ _{\beta}\otimes_{\alpha}}_{\
\nu}[\eta\surl{\otimes}_{tr}\Lambda_{\nu}(y)]&\mapsto\eta\surl{\otimes}_{tr}\alpha(y)\Xi
\end{aligned}$$
for all $\eta\in L^2(M),\Xi\in L^2(M)\surl{\otimes}_{tr}L^2(M)$ and
$y\in M$. Consequently, by the previous proposition $W^*=I^*I'$.
\end{rema}

\begin{coro}
We can reconstruct the von Neumann algebra thanks to $W$:
$$M'\surl{\otimes}_{Z(M)}
M=<(id*\omega_{\xi,\eta})(W^*)\,|\,\xi\in
D((H\surl{\otimes}_{tr}H)_{\hat{\beta}},\nu^o),\eta\in
D(_{\alpha}(H\surl{\otimes}_{tr}H),\nu)>^{-\textsc{w}}$$
\end{coro}

\begin{proof}
By \ref{appartenance}, we know that:
$$<(id*\omega_{\xi,\eta})(W^*)\,|\,\xi\in
D((H\surl{\otimes}_{tr}H)_{\hat{\beta}},\nu^o,\eta\in
D(_{\alpha}(H\surl{\otimes}_{tr}H ),\nu)>^{-\textsc{w}}\subset
M'\surl{\otimes}_{Z(M)}M$$ Let $\eta,\xi\in H$ and $m,e\in {\mathcal
N}_{\nu}$. Then, for all $\Xi_1,\Xi_2\in H\surl{\otimes}_{tr}H$, we
have by \ref{simex2}:
$$
\begin{aligned}
&\ \quad
((id*\omega_{\eta\surl{\otimes}_{tr}\Lambda_{\nu}(m),\xi\surl{\otimes}_{tr}
J_{\nu}\Lambda_\nu(e)})(W^*)\Xi_1|\Xi_2)\\
&=(W^*(\Xi_1\surl{\ _{\alpha} \otimes_{\hat{\beta}}}_{\ \nu^o}
[\eta\surl{\otimes}_{tr}\Lambda_{\nu}(m)])|\Xi_2\surl{\ _{\beta}
\otimes_{\alpha}}_{\ \nu}[\xi\surl{\otimes}_{tr}J_{\nu}\Lambda_\nu(e)])\\
&=(I^*(\eta\surl{\otimes}_{tr}(1\surl{\otimes}_{Z(M)}m)\Xi_1)|\Xi_2\surl{\
_{\beta}\otimes_{\alpha}}_{\ \nu}[\xi\surl{\otimes}_{tr}J_{\nu}\Lambda_\nu(e)])\\
&=(\eta\surl{\otimes}_{tr}(1\surl{\otimes}_{Z(M)}m)\Xi_1|\xi\surl{\otimes}_{tr}
(J_{\nu}eJ_{\nu}\surl{\otimes}_{Z(M)}1)\Xi_2)\\
&=((<\eta,\xi>_{tr}J_{\nu}e^*J_{\nu}\otimes m)\Xi_1|\Xi_2)
\end{aligned}$$
Consequently, we get the reverse inclusion thanks to the relation:
$$(id*\omega_{\eta\surl{\otimes}_{tr}\Lambda_{\nu}(m),\xi\surl{\otimes}_{tr}
J_{\nu}\Lambda_\nu(e)})(W^*)=<\eta,\xi>_{tr}J_{\nu}e^*J_{\nu}\surl{\otimes}_{Z(M)}
m$$
\end{proof}

Now, we compute $G$ so as to get the antipode.

\begin{prop}
If $F_{\nu}=S_{\nu}^*$ comes from Tomita's theory, then we have:
$$G=\sigma_{tr}\circ (F_{\nu}\surl{\otimes}_{tr}F_{\nu})$$
\end{prop}

\begin{proof}
Let $a=J_{\nu}a_1J_{\nu}\surl{\otimes}_{Z(M)}
a_2,b=J_{\nu}b_1J_{\nu}\surl{\otimes}_{Z(M)}
b_2,c=J_{\nu}c_1J_{\nu}\surl{\otimes}_{Z(M)} c_2$ and
$d=J_{\nu}d_1J_{\nu}\surl{\otimes}_{Z(M)} d_2$ be elements of
$M'\surl{\otimes}_{Z(M)} M$ analytic w.r.t
$\nu'\surl{\star}_{Z(M)}\nu$. Then, by \ref{simex2}, the value of
$(\lambda^{\beta,\alpha}_{\Lambda_{\nu}(\sigma_{i/2}^{\nu}(b_1))\surl{\otimes}_{tr}
\Lambda_{\nu}(\sigma_{-i}^{\nu}(b_2^*))})^*W^*$ on
$$[\Lambda_{\nu'}(J_{\nu}a_1J_{\nu})
\surl{\otimes}_{tr}\Lambda_{\nu}(a_2)]\surl{\ _{\alpha}
\otimes_{\hat{\beta}}}_{\
\nu^o}[\Lambda_{\nu'}(J_{\nu}d_1^*c_1^*J_{\nu})
\surl{\otimes}_{tr}\Lambda_{\nu}(d_2^*c_2^*)]$$ is equal to:
$$
\begin{aligned}
&\ \quad
(\lambda^{\beta,\alpha}_{\Lambda_{\nu}(\sigma_{i/2}^{\nu}(b_1))\surl{\otimes}_{tr}
\Lambda_{\nu}(\sigma_{-i}^{\nu}(b_2^*))})^*I^*(\Lambda_{\nu'}(J_{\nu}d_1^*c_1^*J_{\nu})
\surl{\otimes}_{tr}\Lambda_{\nu'}(J_{\nu}a_1J_{\nu})\surl{\otimes}_{tr}\Lambda_{\nu}(d_2^*c_2^*a_2))\\
&=\!\left[\rho^{tr}_{\Lambda_{\nu}(\sigma_{-i}^{\nu}(b_2^*))}(1\surl{\otimes}_{Z(M)}
\sigma_{i/2}^{\nu}(b_1))\right]^*\!\!\!\!(\Lambda_{\nu'}(J_{\nu}d_1^*c_1^*J_{\nu})
\surl{\otimes}_{tr}\Lambda_{\nu'}(J_{\nu}a_1J_{\nu})\surl{\otimes}_{tr}\Lambda_{\nu}(d_2^*c_2^*a_2))\\
&=<d_2^*c_2^*\Lambda_{\nu}(a_2),\Lambda_{\nu}(\sigma_{-i}^{\nu}(b_2^*))>_{tr}\,
\Lambda_{\nu'}(J_{\nu}d_1^*c_1^*J_{\nu})\surl{\otimes}_{tr}\sigma_{-i/2}^{\nu}(b_1^*)
\Lambda_{\nu'}(J_{\nu}a_1J_{\nu})\\
&=<\Lambda_{\nu}(a_2b_2),\Lambda_{\nu}(c_2d_2)>_{tr}\,
J_{\nu}\Lambda_{\nu}(d_1^*c_1^*)\surl{\otimes}_{tr}
J_{\nu}\Lambda_{\nu}(a_1b_1)
\end{aligned}$$
Consequently, by definition of $G$:
$$G\left[<\Lambda_{\nu}(a_2b_2),\Lambda_{\nu}(c_2d_2)>_{tr}\,
J_{\nu}\Lambda_{\nu}(d_1^*c_1^*)\surl{\otimes}_{tr}
J_{\nu}\Lambda_{\nu}(a_1b_1)\right]$$ is equal to the value of
$G(\lambda^{\beta,\alpha}_{\Lambda_{\nu}(\sigma_{i/2}^{\nu}(b_1))\surl{\otimes}_{tr}
\Lambda_{\nu}(\sigma_{-i}^{\nu}(b_2^*))})^*W^*$ on:
$$[\Lambda_{\nu'}(J_{\nu}a_1J_{\nu})
\surl{\otimes}_{tr}\Lambda_{\nu}(a_2)]\surl{\ _{\alpha}
\otimes_{\hat{\beta}}}_{\
\nu^o}[\Lambda_{\nu'}(J_{\nu}d_1^*c_1^*J_{\nu})
\surl{\otimes}_{tr}\Lambda_{\nu}(d_2^*c_2^*)]$$ which is equal to
the value of
$(\lambda^{\beta,\alpha}_{\Lambda_{\nu}(\sigma_{i/2}^{\nu}(d_1))\surl{\otimes}_{tr}
\Lambda_{\nu}(\sigma_{-i}^{\nu}(d_2^*))})^*W^*$ on:
$$[\Lambda_{\nu'}(J_{\nu}c_1J_{\nu})
\surl{\otimes}_{tr}\Lambda_{\nu}(c_2)]\surl{\ _{\alpha}
\otimes_{\hat{\beta}}}_{\
\nu^o}[\Lambda_{\nu'}(J_{\nu}b_1^*a_1^*J_{\nu})
\surl{\otimes}_{tr}\Lambda_{\nu}(b_2^*a_2^*)]$$ This last vector is
$<\Lambda_{\nu}(c_2d_2),\Lambda_{\nu}(a_2b_2)>_{tr}\,
J_{\nu}\Lambda_{\nu}(b_1^*a_1^*)\surl{\otimes}_{tr}
J_{\nu}\Lambda_{\nu}(c_1d_1)$. Since $G$ is closed, we get:
$$G\left[J_{\nu}\Lambda_{\nu}(d_1^*c_1^*)\surl{\otimes}_{tr}
J_{\nu}\Lambda_{\nu}(a_1b_1)\right]=\left[
J_{\nu}\Lambda_{\nu}(b_1^*a_1^*)\surl{\otimes}_{tr}
J_{\nu}\Lambda_{\nu}(c_1d_1)\right]$$ so that $G$ coincides with
$\sigma_{tr}(F_{\nu}\surl{\otimes}_{tr}F_{\nu})$.
\end{proof}

The polar decomposition of $G=ID^{1/2}$ is such that
$D=\Delta_{\nu}^{-1}\surl{\otimes}_{tr}\Delta_{\nu}^{-1}$ and
$I=\sigma_{tr}(J_{\nu}\surl{\otimes}_{tr} J_{\nu})$ so that the
scaling group is
$\tau_t=\sigma_{-t}^{\nu'}\surl{\star}_{Z(M)}\sigma_t^{\nu}$ for all
$t\in\mathbb{R}$ and the unitary antipode is
$R=\varsigma_{Z(M)}\circ (j\surl{\otimes}_{Z(M)}j)$. We also notice
that $\nu'\surl{\star}_{Z(M)}\nu$ is $\tau$-invariant.

\begin{rema}
If $M$ is the commutative von Neumann algebra $L^{\infty}(X)$, then
the structure coincides with the quantum space $X$.
\end{rema}

\subsection{Dual structure} Here we compute the dual structure and we
observe that this is not of adapted measured quantum groupoid's
type.

\begin{prop}\label{cru2}
For all $e,y\in {\mathcal N}_{\nu}$ and $\eta,\zeta\in H$, we have:
$$(\omega_{\Lambda_{\nu}(y)\surl{\otimes}_{tr}\eta,\zeta\surl{\otimes}_{tr}
J_{\nu}\Lambda_{\nu}(e)}* id)(W)=1\surl{\otimes}_{Z(M)}
J_{\nu}e^*J_{\nu}(\rho^{tr}_{\zeta})^*\sigma_{tr}\rho^{tr}_{\eta}y$$
\end{prop}

\begin{proof}
For all $\Xi\in H\surl{\otimes}_{tr}H,\xi\in H$ and $m\in {\mathcal
N}_{\nu}$, we have:
$$
\begin{aligned}
&\ \quad
((\omega_{\Lambda_{\nu}(y)\surl{\otimes}_{tr}\eta,\zeta\surl{\otimes}_{tr}
J_{\nu}\Lambda_{\nu}(e)}* id)(W)\Xi|\xi\surl{\otimes}_{tr}
\Lambda_{\nu}(m))\\&=([\Lambda_{\nu}(y)\surl{\otimes}_{tr}\eta]\surl{\
_{\beta} \otimes_{\alpha}}_{\ \nu}\Xi|W^*([\zeta\surl{\otimes}_{tr}
J_{\nu}\Lambda_{\nu}(e)]\surl{\ _{\alpha}
\otimes_{\hat{\beta}}}_{\ \nu^o}[\xi\surl{\otimes}_{tr}\Lambda_{\nu}(m)]))\\
&=((1\surl{\otimes}_{Z(M)}y)\Xi\surl{\otimes}_{tr}\eta
|\xi\surl{\otimes}_{tr}\zeta\surl{\otimes}_{tr}
mJ_{\nu}\Lambda_{\nu}(e))=(\Xi\surl{\otimes}_{tr}\eta
|\xi\surl{\otimes}_{tr}
y^*\zeta\surl{\otimes}_{tr}J_{\nu}eJ_{\nu}\Lambda_{\nu}(m))\\
&=((1\surl{\otimes}_{Z(M)}\rho^{tr}_{\eta})\Xi
|(1\surl{\otimes}_{Z(M)}\sigma_{tr}\rho^{tr}_{y^*\zeta}
J_{\nu}eJ_{\nu})(\xi\surl{\otimes}_{tr}\Lambda_{\nu}(m)))\\
&=((1\surl{\otimes}_{Z(M)}
J_{\nu}e^*J_{\nu}(\rho^{tr}_{\zeta})^*\sigma_{tr}\rho^{tr}_{\eta}y)\Xi
|\xi\surl{\otimes}_{tr}\Lambda_{\nu}(m))
\end{aligned}$$
\end{proof}

\begin{coro} We have $\widehat{M'\surl{\otimes}_{Z(M)}M}=1\surl{\otimes}_{Z(M)}Z(M)'$
which is identified with $Z(M)'$.
\end{coro}

\begin{proof}
We already know that $\alpha(M)\cup\hat{\beta}(M)\subset
\widehat{M'\surl{\otimes}_{Z(M)}M}$ so that
$1\surl{\otimes}_{Z(M)}Z(M)'\subset\widehat{M'\surl{\otimes}_{Z(M)}M}$.
The reverse inclusion comes from the previous proposition..
\end{proof}

With this identification between $1\surl{\otimes}_{Z(M)}Z(M)'$ and
$Z(M)'$, the dual structure admits $M$ for basis, $id$ for
representation and $j$ for anti-representation. The dual co-product
necessarily satisfies $\widehat{\Gamma}(mn)=m\surl{\ _{j}
\otimes_{id}}_{\ \nu}n$ for all $m\in M$ and $n\in M'$. If $I_{\nu}$
is the canonical isomorphism from $L^2(M)\surl{\
_{j}\otimes_{id}}_{\ \nu}L^2(M)$ onto $L^2(M)$ given by
$I_{\nu}(\Lambda_{\nu}(x)\surl{\ _{j}\otimes_{id}}_{\
\nu}\eta)=\alpha(x)\eta$ for all $x\in {\mathcal N}_{\nu}$ and
$\eta\in L^2(M)$, then we have $I_{\nu}(m\surl{\ _{\beta}
\otimes_{\alpha}}_{\ \nu}n)=mnI_{\nu}$ and we can identify the von
Neumann algebra $M'\surl{\ _{\beta}\star_{\alpha}}_{\ M}M$ with
$Z(M)$ and the von Neumann algebra $Z(M)'\surl{\
_{\beta}\star_{\alpha}}_{\ M}Z(M)'$ with $Z(M)'$. The dual
co-product is then identity through this identification.

\begin{lemm}
We have
$\hat{\Lambda}((\omega_{\Xi,\Lambda_{\nu}(m)\surl{\otimes}_{tr}
J_{\nu}\Lambda_{\nu}(e)}*id)(W))=(m^*\surl{\otimes}_{Z(M)}
J_{\nu}e^*J_{\nu})\Xi$ for all $m,e\in {\mathcal N}_{\nu}$ et
$\Xi\in D((H\surl{\otimes}_{tr} H)_{\beta},\nu^o)$.
\end{lemm}

\begin{proof}
Let $m_1,m_2\in {\mathcal N}_{\nu}$. Then, we have:
$$
\begin{aligned}
&\ \quad
(\hat{\Lambda}((\omega_{\Xi,\Lambda_{\nu}(m)\surl{\otimes}_{tr}
J_{\nu}\Lambda_{\nu}(e)}*id)(W))|\Lambda_{\nu'}(J_{\nu}m_1J_{\nu})\surl{\otimes}_{Z(M)}
\Lambda_{\nu}(m_2))\\
&=\omega_{\Xi,\Lambda_{\nu}(m)\surl{\otimes}_{tr}
J_{\nu}\Lambda_{\nu}(e)}(J_{\nu}m_1^*J_{\nu}\surl{\otimes}_{Z(M)}
m_2^*)\\
&=((J_{\nu}m_1^*J_{\nu}\surl{\otimes}_{Z(M)} m_2^*)\Xi
|\Lambda_{\nu}(m)\surl{\otimes}_{tr}
J_{\nu}\Lambda_{\nu}(e))\\
&=(\Xi |mJ_{\nu}\Lambda_{\nu}(m_1)\surl{\otimes}_{tr}
J_{\nu}eJ_{\nu}\Lambda_{\nu}(m_2))\\
&=((m^*\surl{\otimes}_{Z(M)}J_{\nu}e^*J_{\nu})\Xi
|\Lambda_{\nu'}(J_{\nu}m_1J_{\nu})\surl{\otimes}_{tr}
\Lambda_{\nu}(m_2))
\end{aligned}$$

\end{proof}

\begin{prop}
The dual operator-valued weight $\widehat{T_R}$ coincide with
$T^{-1}$ in sense of proposition 12.11 of \cite{St}. Also, the dual
operator-valued weight $\widehat{T_L}$ coincide with $j\circ
T^{-1}\circ j$.
\end{prop}

\begin{proof}
Via the identification between $\widehat{M'\surl{\otimes}_{Z(M)}M}$
and $Z(M)'$, we have, by proposition \ref{cru2}:
$$(\omega_{\Xi,\Lambda_{\nu}(m)\surl{\otimes}_{tr}
J_{\nu}\Lambda_{\nu}(e)}*id)(W)=
J_{\nu}e^*J_{\nu}[(\rho_{\zeta}^{tr})^*\sigma_{tr}\rho_{\eta}^{tr}]y$$

Let $m,e,y\in {\mathcal N}_{\nu}$ and $\eta\in H$. On one hand, we
compute:
$$
\begin{aligned}
&\ \quad
||\hat{\Lambda}((\omega_{\Lambda_{\nu}(y)\surl{\otimes}_{tr}\eta,\Lambda_{\nu}(m)\surl{\otimes}_{tr}
J_{\nu}\Lambda_{\nu}(e)}*id)(W))||^2=||m^*\Lambda_{\nu}(y)\surl{\otimes}_{tr}
J_{\nu}e^*J_{\nu}\eta ||^2\\
&=(<J_{\nu}e^*J_{\nu}\eta,J_{\nu}e^*J_{\nu}\eta>_{tr}\Lambda_{\nu}(m^*y)|\Lambda_{\nu}(m^*y))
\end{aligned}$$
On the other hand, we have:
$$
\begin{aligned}
&\ \quad
||\hat{\Lambda}((\omega_{\Lambda_{\nu}(y)\surl{\otimes}_{tr}\eta,\Lambda_{\nu}(m\surl{\otimes}_{tr}
J_{\nu}\Lambda_{\nu}(e)}*id)(W)))||^2\\
&=\hat{\Phi}((\rho^{tr}_{J_{\nu}e^*J_{\nu}\eta})^*\sigma_{tr}\rho^{tr}_{\Lambda_{\nu}(y^*m)}
(\rho^{tr}_{\Lambda_{\nu}(y^*m)})^*\sigma_{tr}\rho^{tr}_{J_{\nu}e^*J_{\nu}\eta})\\
&=\hat{\Phi}((\rho^{tr}_{J_{\nu}e^*J_{\nu}\eta})^*[\theta^{tr}(\Lambda_{\nu}(y^*m),
\Lambda_{\nu}(y^*m))\surl{\otimes}_{Z(M)}1]\rho^{tr}_{J_{\nu}e^*J_{\nu}\eta})\\
&=\hat{\Phi}(<J_{\nu}e^*J_{\nu}\eta,J_{\nu}e^*J_{\nu}\eta>_{tr}\theta^{tr}(\Lambda_{\nu}(y^*m),
\Lambda_{\nu}(y^*m)))
\end{aligned}$$
Then we conclude that, for all $m,y\in {\mathcal N}_{\nu}$, we have:
$$
\begin{aligned}
&\ \quad\hat{\Phi}(\theta^{tr}(\Lambda_{\nu}(y^*m),
\Lambda_{\nu}(y^*m)))=||\Lambda_{\nu}(m^*y)||^2=||\Delta_{\nu}^{-1/2}J_{\nu}\Lambda_{\nu}(y^*m)||^2\\
&=\nu'(\theta^{\nu}(J_{\nu}\Lambda_{\nu}(y^*m),J_{\nu}\Lambda_{\nu}(y^*m)))
=\nu'\circ T^{-1}(\theta^{tr}(J_{\nu}\Lambda_{\nu}(y^*m),J_{\nu}\Lambda_{\nu}(y^*m)))\\
&=\nu\circ j\circ T^{-1}\circ j(\theta^{tr}(\Lambda_{\nu}(y^*m)),\Lambda_{\nu}(y^*m)))\\
\end{aligned}$$
Therefore $\widehat{T_L}=j\circ T^{-1}\circ j$ and we get the
proposition.
\end{proof}

\begin{prop}
The dual quantum space quantum groupoid can be identify with
$(M,Z(M)',id,j,\nu,id,j\circ T^{-1}\circ j,T^{-1})$ which is a
measured quantum groupoid but not a adapted measured quantum
groupoid. Moreover, expressions for co-involution and scaling group
are given, for all $x\in Z(M)'$ and $t\in\mathbb{R}$:
$$\hat{R}(x)=J_{\nu}x^*J_{\nu}\quad\text{ and
}\quad\hat{\tau}_t(x)=\Delta_{\nu}^{it}x\Delta_{\nu}^{-it}$$
\end{prop}

\begin{proof}
The proposition gathers results of the section. Nevertheless we lay
stress on the following point. We have, for all $t\in\mathbb{R}$ and
$m\in M$:
$$\sigma_{t}^{T^{-1}}(m)=\sigma_{-t}^T(m)=\sigma_{-t}^{\nu}(m)$$
instead of $\sigma_t^{\nu}(m)$ to have a adapted measured quantum
groupoid.
\end{proof}

\begin{rema}
If $M$ is a factor, $\mathcal{L}(H)$ is the von Neumann algebra
underlying the structure of quantum space quantum groupoid whereas
$M'\otimes M$ is the underlying von Neumann algebra of the dual
structure. In general, they are not isomorphic. Nevertheless, if $M$
is abelian or if $M$ is a type $I$ factor (and henceforth a sum of
type $I$ factors cf. paragraph \ref{ogqm}), the structure is
self-dual. In the abelian case $M=L^{\infty}(X)$, we recover the
space groupoid $X$. This example comes from the inclusion of von
Neumann algebras (\cite{E1}):
$$Z(M)\subset M \subset Z(M)'\subset\ldots$$
\end{rema}

\section{Pairs quantum groupoid}\label{pqg}
\subsection{Definition} Let $M$ be a von Neumann algebra. $M$ acts on
$H=L^2(M)=L^2_{\nu}(M)$ where $\nu$ is a n.s.f weight on $M$. We
denote by $M'$ the commutant of $M$ in $\mathcal{L}(L^2(M))$.
$M'\otimes M$ acts on $L^2(M)\otimes L^2(M)$.

Let $\alpha$ (resp. $\beta$) be the (resp. anti-) representation of
$M$ to $M'\otimes M$ such that $\alpha(m)=1\otimes m$ (resp.
$\beta(m)=j(m)\otimes 1$) where $j(x)=J_{\nu}x^*J_{\nu}$ for all
$x\in \mathcal{L}(L^2_{\nu}(M))$.

\begin{prop}\label{renv}
The following formula:
$$
\begin{aligned}
I: [L^2(M)\otimes L^2(M)]\surl{\ _{\beta}
  \otimes_{\alpha}}_{\ \nu}[L^2(M)\otimes L^2(M)]&\rightarrow L^2(M)\otimes L^2(M)
  \otimes L^2(M)\\
[\Lambda_{\nu}(y)\otimes\eta]\surl{\ _{\beta}
  \otimes_{\alpha}}_{\ \nu}\Xi&\mapsto \alpha(y)\Xi\otimes\eta
\end{aligned}$$
for all $\eta\in L^2(M),\Xi\in L^2(M)\otimes L^2(M)$ and $y\in M$,
defines a canonical isomorphism such that we have $I([m\otimes
x]\surl{\ _{\beta}\otimes_{\alpha}}_{\ \nu}[y\otimes n])=(y\otimes
mn\otimes x)I$, for all $m\in M, n\in M'$ and $x,y\in
\mathcal{L}(L^2(M))$.
\end{prop}

\begin{proof}
Straightforward.
\end{proof}
Then, we can identify $(M'\otimes M)\surl{\
_{\beta}\star_{\alpha}}_{\ M}(M'\otimes M)$ with $M'\otimes
Z(M)\otimes M$. We define a normal *-homomorphism $\Gamma$ by:
$$
\begin{aligned}
M'\otimes M &\rightarrow (M'\otimes M)
\surl{\ _{\beta}\star_{\alpha}}_{\ \nu} (M'\otimes M)\\
n\otimes m&\mapsto I^*(n\otimes 1\otimes m)I=[1\otimes m] \surl{\
_{\beta}\otimes_{\alpha}}_{\ \nu}[n\otimes 1]
\end{aligned}$$

\begin{theo}\label{pair}
$(M,M'\surl{\otimes}_{Z(M)}M,\alpha,\beta,\Gamma,\nu,\nu'\otimes
id,id\otimes\nu)$ is a adapted measured quantum groupoid w.r.t $\nu$
called \textbf{pairs quantum groupoid}.
\end{theo}

\begin{proof}
By definition, $\Gamma$ is a morphism of Hopf bimodule. We have to
prove co-product relation. For all $m\in M$ and $n\in M'$, we have:
$$
\begin{aligned}
(\Gamma\surl{\ _{\beta}\star_{\alpha}}_{\
\nu}id)\circ\Gamma(n\otimes m )&=[1\otimes m]\surl{\
_{\beta}\otimes_{\alpha}}_{\ \nu}[1\otimes
1]\surl{\ _{\beta}\otimes_{\alpha}}_{\ \nu}[n\otimes 1]\\
&=(id\surl{\ _{\beta}\star_{\alpha}}_{\
\nu}\Gamma)\circ\Gamma(n\otimes m)
\end{aligned}$$
$R=\varsigma\circ (\beta_{\nu}\otimes\beta_{\nu})$, where
$\varsigma: M'\otimes M\rightarrow M\otimes M'$ is the flip, is a
co-involution so it is sufficient to show that $T_L=\nu'\otimes id$
is left invariant and $\beta$-adapted w.r.t $\nu$. Let $m\in M,n\in
M'$ and $\xi\in D((L^2(M)\otimes L^2(M))_{\beta,\nu^o})$. We put
$\Phi=\nu\circ\alpha^{-1}\circ T_L$ and we compute:
$$
\begin{aligned}
\omega_{\xi}((id\surl{\ _{\beta}\star_{\alpha}}_{\
\nu}\Phi)\Gamma(n\otimes m))&=\Phi((\omega_{\xi}\surl{\
_{\beta}\star_{\alpha}}_{\ \nu}id)([1\otimes m]\surl{\ _{\beta}
  \otimes_{\alpha}}_{\ \nu}[n\otimes 1]))\\
&=\Phi([n\otimes 1]\alpha(<[1\otimes
m]\xi,\xi>_{\beta,\nu^o}))\\
&=\nu'(n)\nu(<[1\otimes
m]\xi,\xi>_{\beta,\nu^o})\\
&=\nu'(n)\omega_{\xi}(1\otimes m)=\omega_{\xi}(T_L(n\otimes m))
\end{aligned}$$ Finally, we prove that
$T_R=R\circ T_L\circ R=id\otimes \nu$ is $\alpha$-adapted w.r.t
$\nu$. For all $t\in\mathbb{R}$, we have:
$$\sigma_t^{T_R}=\sigma_t^{\nu'\otimes\nu}\ _{|(M'\otimes
M)\cap \beta(M)'}=\sigma_t^{\nu'\otimes\nu}\ _{|Z(M)\otimes
M}=id\otimes\sigma_t^{\nu}\ _{|Z(M)\otimes M}$$ so that we have for
all $t\in\mathbb{R}$ and $m\in M$:
$$\sigma_t^{T_R}\circ\alpha(m)=1\otimes\sigma_t^{\nu}(m)=\alpha(\sigma_t^{\nu}(m))$$
\end{proof}

\begin{rema}
If $M=L^{\infty}(X)$, we find the structure of pairs groupoid
$X\times X$.
\end{rema}

\subsection{Fundamental elements}

By \ref{raccourci}, we can compute the pseudo-multiplicative
unitary. Let first notice that $\Phi=\nu'\otimes\nu =\Psi$ so that
$\lambda=\delta=1$ and:
$$\alpha=1\otimes id, \hat{\alpha}=id\otimes 1, \beta=\beta_{\nu}\otimes 1
\text{ and } \hat{\beta}=1\otimes \beta_{\nu}$$ For example, we have
$D((H\otimes H)_{\hat{\beta},\nu^o})\supset H\otimes
D(H_{\beta_{\nu}},\nu^o)=H\otimes \Lambda_{\nu}({\mathcal N}_{\nu})$
and for all $\eta\in H$ and $y\in {\mathcal N}_{\nu}$, we have
$R^{\hat{\beta},\nu^o}(\eta\otimes\Lambda_{\nu}(y))=\lambda_{\eta}
R^{\beta_{\nu},\nu^o}(\Lambda_{\nu}(y))=\lambda_{\eta}y$.

\begin{lemm}\label{simex}
We have, for all $\eta\in H$ and $e\in {\mathcal N}_{\nu}$:
$$I\rho^{\beta,\alpha}_{\eta\otimes J_{\nu}\Lambda_{\nu}(e)}=
\lambda_{\eta}J_{\nu}eJ_{\nu} \otimes 1\text{ and }
I\lambda_{\Lambda_{\nu}(y)\otimes\eta}^{\beta,\alpha}=\rho_{\eta}(1\otimes
y)$$
\end{lemm}

\begin{proof}
Straightforward.
\end{proof}

\begin{prop}
We have, for all $\Xi\in H\otimes H,\eta\in H$ and $m\in {\mathcal
N}_{\nu}$:
$$W^*(\Xi\surl{\ _{\alpha} \otimes_{\hat{\beta}}}_{\
\nu^o}(\eta\otimes\Lambda_{\nu}(m))) =I^*(\eta\otimes(1\otimes
m)\Xi)$$
\end{prop}

\begin{proof}
For all $m,e\in {\mathcal N}_{\nu}$ and $m',e'\in {\mathcal
N}_{\nu'}$, we have by the previous lemma:
$$
\begin{aligned}
I\Gamma(m'\otimes
m)\rho^{\beta,\alpha}_{J_{\nu'}\Lambda_{\nu'}(e')\otimes
J_{\nu}\Lambda_{\nu}(e)}&=(m'\otimes 1\otimes m)I
\rho^{\beta,\alpha}_{J_{\nu'}\Lambda_{\nu'}(e')\otimes
J_{\nu}\Lambda_{\nu}(e)}\\
&=(m'\otimes 1\otimes
m)\lambda_{J_{\nu'}\Lambda_{\nu'}(e')}J_{\nu}eJ_{\nu}\otimes 1\\
&=\lambda_{J_{\nu'}e'J_{\nu'}\Lambda_{\nu'}(m')}
J_{\nu}eJ_{\nu}\otimes m
\end{aligned}$$
On the other hand, we have by \ref{renv}:
$$
\begin{aligned}
&\ \quad I([1\otimes 1]\surl{\ _{\beta} \otimes_{\alpha}}_{\ \nu}
[J_{\nu'}e'J_{\nu'}\otimes J_{\nu}eJ_{\nu}])W^*
\rho^{\alpha,\hat{\beta}}_{\Lambda_{\nu'}(m')\otimes\Lambda_{\nu'}(m')}\\
&=(J_{\nu'}e'J_{\nu'}\otimes J_{\nu}eJ_{\nu}\otimes 1)
IW^*\rho^{\alpha,\hat{\beta}}_{\Lambda_{\nu'}(m')\otimes\Lambda_{\nu'}(m')}
\end{aligned}$$
Then by \ref{raccourci} and taking the limit over $e$ and $e'$ which
go to $1$, we get for all $\Xi\in H\otimes H$:
$$W^*(\Xi\surl{\ _{\alpha} \otimes_{\hat{\beta}}}_{\
\nu^o}(\Lambda_{\nu'}(m')\otimes\Lambda_{\nu}(m)))
=I^*(\Lambda_{\nu'}(m')\otimes(1\otimes m)\Xi)$$ Now, if $\Xi\in
D(_{\alpha}(H\otimes H),\nu)$, by continuity and density of
$\Lambda_{\nu'}({\mathcal N}_{\nu'})$, we have for all $\Xi\in
D(_{\alpha}(H\otimes H),\nu)$:
$$W^*(\Xi\surl{\ _{\alpha} \otimes_{\hat{\beta}}}_{\
\nu^o}(\eta\otimes\Lambda_{\nu}(m))) =I^*(\eta\otimes(1\otimes
m)\Xi)$$ Since $\eta\otimes\Lambda_{\nu}(m)\in D((H\otimes
H)_{\hat{\beta},\nu^o})$, the previous relation holds by continuity
for all $\Xi\in H\otimes H$.
\end{proof}

\begin{rema}
If $\sigma$ denotes the flip of $L^2(M)\otimes L^2(M)$, then
$\sigma\circ\hat{\beta}=\beta\circ\sigma$ and if
$I'=(1\otimes\sigma)I (\sigma\surl{\ _{\hat{\beta}}
\otimes_{\alpha}}_{\ \nu} [1\otimes 1])\sigma_{\nu^o}$, then $I'$ is
the identification:
$$
\begin{aligned}
I': [L^2(M)\otimes L^2(M)]\surl{\ _{\alpha}
\otimes_{\hat{\beta}}}_{\ \nu^o}[L^2(M)\otimes L^2(M)]&\rightarrow
L^2(M)\otimes L^2(M)
  \otimes L^2(M)\\
\Xi\surl{\ _{\beta}
  \otimes_{\alpha}}_{\ \nu}[\eta\otimes\Lambda_{\nu}(y)]&\mapsto
  \eta\otimes\alpha(y)\Xi
\end{aligned}$$
for all $\eta\in L^2(M),\Xi\in L^2(M)\otimes L^2(M)$ and $y\in M$.
Consequently, by the previous proposition $W^*=I^*I'$.
\end{rema}

\begin{coro}
We can re-construct the underlying von Neumann algebra thanks to
$W$:
$$M'\otimes M=<(id*\omega_{\xi,\eta})(W^*)\,|\,\xi\in D((H\otimes
H)_{\hat{\beta}},\nu^o),\eta\in D(_{\alpha}(H\otimes
H),\nu)>^{-\textsc{w}}$$
\end{coro}

\begin{proof}
By \ref{appartenance}, we know that:
$$<(id*\omega_{\xi,\eta})(W^*)\,|\,\xi\in
D((H\otimes H)_{\hat{\beta}},\nu^o,\eta\in D(_{\alpha}(H\otimes H
),\nu)>^{-w}\subset M'\otimes M$$ Let $\eta,\xi\in H$ and $m,e\in
{\mathcal N}_{\nu}$. Then, for all $\Xi_1,\Xi_2\in H\otimes H$, we
have, by \ref{simex}:
$$
\begin{aligned}
&\ \quad((id*\omega_{\eta\otimes\Lambda_{\nu}(m),\xi\otimes
J_{\nu}\Lambda_\nu(e)})(W^*)\Xi_1|\Xi_2)\\
&=(W^*(\Xi_1\surl{\ _{\alpha} \otimes_{\hat{\beta}}}_{\ \nu^o}
\eta\otimes\Lambda_{\nu}(m))|\Xi_2\surl{\ _{\beta}
\otimes_{\alpha}}_{\ \nu}\xi\otimes J_{\nu}\Lambda_\nu(e))\\
&=(I^*(\eta\otimes(1\otimes m)\Xi_1)|\Xi_2\surl{\ _{\beta}
\otimes_{\alpha}}_{\ \nu}\xi\otimes J_{\nu}\Lambda_\nu(e))\\
&=(\eta\otimes(1\otimes m)\Xi_1|\xi\otimes (J_{\nu}eJ_{\nu}\otimes
1)\Xi_2)\\
&=(\eta|\xi)((J_{\nu}e^*J_{\nu}\otimes m)\Xi_1|\Xi_2)
\end{aligned}$$
Consequently, we get the reverse inclusion thanks to the relation:
\begin{equation}\label{equ1}
(id*\omega_{\eta\otimes\Lambda_{\nu}(m),\xi\otimes
J_{\nu}\Lambda_\nu(e)})(W^*)=(\eta|\xi)(J_{\nu}e^*J_{\nu}\otimes m)
\end{equation}
\end{proof}

Now, we compute $G$ so as to get the antipode.

\begin{prop}
If $F_{\nu}=S_{\nu}^*$ comes from Tomita's theory, we have:
$$G=\sigma(F_{\nu}\otimes F_{\nu})$$
\end{prop}

\begin{proof}
For all $a=J_{\nu}a_1J_{\nu}\otimes a_2,b=J_{\nu}b_1J_{\nu}\otimes
b_2,c=J_{\nu}c_1J_{\nu}\otimes c_2$ and $d=J_{\nu}d_1J_{\nu}\otimes
d_2$ be analytic elements of $M'\otimes M$ w.r.t $\nu'\otimes\nu$.
Then, by \ref{simex}, we have:
$$
\begin{aligned}
&\
\quad(\lambda^{\beta,\alpha}_{\Lambda_{\nu}(\sigma_{i/2}^{\nu}(b_1))\otimes
\Lambda_{\nu}(\sigma_{-i}^{\nu}(b_2^*))})^*W^*(\Lambda_{\nu'\otimes\nu}(a)\surl{\
_{\alpha} \otimes_{\hat{\beta}}}_{\
\nu^o}\Lambda_{\nu'\otimes\nu}((J_{\nu}d_1^*J_{\nu}\otimes
d^*_2)c^*))\\
&=(\lambda^{\beta,\alpha}_{\Lambda_{\nu}(\sigma_{i/2}^{\nu}(b_1))\otimes
\Lambda_{\nu}(\sigma_{-i}^{\nu}(b_2^*))})^*I^*(\Lambda_{\nu'}(J_{\nu}d_1^*c_1^*J_{\nu})
\otimes (1\otimes d_2^*c_2^*)\Lambda_{\nu'\otimes\nu}(a))\\
&=\left[\rho_{\Lambda_{\nu}(\sigma_{-i}^{\nu}(b_2^*))}(1\otimes
\sigma_{i/2}^{\nu}(b_1))\right]^*(\Lambda_{\nu'}(J_{\nu}d_1^*c_1^*J_{\nu})
\otimes (1\otimes d_2^*c_2^*)\Lambda_{\nu'\otimes\nu}(a))\\
&=(d_2^*c_2^*\Lambda_{\nu}(a_2)|\Lambda_{\nu}(\sigma_{-i}^{\nu}(b_2^*)))\,
\Lambda_{\nu'}(J_{\nu}d_1^*c_1^*J_{\nu})\otimes\sigma_{-i/2}^{\nu}(b_1^*)
\Lambda_{\nu'}(J_{\nu}a_1J_{\nu})\\
&=\nu(d_2^*c_2^*a_2b_2)\, J_{\nu}\Lambda_{\nu}(d_1^*c_1^*)\otimes
J_{\nu}\Lambda_{\nu}(a_1b_1)
\end{aligned}$$
Consequently, by definition of $G$, we have:
$$
\begin{aligned}
&\quad G\left[\nu(d_2^*c_2^*a_2b_2)\,
J_{\nu}\Lambda_{\nu}(d_1^*c_1^*)\otimes
J_{\nu}\Lambda_{\nu}(a_1b_1)\right]\\
&=G(\lambda^{\beta,\alpha}_{\Lambda_{\nu}(\sigma_{i/2}^{\nu}(b_1))\otimes
\Lambda_{\nu}(\sigma_{-i}^{\nu}(b_2^*))})^*W^*(\Lambda_{\nu'\otimes\nu}(a)\surl{\
_{\alpha} \otimes_{\hat{\beta}}}_{\
\nu^o}\Lambda_{\nu'\otimes\nu}((J_{\nu}d_1^*J_{\nu}\otimes
d^*_2)c^*))\\
&=(\lambda^{\beta,\alpha}_{\Lambda_{\nu}(\sigma_{i/2}^{\nu}(d_1))\otimes
\Lambda_{\nu}(\sigma_{-i}^{\nu}(d_2^*))})^*W^*(\Lambda_{\nu'\otimes\nu}(c)\surl{\
_{\alpha} \otimes_{\hat{\beta}}}_{\
\nu^o}\Lambda_{\nu'\otimes\nu}((J_{\nu}b_1^*J_{\nu}\otimes
b^*_2)a^*))\\
&=\nu(b_2^*a_2^*c_2d_2)\, J_{\nu}\Lambda_{\nu}(b_1^*a_1^*)\otimes
J_{\nu}\Lambda_{\nu}(c_1d_1)
\end{aligned}$$
Since $G$ is anti-linear, we get:
$$G\left[J_{\nu}\Lambda_{\nu}(d_1^*c_1^*)\otimes
J_{\nu}\Lambda_{\nu}(a_1b_1)\right]=\left[
J_{\nu}\Lambda_{\nu}(b_1^*a_1^*)\otimes
J_{\nu}\Lambda_{\nu}(c_1d_1)\right]$$ so that $G$ coincides with
$\sigma(F_{\nu}\otimes F_{\nu})$.
\end{proof}

The polar decomposition of $G=ID^{1/2}$ is such that
$D=\Delta_{\nu}^{-1}\otimes\Delta_{\nu}^{-1}$ and
$I=\Sigma(J_{\nu}\otimes J_{\nu})$ so that the scaling group is
$\tau_t=\sigma_{-t}^{\nu'}\otimes\sigma_t^{\nu}$ for all
$t\in\mathbb{R}$ and the unitary antipode is $R=\varsigma\circ
(\beta_{\nu}\otimes\beta_{\nu})$. We also notice that
$\nu'\otimes\nu$ is $\tau$-invariant.

\begin{coro}
We have ${\mathcal D}(S)={\mathcal D}(\sigma_{i/2}^{\nu'})\otimes
{\mathcal D}(\sigma_{-i/2}^{\nu})$ and we have
$$S(J_{\nu}eJ_{\nu}\otimes
m^*)=J_{\nu}\sigma_{i/2}^{\nu}(m)J_{\nu}\otimes\sigma_{-i/2}^{\nu}(e^*)$$
for all $e,m\in {\mathcal D}(\sigma_{i/2}^{\nu})$. Moreover
$(id*\omega_{\xi,\eta})(W)\in {\mathcal D}(S)$ and:
$$S((id*\omega_{\xi,\eta})(W))=(id*\omega_{\xi,\eta})(W^*)$$ for
all $\xi,\eta\in D(_{\alpha}(H\otimes H),\nu)\cap D((H\otimes
H)_{\hat{\beta}},\nu^o)$.
\end{coro}

\begin{proof}
The first part of the corollary is straightforward by what precedes.
Let $\zeta,\eta\in H$ and $e,m\in {\mathcal D}(\sigma_{i/2}^{\nu})$.
By \ref{equ1}, we have:
$$\begin{aligned}
S((id*\omega_{\zeta\otimes J_{\nu}\Lambda_{\nu}(e),\eta\otimes
\Lambda_{\nu}(m)})(W))&=S((\zeta|\eta)J_{\nu}eJ_{\nu}\otimes
m^*)\\
&=(\zeta|\eta)J\sigma_{i/2}^{\nu}(m)J\otimes\sigma_{-i/2}^{\nu}(e^*)\\
&=(id*\omega_{\zeta\otimes J_{\nu}\Lambda_{\nu}(e),\eta\otimes
\Lambda_{\nu}(m)})(W^*)
\end{aligned}$$
Since $S$ is closed, we can conclude.
\end{proof}

\subsection{Dual structure}

We are now computing the dual structure.

\begin{prop}\label{cru}
For all $e,y\in {\mathcal N}_{\nu}$ and $\eta,\zeta\in H$, we have:
$$(\omega_{\Lambda_{\nu}(y)\otimes\eta,\zeta\otimes
J_{\nu}\Lambda_{\nu}(e)}* id)(W)=1\otimes
J_{\nu}e^*J_{\nu}\rho^*_{\zeta}\Sigma\rho_{\eta}y$$
\end{prop}

\begin{proof}
For all $\Xi\in H\otimes H,\xi\in H$ and $m\in {\mathcal N}_{\nu}$,
we have:
$$
\begin{aligned}
&\ \quad ((\omega_{\Lambda_{\nu}(y)\otimes\eta,\zeta\otimes
J_{\nu}\Lambda_{\nu}(e)}* id)(W)\Xi|\xi\otimes
\Lambda_{\nu}(m))\\
&=([\Lambda_{\nu}(y)\otimes\eta]\surl{\ _{\beta}
\otimes_{\alpha}}_{\ \nu}\Xi|W^*([\zeta\otimes
J_{\nu}\Lambda_{\nu}(e)]\surl{\ _{\alpha}
\otimes_{\hat{\beta}}}_{\ \nu^o}[\xi\otimes \Lambda_{\nu}(m)]))\\
&=((1\otimes y)\Xi \otimes\eta |\xi\otimes \zeta\otimes
mJ_{\nu}\Lambda_{\nu}(e))=(\Xi\otimes\eta |\xi\otimes
y^*\zeta\otimes J_{\nu}eJ_{\nu}\Lambda_{\nu}(m))\\
&=((1\otimes\rho_{\eta})\Xi |(1\otimes \Sigma\rho_{y^*\zeta}
J_{\nu}eJ_{\nu})(\xi\otimes \Lambda_{\nu}(m)))\\
&=((1\otimes J_{\nu}e^*J_{\nu}\rho^*_{\zeta}\Sigma\rho_{\eta}y)\Xi
|\xi\otimes \Lambda_{\nu}(m))
\end{aligned}$$
\end{proof}

\begin{coro}
We have $\widehat{M'\otimes M}=1\otimes \mathcal{L}(H)$.
\end{coro}

\begin{proof}
By definition, we recall that: $$\widehat{M'\otimes
M}=<(\omega_{\Xi_1,\Xi_2}*id)(W)\,|\, \Xi_1\in D((H\otimes
H)_{\beta},\nu^o),\Xi_2\in D(_{\alpha}(H\otimes H),\nu)>^{-w}$$ and
we notice that $\mathcal{L}(H)\otimes 1\subset
\alpha(M)'\cap\beta(M)'$. First of all, we prove that
$\widehat{M'\otimes M}\subset 1\otimes\mathcal{L}(H)$. Let $\Xi\in
H,\eta\in H$ and $m\in{\mathcal N}_{\nu}$. For all $x\in
\mathcal{L}(H)$, we have:
$$
\begin{aligned}
&\quad\,([1\otimes 1]\surl{\ _{\beta} \otimes_{\alpha}}_{\
\nu}[x\otimes 1])W^*(\Xi\surl{\ _{\alpha} \otimes_{\hat{\beta}}}_{\
\nu}[\eta\otimes J_{\nu}\Lambda_{\nu}(m)])\\
&=([1\otimes 1]\surl{\ _{\beta} \otimes_{\alpha}}_{\ \nu}[x\otimes
1])I^*(\eta\otimes
(1\otimes m)\Xi)=I^*(x\eta\otimes (1\otimes m)\Xi)\\
&=W^*(\Xi\surl{\ _{\alpha} \otimes_{\hat{\beta}}}_{\
\nu}[x\eta\otimes J_{\nu}\Lambda_{\nu}(m)])\\
&=W^*([1\otimes 1]\surl{\ _{\beta} \otimes_{\alpha}}_{\
\nu}[x\otimes 1])(\Xi\surl{\ _{\alpha} \otimes_{\hat{\beta}}}_{\
\nu}[\eta\otimes J_{\nu}\Lambda_{\nu}(m)])
\end{aligned}$$
Therefore, we get that $([1\otimes 1]\surl{\ _{\beta}
\otimes_{\alpha}}_{\ \nu}[x\otimes 1])W^*=W^*([1\otimes 1]\surl{\
_{\beta} \otimes_{\alpha}}_{\ \nu}[x\otimes 1])$ for all $x\in
\mathcal{L}(H)$. Thus, we get:
$$\widehat{M'\otimes M}\subset (\mathcal{L}(H)\otimes 1)'=1\otimes\mathcal{L}(H)$$

Then, we prove the reverse inclusion. By the previous proposition,
we state that, for all $v,w\in H$:
$$\begin{aligned}
(\omega_{\Lambda_{\nu}(y)\otimes\eta,\zeta\otimes
J_{\nu}\Lambda_{\nu}(e)}* id)(W)(v\otimes w)&=(1\otimes
J_{\nu}e^*J_{\nu})\rho^*_{y^*\zeta}(1\otimes\Sigma)\rho_{\eta}(v\otimes
w)\\
&=v\otimes(w|y^*\zeta)J_{\nu}e^*J_{\nu}\eta \end{aligned}$$ and
therefore we have:
$$
\begin{aligned}
\widehat{M'\otimes M}&\supset
<\omega_{\Lambda_{\nu}(y)\otimes\eta,\zeta\otimes
J_{\nu}\Lambda_{\nu}(e)}* id)(W)\,|\, \eta,\zeta\in H,e,y\in
{\mathcal N}_{\nu}>^{-w}\\
&= <1\otimes p\,|\, p \text{ rank 1 projection } >^{-w}= 1\otimes
\mathcal{L}(H)
\end{aligned}$$

\end{proof}

We verify that $(M'\otimes M)\cap\widehat{M'\otimes M}=1\otimes
M=\alpha(M)$. The dual co-product is given by
$\hat{\Gamma}(\widehat{x})=\sigma_{\nu^o}W(\widehat{x}\surl{\
_{\beta} \otimes_{\alpha}}_{\ \nu} 1)W^*\sigma_{\nu}$ for all
$\widehat{x}\in\widehat{M'\otimes M}$. We us the following
identification:
$$
\begin{aligned}
\Phi : 1\otimes \mathcal{L}(H) &\mapsto\mathcal{L}(H)\\
1\otimes x &\rightarrow x
\end{aligned}$$
which is implemented by $\lambda_e$ where $e\in H$ is a normalized
vector.Then, we know that $\Phi\surl{\
_{\hat{\beta}}\star_{\alpha}}_{\ \nu}\Phi$ is the identification
between $[1\otimes\mathcal{L}(H)]\surl{\
_{\hat{\beta}}\otimes_{\alpha}}_{\ \nu}[1\otimes\mathcal{L}(H)]$ and
$\mathcal{L}(H)\surl{\ _{\beta_{\nu}}\star_{id}}_{\
\nu}\mathcal{L}(H)\simeq\mathcal{L}(H)$.

\begin{prop}
We have $W^*\sigma_{\nu}(\lambda_e\surl{\ _{\beta_{\nu}}
\otimes_{id}}_{\
\nu}\lambda_e)=I^*(\lambda_e\otimes\lambda_e)I_{\nu}$ for all vector
$e\in H$ of norm $1$.
\end{prop}

\begin{proof}
Let $m\in {\mathcal N}_{\nu}$ and $\eta\in H$. We have:
$$
\begin{aligned}
&\ \quad W^*\sigma_{\nu}(\lambda_e\surl{\ _{\beta_{\nu}}
\otimes_{id}}_{\ \nu}\lambda_e)(\Lambda_{\nu}(m)\surl{\
_{\beta_{\nu}} \otimes_{id}}_{\ \nu}\eta)
=W^*\sigma_{\nu}([e\otimes\Lambda_{\nu}(m)]\surl{\ _{\hat{\beta}}
\otimes_{\alpha}}_{\ \nu} [e\otimes\eta])\\
&=W^*([e\otimes\eta]\surl{\ _{\alpha} \otimes_{\hat{\beta}}}_{\
\nu^o}[e\otimes\Lambda_{\nu}(m)])=I^*(e\otimes e\otimes
m\eta)\\
&=I^*(\lambda_e\otimes\lambda_e)I_{\nu}(\Lambda_{\nu}(m)\surl{\
_{\beta_{\nu}} \otimes_{id}}_{\ \nu}\eta)
\end{aligned}$$
\end{proof}

\begin{coro}
For all $x\in\mathcal{L}(H)$, we have:
$$(\Phi\surl{\
_{\hat{\beta}}\star_{\alpha}}_{\
\nu}\Phi)\circ\hat{\Gamma}\circ\Phi^{-1}(x)=I_{\nu}^*xI_{\nu}$$
\end{coro}

\begin{proof}
Let $x\in\mathcal{L}(H)$ and $e\in H$ be a vector of norm $1$. We
have:
$$
\begin{aligned}
&\ \quad (\Phi\surl{\ _{\hat{\beta}}\star_{\alpha}}_{\
\nu}\Phi)\circ\hat{\Gamma}\circ\Phi^{-1}(x)\\
&=(\lambda_e^*\surl{\ _{\beta_{\nu}} \otimes_{id}}_{\
\nu}\lambda_e^*)\sigma_{\nu^o}W([1\otimes x]\surl{\
_{\hat{\beta}}\otimes_{\alpha}}_{\ \nu}[1\otimes 1
])W^*\sigma_{\nu}(\lambda_e\surl{\ _{\beta_{\nu}}
\otimes_{id}}_{\ \nu}\lambda_e)\\
&=I^*_{\nu}(\lambda_e^*\otimes\lambda_e^*)I([1\otimes x]\surl{\
_{\hat{\beta}}\otimes_{\alpha}}_{\
\nu}[1\otimes 1])I^*(\lambda_e\otimes\lambda_e)I_{\nu}\\
&=I^*_{\nu}(\lambda_e^*\otimes\lambda_e^*)(1\otimes 1\otimes x)
(\lambda_e\otimes\lambda_e)I_{\nu}=I^*_{\nu}xI_{\nu}
\end{aligned}$$
\end{proof}

Now, we are computing the dual operator-valued weight.

\begin{lemm}
We have $\hat{\Lambda}((\omega_{\Xi,\Lambda_{\nu}(m)\otimes
J_{\nu}\Lambda_{\nu}(e)}*id)(W))=(m^*\otimes J_{\nu}e^*J_{\nu})\Xi$
for all $m,e\in {\mathcal N}_{\nu}$ et $\Xi\in D((H\otimes
H)_{\beta},\nu^o)$.
\end{lemm}

\begin{proof}
Let $m_1,m_2\in {\mathcal N}_{\nu}$. We have:
$$
\begin{aligned}
&\quad\, (\hat{\Lambda}((\omega_{\Xi,\Lambda_{\nu}(m)\otimes
J_{\nu}\Lambda_{\nu}(e)}*id)(W))|\Lambda_{\nu'\otimes\nu}(J_{\nu}m_1J_{\nu}\otimes
m_2))\\
&=\omega_{\Xi,\Lambda_{\nu}(m)\otimes
J_{\nu}\Lambda_{\nu}(e)}(J_{\nu}m_1^*J_{\nu}\otimes
m_2^*)=((J_{\nu}m_1^*J_{\nu}\otimes m_2^*)\Xi
|\Lambda_{\nu}(m)\otimes
J_{\nu}\Lambda_{\nu}(e))\\
&=(\Xi |mJ_{\nu}\Lambda_{\nu}(m_1)\otimes
J_{\nu}eJ_{\nu}\Lambda_{\nu}(m_2))=((m^*\otimes
J_{\nu}e^*J_{\nu})\Xi
|\Lambda_{\nu'\otimes\nu}(J_{\nu}m_1J_{\nu}\otimes m_2))
\end{aligned}$$
\end{proof}

\begin{prop}
We have $\hat{T_L}=id\otimes E_{\nu'}$ where $E_{\nu'}$ is the
operator-valued weight from $\mathcal{L}(H)$ to $M$ obtained from
the weight $\nu'$.
\end{prop}

\begin{proof}
Let $m,e,y\in {\mathcal N}_{\nu}$ and $\eta\in H$. On one hand, we
compute:
$$\begin{aligned}
||\hat{\Lambda}((\omega_{\Xi,\Lambda_{\nu}(m)\otimes
J_{\nu}\Lambda_{\nu}(e)}*id)(W))||^2&=||m^*\Lambda_{\nu}(y)\otimes
J_{\nu}e^*J_{\nu}\eta ||^2\\
&=||J_{\nu}e^*J_{\nu}\eta ||^2||\Lambda_{\nu}(m^*y)||^2
\end{aligned}$$ On the other hand, by proposition \ref{cru}, we have:
$$
\begin{aligned}
&\ \quad ||\hat{\Lambda}((\omega_{\Xi,\Lambda_{\nu}(m)\otimes
J_{\nu}\Lambda_{\nu}(e)}*id)(W))||^2\\
&=\hat{\Phi}(\rho_{\eta}^*(1\otimes\Sigma)\rho_{y^*\Lambda_{\nu}(m)}(1\otimes
J_{\nu}eJ_{\nu})(1\otimes
J_{\nu}e^*J_{\nu})\rho^*_{y^*\Lambda_{\nu}(m)}(1\otimes\Sigma)\rho_{\eta})\\
&=||J_{\nu}e^*J_{\nu}\eta ||^2\hat{\Phi}(1\otimes
(\Lambda_{\nu}(m^*y)\otimes \Lambda_{\nu}(m^*y)))
\end{aligned}$$
where $\xi\otimes\xi$ is the operator of $\mathcal{L}(H)$ such that
$(\xi\otimes\xi)v=(v|\xi)\xi$. Then, if $\xi\in {\mathcal
D}(S_{\nu})$, then we have:
$$\begin{aligned}
\hat{\Phi}(1\otimes (\xi\otimes\xi))&=||S_{\nu}\xi ||^2
=(\Delta_{\nu}\xi|\xi)=(\frac{\text{d}\nu}{\text{d}\nu'}\xi
|\xi)=\nu(\theta^{\nu'}(\xi,\xi))\\
&=\nu\circ E_{\nu'}(\xi\otimes\xi)=\nu\circ\alpha^{-1}\circ
\hat{T_L}(1\otimes (\xi\otimes\xi))
\end{aligned}$$
\end{proof}
We also have the following formulas, for all $x\in \mathcal{L}(H)$
and $t\in\mathbb{R}$:
$$\hat{R}(1\otimes x)=1\otimes J_{\nu}x^*J_{\nu}\quad\text{ and
}\quad\hat{\tau}_t(1\otimes x)=1\otimes
\Delta_{\nu}^{it}x\Delta_{\nu}^{-it}$$ The right invariant
operator-valued weight is given by:
$\hat{T_R}=\hat{R}\circ\hat{T_L}\circ\hat{R}=(id\otimes E_{\nu})$.

\begin{prop}
The dual pairs quantum groupoid can be identify with
$(M,\mathcal{L}(H),id,j,id,\nu,E_{\nu'},E_{\nu})$ which is a
measured quantum groupoid but not a adapted measured quantum
groupoid. Moreover, expressions for co-involution and scaling group
are given, for all $x\in \mathcal{L}(H)$ and $t\in\mathbb{R}$:
$$\hat{R}(x)=J_{\nu}x^*J_{\nu}\quad\text{ and
}\quad\hat{\tau}_t(x)=\Delta_{\nu}^{it}x\Delta_{\nu}^{-it}$$
\end{prop}

\begin{proof}
The proposition gathers results of the section. Nevertheless we lay
stress on the following point. We have, for all $t\in\mathbb{R}$ and
$m\in M$:
$$\sigma_{t}^{E_{\nu}}(m)=\sigma_{-t}^{\nu}(m)$$
instead of $\sigma_t^{\nu}(m)$ to have a adapted measured quantum
groupoid.
\end{proof}

\begin{rema}
If $M=L^{\infty}(X)$, we find the structure of pairs groupoid
$X\times X$. This example comes from the inclusion of von Neumann
algebras (\cite{E1}):
$$\mathbb{C}\subset M \subset \mathcal{L}(L^2(M))\subset
\mathcal{L}(L^2(M))\otimes M\subset\ldots$$
\end{rema}

\section{Inclusions of von Neumann algebras}\label{ivnasur}

Let $M_0\subseteq M_1$ be an inclusion of von Neumann algebras. We
call \textbf{basis construction} the following inclusions:

$$M_0\subseteq M_1\subseteq M_2=J_1M'_0J_1=End_{M_0^o}(L^2(M_1))$$

By iteration, we construct Jones' tower $M_0\subseteq M_1\subseteq
M_2\subseteq M_3\subseteq\cdots$

\begin{defi}If $M'_0\cap M_1\subseteq M'_0\cap M_2\subseteq M'_0\cap
M_3$ is a basis construction, then the inclusion is said to be
\textbf{of depth $2$}.
\end{defi}

Let $T_1$ be a n.s.f operator-valued weight from $M_1$ to $M_0$. By
Haagerup's construction \cite{St} (12.11) and \cite{EN} (10.1), it
is possible to define a canonical n.s.f operator-valued weight $T_2$
from $M_2$ to $M_1$ such that, for all $x,y\in {\mathcal N}_{T_1}$,
we have:
$$T_2(\Lambda_{T_1}(x)\Lambda_{T_1}(y)^*)=xy^*$$ By iteration, we define,
for all $i\geq 1$, a n.s.f operator-valued weight $T_i$ from $M_i$
to $M_{i-1}$. If $\psi_0$ is n.s.f weight sur $M_0$, we put
$\psi_i=\psi_{i-1}\circ T_i$.

\begin{defi}
\cite{EN} (11.12), \cite{EV} (3.6). $T_1$ is said to be regular if
restrictions of $T_2$ to $M_0'\cap M_2$ and of $T_3$ to $M_1'\cap
M_3$ are semifinite.
\end{defi}

\begin{prop}\cite{EV} (3.2, 3.8, 3.10).
If $M_0\subset M_1$ is an inclusion with a regular n.s.f
operator-valued weight $T_1$ from $M_1$ to $M_0$, then there exists
a natural *-representation $\pi$ of $M_0'\cap M_3$ on $L^2(M_0'\cap
M_2)$ whose restriction to $M_0'\cap M_2$ is the standard
representation of $M_0'\cap M_2$. Moreover, the inclusion is of
depth $2$ if, and only if $\pi$ is faithful.
\end{prop}

The following theorem exhibits a structure of measured quantum
groupoid coming from inclusion of von Neumann algebras.

\begin{theo}\label{ex2}
Let $M_0\subset M_1$ be a depth 2 inclusion of $\sigma$-finite von
Neumann algebras, equipped with a regular nsf operator-valued weight
$T_1$ in the sense of \cite{E1} and \cite{E3}. Moreover, assume
there exists on $M_0'\cap M_1$ a nsf weight $\chi$ invariant under
the modular automorphism group $T_1$. Then, by theorem 8.2 of
\cite{E3}, we have:
\begin{enumerate}
\item there exists an application $\widetilde{\Gamma}$ from
$M_0'\cap M_2$ to $$(M_0'\cap M_2)\surl{\ _{j_1}
\star_{id}}_{M_0'\cap M_1}(M_0'\cap M_2)$$ such that $(M_0'\cap
M_1,M_0'\cap M_2,id,j_1,\widetilde{\Gamma})$ is a Hopf-bimodule,
(where id means here the injection of $M_0'\cap M_1$ into $M_0'\cap
M_2$, and $j_1$ means here the restriction of $j_1$ coming from
Tomita's theory to $M_0'\cap M_1$, considered then as an
anti-representation of $M_0'\cap M_1$ into $M_0'\cap M_2$).
Moreover, the anti-automorphism $j_1$ of $M_0'\cap M_2$ is a
co-involution for this Hopf-bimodule structure.
\item the nsf operator-valued weight $\widetilde{T_2}$ from $M_0'\cap
M_2$ to $M_0'\cap M_1$, restriction of the second canonical weight
construct from Jones' tower and $T_1$, is left invariant.
\item Let $\chi_2$ be the weight $\chi\circ\widetilde{T_2}$; there
exist a one-parameter group of automorphisms $\widetilde{\tau_t}$ of
$M_0'\cap M_2$, commuting with the modular automorphism group
$\sigma^{\chi_2}$, such that, for all $t\in\mathbb R$, we have:
$$\widetilde{\Gamma}\circ\sigma^{\chi_2}_t=(\widetilde{\tau_t}\surl{\ _{j_1}
\star_{id}}_{\chi}\sigma^{\chi_2}_t)\circ\widetilde{\Gamma}$$
Moreover, we have
$j_1\circ\widetilde{\tau_t}=\widetilde{\tau_t}\circ j_1$.
\end{enumerate}
Then, $(M_0'\cap M_1,M_0'\cap
M_2,id,j_1,\widetilde{\Gamma},\widetilde{T_2},j_1,\widetilde{\tau},\chi)$
is a measured quantum groupoid.
\end{theo}

\begin{proof}
We have $\Phi=\chi_2$. Then, by proposition 6.6 of \cite{E3}, we
have the relation between $R$ and $\Gamma$. Also, we notice that
$\widetilde{\tau}$ coincide with $\sigma^{\chi}$ on $M_0'\cap M_1$
by theorem 5.10 of \cite{E3} and that we have, for all $n\in
M_0'\cap M_1$ and $t\in\mathbb R$:
$$\sigma_t^{\widetilde{T_2}}(j_1(n))=\sigma_t^{T_2}(j_1(n))=j_1(\sigma_t^{T_1}(n))$$
by corollary 4.8 and by 4.1 of \cite{E3}. So that,
$\gamma=\sigma^{T_1}$ leaves $\chi$ invariant by hypothesis.
\end{proof}

Then we can show that the dual structure coincide with the natural
one on the second relative commutant of Jones' tower.

\begin{theo}
Let $M_0\subset M_1$ be a depth 2 inclusion of $\sigma$-finite von
Neumann algebras, equipped with a regular nsf operator-valued weight
$T_1$ in the sense of \cite{E1} and \cite{E3}. Moreover, assume
there exists on $M_0'\cap M_1$ a nsf weight $\chi$ invariant under
the modular automorphism group $T_1$.
\begin{enumerate}
\item there exists an application $\widetilde{\Gamma}'$ from
$M_1'\cap M_3$ to $$(M_1'\cap M_3)\surl{\ _{j_1} \star_{j_2\circ
j_1}}_{M_0'\cap M_1}(M_1'\cap M_3)$$ such that $(M_0'\cap
M_1,M_1'\cap M_3,j_2\circ j_1,j_1,\widetilde{\Gamma}')$ is a
Hopf-bimodule, where $j_1,j_2$ come from Tomita's theory. Moreover,
the anti-automorphism $j_2$ of $M_1'\cap M_3$ is a co-involution for
this Hopf-bimodule structure.
\item the nsf operator-valued weight $\widetilde{T_3}$ from $M_1'\cap
M_3$ to $M_1'\cap M_2=j_1(M_0'\cap M_1)$, restriction of the second
canonical weight construct from Jones' tower and $T_1$, is left
invariant.
\item Let $\chi_3$ be the weight $\chi\circ\widetilde{T_3}$; there
exist a one-parameter group of automorphisms $\widetilde{\tau_t}$ of
$M_1'\cap M_3$, commuting with the modular automorphism group
$\sigma^{\chi_3}$, such that, for all $t\in\mathbb R$, we have:
$$\widetilde{\Gamma}\circ\sigma^{\chi_3}_t=(\widetilde{\tau_t}\surl{\ _{j_1}
\star_{j_2\circ
j_1}}_{\chi}\sigma^{\chi_3}_t)\circ\widetilde{\Gamma}$$ Moreover, we
have $j_2\circ\widetilde{\tau_t}=\widetilde{\tau_t}\circ j_2$.
\end{enumerate}
Then, $(M_0'\cap M_1,M_1'\cap M_3,j_2\circ j_1
,j_1,\widetilde{\Gamma},\widetilde{T_3},j_2,\widetilde{\tau},\chi)$
is the dual measured quantum groupoid of $M_0'\cap M_2$ (equipped
with the structure described on \ref{ex2}).
\end{theo}

\begin{proof}
All objects are constructed from the fundamental unitary that's why
the Hopf-bimodule structure of the dual coincide with the structure
on the second relative commutant. The uniqueness theorem implies
that the dual operator-valued weight coincide with the restriction
of $T_3$ upto an element of the basis.
\end{proof}

We can't characterize, at this stage, inclusions of von Neumann
algebras among measured quantum groupoids. A way to answer the
question is to know if each measured quantum groupoid acts on a von
Neumann algebra .

\section{Operations on adapted measured quantum groupoids}

\subsection{Elementary operations}\label{ogqm}
\label{operations}

\subsubsection{Sum of adapted measured quantum groupoids}\label{sgqm}

A union of groupoids is still a groupoid. We establish here a
similar result at the quantum level:

\begin{prop}
Let $(N_i,M_i,\alpha_i,\beta_i,\Gamma_i,\nu_i,T_L^i,T_R^i)_{i\in I}$
be a family of adapted measured quantum groupoids. In the von
Neumann algebra level, if we identify $\bigoplus_{i\in I}M_i\surl{\
_{\beta}\star_{\alpha}}_{\ N_i}M_i$ with $\left(\bigoplus_{i\in
I}M_i\right)\surl{\ _{\beta}\star_{\alpha}}_{\ \bigoplus_{i\in
I}N_i}\left(\bigoplus_{i\in I}M_i\right)$, then we get:
$$(\bigoplus_{i\in I}N_i,\bigoplus_{i\in I}M_i,\bigoplus_{i\in I}\alpha_i,\bigoplus_{i\in I}\beta_i,
\bigoplus_{i\in I}\Gamma_i,\bigoplus_{i\in I}\nu_i,\bigoplus_{i\in
I}T_L^i,\bigoplus_{i\in I}T_R^i)$$ a adapted measured quantum
groupoid where operators act on the diagonal.
\end{prop}

\begin{proof}
Straightforward.
\end{proof}

In particular, the sum of two quantum groups with different scaling
constants (\cite{VaV} for examples) produce a adapted measured
quantum groupoid with non scalar scaling operator.

\subsubsection{Tensor product of adapted measured quantum groupoids}
Cartesian product of groups correspond to tensor product of quantum
groups. In the same way, we have:

\begin{prop}
Let $(N_i,M_i,\alpha_i,\beta_i,\Gamma_i,\nu_i,T_L^i,T_R^i)$ be
adapted measured quantum groupoids for $i=1,2$. If we identify
$(M_1\surl{\ _{\beta_1}\star_{\alpha_1}}_{\ N_1}M_1)\otimes
(M_2\surl{\ _{\beta_2}\star_{\alpha_2}}_{\ N_2}M_2)$ with
$(M_1\otimes M_2)\surl{\
_{\beta_1\otimes\beta_2}\star_{\alpha_1\otimes\alpha_2}}_{\
N_1\otimes N_2}(M_1\otimes M_2)$ as von Neumann algebras, then we
have:
$$(N_1\otimes N_2,M_1\otimes M_2,\alpha_1\otimes\alpha_2,\beta_1\otimes\beta_2,
\Gamma_1\otimes\Gamma_2,\nu_1\otimes\nu_2,T_L^1\otimes
T_L^2,T_R^1\otimes T_R^2)$$ is a adapted measured quantum groupoid.
\end{prop}

\begin{proof}
Straightforward.
\end{proof}

\subsubsection{Direct integrals of adapted measured quantum groupoids}

In this section, $X$ denote $\sigma$-compact, locally compact space
and $\mu$ a Borel measure on $X$. Theory of hilbertian integrals is
described in \cite{T}.

\begin{prop}
Let $(N_p,M_p,\alpha_p,\beta_p,\Gamma_p,\nu_p,T_L^p,T_R^p)_{p\in X}$
be a family of adapted measured quantum groupoids. In the von
Neumann algebra level, if we identify $\int_X^{\oplus}\!M_p\surl{\
_{\beta}\star_{\alpha}}_{\ N_p}M_p\,d\mu(p)$ and
$\left(\int_X^{\oplus}\!M_p\,d\mu(p)\right)\surl{\
_{\beta}\star_{\alpha}}_{\
\int_X^{\oplus}\!N_p\,d\mu(p)}\left(\int_X^{\oplus}\!M_p\,d\mu(p)\right)$,
we have:
$$(\int_X^{\oplus}\!N_p\,d\mu(p),\int_X^{\oplus}\!M_p\,d\mu(p),\int_X^{\oplus}\!\alpha_p\,d\mu(p),\int_X^{\oplus}\!\beta_p\,d\mu(p),\cdots$$
$$\cdots\int_X^{\oplus}\!\Gamma_p\,d\mu(p),\int_X^{\oplus}\!\nu_p\,d\mu(p),\int_X^{\oplus}\!T_L^p\,d\mu(p),\int_X^{\oplus}\!T_R^p\,d\mu(p))$$
is a adapted measured quantum groupoid.
\end{prop}

\begin{proof}
Left to the reader.
\end{proof}

\cite{Bl} gives examples. In this case, the basis is $L^{\infty}(X)$
and $\alpha=\beta=\hat{\beta}$. The fundamental unitary comes from a
space onto the same space and then can be viewed as a field of
multiplicative unitaries.

\subsection{Opposite and commutant structures}
\begin{defi}
We call Hopf-bimodule morphism from
$(N,M_1,\alpha_1,\beta_1,\Gamma_1)$ to
$(N,M_2,\alpha_2,\beta_2,\Gamma_2)$ a morphism $\pi$ of von Neumann
algebras from $M_1$ to $M_2$ such that:
\begin{center}
\begin{minipage}{8cm}
\begin{enumerate}[i)]
\item $\pi\circ\alpha_1=\alpha_2$ et $\pi\circ\beta_1=\beta_2$ ;
\item $\Gamma_2\circ\pi=(\pi\surl{\ _{\beta_1} \star_{\alpha_1}}_{\
N}\pi)\circ\Gamma$.
\end{enumerate}
\end{minipage}
\end{center}
Also, we call anti-morphism of Hopf-bimodule from
$(N,M_1,\alpha_1,\beta_1,\Gamma_1)$ to\\
$(N,M_2,\alpha_2,\beta_2,\Gamma_2)$ a morphism $j$ of von Neumann
algebras from $M_1$ to $M_2$ such that:
\begin{center}
\begin{minipage}{8cm}
\begin{enumerate}[i)]
\item $j\circ\alpha_1=\beta_2$ et $j\circ\beta_1=\alpha_2$ ;
\item $\Gamma_2\circ j=(j\surl{\ _{\beta_1} \star_{\alpha_1}}_{\
N}j)\circ\Gamma_1$.
\end{enumerate}
\end{minipage}
\end{center}
\end{defi}

\begin{defi}
For all Hopf-bimodule $(N,M_1,\alpha_1,\beta_1,\Gamma_1)$ and all
1-1 morphism of von Neumann algebras $\pi$ from $M_1$ onto $M_2$,
$M_2$, $(N,M_2,\pi\circ\alpha_1,\pi\circ\beta_1,(\pi\surl{\
_{\beta_1} \star_{\alpha_1}}_{\ N}\pi)\circ\Gamma\circ\pi^{-1})$ is
a Hopf-bimodule called Hopf-bimodule image by $\pi$. Also, if $j$ is
a 1-1 anti-morphism of von Neumann algebras from $M_1$ onto $M_2$,
then $(N^0,M_2,j\circ\alpha_1,j\circ\beta_1,(j\surl{\ _{\beta_1}
\star_{\alpha_1}}_{\ N}j)\circ\Gamma\circ j^{-1})$ is is a
Hopf-bimodule called Hopf-bimodule image by $j$.
\end{defi}

\begin{prop}\label{mor}
Let $\pi$ a 1-1 morphism from $(N,M_1,\alpha_1,\beta_1,\Gamma_1)$
onto $(N,M_2,\alpha_2,\beta_2,\Gamma_2)$. If
$(N,M_1,\alpha_1,\beta_1,\Gamma_1,R,T_L,\tau,\nu)$ is a measured
quantum groupoid, then $(N,M_2,\alpha_2,\beta_2,\Gamma_2,\pi\circ
R\circ\pi^{-1},\pi\circ T_L\circ\pi^{-1},\pi\circ
\tau\circ\pi^{-1},\nu)$ is a measured quantum groupoid such that:
$$\lambda_2=\pi(\lambda_1)\text{ et } \delta_2=\pi(\delta_1)$$ We
denote by $\Phi^1=\nu\circ\alpha_1^{-1}\circ T_L$ and
$\Phi^2=\Phi^1\circ\pi^{-1}$. If $I$ is the unitary from
$H_{\Phi^1}$ onto $H_{\Phi^2}$ such that
$I\Lambda_{\Phi^1}(a)=\Lambda_{\Phi^2}(\pi(a))$ for all $a\in
{\mathcal N}_{\Phi^1}$, then fundamental unitaries are linked by:
$$W_2=(I\surl{\ _{\alpha_1}
  \otimes_{\hat{\beta_1}}}_{\ N^o}I)W_1(I^*\surl{\ _{\beta_2}
  \otimes_{\alpha_2}}_{\ N}I^*)$$
\end{prop}

\begin{proof}
It is easy to state that $(N,M_2,\alpha_2,\beta_2,\Gamma_2,\pi\circ
R\circ\pi^{-1},\pi\circ T_L\circ\pi^{-1},\pi\circ
\tau\circ\pi^{-1},\nu)$ is a measured quantum groupoid. For all
$v\in D((H_{\Phi^1})_{\beta_1},\nu^o)$, $a\in {\mathcal N}_{T_L}\cap
{\mathcal N}_{\Phi^1}$ and $(N^o,\nu^o)$-basis $(\xi_i)_{i\in I}$ of
$(H_{\Phi^1})_{\beta_1}$, we have:
$$
\begin{aligned}
&\ \quad W_1^*(I^*\surl{\ _{\alpha_2}\otimes_{\hat{\beta_2}}}_{\ \
N^o}I^*)(Iv\surl{\ _{\alpha_2}\otimes_{\hat{\beta_2}}}_{\ \
\nu^o}\Lambda_{\Psi^2}(\pi(a))=W_1^*(v\surl{\ _{\alpha_1}
\otimes_{\hat{\beta_1}}}_{\ \ \nu^o}\Lambda_{\Psi^1}(a))\\
&=\sum_{i\in I}\xi_i\surl{\ _{\beta_1}
    \otimes_{\alpha_1}}_{\ \nu}\Lambda_{\Phi^1}((\omega_{v,\xi_i}\surl{\
_{\beta_1}\star_{\alpha_1}}_{\ \nu}id)\Gamma_1(a))\\
&=\sum_{i\in I}\xi_i\surl{\ _{\beta_1}
    \otimes_{\alpha_1}}_{\ \nu}
    \Lambda_{\Phi^1}((\omega_{v,\xi_i}\surl{\
_{\beta_1}\star_{\alpha}}_{\ \nu}id)(\pi^{-1}\surl{\ _{\beta_1}
\star_{\alpha_1}}_{\
N}\pi^{-1})\Gamma_2(\pi(a))))\\
&=(I^*\surl{\ _{\beta_2}\otimes_{\alpha_2}}_{\ N}I^*)\sum_{i\in
I}I\xi_i\surl{\ _{\beta_2} \otimes_{\alpha_2}}_{\
\nu}\Lambda_{\Phi^2}((\omega_{Iv,I\xi_i}\surl{\
_{\beta_2}\star_{\alpha_2}}_{\ \nu}id)\Gamma_2(a))\\
&=(I^*\surl{\ _{\beta_2}\otimes_{\alpha_2}}_{\ N}I^*)W_2^*(Iv\surl{\
_{\alpha_2}
    \otimes_{\hat{\beta_2}}}_{\ \ \nu^o}\Lambda_{\Psi^2}(\pi(a)))
\end{aligned}$$ Then, we have proved that $W_2=(I\surl{\ _{\alpha_1}
\otimes_{\hat{\beta_1}}}_{\ N^o}I)W_1(I^*\surl{\ _{\beta_2}
\otimes_{\alpha_2}}_{\ N}I^*)$. For all $a,b\in {\mathcal
N}_{jT_Lj}\cap {\mathcal N}_{\Phi\circ j}$, we have:
$$
\begin{aligned}
&\ \quad R_2((id\surl{\ _{\beta_2}\otimes_{\alpha_2}}_{\
\nu}\omega_{J_{\Phi^2}\Lambda_{\Phi^2}(\pi(a))})\Gamma_2(\pi(b^*)\pi(b)))\\
&=(id\surl{\ _{\beta_2}\otimes_{\alpha_2}}_{\
\nu}\omega_{J_{\Phi^2}\Lambda_{\Phi^2}(\pi(b))})\Gamma_2(\pi(a^*a))\\
&=(id\surl{\ _{\beta_2}\otimes_{\alpha_2}}_{\
\nu}\omega_{J_{\Phi^2}\Lambda_{\Phi^2}(\pi(b))})(\pi\surl{\
_{\beta_1}\star_{\alpha_1}}_{\ N}\pi)\Gamma_1(a^*a)\\
&=\pi((id\surl{\ _{\beta_1}\otimes_{\alpha_1}}_{\ \
\nu^o}\omega_{J_{\Phi^1}\Lambda_{\Phi^1}(b)})\Gamma_1(a^*a))=\pi
R_1((id\surl{\ _{\beta_1}\otimes_{\alpha_1}}_{\ \
\nu^o}\omega_{J_{\Phi^1}\Lambda_{\Phi^1}(a)})\Gamma_1(b^*b))\\
&=\pi R_1\pi^{-1}((id\surl{\ _{\beta_2}\otimes_{\alpha_2}}_{\
\nu}\omega_{J_{\Phi^2}\Lambda_{\Phi^2}(\pi(a))})\Gamma_2(\pi(b^*)\pi(b)))
\end{aligned}$$ from which we get that $R_2=\pi\circ
R_1\circ\pi^{-1}$ and then $S_2=\pi\circ S_1\circ\pi^{-1}$ and
$\tau_2=\pi\circ\tau_1\circ \pi^{-1}$ thanks to fundamental
unitaries. Finally, we have for all $t\in\mathbb{R}$:
$$\begin{aligned}
\ [D\Phi^2\circ R_2:D\Phi^2]_t&=[D\Phi^1\circ
R_1\circ\pi^{-1}:D\Phi^1\circ\pi^{-1}]_t
\\&=\pi([D\Phi^1\circ
R_1:D\Phi^1]_t)=\pi(\lambda_1)^{\frac{it^2}{2}}\pi(\delta_1)^{it}
\end{aligned}$$ and, so we have $\delta_2=\pi(\delta_1)$ and
$\lambda_2=\pi(\lambda_1)$.
\end{proof}

\begin{prop}\label{antimor}
Let $j$ a 1-1 anti-morphism from $(N,M_1,\alpha_1,\beta_1,\Gamma_1)$
onto $(N^0,M_2,\alpha_2,\beta_2,\Gamma_2)$. If
$(N,M_1,\alpha_1,\beta_1,\Gamma_1,R,T_L,\tau,\nu)$ is a measured
quantum groupoid, then $(N^0,M_2,\beta_2,\alpha_2,\Gamma_2,j\circ
R\circ j^{-1},j\circ T_L\circ j^{-1},j\circ \tau_{-t}\circ
j^{-1},\nu^0)$ is a measured quantum groupoid such that:
$$\lambda_2=j(\lambda_1^{-1})\text{ et } \delta_2=j(\delta_1)$$ We
denote by $\Phi^1=\nu\circ\alpha_1^{-1}\circ T_L$ and
$\Phi^2=\Phi^1\circ j^{-1}$. If $J$ is the unitary from $H_{\Phi^1}$
onto $H_{\Phi^2}$ such that
$I\Lambda_{\Phi^1}(a)=J_{\Phi^2}\Lambda_{\Phi^2}(j(a^*))$ for all
$a\in {\mathcal N}_{\Phi^1}$, then fundamental unitaries are linked
by:
$$W_2=(J\surl{\ _{\alpha_1}
  \otimes_{\hat{\beta_1}}}_{\ N^o}J)W_1(J^*\surl{\ _{\alpha_2}
  \otimes_{\beta_2}}_{\ N}J^*)$$
\end{prop}

\begin{proof}
The proof is very similar to the previous one.
\end{proof}

\begin{defi}
We call opposite quantum groupoid the image by the co-involution $R$
of the Hopf-bimodule, denoted by
$(N,M,\alpha,\beta,\Gamma,R,T_L,\tau,\nu)^{\text{op}}$. THe
Hopf-bimodule is then the symmetrized one
$(N^o,M,\beta,\alpha,\varsigma_N\circ\Gamma)$.
\end{defi}

\begin{rema}
If $N$ is abelian, $\alpha=\beta$, $\Gamma=\varsigma_N\circ\Gamma$
then the measured quantum groupoid is equal to its opposite : we
speak about symmetric quantum groupoid.
\end{rema}

We put $j$ the canonical *-anti-isomorphism from $M$ onto $M'$
coming from Tomita's theory, trivial on the center of $M$ and given
by $j(x)=J_{\Phi}x^*J_{\Phi}$. Then, we have
$j\circ\alpha=\widehat{\beta}$ and we put $\varrho=j\circ\beta$. We
can construct a unitary $j\surl{\ _{\varrho} \star_{\hat{\beta}}}_{\
N^o}j$ from $M'\surl{\ _{\varrho} \star_{\hat{\beta}}}_{\ N^o}M'$
onto $M\surl{\ _{\beta} \star_{\alpha}}_{\ N}M$ the adjoint of which
is $j\surl{\ _{\beta} \star_{\alpha}}_{\ N}j$.

\begin{defi}
We call commutant quantum groupoid the image by $j$ of the Hopf
bimodule. It is denoted by
$(N,M,\alpha,\beta,\Gamma,R,T_L,\tau,\nu)^c$. The Hopf-bimodule is
equal to $(N^o,M',\hat{\beta},\varrho,(j\surl{\
_{\beta}\star_{\alpha}}_{\ N}j)\circ\Gamma\circ j)$. We put
$\Gamma^c=(j\surl{\ _{\beta}\star_{\alpha}}_{\ N}j)\circ\Gamma\circ
j$.
\end{defi}

We describe fundamental objects of the structures.

\begin{prop}
We have the following formulas:
\begin{enumerate}[i)]
\item $W^{\text{op}}=\sigma_{\nu^o}W'^*\sigma_{\nu^o}$,
$R^{\text{op}}=R$, $\tau^{\text{op}}_t=\tau_{-t}$,
$\delta^{\text{op}}=\delta^{-1}$ et
$\lambda^{\text{op}}=\lambda^{-1}$ ;
\item $W^c=(J_{\Phi}\surl{\ _{\beta}
    \otimes_{\alpha}}_{\ N}J_{\Phi})W(J_{\Phi}\surl{\ _{\varrho}
    \otimes_{\hat{\beta}}}_{\ N^o}J_{\Phi})$,
$R^c=jRj$, $\tau^c_t=j\tau_{-t}j$, $\delta^c=j(\delta)$ et
$\lambda^c=\lambda^{-1}$.
\end{enumerate}
\end{prop}

\begin{proof}
It is an easy consequence of propositions \ref{mor} and
\ref{antimor} except for the relation between $W^{\text{op}}$ and
$W'$. For all $v\in D(_{\alpha}H_{\Psi},\nu)$, $a\in {\mathcal
N}_{T_R}\cap {\mathcal N}_{\Psi_R}$ and $(N,\nu)$-basis
$(\eta_i)_{i\in I}$ of $\ _{\alpha}H_{\Psi}$, we have:
$$
\begin{aligned}
(W^{\text{op}})^*\sigma_{\nu^o}(\Lambda_{\Psi}(a)\surl{\
_{\hat{\alpha}}
    \otimes_{\beta}}_{\ \ \nu^o}v)&=(W^{\text{op}})^*(v\surl{\
_{\hat{\alpha}}
    \otimes_{\beta}}_{\ \ \nu^o}\Lambda_{\Psi}(a))\\
&=\sum_{i\in I}\eta_i\surl{\ _{\alpha}\otimes_{\beta}}_{\ \
\nu^o}\Lambda_{\Psi}((\omega_{v,\eta_i}\surl{\ _{\alpha}
\star_{\beta}}_{\ \ \nu^o}id)(\varsigma_N\circ\Gamma(a)))\\
&=\sigma_{\nu}\sum_{i\in I}\Lambda_{\Psi}((id\surl{\ _{\alpha}
\star_{\beta}}_{\ \ \nu^o}\omega_{v,\eta_i})\Gamma(a))\surl{\
_{\beta}\otimes_{\alpha}}_{\ \nu}\eta_i\\
&=\sigma_{\nu}W'(\Lambda_{\Psi}(a)\surl{\ _{\hat{\alpha}}
    \otimes_{\beta}}_{\ \ \nu^o}v)
\end{aligned}$$ Then, we have proved that $W^{\text{op}}=\sigma_{\nu^o}W'^*\sigma_{\nu^o}$.
\end{proof}

\begin{coro}
We have $W'=(J_{\widehat{\Phi}}\surl{\ _{\alpha}\otimes_{\beta}}_{\
\
\nu^o}J_{\widehat{\Phi}})\sigma_{\nu}W^*\sigma_{\nu}(J_{\widehat{\Phi}}^*\surl{\
_{\hat{\alpha}}\otimes_{\beta}}_{\ \ \nu^o}J_{\widehat{\Phi}}^*)$.
\end{coro}

\begin{proof}
It is a consequence of the previous proposition and proposition
\ref{antimor}.
\end{proof}

\begin{rema}
The application $j\circ R$, implemented by $J_{\Phi}\hat{J}$, gives
an isomorphism between the measured quantum groupoid and the
opposite of the commutant one.
\end{rema}

\begin{prop}
We have the following equalities:
\begin{center}
\begin{minipage}{12cm}
\begin{enumerate}[i)]
\item
$(N,M,\alpha,\beta,\Gamma,\nu,T_L,T_R)^{\text{op}\wedge}=(N,M,\alpha,\beta,\Gamma,\nu,T_L,T_R)^{\wedge
c }$
\item
$(N,M,\alpha,\beta,\Gamma,\nu,T_L,T_R)^{c\wedge}=(N,M,\alpha,\beta,\Gamma,\nu,T_L,T_R)^{\wedge\text{op}}$
\item
$(N,M,\alpha,\beta,\Gamma,\nu,T_L,T_R)^{c\
\text{op}}=(N,M,\alpha,\beta,\Gamma,\nu,T_L,T_R)^{\text{op}\ c}$
\end{enumerate}
\end{minipage}
\end{center}
\end{prop}

\begin{proof}
The dual of the opposite and the commutant of the dual have the same
basis $N^0$. The von Neumann algebra of the first one is generated
by the operators $(\omega *id)(W^{\text{op}})$ so is equal to
$J_{\widehat{\Phi}}\{(\omega
*id)(W)\}''J_{\widehat{\Phi}}^*=\widehat{M}'$. The representation
and the anti-representation over $N^0$ are both given by
$\widehat{\beta}$ and $\alpha$. Finally, for all $x\in\widehat{M}'$,
we have:
$$\Gamma^{\text{op}\wedge}(x)=\sigma_{\nu} W^{\text{op}}(x\surl{\
_{\alpha}\otimes_{\beta}}_{\
N^o}1)(W^{\text{op}})^*\sigma_{\nu^o}=W'^*(1\surl{\
_{\beta}\otimes_{\alpha}}_{\ N}x)W'$$ By the previous corollary, we
have:$$\Gamma^{\text{op}\wedge}(x)=\sigma_{\nu}(J_{\widehat{\Phi}}\surl{\
_{\alpha}\otimes_{\hat{\beta}}}_{\ \
\nu^o}J_{\widehat{\Phi}})W(J_{\widehat{\Phi}}xJ_{\widehat{\Phi}}\surl{\
_{\beta}\otimes_{\alpha}}_{\ N}1)W^*({\mathcal J}\surl{\
_{\beta}\otimes_{\hat{\alpha}}}_{\
\nu}J_{\widehat{\Phi}})\sigma_{\nu^o}$$
$$=(J_{\widehat{\Phi}}\surl{\ _{\hat{\beta}}\otimes_{\alpha}}_{\ \nu}J_{\widehat{\Phi}})
\hat{\Gamma}(J_{\widehat{\Phi}}xJ_{\widehat{\Phi}})(J_{\widehat{\Phi}}\surl{\
_{\hat{\alpha}}\otimes_{\beta}}_{\
\nu^o}J_{\widehat{\Phi}})=\hat{\Gamma}^{\text{c}}(x)$$ So i) is
proved. ii) comes from i) and the bi-duality theorem.

The opposite of the commutant and the commutant of the opposite have
the same basis $N^0$ and the same von Neumann algebra $M'$. The
representation and the anti-representation are both given by
$\varrho$ and $\widehat{\beta}$. By \cite{Vae}, we have
$J_{\Psi}=\lambda^{i/4}J_{\Phi}$. Then we get, for all $x\in M'$:
$$
\begin{aligned}
\Gamma^{\text{op c}}(x)&=(J_{\Psi}\surl{\
_{\alpha}\otimes_{\hat{\beta}}}_{\ \
\nu^o}J_{\Psi})\varsigma_N\Gamma(J_{\Psi}xJ_{\Psi})(J_{\Psi}\surl{\
_{\alpha}\otimes_{\hat{\beta}}}_{\ \ \nu^o}J_{\Psi})\\
&=\varsigma_{N^o}(J_{\Phi}\surl{\ _{\alpha}\otimes_{\hat{\beta}}}_{\
\ \nu^o}J_{\Phi})\Gamma(J_{\Phi}xJ_{\Phi})(J_{\Phi}\surl{\
_{\alpha}\otimes_{\hat{\beta}}}_{\ \ \nu^o}J_{\Phi})\Gamma^{\text{c
op}}(x)
\end{aligned}$$ because $\lambda^{i/4}\in
Z(M)\cap\alpha(N)\cap\beta(N)$ So iii) is proved.

\end{proof}

\backmatter

\end{document}